\documentclass[11pt]{amsbook}
\usepackage{amscd,amssymb}
\usepackage[mathscr]{eucal}
\usepackage[english]{babel}
\input pstricks
\input pst-node

\usepackage{pstricks, pst-node}
\input xy
\usepackage[all]{xy}

\usepackage{amsmath}

\usepackage[latin1]{inputenc}%
\usepackage{a4}

\newlength{\itemlaenge}

\begin{document}


\newtheoremstyle{mytheorem}
  {}
  {}
  {\slshape}
  {}
  {\scshape}
  {.}
  { }
  {}

\newtheoremstyle{mydefinition}
  {}
  {}
  {\upshape}
  {}
  {\scshape}
  {.}
  { }
  {}

\theoremstyle{mytheorem}
\newtheorem{lemma}{Lemma}[section]
\newtheorem{prop}[lemma]{Proposition}
\newtheorem*{prop*}{Proposition}
\newtheorem{prop_intro}{Proposition}
\newtheorem{cor}[lemma]{Corollary}
\newtheorem{cor_intro}[prop_intro]{Corollary}
\newtheorem{thm}[lemma]{Theorem}
\newtheorem{thm_intro}[prop_intro]{Theorem}
\newtheorem*{thm*}{Theorem}
\theoremstyle{mydefinition}
\newtheorem{rem}[lemma]{Remark}
\newtheorem*{rem*}{Remark}
\newtheorem{rem_intro}[prop_intro]{Remark}
\newtheorem{rems_intro}[prop_intro]{Remarks}
\newtheorem*{notation*}{Notation}
\newtheorem*{warning*}{Warning}
\newtheorem{rems}[lemma]{Remarks}
\newtheorem{defi}[lemma]{Definition}
\newtheorem*{defi*}{Definition}
\newtheorem{defi_intro}[prop_intro]{Definition}
\newtheorem{defis}[lemma]{Definitions}
\newtheorem{exo}[lemma]{Example}
\newtheorem{exo_intro}[prop_intro]{Example}
\newtheorem{exos_intro}[prop_intro]{Examples}

\numberwithin{equation}{section}

\newcommand{\bqn}{\begin{eqnarray*}}
\newcommand{\eqn}{\end{eqnarray*}}

\newcommand{\bibURL}[1]{{\unskip\nobreak\hfil\penalty50{\tt#1}}}

\def\ti{-\allowhyphens}
\newcommand{\thismonth}{\ifcase\month 
  \or January\or February\or March\or April\or May\or June%
  \or July\or August\or September\or October\or November%
  \or December\fi}
\newcommand{\thismonthyear}{{\thismonth} {\number\year}}
\newcommand{\thisdaymonthyear}{{\number\day} {\thismonth} {\number\year}}
%
%
%

\newcommand{\Stab}{\operatorname{Stab}}
\newcommand{\sign}{\operatorname{sign}}
\newcommand{\aut}{\operatorname{Aut}}
\newcommand{\Is}{\operatorname{Is}}
\newcommand{\codim}{\operatorname{codim}}
\newcommand{\esssup}{\operatorname{ess\,sup}}
\newcommand{\supp}{\operatorname{supp}}
\newcommand{\BB}{{\mathbb B}}
\newcommand{\CC}{{\mathbb C}}
\newcommand{\DD}{{\mathbb D}}
\newcommand{\FF}{{\mathbb F}}
\newcommand{\HH}{{\mathbb H}}
\newcommand{\GG}{{\mathbb G}}
\newcommand{\KK}{{\mathbb K}}
\newcommand{\NN}{{\mathbb N}}
\newcommand{\PP}{{\mathbb P}}
\newcommand{\QQ}{{\mathbb Q}}
\newcommand{\RR}{{\mathbb R}}
\newcommand{\TT}{{\mathbb T}}
\newcommand{\ZZ}{{\mathbb Z}}

\newcommand{\Aa}{{\mathcal A}}
\newcommand{\Bb}{{\mathcal B}}
\newcommand{\Cc}{{\mathcal C}}
\newcommand{\Dd}{{\mathcal D}}
\newcommand{\Ee}{{\mathcal E}}
\newcommand{\Ff}{{\mathcal F}}
\newcommand{\Hh}{{\mathcal H}}
\newcommand{\Jj}{{\mathcal J}}
\newcommand{\Ll}{{\mathcal L}}
\newcommand{\Mm}{{\mathcal M}}
\newcommand{\Nn}{{\mathcal N}}
\newcommand{\Oo}{{\mathcal O}}
\newcommand{\Qq}{{\mathcal Q}}
\newcommand{\Pp}{{\mathcal P}}
\newcommand{\Rr}{{\mathcal R}}
\newcommand{\Ss}{{\mathcal S}}
\newcommand{\Tt}{{\mathcal T}}
\newcommand{\Vv}{{\mathcal V}}
\newcommand{\Xx}{{\mathcal X}}
\newcommand{\Zz}{{\mathcal Z}}

\newcommand{\frakg}{{\mathfrak g}}
\newcommand{\frakk}{{\mathfrak k}}
\newcommand{\frakh}{{\mathfrak h}}
\newcommand{\fraka}{{\mathfrak a}}
\newcommand{\frake}{{\mathfrak e}}
\newcommand{\frakp}{{\mathfrak p}}
\newcommand{\frako}{{\mathfrak o}}
\newcommand{\fraks}{{\mathfrak s}}
\newcommand{\frakl}{{\mathfrak l}}
\newcommand{\frakm}{{\mathfrak m}}
\newcommand{\frakn}{{\mathfrak n}}
\newcommand{\fraku}{{\mathfrak u}}
\newcommand{\frakc}{{\mathfrak c}}
\newcommand{\frakr}{{\mathfrak r}}
\newcommand{\frakpp}{{\mathfrak p_+}}
\newcommand{\frakpm}{{\mathfrak p_-}}
\newcommand{\frakkc}{{{\mathfrak k}_\CC}}
\newcommand{\frakgc}{{{\mathfrak g}_\CC}}
\newcommand{\frakpc}{{{\mathfrak p}_\CC}}
\newcommand{\frakB}{{\mathfrak B}}

\renewcommand{\a}{\alpha}
\newcommand{\e}{\epsilon}
\newcommand{\eps}{\epsilon}
\renewcommand{\b}{\beta}
\newcommand{\g}{\gamma}
\newcommand{\G}{\Gamma}
\renewcommand{\L}{\Lambda}
\renewcommand{\l}{\lambda}

\newcommand{\dD}{{\mathbf D}}
\newcommand{\gG}{{\mathbf G}}
\newcommand{\hH}{{\mathbf H}}
\newcommand{\pP}{{\mathbf P}}
\newcommand{\lL}{{\mathbf L}}
\newcommand{\qQ}{{\mathbf Q}}
\newcommand{\nN}{{\mathbf N}}
\newcommand{\wW}{{\mathbf W}}
\newcommand{\uU}{{\mathbf U}}

\newcommand{\<}{\langle}
\renewcommand{\>}{\rangle}

\def\ol{\overline}

\def\h{{\rm H}}
\def\hb{{\rm H}_{\rm b}}
\def\ehb{{\rm EH}_{\rm b}}
\def\ha{{\rm H}_{(G,K)}}
\def\hc{{\rm H}_{\rm c}}
\def\hcb{{\rm H}_{\rm cb}}
\def\ehbc{{\rm EH}_{\rm cb}}
\def\linfty{L^\infty}
\def\linftyw{L^\infty_{\rm w*}}
\def\linftya{L^\infty_{\mathrm{w*,alt}}}
\def\la{L^\infty_{\mathrm{alt}}}
\def\cb{{\rm C}_{\rm b}}
\def\binfty{\mathcal B^\infty_{\mathrm alt}}

\def\one{\mathbf{1\kern-1.6mm 1}}
\def\sous#1#2{{\raisebox{-1.5mm}{$#1$}\backslash \raisebox{.5mm}{$#2$}}}
\def\rest#1{{\raisebox{-.95mm}{$\big|$}\raisebox{-2mm}{$#1$}}}
\def\homeo#1{{\sl H\!omeo}^+\!\left(#1\right)}
\def\thomeo#1{\widetilde{{\sl \!H}\!omeo}^+\!\left(#1\right)}
\def\bu{\bullet}
\def\weak{weak-* }
\def\property{\textbf{\rm\textbf A}}
\def\cont{\mathcal{C}}
\def\id{{\it I\! d}}
\def\opposite{^{\rm op}}
\def\oddex#1#2{\left\{#1\right\}_{o}^{#2}}
\def\comp#1{{\rm C}^{(#1)}}
\def\ro{\varrho}
\def\ti{-\allowhyphens}
\def\lra{\longrightarrow}

\def\fix{{\operatorname{Fix}}}
\def\Mat{{\operatorname{Mat}}}
\def\zent{{\operatorname{Zent}}}
\def\C{{\operatorname{C}}}
\def\sym{{\operatorname{Sym}}}
\def\stab{{\operatorname{Stab}}}
\def\arg{{\operatorname{arg}}}
\def\HTP{{\operatorname{HTP}}}
\def\h2{{\operatorname{H_2}}}
\def\h1{{\operatorname{H_1}}}
\def\pr{{\operatorname{pr}}}
\def\rk{{\operatorname{rank}}}
\def\tr{{\operatorname{tr}}}
\def\codim{{\operatorname{codim}}}
\def\nt{{\operatorname{nt}}}
\def\d{{\operatorname{d}}}
\def\Gr{{\operatorname{Gr}}}
\def\id{{\operatorname{Id}}}
\def\ker{{\operatorname{ker}}}
\def\im{{\operatorname{im}}}

\def\ad{\operatorname{ad}}
\def\Ad{\operatorname{Ad}}
\def\adg{\operatorname{ad}_\frakg}
\def\adp{\operatorname{ad}_\frakp}
\def\bg{B_\frakg}
\def\creg{C_{\rm reg}}
\def\cs{\check S}
\def\cst{{\check S}^{(3)}}
\def\det{{\operatorname{det}}}
\def\deta{{\operatorname{det}_A}}
\def\diag{{\operatorname{diag}}}
\def\gmodp{\gG(\RR)/\pP(\RR)}
\def\gmodq{\gG(\RR)/\qQ(\RR)}
\def\hom{\operatorname{Hom}}
\def\isp{\operatorname{Is}_{\<\cdot,\cdot\>}}
\def\isptwo{\operatorname{Is}_{\<\cdot,\cdot\>}^{(2)}}
\def\ispth{\operatorname{Is}_{\<\cdot,\cdot\>}^{(3)}}
\def\isf{\operatorname{Is}_F}
\def\isfi{\operatorname{Is}_{F_i}}
\def\isft{\operatorname{Is}_F^{(3)}}
\def\isfit{\operatorname{Is}_{F_i}^{(3)}}
\def\isftwo{\operatorname{Is}_F^{(2)}}
\def\lin{\operatorname{Lin}(L_+,L_-)}
\def\ll{{\Ll_1,\Ll_2}}
\def\kahler{ K\"ahler }
\def\kg{\kappa_G}
\def\kgb{\kappa_G^{\rm b}}
\def\kib{\kappa_i^{\rm b}}
\def\kibt{\tilde\kappa_i^{\rm b}}
\def\oc{\overline c}
\def\om{\overline m}
\def\pii{\pi_i}
\def\piit{\tilde\pi_i}
\def\psupq{{\rm PSU}(p,q)}
\def\psuvi{{\rm PSU}\big(V,\<\cdot,\cdot\>_i\big)}
\def\rg{r_G}
\def\slv{{\rm SL}(V)}
\def\supq{{\rm SU}(p,q)}
\def\suv{{\rm SU}\big(V,\<\cdot,\cdot\>\big)}
\def\suvi{{\rm SU}\big(V,\<\cdot,\cdot\>_i\big)}
\def\jm{J^{(m)}}
\def\val{\operatorname{Val}}
\def\vc{V_\CC}
\def\vconj{\overline v}
\def\wconj{\overline w}

\def\binfty{\mathcal B^\infty_{\mathrm {alt}}}
\def\hcb{{\rm H}_{\rm cb}}
\def\to{\rightarrow}
\def\bu{\bullet}
\def\la{L^\infty_{\mathrm{alt}}}
\def\hb{{\rm H}_{\rm b}}
\def\hc{{\rm H}_{\rm c}}
\def\h{{\rm H}}
\def\kxb{{\kappa_X^{\rm b}}}
\def\sk{\mathsf k}
\renewcommand{\phi}{\varphi}

\def\eI{{^{I}\!\mathrm{E}}}
\def\eII{{^{I\!I}\!\mathrm{E}}}
\def\dI{{^{I}\!d}}
\def\dII{{^{I\!I}\!d}}
\def\fI{{^{I}\!F}}
\def\fII{{^{I\!I}\!F}}

\def\ind{\mathbf{i}}
\def\Ind{\mathrm{Ind}}
\def\No{N\raise4pt\hbox{\tiny o}\kern+.2em}
\def\no{n\raise4pt\hbox{\tiny o}\kern+.2em}
\def\bsl{\backslash}
%
%
%
%

\newcommand{\mat}[4]{
\left(
\begin{array}{cc}
{\scriptstyle #1} & {\scriptstyle #2} \\ {\scriptstyle #3} & {\scriptstyle
#4}
\end{array}\right)}

\newcommand{\matt}[6]{
\left(
\begin{array}{cc}
{\scriptstyle #1} & {\scriptstyle #2} \\ {\scriptstyle #3} & {\scriptstyle
#4}\\{\scriptstyle #5} & {\scriptstyle #6}
\end{array}\right)}

\newcommand{\be}[3]{{\beta(y_{#1},y_{#2},y_{#3})}}
\newcommand{\db}[4]{0=d\beta(y_{#1},y_{#2},y_{#3},y_{#4})
= \beta(y_{#1},y_{#2},y_{#3})+\beta(y_{#1},y_{#3},y_{#4})-\beta(y_{#1},y_{#2},y_{#4})
-\beta(y_{#2},y_{#3},y_{#4})}
%
%
%
%
%

\title[Bounded cohomology and geometry]{Bounded Cohomology and Geometry}
\author{Anna Wienhard}
\email{wienhard@math.uni-bonn.de}
\address{Rheinische Friedrich-Wilhelms-Universit\"at Bonn, Beringstrasse 1, D-53115 Bonn, Germany}




%
%
\maketitle
%
%
%
%
%


\tableofcontents

\chapter*{Preface and Acknowledgement}
My first encounter with bounded cohomology dates back several years ago, 
when I participated in a student seminar organized by Prof. Ballmann, where 
we discussed Thurston's proof of Mostow rigidity for hyperbolic manifolds 
using Gromov's notion of simplicial volume. Since then 
the links between bounded cohomology and geometry kept on fascinating me. 
Thanks to my advisor, I had the opportunity to talk to Nicolas Monod during 
a conference in Oberwolfach, getting to know 
about the functorial approach to bounded cohomology developed by 
Burger and Monod \cite{Burger_Monod_GAFA}.
After trying to solve some other problems, 
I happily 
ended up thinking about bounded cohomology and Hermitian symmetric spaces, 
finding the mixture of methods from bounded cohomology and 
from the geometry of symmetric spaces of higher rank which I like a lot.

I am grateful to my advisor Prof. Werner Ballmann 
who taught me a lot of mathematics and mathematical thinking since my first 
semester at university. Without his ideas, the freedom and the advise 
he gave me and his support, I would never have started or finished this work.
I owe special thanks to Prof. Marc Burger and Prof. Alessandra Iozzi, 
not only for all what I learned from them during my stays in Z\"urich, but 
especially for the time we spent working together.
Several results are actually joint work with them.\footnote{Some of the results have been announced in 
\cite{Burger_Iozzi_Wienhard_ann}.}  
This work owes much to ideas of Prof. Marc Burger and I would like to express 
my thanks to him for becoming a 
second advisor for me and for accepting to be the
 second referee of this thesis.
I am happy that Prof. Flume and Prof. Albeverio accepted to
participate in my ``Promotionskommission''.
I thank Theo B\"uhler, who patiently corrected a lot of mistakes and
the proof of Lemma~\ref{lem:dichotomy} (Chapter 8) in an earlier
version of this work. 
The remaining mistakes remain my responsibility.

Several instutions made it possible for me to travel a lot during the 
time of my Ph.D. My thanks go to all of them.
The fellowship I received from the ``Studienstiftung des deutschen Volkes'' 
for my Ph.D. studies gave me the financial independence to leave the 
University of Bonn to visit different places. 
I thank the ``Forschungsinstitut f\"ur Mathematik'' 
at the ETH-Z\"urich for its 
hospitality and financial support during a six month stay 
which was very important for my studies. The 
``Institute for Pure and Applied Mathematics'' at UCLA 
gave me the opportunity to participate 
in a three month program on Symplectic Geometry. I enjoyed my one month stays 
at ``Mathematical Science Research Institute'' in Berkeley and at the 
``Institut Joseph Fourier'' in Grenoble. 
The financial support from the ``Bonner International Graduate School'' 
and the 
``Sonderforschungsbereiche 256/611'' were essential for the realization 
of all my visits of other institutes and my participation in 
several conferences.

I would not have survived without the kindness, encouragement and company 
of colleagues and friends I had and have 
in Bonn, Z\"urich, Grenoble and several other 
places of the world. Thank you.


\vskip1cm
\chapter*{Introduction}
This work is motivated by interesting relations between bounded cohomology 
and geometric properties of spaces and groups. 
We want to give a contribution to the understanding 
of these relations in the context of Hermitian 
symmetric spaces of noncompact type.

The first who introduced bounded cohomology in the realm of geometry 
was Gromov with his celebrated paper ``Volume and bounded cohomology'' 
\cite{Gromov_82}, where he defined  
the simplicial volume to give a lower bound for the 
minimal volume of a Riemannian manifold. 
Several people followed him relating bounded cohomology to geometric 
properties of manifolds and groups. We recall some results where 
bounded cohomology is used in the context of geometric problems. 
Thurston \cite{Thurston} gave a proof of Mostow's rigidity theorem 
for hyperbolic manifolds using bounded cohomology. Ghys
\cite{Ghys_87} 
has proven that the bounded Euler class classifies actions by orientation 
preserving homeomorphism 
on a circle up to semi-conjugacy. In a different way Soma 
\cite{Soma_tame, Soma_third} used bounded cohomology to exhibit 
rigidity phenomena of hyperbolic three-manifolds. 
An important aspect of how much geometric information bounded 
cohomology might contain is given by the characterization 
of hyperbolic 
groups in terms of bounded cohomology obtained by Mineyev 
\cite{Mineyev_straighten, Mineyev_bounded}.
New geometric applications of bounded cohomology emerged from the
functorial approach to bounded cohomology developed by Burger and Monod 
in \cite{Burger_Monod_JEMS, Burger_Monod_GAFA, Burger_Iozzi_App, Monod_book} 
containing 
already some applications. 
Further applications to measure and orbit equivalence of groups using 
this approach were given 
by Monod and Shalom \cite{Monod_Shalom_1, Monod_Shalom_2}. 
Making use of this functorial approach as well, but exploring it in another 
direction Burger and Iozzi studied several rigidity questions about  
Hermitian symmetric spaces \cite{Burger_Iozzi_00, Burger_Iozzi_def,
Burger_Iozzi_supq}, in particular about complex hyperbolic space, 
relating them to bounded cohomology. 
Our work lies very much in the direction of the last mentioned results 
relating the geometry of Hermitian symmetric spaces and properties 
of groups acting isometrically on them to bounded cohomology classes. 
The restriction to Hermitian symmetric spaces is not superficial 
but based on the still very restricted knowledge of bounded cohomology in the 
context of arbitrary symmetric spaces.

Some results presented in this thesis are 
concerned with 
Hermitian symmetric spaces, their Shilov boundaries and a particular class 
of totally geodesic embeddings of Hermitian symmetric spaces which is defined 
using bounded cohomology. The other results concern finitely generated groups 
acting isometrically on Hermitian symmetric spaces. Several of the
latter 
results are joint work with Marc Burger and Alessandra Iozzi \cite{Burger_Iozzi_Wienhard_ann}.

A Hermitian symmetric space $X$ of noncompact type 
is a symmetric space $X=G/K$ 
of noncompact type endowed with a 
$G$-invariant complex structure. 
The imaginary part of the $G$-invariant Hermitian metric on $X$ defines 
a $G$-invariant closed two form, 
the K\"ahler form $\omega \in \Omega^2(X)^G$.
Let $c_G:G\times G\times G \to \RR$ be the map defined by integration of $\omega$ over geodesic triangles in $X$:
\bqn
c_G(g_1,g_2,g_3):= \int_{\Delta(g_1 x,g_2 x,g_3 x)} \omega,
\eqn
where $\Delta(g_1 x,g_2 x,g_3 x)$ is a geodesic triangle with vertices in  
$g_1 x,g_2 x,g_3 x$ for $x$ in $X$ arbitrary.  
Then $c_G$ is a bounded $2$-cocycle on the group $G$ 
\cite{Gromov_82,Dupont_Guichardet,Domic_Toledo,Clerc_Orsted_2}, i.e. 
a continuous uniformly bounded function $c_G:G\times G\times G \to \RR$ 
satisfying the cocycle identity 
$c_G(g_1,g_2,g_3)-
c_G(g_1,g_2,g_4)+
c_G(g_1,g_3,g_4)-
c_G(g_2,g_3,g_4)=0$.

Every Hermitian symmetric space of noncompact type admits a
realization as a {\em bounded symmetric domain} $\Dd$ in $\CC^N$. 
The {\em Shilov boundary} $\cs$ is a special part of the boundary
$\partial \Dd$ of $\Dd$. 
The geometric applications of bounded cohomology discussed in this
work are based on understanding 
relations of this bounded cocycle or ``its brother'' which is defined on 
triples of points on the {\em Shilov boundary} $\cs$ and the geometry
of $\Dd$.
Define $\b_\Dd: \Dd^3 \to \RR$ as above by 
$\b_\Dd(x_1,x_2,x_3):= \int_{\Delta(x_1,x_2,x_3)} \omega$,
where $\Delta(x_1,x_2 ,x_3 )$ is the geodesic triangle in $\Dd$ 
with vertices in 
$x_1,x_2 ,x_3 $. The cocycle $\b_\Dd$ extends continuously to the space of 
 triples of pairwise transverse points $\cs^{(3)}$ in the Shilov boundary $\cs$ 
of $\Dd$. Let $\b$ 
denote the restriction of $\b_\Dd$ to $\cs^{(3)}$.
Both cocycles, $c_G$ and $\b$ determine the same bounded cohomology class 
$k_X^b \in \hcb^2(G)$, called the {\em bounded K\"ahler class}. This class 
was already used in the work of Burger and Iozzi  
(\cite{Burger_Iozzi_00, Burger_Iozzi_def,Burger_Iozzi_supq}). 
It is an important ingredient for our results.

\subsection*{Hermitian symmetric spaces, their Shilov boundaries and totally geodesic embeddings}
\subsubsection*{Tube type and not}
Any 
Hermitian symmetric space is biholomorphically equivalent to a bounded symmetric 
domain $\Dd\subset \CC^N$, generalizing the Poincar\'e-disc model of the 
hyperbolic plane $\HH^2$. 
The generalization of the upper half space model is not available for all 
Hermitian symmetric spaces, but only for a subclass of them. 
Namely, a Hermitian symmetric space $X$ is said to be of {\em tube type}, 
if $X$ is biholomorphically equivalent to a tube domain 
$T_\Omega=V+ i\Omega$, where $V$ is a 
real vectorspace and $\Omega\subset V$ is an open cone.
The difference of being of tube type or not is reflected in the root 
system and in properties of the Shilov boundary $\cs$. 
We describe different models for the Shilov boundary, 
giving explicit parametrizations of the Shilov boundary of Hermitian 
symmetric spaces of tube type.
The Shilov boundary $\cs$ may in particular be seen 
as a generalized flag variety, $\cs=G/Q$, where 
$Q$ is a specific maximal parabolic subgroup of $G$.
As a consequence of a detailed study of the generalized Bruhat decomposition
 of $\cs$ 
we obtain a new characterization 
of the Hermitian symmetric spaces which are not 
of tube type in terms of triples on the Shilov boundary. 

\begin{thm*}[Theorem~\ref{connected} in Chapter 4]
Denote by $\cs^{(3)}$ the space of triples of pairwise transverse points in 
the Shilov boundary $\cs$ of $X$. 
Then $X$ is not of tube type if and only if $\cs^{(3)}$ is connected. 
\end{thm*}

As a corollary of the new classification of non-tube type 
Hermitian symmetric spaces, we obtain 
\begin{prop*}[Proposition~\ref{infinite_values} in Chapter 5]
Denote by $\cs^{(3)}$ the space of triples of pairwise transverse points in 
the Shilov boundary $\cs$ of $X$. 
Then there exists a normalization constant 
$c\in \RR$
such that\\
(1) $\b(\cs^{(3)})= [-c r_X, +c r_X]$  if $X$ is not of tube type\\
or\\
(2) $\b(\cs^{(3)})= c(r_X+ 2\ZZ) \cap [-c r_X, +c r_X]$, if $X$ is of tube type,\\ 
where $r_X$ is the rank of $X$.
\end{prop*}

Another consequence of our study of the generalized Bruhat decomposition 
of $\cs$ is the definition of a  Maslov index 
for all Hermitian symmetric 
spaces of tube type generalizing the usual Maslov index for the
symplectic group. The set of points $x\in \cs$ which are not 
transverse to $z\in \cs$ 
define a codimension one cycle in $\cs$. 
The intersection number of a loop with this cycle defines the generalized Maslov index 
$\mu \in \h^1(\cs, \ZZ)$.

\subsubsection*{Tight embeddings and homomorphisms}
Let $H, G$ be locally compact 
groups. A continuous homomorphism $\pi: H\to G$ 
induces canonical pull-back maps $\pi^*$ in continuous cohomology and 
$\pi_b^*$ in continuous bounded cohomology, such that the following diagram commutes:
\bqn
\xymatrix{
\hcb^*(G) \ar[d]^\kappa \ar[r]^{\pi_b^*} & \hcb^*(H) \ar[d]^\kappa\\
\hc^*(G) \ar[r]^{\pi^*} & \hc^*(H),
}
\eqn
where $\kappa$ is the natural comparison map between continuous bounded 
cohomology and continuous cohomology.

The continuous bounded cohomology groups come equipped with a canonical 
seminorm $||\cdot||$, with respect to which $\pi_b^*$ is norm decreasing, 
i.e. 
\bqn
||\pi_b^*(\a)|| \leq ||\a|| \text{ for all } \a\in \hcb^*(G).
\eqn

Requiring that a homomorphism preserves the norm of some classes in 
$\hcb^*(G)$ imposes restrictions on it. We study this in a 
specific situation.
When $G$ is a connected semisimple Lie group 
and its associated symmetric space $X$ is Hermitian, we call 
a continuous homomorphism $\pi:H\to G$ {\em tight} if 
it preserves the norm of the bounded K\"ahler class $k_X^b\in
\hcb^2(G)$, i.e. if $||\pi_b^*(k_X^b)||= ||k_X^b||$.
One can show that the notion of tightness only depends on the $G$-invariant complex structure on $X$ and 
not on the particular choice of a compatible Hermitian metric. 
An important point is the following 
property of tight homomorphisms.

\begin{thm*}[Theorem~\ref{thm:tightness} in Chapter 6]
Let $H$ be a locally compact group and $G$ a connected semisimple Lie
group with finite center, which is of Hermitian type. 
Suppose $\pi:H\to G$ is a continuous tight homomorphism. 
Then the Zariski closure of the image $L=\overline{\pi(H)}^Z$ 
is reductive with compact center. The symmetric space associated to 
the semisimple part of $L$ is Hermitian symmetric.
\end{thm*}

Clearly, there is a parallel notion for symmetric spaces. 
Namely, let $Y, X$ be symmetric spaces of noncompact type 
with $X$ Hermitian. 
A totally geodesic embedding $f: Y\to X$ is tight if the corresponding homomorphism, 
$\pi: H_Y\to {\rm Is}(X)^\circ$ is tight, where $H_Y$ is an appropriate finite extension of ${\rm Is}(Y)^\circ$.
We can give a more geometric definition of tightness of totally
geodesic embeddings. 
A totally geodesic embedding $f:Y\to X$ is tight if and only if 
\bqn
\sup_{\Delta\subset Y} \int_\Delta f^*\omega_X = \sup_{\Delta\subset X} \int_\Delta \omega_X
\eqn
where the suprema are taken over all geodesic triangles $\Delta$ in
$Y$ or $X$, respectively. 

We give some examples of tight embeddings:\\
1) Every tight holomorphic embedding from the disc $\DD$ into a 
Hermitian symmetric space $X$ is given by a diagonal embedding into a maximal 
polydisc in $X$.\\
2) The $2n$-dimensional irreducible representation ${\rm SL}(2,\RR) \to {\rm Sp}(2n, \RR)$, 
is tight for every $n\geq 1$. However, the associated tight embedding $\DD\to \Hh_n$ is holomorphic if and only if $n=1$.

This shows that tight embeddings are not necessarily 
holomorphic. Still, 
they behave nicely with respect to some structures of
 Hermitian symmetric spaces.

\begin{prop*}[Proposition~\ref{tight_embedding} in Chapter 6]
A tight embedding $f:Y\to X$ extends uniquely to a continuous
equivariant map $\ol{f}: \cs_Y\to \cs_X$ between 
the Shilov boundaries.
\end{prop*} 

As a consequence there do not exist tight embeddings of 
non-tube type Hermitian symmetric spaces 
into Hermitian symmetric spaces of tube type. 
Methods from bounded cohomology (see \cite{Burger_Monod_GAFA}) and 
results of \cite{Clerc_Orsted_2} 
are the main tools used to obtain these properties of tight embeddings and 
homomorphisms.

For totally geodesic embeddings into 
Hermitian symmetric spaces of tube type we exhibit a simple criterion 
of tightness in terms of the corresponding 
Lie algebra homomorphisms. 
Using this criterion we classify tight embeddings of the Poincar\'e
disc. All of them are products of embeddings obtained from irreducible representations 
of $\fraks\frakl(2,\RR)$ into $\fraks\frakp(2n, \RR)$. 
More precisely: 

\begin{thm*}[Theorem~\ref{thm:hull} in Chapter 6]
Suppose that $X$ is a Hermitian symmetric space of noncompact type and
 $f:\DD\to X$ is a tight totally geodesic embedding. 
Then the smallest Hermitian symmetric subspace $Y\subset X$ 
containing $f(\DD)$ 
is a product of Hermitian symmetric subspaces $Y_i$ of $X$, 
$Y= \Pi_{i=1}^k Y_i$,
where $Y_i$ are Hermitian symmetric spaces 
associated to symplectic groups ${\rm Sp}(2n_i,\RR)$. 
Moreover, $\sum_{i=1}^k n_i \leq r_X$ and the embeddings  $f_i:\DD\to Y_i$ 
correspond to irreducible representations 
$\fraks\frakl(2,\RR)\to \fraks\frakp(2n_i, \RR)$.
\end{thm*}

\subsection*{Applications to group actions on Hermitian symmetric spaces}
Let $\G$ be a finitely generated group. An action of $\G$ by isometries on a Hermitian 
symmetric space $X=G/K$ is given by a representation $\rho:\G\to G$. 
The bounded K\"ahler class $k_X^b\in \hcb^2(G)$ gives rise to the 
{\em bounded K\"ahler invariant} $\rho^*(k_X^b) \in \hb^2(\G)$ 
of a representation $\rho:\G\to G$ by taking the pull-back in 
bounded cohomology. Our results are based on this invariant, 
first defined in \cite{Burger_Iozzi_supq}. 
Note that tight embeddings where defined by using the weaker numerical 
invariant obtained from the bounded K\"ahler invariant by taking its norm. 
\subsubsection*{Zariski dense representations}
The information encoded by the bounded K\"ahler invariant depends quite 
a lot on the situation. 
In the case where the Hermitian symmetric space is not of tube type 
and the action of $\G$ does not preserve any proper subspace of $X$, 
the bounded K\"ahler invariant  
$\rho^*(k_X^b) \in H^2_b(\G)$  determines the 
action of $\G$ on $X$. 
More precisely we prove the following generalization of the 
result in \cite{Burger_Iozzi_supq} to all Hermitian symmetric spaces.
\begin{thm*}[Theorem~\ref{thm:z_dense} in Chapter 7]
Assume that $X=G/K$ is a Hermitian symmetric space which is not of 
tube type. 
Let $\rho:\G\to G$ be a representation with Zariski dense image. 
Then 
\begin{enumerate}
\item{$\rho^*(k_X^b) \neq 0$.}
\item{$\rho^*(k_X^b)$ determines $\rho$ up to conjugation with
    elements of $G$.}
\end{enumerate}
\end{thm*}
Unfortunately we do not know how to extract geometric information 
about the group action from the bounded K\"ahler invariant. 
We first have to form a numerical invariant out of it.
\subsubsection*{A generalization of Teichm\"uller space}
In the special case when $\G=\G_g$ is the fundamental group of a closed 
Riemann surface $\Sigma_g$ of genus $g\geq 2$, 
we may evaluate the (bounded) K\"ahler 
invariant on the fundamental class $[\Sigma_g]$ of the surface 
to obtain a numerical invariant of the representation, known as the {\em Toledo invariant}. 
Thus from the bounded K\"ahler invariant $\rho^*(k_X^b)$ we obtain a 
function 
\bqn
{\rm Tol}: \hom(\G_g, G)\to \RR,
\eqn
which is constant on connected components and satisfies  
(with appropriate normalizations) 
\bqn\label{eqn:tol}
|{\rm Tol}(\rho)|\leq 4\pi (g-1)r_X,
\eqn
where $r_X$ is the rank of $X$. 
We call a representation {\em maximal} if $|{\rm Tol}(\rho)|= 4\pi (g-1)r_X$.

In the case $G={\rm PSL}(2,\RR)$ the Toledo invariant is called Euler
number and  
takes values in $2\pi \ZZ$; 
Goldman \cite{Goldman_thesis,Goldman_88} 
has shown that it separates connected 
components of the representation variety and 
that the set of maximal representations consists of the two Teichm\"uller components.
Toledo \cite{Toledo_89} 
showed that a maximal representation $\rho: \G_g\to {\rm PU}(n,1)$ 
stabilizes a complex geodesic in complex hyperbolic $n$-space $\HH^n_\CC$.
Generalizations of Toledo's result for $G={\rm PU}(n,m)$ were proven in 
\cite{Hernandez} and \cite{Bradlow_GarciaPrada_Gothen}. 
We prove 
\begin{thm*}[Theorem~\ref{thm:main} in Chapter 8, see also
  \cite{Burger_Iozzi_Wienhard_ann}]
Let $G$ be a connected semisimple real algebraic group and assume that its associated 
symmetric space is Hermitian. 
Let $\rho :\G_g \rightarrow \rm G$
be a maximal representation. Then:\\
(1) $\rho(\G)$ stabilizes a maximal Hermitian symmetric subspace of
tube type $T$ in $X$.\\ 
(2) The Zariski closure $L$ of $\rho (\G_g)$ is reductive.\\
(3) The symmetric subspace associated to $L$ is a Hermitian symmetric space of tube type and the inclusion $Y\to T$ is tight. \\
(4) The $\G_g$-action on $Y$ (via $\rho$) is properly discontinuous without fixed points.\\
(5) There exists a Hermitian symmetric subspace of tube type $T_Y$
containing $Y$ on which $\rho(\G_g)$ acts properly discountinuously
without fixed points.
\end{thm*}
\begin{rem*}
By a Hermitian symmetric subspace we understand a symmetric space
which is totally geodesically and holomorphically embedded into $X$.
\end{rem*}
Concerning the optimality of this theorem, recall that a maximal polydisc $P$ 
in a Hermitian symmetric space $X$ is a Hermitian symmetric subspace isomorphic to $\DD^r$, 
$r=r_X$. Maximal polydiscs exist and are ${\rm Is}(X)^\circ$-conjugate. 
Denote by 
\bqn
h_P: {\rm SU}(1,1)^r \to {\rm Is}(X)^\circ
\eqn
the homomorphism associated to the inclusion $P\to X$. 

\begin{thm*}[Theorem~\ref{Zariski_dense} in Chapter 8, see also
  \cite{Burger_Iozzi_Wienhard_ann}]
Let $X$ be a Hermitian symmetric space of tube type and let $P\subset
X$ be a maximal polydisc. 
Let $\rho_0: \G_g\to {\rm Is}(X)^\circ$ be 
the maximal representation obtained from 
composing a diagonal discrete injective embedding 
$\G_g\to {\rm SU}(1,1)^r$ with $h_P$. 

Then $\rho_0$ admits a continuous deformation 
\bqn
\rho_t:\G_g\to {\rm Is}(X)^\circ\, , \, t\geq 0
\eqn
with $\rho_t$ maximal and $\rho_t(\G_g)$ Zariski dense in ${\rm
  Is}(X)^\circ$, for all $t>0$.
\end{thm*}

Maximal representations of surface groups give rise to Kleinian groups 
on tube type domains, so one might turn the attention to 
questions regarding limit sets. 
Let $\Dd\subset \CC^n$ be the Harish-Chandra realization 
of a Hermitian symmetric space as a bounded symmetric domain. Then $\aut(\Dd)$ acts on $\partial \Dd$ 
with finitely many orbits and a unique closed one, which coincides with the Shilov boundary $\cs$ of $\Dd$.

For a discrete subgroup $\Lambda <\aut(\Dd)$ and $x\in \Dd$, define  
$\Ll_\Lambda := \overline{\Lambda x} \cap \cs$, 
where  $\overline{\Lambda x}$ is the topological closure in $\CC^n$. 
We will call $\Ll_\Lambda$ the Shilov limit set of $\Lambda$.
While $ \overline{\Lambda x} \cap \partial \Dd$ may depend on the base point $x$, $\Ll_\Lambda$ does not. 
It can however happen that 
$\Ll_\Lambda$ is void even if $\Lambda$ is unbounded. 
\begin{thm*}[Theorem~\ref{thm:limit_set} in Chapter 8]
Let $\Dd$ be a bounded symmetric 
domain of tube type and let $\rho: \G_g\to \aut(\Dd)$ 
be a maximal representation. 
Then, the Shilov limit set $\Ll_\rho \subset \cs$ 
of $\rho(\G_g)$ is a rectifiable circle.
\end{thm*}
Maximal representations of surface groups into isometry groups of Hermitian symmetric spaces of 
tube type provide (following Goldman's theorem) 
a generalization of Teichm\"uller space in the context of Hermitian 
symmetric spaces. 
This generalization of Teichm\"uller space fits together with another approach 
generalizing Teichm\"uller space in the context of split real forms 
of complex Lie groups.
Hitchin \cite{Hitchin} 
singled out a component in the space $\hom^{red}(\G_g, G)/G$ of $G$-conjugacy 
classes of reductive representations, which is homeomorphic to 
$\RR^{(2g-2)\dim(G)}$. 
For $G={\rm PSL}(n,\RR)$ this Hitchin component is the one containing 
the homomorphisms obtained by composing a uniformization $\G_g\to {\rm PSL}(2,\RR)$ 
with the $n$-dimensional irreducible representation of ${\rm PSL}(2,\RR)$. 
Labourie has recently shown \cite{Labourie_anosov} that representations in 
this Hitchin component are faithful, irreducible, and have discrete image. 
He defined the notion of an Anosov representation and showed that all 
representations in this Hitchin component are Anosov. 

For the group $G$ common to Hitchin's and our setting, i.e. 
$G={\rm Sp}(2n, \RR)$, 
the Hitchin component is (properly) contained in the set of 
maximal representations.
We show that maximal representations also share the property 
of being Anosov.
In order to convey the fundamental idea, consider a 
representation $\rho: \G_g\to {\rm Sp}(V)$ into the symplectic group of a 
nondegenerate symplectic form on a real vector space $V$. 
Fix a hyperbolization $\G_g<{\rm PU}(1,1)$, let $T^1\DD$ be 
the unit tangent bundle of $\DD$ and $(g_t)$ 
the geodesic flow. On the vector bundle
\bqn
E_\rho= \G_g\backslash (T^1\DD\times V)
\eqn
with basis $T^1 \Sigma_g$, where $\Sigma_g=\G_g\backslash\DD$, 
we define the flow $\tilde{\psi}_t(u,v)=(g_t(u), v)$.
Fix a continuous scalar product on $E_\rho$.

\begin{thm*}[Theorem~\ref{thm:anosov_splitting} in Chapter 8]
Assume that $\rho:\G_g\to {\rm Sp}(V)$ is a maximal representation. Then there is a continuous 
$\tilde{\psi}_t$-invariant splitting 
\bqn
E_\rho= E^+_\rho \oplus E^-_\rho
\eqn
such that \\
(1) The splitting $ E_{\rho, u}= E^+_{\rho , u} \oplus E^-_{\rho, u}$ is a decomposition into complementary Lagrangian subspaces.\\
(2) The flow $(\tilde{\psi}_t)_{t\geq 0}$ acts contracting on $E^-_\rho$ and $(\tilde{\psi}_t)_{t\leq 0}$ acts contracting on $E^+_\rho$.
\end{thm*}

\subsection*{Structure of the work}
The first three chapters do not contain any new results 
but provide the reader with the necessary background from  
bounded cohomology (Chapter 1), symmetric spaces (Chapter 2) 
and Hermitian symmetric spaces (Chapter 3).
Starting from Chapter 4, the chapters contain new results. 
Chapter 4 itself contains 
the new characterization of non-tube type  Hermitian symmetric spaces.
It also gives a description of various different models of the Shilov 
boundary, some of which had not been worked out before. 
In Chapter 5 we discuss the different ways of defining a two-cocycle on 
a Hermitian symmetric space. The realization of the cocycle by an 
algebraic function in the general case 
(see \cite{Burger_Iozzi_supq} for $G={\rm SU}(p,q)$) 
is essential for some of the 
applications given in Chapter 7 and 8. 
One of the most important parts of this work is in Chapter 6, where 
tight homomorphisms and tight embeddings are defined and studied.
Chapter 7 contains applications of results from the previous chapters to
isometric actions of arbitrary finitely 
generated groups on Hermitian symmetric spaces.
The characterization of maximal representations of surface groups
is obtained in Chapter 8.
The latter results are joint work with Marc Burger and 
Alessandra Iozzi and where partly announced in  
\cite{Burger_Iozzi_Wienhard_ann}.


\vskip1cm
\chapter{A short introduction to bounded cohomology}\label{chap:bounded}
The concept of bounded cohomology and certain methods developed in its context 
play an important role in the definitions and proofs of the results 
mentioned in 
the introduction. 
We give a short synopsis of bounded cohomology 
including its definition, several properties
and a description of the methods used later on. 
We restrict our attention to cohomology with trivial coefficients. 
Even if we state several properties for arbitrary locally compact groups,  we will mostly work 
in the context
of connected semisimple Lie groups with finite center and their lattices.
For an extensive 
account on bounded cohomology we refer to \cite{Monod_book} and \cite{Burger_Monod_GAFA}, 
for the methods reviewed in Section~\ref{sec:implement} below we refer to \cite{Burger_Iozzi_App}.
 
\section{The definition}
\subsection{The homogeneous standard resolution}
Let $G$ be a locally compact topological group. 
We view abstract groups as topological groups with respect to the discrete topology. 
Denote by $C(G^k)$ (respectively $C_b(G^k)$) the space of continuous
(bounded) real valued  
functions on $G^k$,
\bqn
C(G^k)&:=&\{ f:G^k\to \RR \,|\, f \text{ continuous}\}\\
C_b(G^k)&:=&\{ f:G^k\to \RR \,|\, f \text{ continuous} ,\, \sup_{x\in G^k} |f(x)| <\infty\}.
\eqn

We endow $G^k$ with the diagonal left $G$-action and $C(G^k)$ (respectively $C_b(G^k)$) with the induced 
left action 
\bqn
(hf)(x):= f(h^{-1} x) \quad \forall h\in G, \, \forall x\in G^k.
\eqn

Denote by $C(G^k)^G$ (respectively $C_b(G^k)^G$) the space of $G$-invariant functions in $C(G^k)$ (respectively $C_b(G^k)$), 
\bqn
C(G^k)^G&:=&\{f\in C(G^k)\,|\, (hf)(x)=f(x) \,\forall h\in G, \, \forall x\in G^k\} \\
C_b(G^k)^G&:=&\{f\in C_b(G^k)\,|\, (hf)(x)=f(x) \,\forall h\in G, \, \forall x\in G^k\}.
\eqn

Define $G$-equivariant (coboundary) maps 
\bqn
d: C(G^k) \to C(G^{k+1}), \,\, d: C_b(G^k) \to C_b(G^{k+1})
\eqn
by 
\bqn
(df)(g_0, \dots,g_{k}):= \sum_{i=0}^k (-1)^i f(g_0, \dots,\hat{g_i}, \dots g_{k}).
\eqn
The complexes $(C(G^{\bullet+1}), d)$ and $(C_b(G^{\bullet+1}), d)$ are
called {\em (bounded) homogeneous standard resolution}. 

The maps 
$d$ induce coboundary maps $d: C_{(b)}(G^k)^G \to C_{(b)}(G^{k+1})^G$, 
and give rise to two complexes
\bqn
C(G^\bullet)^G &:& \, 0 \to  C(G)^G \to \dots  \to C(G^k)^G \to C(G^{k+1})^G\to \dots\\
C_b(G^\bullet)^G &:& \, 0 \to C_b(G)^G \to \dots \to C_b(G^k)^G \to C_b(G^{k+1})^G\to \dots.
\eqn

The {\em continuous cohomology} $\hc(G)$ of $G$ (with real
coefficients) 
is defined as the cohomology of the complex $C(G^\bullet)^G$, i.e.  
\bqn
\hc^k(G):= \frac{\ker(d: C(G^{k+1})^G \to C(G^{k+2})^G)}{\im(d: C(G^k)^G \to C(G^{k+1})^G)}.
\eqn

The {\em continuous bounded cohomology} $\hcb(G)$ of $G$ (with real coefficients) is defined as the cohomology of the complex $C_b(G^\bullet)^G$, 
\bqn
\hcb^k(G):=\frac{\ker(d: C_b(G^{k+1})^G \to C_b(G^{k+2})^G)}{\im(d: C_b(G^k)^G \to C_b(G^{k+1})^G)}.
\eqn

The supremum norm on $C_b(G^{k+1})^G$ induces a 
seminorm $||\cdot||$ on $\hcb^k(G)$, 
given by
\bqn
||\a||:= \inf_{f  \in \a} \, \, (\sup_{x\in G^{k+1}} |f(x)|). 
\eqn
Since the image of $d$ is not necessarily closed, $||\cdot||$ is in
general not a norm.

\subsection{The comparison map}
The inclusion $C_b(G^\bullet)^G\subset C(G^\bullet)^G$ as a subcomplex induces natural linear comparison 
maps in cohomology 
\bqn
\kappa: \hcb^*(G) \to \hc^*(G).
\eqn

These maps are in general neither surjective nor injective. Their 
behaviour contains in certain cases information about the group $G$ 
(see e.g. \cite{Mineyev_bounded, Mineyev_straighten}). 
For semisimple Lie groups the comparison map $\kappa$ is an
isomorphism in degree two (see Proposition~\ref{prop:comparison_two}).

\subsection{Group homomorphisms}
Given two locally compact 
topological groups $G,H$ and a continuous group homomorphism $\rho: H\to G$, 
precomposition with $\rho$ induces maps 
$\rho^*:\hc^*(G) \to \hc^*(H)$ 
and 
$\rho_b^*:\hcb^*(G) \to \hcb^*(H)$, 
by $ \rho^*(f)(h_0, \dots,h_{k})= f(\rho(h_0), \dots,\rho(h_{k}))$,
such that the diagram
\bqn
\xymatrix{
\hcb^*(G) \ar[d]^{\kappa^*} \ar[r]^{\rho_b^*} & \hcb^*(H) \ar[d]^{\kappa^*}\\
\hc^*(G) \ar[r]^{\rho^*} & \hc^*(H)
}
\eqn
commutes.

The map $\rho_b^*:\hcb^*(G) \to \hcb^*(H)$ is norm decreasing in each
degree, i.e. for all $\a \in  \hcb^k(G)$, we have
 $||\rho_b^*(\a)||\leq ||\a||$.

Later we will study group homomorphisms $\rho$ satisfying 
\bqn
||\rho_b^*(\a)||= ||\a||
\eqn
for specific $\a\in \hcb^2(G)$. These are the tight homomorphisms which
were discussed in the introduction.
\section{Functorial approach}
We sketch the idea of the functorial approach to 
continuous bounded cohomology developed in \cite{Monod_book} and \cite{Burger_Monod_GAFA}. 
For the corresponding functorial description of the usual continuous cohomology see \cite{Guichardet}. 
An important advantage of the functorial approach is that one obtains 
complexes realizing the cohomology which are easier to compute
 than the homogeneous standard resolution described above.

\subsection{Resolutions by relatively injective $G$-modules}
For the functorial approach one defines a class of so called {\em relatively injective $G$-modules}. 
One considers {\em strong resolutions} of $\RR$ 
by relatively injective $G$-modules, these are resolutions admitting 
special contraction homotopies. 
Finally the cohomology of the (unaugmented) 
complex of $G$-invariants of such resolutions is canonically isomorphic to the cohomology of $G$ with 
real coefficients.
 Relative injectivity is an extension 
property of the modules $E_i$ that guarantees that, 
given two strong resolutions of e.g. $\RR$ by relatively injective $G$-modules $E_i, F_i$, 
the identity map on the coefficients extends to a 
$G$-morphism of the resolutions, 
thus to an isomorphism of the 
cohomology groups defined by the corresponding complexes of $G$-invariants.

\subsection{Bounded cohomology via resolutions}
We review the notions for bounded continuous cohomology, for the
corresponding notions for usual continuous cohomology see \cite{Guichardet}.

\begin{defi}
A Banach $G$-module $E$ is a Banach space $E$ together with an
isometric $G$-action. 
\end{defi}

\begin{defi}
The maximal continuous submodule $\Cc E$ of a Banach $G$-module $E$ is
the set of vectors $e\in E$ such that the map $G\to E$, $g\mapsto ge$
is continuous. 
\end{defi}
%
%
\begin{defi}
A Banach $G$-module $E$ is called 
{\em $G$-relatively injective} (in the category of 
Banach $G$-modules) if it satisfies the
following extension property:
given an injective $G$-morphism $i:A \rightarrow B$ of Banach
$G$-modules $A$ and $B$, a linear norm one operator 
(not necessarily a $G$-morphism) $\sigma: B \rightarrow A$ 
such that 
$\sigma \circ i= id_A$ and $\alpha:A \rightarrow E$ a $G$-morphism, 
there exists a $G$-morphism $\beta: B \rightarrow E$, satisfying
$\beta\circ i=\alpha$, 
$\vert\vert \beta\vert\vert \leq \vert\vert \alpha \vert\vert$.
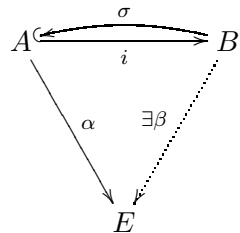
\[
\xymatrix
{
A \ar@{^{(}->}[rr]_i  \ar[ddr]^{\alpha}
& &
B \ar@/_/[ll]_{\sigma} \ar@{.>}[ddl]_{\exists\beta}
\\
\\
& E &
}
\]
\end{defi}

\begin{defi}
A {\em resolution} $E_\bullet$ of $\RR$ is 
an exact complex of Banach $G$-modules $E_i$, with $E_{-1}=\RR$.
A resolution is called a resolution by relatively injective modules, 
if all modules $E_i$, $i\geq 0$ are relatively injective.
\end{defi}
\begin{defi}
The complex $E_\bullet$ is called
{\em strong}, if the complex $\Cc E_\bullet$ of maximal continuous
submodules admits linear operators 
$K_i:\Cc E_i \rightarrow \Cc E_{i-1}$ with 
$\vert\vert K_i \vert\vert \leq 1$, 
$K_{i+1} d_{i} + d_{i-1} K_i =1$ for all $i\geq 0$
and $K_0 d_{-1} =id$.
\end{defi}

\begin{thm}[\cite{Burger_Monod_GAFA,Monod_book}]
Let $E_\bullet$ be a strong resolution by relatively injective Banach
$G$-modules.
Then the cohomology of the corresponding 
unaugmented complex of $G$-invariants
${E_\bullet}^G$ is topologically isomorphic to the continuous bounded cohomology of $G$
 with coefficients in $\RR$, 
i.e. $H^*({E_\bullet}^G)\cong \hcb^*(G)$.
\end{thm}
The isomorphism  $H^*({E_\bullet}^G)\cong \hcb^*(G)$ is not necessarily an isometric isomorphism, 
but the seminorm on $\hcb^*(G)$ defined as the infimum of the
induced seminorms taken over all
strong relatively injective resolutions concides with the seminorm defined by the homogeneous standard resolution above.

\subsection{A preferred resolution for continuous cohomology}
A consequence of the functorial approach to continuous cohomology is, that the continuous cohomology of a 
semisimple connected Lie group $G$ with finite center 
is given by (the cohomology of the complex of) $G$-invariant
differential forms on its associated symmetric space \cite{Guichardet}. 
Let $K$ be a maximal compact subgroup of $G$ and $X=G/K$ the
associated symmetric space. 
Denote by $\Omega^k(X)$ the space of (smooth) differential forms on $X$ and by $d$ the exterior differential. 
Then $(\Omega^\bullet(X),d)$ is a strong relatively injective resolution of $G$-modules. 
The exterior differential restricts to the complex of $G$-invariant differential forms $\Omega^\bullet(X)^G$ as the 
zero mapping. 
Hence there are no coboundaries which have to be taken into account, in particular 
\bqn
\hc^k(G) \cong \frac{\ker(d: \Omega^k(X)^G \to\Omega^{k+1}(X)^G)}{\im(d: \Omega^{k-1}(X)^G \to\Omega^{k}(X)^G)} \cong \Omega^k(X)^G.
\eqn

For a lattice $\G<G$ one obtains
\bqn
\hc^k(\G) \cong \frac{\ker(d: \Omega^k(X)^\G \to\Omega^{k+1}(X)^\G)}{\im(d: \Omega^{k-1}(X)^\G \to\Omega^{k}(X)^\G)} \cong \h_{dR}^k(\G\backslash X).
\eqn

This specific realization of the continuous cohomology allows to make computations 
and shows for example 
that the continuous cohomology of $G$ and $\G$ vanishes in degrees
above the dimension of $X$ what is not clear from the definition using
the complex $C(G^\bullet)^G$.

A (co)chain morphism $F:\Omega^k(X) \to C(G^{k+1})^G$ is given by integration over geodesic simplices
\bqn
F(\omega)(g_0, \dots, g_k) := \int_{\Delta(g_0 K, \dots, g_k K)} \omega, 
\eqn
where $\Delta(g_0 K, \dots, g_k K)$ denotes the geodesic simplex with (ordered) vertices $(g_0 K, \dots, g_k K)$ 
built successively by geodesic cones. 

\subsection{A preferred resolution for bounded continuous cohomology}
Suppose that $G$ is a semisimple Lie group with finite center and that
$\G<G$ is a lattice in $G$.
Let $P< G$ be a minimal parabolic subgroup in $G$. 
The homogeneous space $G/P$, called the Furstenberg boundary of $G$, is compact. Let $\mu$ be the Lebesgue measure on 
$G/P$. 
Denote by $L_{alt}^\infty((G/P)^k)$ the space of essentially bounded
measurable alternating functions on $(G/P)^k$.  
The natural diagonal left action of $G$ on $(G/P)^k$ induces a left
action on $L_{alt}^\infty((G/P)^k)$ which turns it into a Banach
$G$-module which is in fact 
relatively injective since $P$ is amenable.
The restriction of the action of $G$ to $\G$ turns 
$L_{alt}^\infty((G/P)^k)$ into a relatively injective Banach $\G$-module.

Denote by $L_{alt}^\infty((G/P)^k)^G$ the space of essentially invariant functions in $L_{alt}^\infty((G/P)^k)$. 
Define a $G$-morphism $d: L_{alt}^\infty((G/P)^k)\to L_{alt}^\infty((G/P)^{k+1})$  as above by 
\bqn
(df)(g_0 P, \dots,g_{k} P):= \sum_{i=1}^k (-1)^i f(g_0 P, \dots,\widehat{g_i P}, \dots g_{k} P).
\eqn

The complex $(L_{alt}^\infty((G/P)^\bullet, d)$ is a strong $G$-resolution of relatively injective 
Banach $G$-modules, and thus 
\bqn
\hcb^k(G) \cong \frac{\ker (d: L_{alt}^\infty((G/P)^{k+1})^G\to L_{alt}^\infty((G/P)^{k+2})^G)}{\im(d: L_{alt}^\infty((G/P)^k)^G\to L_{alt}^\infty((G/P)^{k+1})^G)}. 
\eqn

Similarly 
\bqn
\hcb^k(\G) \cong \frac{\ker (d: L_{alt}^\infty((G/P)^{k+1})^\G\to L_{alt}^\infty((G/P)^{k+2})^\G)}{\im(d: L_{alt}^\infty((G/P)^k)^\G\to L_{alt}^\infty((G/P)^{k+1})^\G)}. 
\eqn

The seminorm induced from the norm on $L_{alt}^\infty((G/P)^{k+1})^G$ realizes the seminorm on $\hcb^k(G)$.

This resolution is a preferred resolution, because of its geometric
origin and because it is ``very small'' in the following sense:  
\begin{prop}\label{prop:L_alt}\cite{Burger_Monod_GAFA}
The actions of $G$ and $\G$ are amenable on $G/P$ and ergodic on
$G/P\times G/P$, in particular $L_{alt}^\infty((G/P)^2)^G =
L_{alt}^\infty((G/P)^2)^\G= 0$ and there are isometric isomorphisms 
\bqn
\hcb^2(G) &\cong& \ker (d: L_{alt}^\infty((G/P)^3)^G\to L_{alt}^\infty((G/P)^{4})^G)\\
\hcb^2(\G) &\cong& \ker (d: L_{alt}^\infty((G/P)^3)^\G\to L_{alt}^\infty((G/P)^{4})^\G).
\eqn
\end{prop} 

Recall the notion of amenability of a $G$-action (\cite{Zimmer_book}): 
\begin{defi}
Let $S$ be a standard Borel space equipped with a probability
measure. Assume that $S$ is endowed with a measure class preserving
Borel $G$-action.
The $G$-action on $S$ is said to be {\em amenable} if for every separable
Banach space $E$ and every Borel right cocycle $\a: S\times G\to {\rm
  Iso}(E)$  the dual $\a^*$-twisted action on $E^*$ satisfies that any
$\a^*$-invariant Borel field $\{A_s\}_{s\in S} \subset E_1^*$ of
nonempty convex weak-* compact subsets $A_s$ admits an
$\a^*$-invariant section. 
\end{defi}

When $S$ is a regular $G$-space the action of $G$ being amenable
relates to certain modules being relatively injective.
\begin{prop}\cite[Theroem~5.7.1.]{Monod_book}
Let $S$ be a regular $G$-space. The following are equivalent:\\
\begin{enumerate}
\item{ The $G$-action on $S$ is amenable.}
\item{ The $G$-module $L^\infty(S^{n+1})$ is relatively
    injective for all $n\geq 0$.}
\end{enumerate}
\end{prop}

\begin{defi}
The $G$-action on $S$ is said to be {\em doubly ergodic} if $G$ acts 
ergodically on $S\times S$.
\end{defi}
This property is usually called mixing, but we adapt the notion of
doubly ergodicity which applies to more general $G$-modules.

For a finitely generated group $\Lambda$
let $(S,\nu)$ be a Poisson boundary for $\Lambda$ such that $\Lambda$ acts 
amenably and doubly ergodically on $(S,\nu)$. 
(For the construction see for example 
\cite{Burger_Monod_GAFA, Burger_Iozzi_supq}.)
The action of $\G$ on $S$ is amenable and doubly ergodic, hence 
the modules $L_{alt}^\infty((S)^{n+1})$ are relatively injective and 
$L_{alt}^\infty((S)^2)^\Lambda= \{0\}$, thus
 \bqn
\hcb^2(\Lambda) &=& \ker (d: L_{alt}^\infty((S)^3)^\Lambda\to L_{alt}^\infty((S)^{4})^\Lambda).
\eqn

\begin{rem}
The bounded cohomology of a surface group in degree $3$ is not a
Banach space \cite{Soma_third}, 
therefore we cannot hope for a generalization of this result to higher degrees 
for discrete groups. 
Nevertheless, as far as the author knows, there is no example of a
continuous bounded cohomology 
space of a semisimple connected Lie group with finite center not being a Banach space.  
\end{rem}

\section{Implementing pull-back maps}\label{sec:implement}
Any (continuous) group homorphism $\rho:\Lambda\to G$ 
induces a natural pull-back map in cohomology as described above with
the homogeneous bar resolution.  
Let $H, G$ be connected semisimple Lie groups with finite center and $\Lambda$ a finitely generated group. 
Suppose that $\rho:\Lambda\to G$ is a homomorphism. We want to find an implementation of the pull-back maps 
$\rho^*$ respectively $\rho_b^*$ 
in our preferred resolutions $(\Omega^\bullet(X),d)$ respectively $(L_{alt}^\infty((G/P)^\bullet, d)$

\subsection{The pull-back in continuous cohomology}
Denote by $Y, X$ the symmetric spaces of noncompact type associated to $H$ and $G$. Assume that 
$\Lambda$ is a lattice in $H$ and denote by $M= \Lambda\backslash Y$ the quotient by $\Lambda$.
The homomorphism induces a flat bundle over $M$
\bqn
X_\rho:= X\times_\rho Y,  
\eqn
where $X\times_\rho Y$ is obtained by dividing $X\times Y$ 
by the equivalence relation  $(x,y) \sim (\rho(h)x, h y)$, $h\in
\Lambda$. 
Since the bundle is flat and the space $X$ is contractible, 
we obtain a (smooth) section, hence a $\rho$-equivariant (smooth) map 
$f: Y\to X$. 

Given a $G$-invariant differential form $\omega\in \Omega^k(X)^G$, 
the pull-back $f^*\omega$ is a well defined (smooth) closed differential form 
 $f^*\omega \in \Omega^k(Y)$, which, by equivariance of $f$, is $\Lambda$-invariant, 
hence $\rho^*([\omega])=[f^*\omega] \in \hc^k(\Lambda)$.  

\subsection{The pull-back in bounded continuous cohomology}
For continuous bounded cohomology, we would like to implement the 
pull-back 
via a $\rho$-equivariant measurable boundary map 
$\phi: S \to \Mm^1(G/P)$ of a doubly 
ergodic Poisson boundary $(S,\nu)$ of $\Lambda$ into the space of probability
measures on the Furstenberg boundary $G/P$ of $G$  
which always exists in the above situation, since the action of
$\Lambda$ on $S$ is amenable and $G/P$ is a compact metrizable space.
\begin{prop}\cite[Theorem~4.3.9.]{Zimmer_book}
Let $\rho:\Lambda \to G$ be a representation. 
Assume that the action of $\Lambda$ on $S$ is amenable. 
There exists a measurable almost everywhere $\rho$-equivariant 
map $\psi:S \to \Mm^1(G/P)$.
\end{prop}
\begin{proof}
We consider the special situation 
where $E=C(G/P)$ is the space of continuous functions on $G/P$. Then 
the space of probability measure 
\bqn
\Mm^1(G/P)\subset E_1^*
\eqn
is compact
and convex. Consider the field $A_s= \Mm^1(G/P)$ for all $s\in S$ 
and apply the 
definition of amenability of a $G$-action. 
Then $\Lambda$ acts on the affine $\Lambda$-space $F(S,\Mm^1(G/P))$ of measurable maps
from $S$ to $\Mm^1(G/P)$ by $(g\phi)(s) = \rho(g) \phi(g^{-1} s)$. 
Since $\Lambda$ acts amenably on $S$, 
there is a fixed point $\phi$ in $F(S,\Mm^1(G/P))$, i.e. 
a measurable map $\phi: S\to \Mm^1(G/P)$ such that 
$\rho(g) \phi(g^{-1} s)=(g\phi)(s)=\phi(s)$.
\end{proof}
Even when the boundary map takes values in the space of Dirac measures 
on $G/P$, which we identify with $G/P$, 
we encounter the problem that the pushforward measure $\phi_*(\nu)$ on $G/P$ 
is in general not absolutely continuous with respect to the Lebesgue measure 
on $G/P$. Thus the precomposition 
with $\phi$ does not give a map 
$\phi^*:L_{alt}^\infty((G/P)^k)^G \to
L_{alt}^\infty((S)^k)^\Lambda$.

Fortunately in the geometric context, 
the specific representatives of bounded cohomology classes are generally  
not only defined almost everywhere 
essentially bounded measurable functions on the Furstenberg boundary, 
but they are {\em strict} (i.e. everywhere defined) Borel cocycles. 
Using the functorial approach it is shown in \cite{Burger_Iozzi_App}, that the pull--back of these cocycles 
can be implemented via the $\rho$-equivariant measurable boundary map $\phi: S \to \Mm^1(G/P)$ of the boundaries. 
More precisely, let $Z$ be a compact metrizable space with a continuous $G$ action. 
Denote by  $\Bb^\infty(Z^n)$ the Banach space of bounded Borel
measurable functions defined everywhere on $Z^n$, i.e.
\bqn
\Bb^\infty(Z^n):=\{f: Z^n\to \RR \,|\, \text{f Borel measurable}, \, \sup_{x\in Z^n}|f(x)| <\infty\}
\eqn 
and by $\Bb^\infty_{alt}(Z^n)$ the subspace of alternating functions in $\Bb_\infty(Z^n)$.
If we endow $\Bb^\infty_{alt}(Z^n)$  with the induced $G$-action and with the homogeneous boundary operator, 
we obtain a strong resolution of Banach $G$-modules, but the modules
are not known to be relatively injective. 
Nevertheless we have the following results from \cite{Burger_Iozzi_App}: 

\begin{prop}\label{prop:borel}
There is a canonical map 
\bqn
\h^*(\Bb^\infty_{alt}(Z^{*+1})^G) \to \hcb^*(G).
\eqn
In particular, every alternating 
Borel measurable $G$-invariant cocycle 
\bqn
c: Z^{n+1}\to \RR
\eqn
canonically
determines 
a class $[c]\in \hcb^{n}(G)$.
\end{prop}

\begin{prop}
Suppose that $\phi:S \to \Mm^1(Z)$ is a measurable almost everywhere 
$\rho$-equivariant map. 
Let $c: Z^{n+1}\to \RR$ be an alternating Borel measurable $G$-invariant cocycle and $[c]\in \hcb^n(G)$ 
its associated bounded cohomology class. 
Then the pull-back $\phi^*(c)$ defined by 
\bqn
\phi^*(c) (y_0, \cdots, y_{n}) := \phi(y_0) \otimes \dots \otimes \phi(y_{n}) (c)
\eqn
defines a cocycle in $L_{alt}^\infty((S)^{n+1})^\Lambda$ representing 
the bounded cohomology class $\rho^*([c]) \in \hcb^n(\Lambda)$.
\end{prop}
\begin{rem}
We described the method to implement pull-back maps for homomorphisms
of a discrete group $\Lambda$ into $G$. 
The same methods allow to pull-back cocycles for homomorphisms $\rho:
H\to G$, where $H$ is for example a connected semisimple Lie group
with finite center. 
\end{rem}
\section{Properties}\label{sec:hcb_prop}
We collect some properties of continuous bounded cohomology which will be used later on. 
\subsection{Subgroups and quotients}
\begin{prop}\cite[Corollary~8.8.5 and 7.5.10]{Monod_book}\label{prop:subgroups}
Let $G$ be a locally compact second countable topological group.
\begin{enumerate}
\item{Suppose $G_0\lhd G$ is a normal subgroup of finite index.\\
        Then $\hcb^*(G)\cong\hcb^*(G_0)^{G/G_0}$ is an isometric isomorphism.}
\item{Suppose $R \lhd G$ is a normal closed amenable subgroup. 
        Then there is a isometric isomorphism $\hcb^*(G)\cong \hcb^*(G/R)$.} 
\end{enumerate}
\end{prop}

\subsection{Degree $2$}
In degree two one has more insight in the bounded cohomology than in other degrees.
As a direct consequence of Proposition~\ref{prop:L_alt}, 
the continuous bounded 
cohomology of $G$ and $\G$ in degree two are Banach spaces. 
This is special for degree two.  
Furthermore there is a K\"unneth type formula in degree two 
(see 
\cite[Corollary 4.4.1.]{Burger_Monod_GAFA}). 
\begin{prop}\cite[Corollary~12.0.4]{Monod_book}\label{prop:degree_two}

Let $G=G_1\times \cdots\times G_k$ be a locally compact second countable topological group
Then 
\bqn
\hcb^2(G,\RR) &\to& \bigoplus_{i=1}^k \hcb^2(G_i, \RR),\\
 \a &\mapsto& \oplus \a|_{G_i}
\eqn
is an isometric isomorphism with respect to the induced $l^1$-norm, 
i.e. $||\a||= \sum_{i=1}^k || \a|_{G_i}||$.
\end{prop}

For semisimple Lie groups with finite center, the comparison map $\kappa$ 
in degree $2$ is completely understood. 

\begin{prop}\cite{Guichardet_Wigner,Dupont_Guichardet}\cite[Lemma~6.1]{Burger_Monod_JEMS}\label{prop:comparison_two}
Let $G$ be a connected semisimple Lie group with finite center, then 
\bqn
\kappa: \hcb^2(G)\to \hc^2(G)  
\eqn
is an isomorphism. 
\end{prop}


\vskip1cm
\chapter{Symmetric spaces}
We recall some basic facts on symmetric spaces. 
For a thorough account of symmetric spaces 
including proofs of the properties given below we refer the reader to the 
books of Helgason \cite{Helgason} and Warner \cite{GWarner}. 
A treatment of symmetric spaces of noncompact type can be found in Eberlein's book 
\cite{Eberlein}.

\section{Definitions}
\begin{defi}
A simply connected Riemmannian manifold $X$ is called a {\em symmetric space} 
if every point $p\in X$ is an isolated fixed point of an involutive isometry 
$s_p: X\to X$.

A symmetric space is said to be of {\em noncompact type} if the 
Ricci curvature of $X$ is negative.
\end{defi}
The involutive isometry $s_p$ is the geodesic symmetry that is defined for any Riemannian manifold $M$ as follows. 
Let $U_0\subset T_p M$ be a (small enough) symmetric neighbourhood of $0 \in T_p M$, 
let $N_p = \exp_p(U_0)\subset M$ be the image of $U_0$ under the Riemannian exponential mapping. 
Then $s_p: N_p \to N_p, \, \, q\mapsto \exp_p(-\exp_p^{-1}(q))$ is called {\em geodesic symmetry}.

The connected component of the group of isometries of $X$ containing the identity, $G={\rm Is}^\circ(X)$ acts transitively on $X$.
Thus choosing a base point $x_0 \in X$, we may realize $X$ as a homogeneous space
$X= G/K$, whith $K=\stab_G(x_0)$.
If furthermore $X$ is of noncompact type, then $G$ is a connected 
semisimple Lie group 
without compact factors and without center. 
On the other hand, if $G$ is a connected 
semisimple Lie group with finite center, 
then we can associate to $G$ a symmetric space of noncompact type: 
Namely, if  $K<G$ is a maximal compact subgroup, 
then $X= G/K$ is a symmetric space of noncompact type. 
If $G$ is compact, then $X$ is a point.

\begin{defi}
Let $G$ be a reductive Lie group. 
The symmetric space associated to $G$ is the symmetric space associated to its semisimple part.\\
A symmetric space $X$ is called {\em irreducible} if $X$ does not split 
as a product $ X= X_1\times X_2$. 
\end{defi}
If $X$ is of noncompact type the irreducibilty of $X$ 
is equivalent to $G={\rm Is}(X)^\circ$ being a simple group.
\subsection{Cartan decomposition}
The involutive isometry $s_p$ at a point $p\in X$ defines 
an involution $\sigma_p$ of $G$, $\sigma_p(g) := s_p\circ g\circ s_p$. 
Under the identification $\frakg\cong T_e G$ the differential $(\sigma_{x_0})_*$ defines the {\em Cartan involution} of the Lie algebra $\frakg$ of $G$ 
\bqn
\theta = (\sigma_{x_0})_*: \frakg\to \frakg.
\eqn
The Lie algebra $\frakg$ decomposes into the $(\pm1)$-eigenspaces of $\theta$, 
$\frakg= \frakk\oplus \frakp$ with $\theta|_{\frakk} = \id_\frakk$, $\theta|_\frakp = -\id_\frakp$, 
satisfying the relations 
$[\frakk, \frakk]\subset \frakk$, 
$[\frakk, \frakp]\subset \frakp$, 
$[\frakp, \frakp]\subset \frakk$.
The $(+1)$-eigenspace $\frakk$ is a Lie subalgebra which is 
identified with the 
Lie algebra of $K$. 
The {\em Cartan decomposition } of $\frakg=\frakk\oplus \frakp$ is orthogonal with respect to 
the Killing form 
\bqn
\frakB: \frakg\times\frakg &\to& \RR\\
(X,Y) &\mapsto& \tr (\ad(X) \circ \ad(Y))
\eqn

The Killing form defines a nondegenerate $\Ad(K)$-invariant symmetric
bilinear form on $\frakg$, 
by 
\bqn
\langle X,Y\rangle := -\frakB(X, \theta(Y)), 
\eqn
whose restriction to $\frakp$ is positive definite.

The differential $\pi_*:\frakg\to T_{x_0} X$ of the 
orbit map $\pi: G\to X, \, g\mapsto g x_0$ has kernel $\frakk$ and 
thus induces an identification $\frakp\cong T_{x_0} X$. 
The restriction of the Killing form to $\frakp$ defines  
an inner product on $T_{x_0} X$. 
Using the transitivity of $G$, we may extend $\langle\,.,.\,\rangle$ to 
a $G$-invariant metric on $X$, which, up to scaling on the factors, coincides 
with the given Riemannian structure on $X$.

\subsection{Maximal abelian subalgebras}
From now on $X$ will be a symmetric space of noncompact type endowed with the metric $q_0= \frac{1}{2} \langle\,.\,,\,.\,\rangle$.

\begin{defi}
The rank $r_X$ of a symmetric space $X$ is 
the maximal dimension of an 
isometrically and totally geodesically embedded Euclidean subspace in $X$, 
\bqn
r_X:= \max \{k \, |\,\,  \RR^k \to X \,\text{isometric and totally geodesic embedding}\}.
\eqn
A $k$-flat, for $k\leq r_X$, in $X$ is the image of an isometric totally geodesic embedding $\RR^k\to X$. 
A maximal flat or a flat is an $r_X$-flat. 
\end{defi}

The rank of a symmetric space on noncompact type is at least $1$, since 
geodesics are isometric embeddings of $\RR$. 
A symmetric space of higher rank contains subspaces on which the curvature is zero. 

Let $\fraka\subset\frakp$ be a maximal abelian subspace and
$A=\exp(\fraka)$.
Then $F_0=Ax_0$ is a maximal flat in $X$.
In particular, $r_X = \max \{\dim(\fraka) \,|\, \fraka\subset\frakp,\, [\fraka, \fraka]=0\}$. 
All maximal flats in $X$ are conjugates by 
 $K$ of $F_0$.

\begin{defi}
A geodesic $\g$ in $X$ is said to be {\em regular} if it is contained in exactly one maximal flat, 
otherwise it is called {\em singular}. 
A vector $X\in\frakp$ is called regular (resp. singular) if the geodesic $g(t) = \exp(tX) x_0$ is regular (resp. singular).
\end{defi}
\begin{defi}
Let $\fraka_{sing}$ be the set of singular vectors in $\fraka$.
\end{defi}

\begin{defi}
A positive Weyl chamber is the choice of a connected component $\fraka^+$ of $\fraka\setminus\fraka_{sing}$. 
\end{defi}
The choice of a maximal compact subgroup $K<G$ and a positive Weyl chamber determine 
the {\em Cartan decomposition} $G= K \exp(\ol{ \fraka^+})K$. This 
decomposition $g=kak'$ is in general not unique, but the element $a\in
exp(\ol{\fraka^+})$ is uniquely determined. 

The Weyl group $W$ of $(\frakg, \fraka)$ is the group of reflections generated by the hyperplanes bounding $\fraka^+$.
The choice of a positive Weyl chamber is related to the choice of a fundamental system of roots.
\begin{defi}
A real {\em root} of $(\frakg, \fraka)$ is a nontrivial linear form $\a\in \fraka^*$ such that 
the corresponding {\em root space} 
\bqn
\frakg_\a:=\{X\in \frakg \, |\, \ad(H)(X)=\a(H)X \}
\eqn 
is nonzero.
\end{defi}

We denote the set of all roots of $(\frakg, \fraka)$ by $\Xi(\frakg, \fraka)$. 
Since for all $H\in \fraka$ the operators 
$\ad(H)$ are commuting and self-adjoint operators, we obtain 
a decomposition into the eigenspaces, i.e. 
into root spaces, $ \frakg = \frakg_0 \oplus_{\a\in \Xi} \frakg_\a$.

The choice of a basis of $\fraka$ induces a lexicographical order on $\Xi(\frakg, \fraka)$ and 
a decomposition of $\Xi=\Xi^+ \cup  \Xi^-$ into the sets of positive roots $\Xi^+$ and negative roots $\Xi^-$.
A positive root $\a$ is called {\em  simple}, if it cannot be written as sum of two positive roots. 
A {\em fundamental system} $\Delta\subset \Xi$ is a set $\Delta= \{\g_1 , \dots , \g_r\}$ of 
positive simple roots such that any root $\a$ can be written as 
$\a = \sum_{i=1}^r n_i \g_i$, where either $n_i\geq 0$ for all $i$ or $n_i\leq0$ for all $i$.

The choice of a fundamental system $\Delta$ is equivalent 
to the choice of a positive Weyl chamber in $\fraka$ 
defined by
\bqn
{\fraka}^+=\{X\in \fraka \,|\,\g_i(X)> 0, \, \forall \g_i \in \Delta\}
\eqn

The Weyl group, generated by elements $s_i$, sending $\g_i$ to $-\g_i$ 
and permuting the other simple positive roots $\g_j$, acts on 
$\Xi (\frakg,\fraka)$.
It contains a unique element $w_{max}$ such that 
$w_{max}(\Xi^+) = \Xi^-$.

The sets of positive and negative roots $\Xi^\pm$ 
determine nilpotent subalgebras $\frakn^\pm= \oplus_{\a\in \Xi^\pm} \frakg_\a$ 
and corresponding nilpotent subgroups $N^\pm = \exp(\frakn^\pm)<G$. 
\begin{prop}
The multiplication map 
\bqn
K\times A\times N^\pm &\to& G, \\
(k,a,n) & \mapsto& kan
\eqn
is a diffeomorphism onto $G$, called the {\em Iwasawa decomposition} of $G$.
\end{prop}
\subsection{Parabolic subgroups}
Let 
\bqn
\fraka^+:=\{X\in \fraka \,|\,\g_i(X)>0, \, \forall \g_i \in \Delta\}
\eqn
be the positive Weyl chamber associated to $\Delta$.

For a set $\theta \subset \Delta$, we define 
\bqn
\fraka_\theta^+:=\{X\in \fraka \,|\,\g_i(X)>0 \,\, \forall \g_i \in
\Delta-\theta \text{ and }\g_i(X)=0 \,\, \forall \g_i \in \theta\}.
\eqn

\begin{defi}
Vectors in $\fraka^+_{\theta_i}$ with $\theta_i=\Delta-\{\gamma_i\}$
are called maximal singular. 
\end{defi}

Let $M_\theta$ denote the centralizer of $\fraka_\theta$ in $K$ and 
$\frakm_\theta$ its Lie algebra. 
To any subset $\theta \subset \Delta$ there is an associated system 
of standard parabolic subgroups $P_\theta= M_\theta A N^+$,
where $P=P_{\emptyset}= M_\emptyset AN^+$ is the minimal parabolic subgroup.
Every parabolic subgroup of $G$ is conjugate via $K$ to 
a standard parabolic subgroup $P_\theta$ for some $\theta \subset \Delta$.
The following inclusion relations hold: $\fraka_{\theta\prime}^+ \subset \overline{\fraka_\theta^+}$ 
if and only if $P_\theta\subset P_{\theta\prime}$.

\section{Boundaries}
There a several nonisomorphic compactifications of a symmetric space $X$ 
of noncompact type. We describe two compactifications 
that are of importance in the later chapters.
For a detailed study of the different compactifications 
see 
\cite{Guivarch_Ji_Taylor} and also the unified approach in \cite{Borel_Ji}. 

\subsection{Geometric boundary and conical compactification}
Define an equivalence relation on the set of all unit speed geodesics in $X$ 
by calling $\g,\eta$ asymptotic 
if $\lim_{t\to +\infty} \d(\g(t),\eta(t)) < \infty$.
The geometric boundary $X(\infty)$ 
is the set of all 
equivalence classes of asymptotic geodesics in $X$.
The \emph{conical compactification} of $X$ 
is defined as 
\bqn
\overline{X}^c= X\cup X(\infty).
\eqn
A geodesic ray is said to be of type $Q$ if $\Stab_G(\g(\infty))=Q$,
where the parabolic subgroup $Q<G$ is minimal
if and only if $\g$ is a regular geodesic ray, maximal 
if and only if $\g$ is a maximal singular geodesic ray.
Since every parabolic subgroup is conjugate by an element of $K$ 
to a standard parabolic subgroup $P_\theta$ determined by a proper 
subset $\theta\subset \Delta$, we 
can define $X_\theta(\infty)$ to be 
the set of equivalence classes of asymptotic rays, where 
$Q$ is conjugate to $P_\theta$. 
A geodesic ray is said to be of type $\theta$ if $\Stab_G(\g(\infty))$ 
is conjugate to the standard parabolic subgroup $P_\theta$. 
In particular we have the following decomposition of the geometric boundary
\bqn
X(\infty)= 
 \amalg_{\theta\subset\Delta}X_\theta(\infty)
\eqn
Representing every equivalence class by a geodesic going through $x_0$ and
associating to every geodesic through $x_0$ its tangent vector at this point 
defines a bijection 
\bqn
X(\infty)\to T^1_{x_0}X \cong \Ad(K) \ol{\fraka^+}
\eqn
which  respects the stratification $X_\theta(\infty) \cong \fraka^+_\theta$ and 
identifies $X(\infty)$ with the unit tangent sphere at $x_0\in X$. 
We endow $X(\infty)$ with the topology of the latter.
The action of $G$ on $X$ extends to a continuous action on $X(\infty)$. 
The boundary strata $X_\theta$ are $G$-invariant and satisfy the 
following closure relations  
$\overline{X_\theta(\infty)}=\amalg_{\theta\subset \theta^\prime} 
X_{\theta^\prime}(\infty)$.

Note that for $\theta_i = \Delta-\{\g_i\}$, 
$X_{\theta_i}(\infty)$ is the $G$-orbit of one point and can be 
identified with $G/P_{\theta_i}$. Since $G=KP_{\theta_i}$, 
$X_{\theta_i}(\infty)$ is also 
the $K$-orbit of a point. 

We define 
\bqn
\Theta_r&:=&\{\theta \subset \Delta|\, \g_r\notin \theta\} \\
\Theta_i&:=&\{\theta \subset \Delta |\, \g_i\notin\theta \, , \theta\notin \Theta_j\, \text{ for } j>i\}.
\eqn
Then $\theta_i:=\Delta-\{\g_i\} \in \Theta_i$.
We denote by $X_{\Theta_i}(\infty)$ the equivalence classes of geodesic rays of type $\theta$ with 
$\theta\in \Theta_i$.

\subsection{The Furstenberg compactifications}
The {\em Furstenberg compactifications} are defined by embedding $X= G/K$ 
into a set of probability measures on a compact space. 
The most popular such compactification is the maximal Furstenberg
compactification 
which embeds $X$ into the space of probability measures $\Mm^1(G/P)$ on the compact homogeneous 
space defined by a minimal parabolic subgroup $P<G$. The closure of $X$ in $\Mm^1(G/P)$ contains the 
set of Dirac measures on $G/P$ which we identify with $G/P$. This part $G/P$ 
of the boundary of the 
maximal Furstenberg compactification is called {\em Furstenberg boundary}.

Furstenberg compactifications can be defined for any proper subset $\theta\subset \Delta$ 
of the set of simple roots.
Given $\theta\subset \Delta$, let $P_\theta$  be the associated standard parabolic subgroup of $G$.
The space  $G/P_\theta$ is a compact homogeneous space on which $K$ acts transitively. 
There exists a unique 
$K$-invariant probability measure $\nu_\theta$ on $G/P_\theta\cong K/M_\theta$. 
The symmetric space $X=G/K$ can be embedded into the space of probability measures $\Mm^1(G/P_\theta)$ 
on $G/P_\theta$ via 
\bqn
F:X &\to& \Mm^1(G/P_\theta), \\
 gK &\mapsto& g_* \nu_\theta.
\eqn
The {\em $\theta$-Furstenberg compactification} of $X$ is the closure of $F(X)$ in  $\Mm^1(G/P_\theta)$.
The boundary of the  {\em $\theta$-Furstenberg compactification} of $X$ contains as 
a distinguished closed subspace the space of Dirac measures on $G/P_\theta$ which we identify with $G/P_\theta$.
This subspace is called the {\em $\theta$-Furstenberg boundary} of $X$.

\begin{rem}
The maximal Furstenberg boundary plays an important role in continuous
bounded cohomology. A parabolic subgroup in $G$ is amenable if and
only if it is a minimal parabolic subgroup. In particular, $G$ acts amenably on $G/P$ and 
thus the resolution by essentially bounded measurable functions on Cartesian products of $G/P$ give a strong 
relatively injective resolution. The action of $G$ on the smaller
Furstenberg boundaries is not amenable, 
but several specific geometric representatives of bounded cocycles
(e.g. the Bergmann cocycle, to be discussed later,  
or an explicit representative of the Euler class) are defined on smaller Furstenberg boundaries. 
The $\theta$-Furstenberg compactification is isomorphic to the
$\theta$-Satake compactification 
defined via irreducible representations of $G$.
The proof of Savage \cite{Savage} 
that the invariant volume form of the 
symmetric space $X$ associated to ${\rm SL}(n, \RR)$ is bounded 
uses a non-maximal Satake compactification of $X$.
\end{rem}

\section{Totally geodesic submanifolds}\label{sec:totally_geodesic}
In contrast to other Riemannian manifolds, symmetric spaces have a lot of totally geodesic submanifolds, 
which can be described in terms of Lie algebras. 
\begin{defi}
Let $M$ be a Riemannian manifold and $S\subset M$ be a connected 
submanifold. 
Then $S$ is said to be {\em totally geodesic} if any geodesic in $S$ is a geodesic in $M$.
\end{defi}

\begin{defi}
Let $\frakg$ be a real Lie algebra. A subspace $\fraks\subset\frakg$ is called a 
{\em Lie triple system} if $[[\fraks, \fraks], \fraks]\subset \fraks$.
\end{defi}


\begin{prop}
There is a one-to-one correspondence between totally geodesic
submanifolds $S$ of $X$ through $x_0$  
and Lie triple systems $\fraks \subset \frakp\subset \frakg$. 
More precisely, 
if $\fraks\subset\frakp$ is a Lie triple system in $\frakg$ then 
the image of $\fraks$ under the Riemannian exponential mapping, 
$S=\exp(\fraks) x_0$, is a totally geodesic submanifold of $X$ through $x_0$. 
If $S$ is a totally geodesic submanifold of $X$ through $x_0$ 
then $\fraks = T_{x_0}S$ is a 
Lie triple system in $\frakg$.  
\end{prop}

In particular, since for simple Lie algebras $[\frakp, \frakp]=\frakk$, 
we obtain a one-to-one correspondence between totally geodesic embeddings 
$f:Y\to X$ of symmetric spaces $Y=H/K_H$, $X= G/K$ and Lie algebra 
homomorphisms $\pi: \frakh\to \frakg$ satisfying 
$\pi(\frakk_H) \subset\frakk, \pi(\frakp_H)\subset \frakp$.

Suppose that $S$ is the image of a totally geodesic embedding 
$f:Y \to X$. 
We would like to associate a homomorphism $\pi : H\to G$ of the isometry groups of $Y$ and $X$.
This is not always possible for the connected components 
of the identity of the isometry group of $Y$, 
but it is possible for some finite extension.
Namely, given any point $p\in S$, the isometry of $X$ 
defined by the geodesic symmetry $s_p$ in $p$ is an isometry of 
$S$, hence of $Y$, but it need not necessarily lie in the connected component 
of the identity of ${\rm Is}(Y)$ even if it lies in ${\rm Is}(X)^\circ$.
Therefore to a totally geodesic embedding $f:Y\to X$ we cannot in
general 
associate a homomorphism $\pi : H\to G$ of the connected 
components of the isometry groups. But there always exists a 
finite extension $H_Y$ of $H$ such that there is a 
homomorphism  $\pi: H_Y\to G$ 
with respect to which $f$ is an equivariant totally geodesic embedding.
We call  $\pi: H_Y\to G$ the {\em corresponding} Lie group homomorphism of $f$.

\section{Real structure and duality}
\subsection{$G$ as real algebraic group}
The Lie algebra $\frakg$ of the connected component of the isometry
group $G$ 
of a symmetric space of noncompact type is 
a real semisimple Lie algebra. 
We want an algebraic structure on $G$. 
Let $\frakg_\CC$ be the 
complexification of $\frakg$, 
then $G_\CC= \aut(\frakg_\CC)^\circ$ is a connected linear algebraic group defined over $\RR$. 
There is a real structure on $G_\CC$ with respect to which
$G=G_\CC(\RR)^\circ$. 
This algebraic structure allows us to work with the Zariski topology on $G$.
\subsection{Duality}
While there is only one compact real form of the simple complex group $G_\CC$, 
there are different noncompact real forms of the same simple complex group $G_\CC$.
For example, all groups ${\rm SU}(p,q)$ with $p+q=n$ are real forms of ${\rm SL}(n,\CC)$, 
whereas the only compact real form is ${\rm SU}(n)$.
Nevertheless, on the level of the associated symmetric spaces, i.e. on the level of the 
Lie algebra together with a Cartan involution, there exists a duality between 
symmetric spaces of compact and noncompact type. 

Let $\frakg$ be a simple real Lie algebra and $\frakg=\frakk\oplus \frakp \subset \frakg_\CC$ 
a Cartan decomposition. 
Then $\frakg^*:= \frakk\oplus i\frakp$ is a simple real Lie algebra. 
The Lie algebra $\frakg$ is noncompact if and only if $\frakg^*$ is compact and vice versa.
Let $G,G^*$ be the corresponding adjoint groups.
If the associated symmetric space $X= G/K$ is of noncompact type, 
the dual symmetric space $X^c= G^*/K$ is of compact type. 
The space $X^c= G^*/K$ is called the {\em compact dual} of $X= G/K$.



\vskip1cm
\chapter{Hermitian symmetric spaces}\label{sec:hermit}

We review some facts about Hermitian symmetric spaces which will be used later. 
We do not give proofs for the classical results. The reader is refered to 
the well written article of A. Kor\'anyi \cite{Koranyi_etal}, 
the article of J.Wolf \cite{Wolf} and the book of I. Satake \cite{Satake_book}.
The results on algebraic groups can be found in A. Borel's book \cite{Borel_lag}.

We use the notations from the previous chapter.
\section{Definitions}
\begin{defi}
A simply connected complex manifold $X$ 
with a Hermitian structure is a {\em Hermitian symmetric space} if 
every $p\in X$ is an isolated fixed point of an involutive holomorphic isometry of $X$.
\end{defi}

\begin{prop}\cite{Guichardet_Wigner}
An irreducible symmetric space $X$ is Hermitian if it satisfies one of the following equivalent properties 
\begin{enumerate} 
\item{$X$ admits a ${\rm Is}(X)^\circ$-invariant  K\"ahler form.}
\item{$X$ carries a ${\rm Is}(X)^\circ$-invariant complex structure.}
\item{The center $\frakc(\frakk)$ of $\frakk$ is nontrivial.}
\item{$\dim (\hc^2(G))\neq 0$.}
\end{enumerate}
\end{prop}
A semisimple Lie group $G$ without compact factors and with finite center is said to be of {\em Hermitian type} if 
its associated symmetric space $X$ is Hermitian symmetric. Fixing a
$G$-invariant complex structure on $X$ 
there exists an unique 
element $Z_\frakg\in\frakc(\frakk)$ in the center of $\frakk$ 
such that $\ad(Z_\frakg)$ is the complex structure of $X$ in $T_{x_0} X\cong \frakp$.
The complexification $\frakg_\CC=\frakg\otimes \CC$ splits under $\ad(Z_\frakg)$ as 
$\frakg_\CC= \frakk_\CC \oplus \frakpp \oplus\frakpm$, where 
$\frakp_\pm$, the 
$(\pm i)$ eigenspaces of $\frakp_\CC$ with respect to  $\ad(Z_\frakg)$, 
are maximal abelian subalgebras of $\frakg_\CC$; 
they can be identified with the holomorphic and antiholomorphic tangent space of $X$ in $x_0$.
Under the complex conjugation $\sigma$ of $\frakg_\CC$ with respect to $\frakg$, respectively 
of $G_\CC$ with respect to $G$, we have $\sigma(\frakp_\pm) = \frakp_\mp$ 
and $\sigma\circ \exp=\exp\circ \sigma$.

Let $\frakh\subset \frakk$ be a maximal abelian subalgebra, in particular $\frakc(\frakk)\subset \frakh$. 
The complexification $\frakh_\CC= \frakh\otimes \CC \subset \frakg_\CC$ is a Cartan subalgebra in $\frakg_\CC$. 
The system $\Psi=\Psi(\frakg_\CC, \frakh_\CC)$ of roots of $\frakg_\CC$ relative to $\frakh_\CC$ induces 
a root space decomposition $\frakg_\CC= \frakh_\CC \oplus \bigoplus_{\a\in \Psi} \frakg_\a$. 
This decomposition is compatible with the decomposition 
$\frakg_\CC=\frakk_\CC \oplus \frakpp \oplus\frakpm$. 
Denote by 
$\Psi^k=\left\{ \a\in \Psi\,|\, \frakg_\a \subset \frakk_\CC\right\} $ the set of compact roots, 
by $\Psi^p_\pm= \left\{ \a\in \Psi\,|\, \frakg_\a \subset \frakp_\pm\right\} $ the set of positive (resp. negative) noncompact roots.
To each root $\a\in \Psi$ we associate a three dimensional simple subalgebra 
$\frakg[\a]=\langle H_\a, E_\a, E_{-\a} \rangle \subset \frakg_\CC$, 
where $H_\a\in \frakh_\CC$ is the element uniquely determined by $\a(H) = 2 \frac{\frakB(H, H_\a)}{\frakB(H_\a, H_\a)}$ for all $H\in \frakh_\CC$,
and $E_\a$ are generators of $\frakg_\a$ determined by requiring that $\tau(E_\a) = -E_{-\a}$, where 
$\tau$ is the complex conjugation of $\frakg_\CC$ with respect to the compact real form $\frakg_U= \frakk+i\frakp$, and that $[E_\a, E_{-\a}]= H_\a$. 
Then $\frakp_\pm = \sum_{\a\in \Psi^p_\pm} \CC E_\a$ and the elements $X_\a= E_\a+E_{-\a}$ and 
$Y_\a= i(E_\a - E_{-\a})$ form a real basis of $\frakp$.

\section{The Borel embedding}
From now on we assume that $G={\rm Is}(X)^\circ$; then  $G_\CC= \aut(\frakg_\CC)^\circ$ 
is a connected linear algebraic group defined over $\RR$ and $G= G_\CC(\RR)^\circ$.
Denote by $P_+, P_-, K_\CC$ the analytic subgroups of $G_\CC$ with Lie algebras 
 $\frakp_+, \frakp_-, \frakk_\CC$, respectively.
Then $K_\CC,\,P_+,\,P_-$ are algebraic subgroups, $K_\CC$ normalizes $P_+$ and $P_-$. 
The semidirect products $K_\CC P_\pm$ have Lie algebras $\sum_{\a\in\Psi^k\cup \Psi^p_\pm} \frakg_\a$.
In particular, $K_\CC P_+$ and $K_\CC P_-$ are opposite parabolic subgroups of $G_\CC$, 
with $P_\pm$ being the unipotent radicals.
Proposition~14.21. in \cite{Borel_lag} implies that the product map induces an isomorphism of varieties of $P_\pm \times K_\CC P_\mp$ onto 
a Zariski open and dense subset of $G_\CC$ which equals $P_\pm K_\CC P_\mp$. 
Thus the decomposition of $g\in P_\pm K_\CC P_\mp$ as  $g= (g)_\pm (g)_{0^\pm} (g)_\mp$ is unique and  
the projection  $P_\pm K_\CC P_\mp\to K_\CC$ defined by this decomposition is a regular map. 

The homogeneous spaces $M_\pm= G_\CC/K_\CC P_\mp$ are compact complex
varieties. 
Let $x_0^\pm$ denote the base points $x_0^\pm=eK_\CC P_\mp\in M_\pm$.
Let $U_\pm= P_\pm x_0^\pm\subset M_\pm$. 
Define maps 
\bqn
\xi_\pm:\frakp_\pm\to U_\pm \subset M_\pm \text{   by   } \xi_\pm(z) = \exp(z) x_0^\pm,
\eqn
with $\exp:\frakg_\CC\to G_\CC$ being 
the exponential map.
\begin{lemma}
The maps $\xi_\pm$ are $K_\CC$-equivariant isomorphisms of varieties.
\end{lemma}
\begin{proof}
The orbit maps $P_\pm \to U_\pm$ are isomorphisms of varieties (see \cite[Theorem~14.12]{Borel_lag}). 
Since $P_\pm$ is unipotent, the exponential maps $\exp:\frakp_\pm \to
P_\pm$ are isomorphisms of varieties (see Appendix). 
The maps $\xi$ are the compositions of the isomorphism $\exp$ 
with the isomorphisms of varieties $P_\pm\to U_\pm$.  
\end{proof}

Since $P_\pm K_\CC P_\mp$ contains $G$ and satisfies $K=G\cap K_\CC P_\pm$, we 
can define holomorphic embeddings, called \emph{Borel embeddings}, 
\bqn
X=G/K\to M_\pm \, \text{ by } \, gK\mapsto g x_0^\pm.
\eqn

\subsection{Kernel functions}
Define $U\subset \frakp_+\times \frakp_-$ to be the subset of $z,w\in \frakp_+\times \frakp_-$ 
such that $\exp(-w) \exp(z)\in P_+ K_\CC P_-$. Define $U'\subset \frakp_+\times \frakp_+$ to be the subset of $z_1,z_2\in \frakp_+\times \frakp_+$ 
such that $(z_1, \sigma(z_2))\in U$. 
%

Let $V_\pm \subset G_\CC\times \frakp_\pm$ be the set of all $g\in  G_\CC,\, z\in \frakp_\pm$ 
such that 
\bqn
g(\exp(z))\in P_\pm K_\CC P_\mp.
\eqn 
We define the
\emph{automorphy factors} 
\bqn
J_\pm :V_\pm\to K_\CC
\eqn 
by $J_\pm (g, z)= (g \exp(z))_{0^\pm}$, and the \emph{complexified automorphy kernel} 
\bqn
K^\CC:U \to K_\CC
\eqn
by  $K^\CC(z , w) = (\exp(-(w)) \exp(z))_{0^+}= J_+ (\exp(-w), z)$.

Define
\bqn
k^\CC:= \det (\Ad_{|_{\frakp_+}} K^\CC):U\to \CC \text{    and    } j_\pm(g, z) :=  \det \Ad_{|_{\frakp_\pm}} J_\pm(g,z): V_\pm \to \CC.
\eqn
\begin{lemma}\label{lem:kernel_prop}
Then $k^\CC$ is a rational function on $\frakp_+\times \frakp_-$ which, on the domain of definition, satisfies
\begin{enumerate}
\item{$k^\CC(gz, gw)= j_-(g,w) k(z,w) j_+(g, z)^{-1}$, where $gz= \log((g \exp(z))_+)$}
\item{$ k^\CC(\sigma(w), \sigma(z))= \overline{k^\CC(z, w)}$}
\end{enumerate}
\end{lemma}
\begin{proof}
The map $k^\CC$ is a regular function on the Zariski open subset $U\subset \frakp_+\times \frakp_-$, 
since it is the composition of two regular maps $K:U\to K_\CC$ and $\det(\Ad|_{\frakp_\pm} ):K_\CC\to \CC$.
Properties (1) and (2) follow by direct calculation. 
\end{proof}

\begin{rem}
In the domains of definition $j_\pm(g,z)$ is the complex Jacobian of the holomorphic action of $G_\CC$ on $\frakp_\pm$, 
 defined by $g(z)= \log((g \exp(z))_+)$
\end{rem}

\section{The Harish-Chandra embedding}
Hermitian symmetric spaces of noncompact type are in one-to-one
correspondence with 
{\em bounded symmetric domains}.

\begin{defi}
A bounded symmetric domain $\Dd\subset \CC^N$ is a bounded open connected subset of $\CC^N$, 
such that every $z\in \Dd$ is the isolated fixed point of an involutive automorphism 
$s_z\in \aut(\Dd)$.
\end{defi}

Let us recall the construction of the Harish-Chandra maps 
$\Phi_\pm:X\to \frakp_\pm$ which realize $X$ as bounded symmetric domains $\Dd_\pm=\Phi_\pm(X) \subset \frakp_\pm \cong \CC^N$. 
The images of $\xi_\pm$ in $M_\pm$ contain the images of  $X$ under the Borel embedding into $M_\pm$. 
The inverses of the maps $\xi_\pm$ restricted to the images of $X$ in $ M_\pm$ 
define the Harish-Chandra maps $\Phi_\pm:X\to \frakp_\pm$.
To take a closer look at the maps $\Phi_\pm$, we need the notion of
{\em strongly orthogonal roots}.

\subsection{Strongly orthogonal roots} 
Two roots $\a,\b \in \Psi$ are called {\em strongly orthogonal} 
if neither $\a-\b$ nor $\a+\b$ is a root.
By a theorem of Harish-Chandra 
(see \cite[p. 582-583]{Harish-Chandra_VI}) 
there exists a (maximal) set $\left\{ \xi_1,\dots,\xi_r\right\} \subset \Psi^p_+$ of $r=r_X$ 
strongly orthogonal positive noncompact roots, 
such that the associated vectors $X_j= X_{\xi_j} \in \frakp, \, j=1,\dots, r$ span a maximal abelian 
subspace $\fraka$ of $\frakp$ over $\RR$.
The set of strongly orthogonal roots determine a fundamental system 
of the root system $\Xi(\frakg,\fraka)$. 

Denote by $\frakh^-$ the real span of $ H_j = H_{\xi_j}$. 
The subalgebra $\frakh^-$ is, via the adjoint action of the Cayley element 
$c=\exp (\frac{\pi}{4} i \sum_{j=1}^r X_j)\in G_\CC$, 
$c\in G_\CC$, conjugate to $\ad(Z_\frakg)\fraka$. 
Under $\ad(Z_\frakg)\Ad(c)$ the root system $\Xi(\frakg, \frakh^-)$ is mapped 
to the root system $\Xi(\frakg, \fraka)$ of 
real roots of $\frakg$ relative to $\fraka$. 
Denote the image of $\xi_j$ by $\mu_j \in \Delta(\frakg, \fraka)$.
Let $\Delta=\left\{ \g_1,\dots\g_r\right\} $ be a fundamental system 
of simple positive roots of $\Xi(\frakg, \fraka)$ 
with respect to the lexicographical order, 
induced by the order on $\Psi$ which was chosen above.

Then, assuming for simplicity that $\frakg$ is simple, 
$\g_j = \mu_j-\mu_{j+1}$ for $j\neq r$ and $\g_r=2\mu_r$ 
if $\Xi$ is of type $C_r$ 
respectively $\g_r=\mu_r$ if $\Xi$ is of type $BC_r$.
In particular, $\g_r$ is distinguished from 
all other positive simple roots 
as the unique simple positive root in $\Delta$ 
which is the restriction of a noncompact positive root.
The positive Weyl chamber in $\fraka$ 
associated to $\Delta$ is given by 
\bqn
\overline{\fraka}^+=\left\{ X\in \fraka \,|\,\g_i(X)\geq 0, \, \forall \g_i \in \Delta\right\} =\left\{ \sum_{i=1}^{r} \lambda_i X_i\, |\, 
\lambda_1\geq \lambda_2 \geq \cdots\geq \lambda_r\geq 0\right\} ,  
\eqn
since $\g_j(\sum_{i=1}^{r} \lambda_i X_i)= \lambda_j-\lambda_{j+1}$ 
for $j=1,\dots, r-1$ and 
$\g_r(\sum_{i=1}^{r} \lambda_i X_i)=\lambda_r$.
\\
\subsection{Polydiscs}
Since the roots $\xi_j$  are strongly orthogonal the associated subalgebras $\frakg[\xi_j]$ 
and its real forms $\frakg(\xi_j) =\frakg[\xi_j] \cap \frakg$  are orthogonal. 
A maximal set of strongly orthogonal roots gives rise to a 
subspace $P\subset X$, $P= \exp(\oplus_{j=1}^r \frakg(\xi_j)) x_0$ which is isomorphic to a product of 
$r$-discs. 

\begin{defi}
A holomorphically embedded Hermitian symmetric subspace $P\subset X$ isomorphic to $\DD^r$ 
is called a \emph{maximal polydisc} in $X$, where $r=\rk_X$. 
\end{defi}
Maximal polydiscs in $X$  
are ${\rm Is}(X)^\circ$-conjugate. 
They are complexifications of maximal flats in $X$
\begin{defi}
A Hermitian symmetric subspace $C\subset P\subset X$ which is isomorphic to a diagonally 
embedded disc $\DD \to \DD^r$ is called a \emph{tight holomorphic disc} in $X$.
\end{defi}

\subsection{The Harish-Chandra realization}
The images of 
\bqn
\fraka_\pm:=\left\{ \sum_{j=1}^{r} \tanh(\lambda_j) E_{\pm j} \, |\,E_{\pm j}=E_{\pm \xi_j},\, \lambda_j\in \RR\right\} \subset \frakp_\pm
\eqn
under the maps $\xi_\pm:\frakp_\pm\to M_\pm$ are 
\begin{equation}\label{eq:flat_hc}
\xi_\pm(\fraka_\pm)= Ax_0^\pm \subset Gx_0^\pm\subset M_\pm,  
\end{equation}
with $A=\exp(\fraka)$.
Since $G=KAK$, $X= KA x_0^\pm \subset M_\pm$ and $[\frakk, \frakp_\pm]\subset \frakp_\pm$ it follows from the $K_\CC$-equivariance of $\xi_\pm$ that 
\bqn
X=\xi_\pm\left( \left\{ \Ad(k) \sum_{j=1}^{r} \tanh(\lambda_j) E_{\pm j}\, |\,k\in K,\, \lambda_j\in \RR\right\} \right).
\eqn
In particular, $X$ is contained in the image of $\xi_\pm$, and 
its preimages 
\bqn
\Dd_\pm= \Phi_\pm(X)= \left\{ \Ad(k) \sum_{i=1}^{r} \tanh(\lambda_j) E_{\pm j}\, |\,k\in K,\, \lambda_j\in \RR\right\}  \subset \frakp_\pm
\eqn 
are bounded symmetric domains, called the \emph{Harish-Chandra realizations} of $X$. 
The Harish-Chandra realization induces an isomorphism ${\rm Is}(X)^\circ \cong \aut(\Dd_\pm)^\circ$, 
which links the Riemannian structure of $X$ with the complex structure of $\Dd_\pm$.
%

\subsection{Structures of a bounded symmetric domain}\label{sec:hermit_domain}
As bounded domain in $\CC^N$, $\Dd$ is 
equipped with a {\em Bergmann kernel} function $k:\Dd\times\Dd\to \CC$ defined by the properties 
\begin{enumerate}
\item{for all $w\in \Dd$: $k_w(z):=k(z,w) \in \Hh^2(\Dd)$}
\item{$k(z,w) = \overline{k(w,z)}$}
\item{$\int_\Dd k(z,w) f(w) d\mu(w) = f(z)$ for all $f\in \Hh^2(\Dd)$.}
\end{enumerate}
where $\Hh^2(\Dd)$ denotes the space of holomorphic square integrable functions 
on $\Dd$.
On a bounded domain one has $k(z,z) >0$ and hence 
$\log k(z,z)$ is well defined. 
Its Hessian $h_{ij}=\frac{\partial^2}{\partial z_i \partial z_j}\log k(z,z)$ defines 
a Hermitian metric called \emph{Bergmann metric}.
The associated K\"ahler form is given by 
\bqn
\omega_\Dd= \frac{1}{2}d d_\CC \log k(z,z).
\eqn
The positive definite real part of the Bergmann metric, 
defines a Riemannian metric $q$ on $\Dd$. 
On a bounded symmetric domain, $q$ is related to the Killing form via $q_0=\frac{1}{2}\frakB$ on $\frakp\cong T_0 \Dd$.

\begin{lemma}\label{lem:kernel}
The Bergmann kernel $k:\Dd\times\Dd\to \CC$ on a bounded symmetric domain $\Dd\subset \frakp_+$ satisfies 
$k(z,w) = k^\CC(z, \sigma(w))$ for all $z,w \in \Dd$.
\end{lemma}
\begin{proof}
The Bergmann kernel $k$ on a bounded symmetric domain arises as $\det(\Ad|_{\frakp_+} K)$ from the automorphy kernel 
$K:U'\subset \frakp_+\times \frakp_+\to K_\CC$ defined by $K(z,w) = (\exp(-\sigma(w))\exp(z))_{0^+}$. 
Hence $K(z,w) = K^\CC(Z, \sigma(w))$ for all $(z,w)\in U'$.
It follows from  $G \subset P_+K_\CC P_-$ that $\Dd\times \Dd \subset
U'$, 
thus $K(z,w) = K^\CC(Z, \sigma(w))$.
\end{proof}

The Bergmann kernel on a bounded symmetric domain is related to the 
polarization of the $K$-invariant polynomial $h$ on $\frakp_+$ which 
is determined by $h(\sum_{i=1}^r \lambda_i E_i) = \Pi_{i=1}^r (1-\lambda_i)$. 
If we normalize the Bergmann kernel on $\Dd$ in such a way that the minimal holomorphic curvature 
of its associated metric is $-1$, it satisfies $k(z,w) = h(z,w)^{-2}$.

\begin{defi}
The  \emph{Shilov boundary} $\cs$ of a bounded domain $\Dd\subset \CC^N$ 
is the unique minimal subset of $\overline{\Dd}$ 
with the property that all functions $f$, continuous on $\overline{\Dd}$ and holomorphic on $\Dd$, satisfy 
$|f(x)|\leq\max_{y\in \cs}|f(y)|$ for all $x\in \Dd$.
\end{defi}

\subsection{The topological boundary and its components}
The compactification of $\Dd$ by taking its closure in $\CC^N$ is
called the topological compactification. Its structure is discussed in 
\cite{Koranyi_Wolf_65_ajm,Wolf,Satake_book}.

The Harish-Chandra realization 
\bqn
\Dd= \Ad (K)\Phi(\exp{\overline\fraka^+} x_0)= 
\Ad (K) \left\{ \sum_{j=1}^r b_j E_j\,|\, 1> b_1\geq \cdots \geq b_r \geq 0\right\} \subset \CC^N
\eqn 
can be compactified by taking the topological closure 
$\overline\Dd= \Dd\cup \partial \Dd $ in $\CC^N$, where 
\bqn
\partial\Dd = \Ad (K) \left\{ \sum_{j=1}^r b_j E_j\,|\, 1\geq b_1\geq \cdots \geq b_r \geq 0\right\} ,
\eqn
with at least one $b_i=1$.
The boundary decomposes into $r=r_X$ $G$-orbits, 
$\partial \Dd = \amalg_{i=1}^r G o_i$, where $o_i = \sum_{j=1}^{i} E_j$.
The $G$-orbits are given by
\bqn
G o_i 
&=& \Ad(K) \left\{ \sum_{j=1}^{r}b_j E_j\,|\, 1= b_1= \cdots=b_i >b_{i+1}\geq b_r \geq 0\right\} \\
&=& \Ad(K)\left(\sum_{j=1}^{i} E_j+ \left\{ \sum_{j=i+1}^{r}b_j E_j\,|\, 1 >b_{i+1}\geq b_r \geq 0\right\} \right)\\
&=& \Ad(K) (o_i +\Dd_i)\\
\eqn
where $\Dd_i$ is a bounded symmetric domain of rank $(r-i)$. 

From the above description it is clear, that 
$\overline{G o_i} = \amalg_{j\geq i} G o_j$, 
and there is a unique closed $G$-orbit in 
$\partial\Dd$, which is $G o_r$. 

\begin{prop}
The unique closed $G$-orbit in $\partial\Dd$ 
is the Shilov boundary $\cs$ of $\Dd$.
\end{prop}

The spaces $\Ad(k)(o_i+ \Dd_i)$, with $ k\in K$ fixed,  
are called \emph{boundary components} of $\Dd$.
Since $\Dd_r$ is reduced to a point, the Shilov boundary 
$\cs=\Ad(K) o_r$ is a $K$-orbit, consisting of point boundary components.

%
%

We want to understand the 
relation between the geometric boundary of the Hermitian symmetric space and its topological boundary.
In every boundary component there exists a unique point which can be joined 
to $0\in \Dd$ by a geodesic. By the transitivity of $G$ on $\Dd$  
this implies that given any $x\in \Dd$ there exists a unique point in each boundary component which 
can be joined to $x$ by a geodesic.
\begin{lemma}(\cite{Wolf})\label{geodesics}
The set of points $\Ad(K)o_i$, $i=1,\dots, r$ is the 
set of end points of geodesics through $0 \in \Dd$.
\end{lemma}
\begin{proof}
Assume that $\g$ is geodesic ray starting in $0\in\Dd$. 
By the $K$-equivariance of the Harish-Chandra map
we may assume that the tangent vector of the geodesic is given by 
 $\g(t)= \Phi(\exp(tH)x_0)$, where $H \in \overline\fraka^+$. 
Let $i\in \NN$ be the number such that $\g_k(H) = 0$ for all $k\geq i$ 
and $\g_i(H)>0$. 
Then 
the image of this geodesic 
under the Harish-Chandra map (see \ref{eq:flat_hc}) is
$\g(t)= \sum \tanh(\l_k t) E_k$, where $\l_k=0$ for all $k>i$ 
and $\l_k>0$ for all $k\leq i$.
Hence $\lim_{t\to\infty} \g(t)= \sum_{k=1}^{i} E_k = o_i$.
\end{proof}
\begin{cor}
Every point $x\in \cs$ can be joined to $0\in \Dd$ by a geodesic. 
The Shilov boundary is the unique $G$-orbit in $\partial\Dd$ with 
this property.
\end{cor}
\begin{lemma}\label{geo_shilov}
The images of geodesic rays of type $\Theta_r$ in $X$ representing 
the same point $\xi \in X(\infty)$ converge in $\Dd$ to the same point 
 in the Shilov boundary.
\end{lemma}
\begin{proof}
We want to show that the images under $\Phi$ of two geodesics 
$\g(t), \eta(t)$ converging to the same point $\xi\in X_{\Theta_r}(\infty)$ 
converge to the same point in $\cs$. 
This is true for $\g(0)=\eta(0)=x_0$ by the calculations made in Lemma~\ref{geodesics}. 
Since the Harish-Chandra map is $K$-equivariant, 
we may assume without loss of generality that 
$\g(t)=\exp(t\sum\l_i X_i) x_0$ and $\eta(t)= \exp(t\sum\mu_i X_i)x$, 
where $\l_i>0$, $\mu_i>0$ and $x=ax_0$ with $a= \exp (\sum a_i X_i)$, $a_i\in\RR$.
The image of the geodesic $\eta(t)= \exp(\sum (t\mu_i+a_i)X_i)x_0$ under $\Phi$ 
is given by $\Phi(\eta)(t) = \sum \tanh(t\mu_i+ a_i) E_i$. 
Since all $\mu_i>0$ and the $a_i$ are bounded 
$\lim_{t\to\infty}\Phi(\eta)(t)= \sum _{i=1}^r E_i = \lim_{t\to\infty} \Phi(\g)(t)$.
\end{proof}

The definition of the topological compactification via the Harish-Chandra embedding 
requires $X$ to be Hermitian symmetric. Nevertheless, the resulting compactification is 
isomorphic to a specific Furstenberg compactification defined for
every 
symmetric space of noncompact type. 
\begin{prop}(\cite{Moore_64_I, Moore_64_II})
The topological compactification of $\Dd\subset \CC^N$ is isomorphic to the 
$\theta_r$-Furstenberg compactification 
of $X$ in  $\Mm(G/P_{\theta_r})$.
\end{prop}

\section{Two classes: tube type and not}
While all Hermitian symmetric spaces have a realization as a bounded symmetric domain 
$\Dd$ in $\CC^N$, generalizing the Poincar\'e disc 
model of the complex hyperbolic line given by the Harish-Chandra realization of $X$ as $\Dd\subset \CC^N$,  
the upper half plane model $\Hh_1=\RR+i\RR$ of the complex hyperbolic line can only be 
generalized for a specific class of Hermitian symmetric spaces.
This class is called {\em ``tube type''}.
\begin{defi}
A Hermitian symmetric space $X$ is said to be of \emph{tube type}
if it is biholomorphically equivalent to a tube type domain 
\bqn
T_\Omega:=\left\{ v+iw|\, v\in V,w\in \Omega\subset V\right\} ,
\eqn 
where $V$ is a real vector space and $\Omega$ is an open cone in $V$.
\end{defi}

Recall that $X_j$ are real vectors in $\frakp$ associated to the strongly orthogonal roots $\xi_j$.
Define the Cayley element to be
\bqn
c:=\exp (\frac{\pi}{4} i \sum_{j=1}^r X_j)\in G_\CC.
\eqn

TIf $X$ is of tube type, then the image of $\Dd_+$ under $c$ is a tube type domain $T_\Omega$ in
$\frakp_+$, 
and thus a generalized upper half
space model of $\Dd$. 
%
%
Irreducible classical domains of tube type are those associated to 
${\rm Sp}(2n, \RR)$, i.e. the Siegel upper half spaces $\Hh_n$, furthermore 
those associated to the groups  
${\rm SO}^*(2n)\, (n \text{ even})$, ${\rm SU} (n,n)$, 
${\rm SO}(2,n)$. The exceptional bounded symmetric domain of rank $3$ is of tube type.
The classical domains associated to 
${\rm SO}^*(2n)\, (n \text{ odd})$, ${\rm SU} (n,m) \, (n\neq m)$, and 
the exceptional domain of rank $2$ are not of tube type.

\subsection{The Cayley transform}
The conjugation of $G < G_\CC$ by 
Cayley element $c=\exp (\frac{\pi}{4} i \sum_{j=1}^r X_j)\in G_\CC$ 
defines a group $G^c= cGc^{-1} \subset G_\CC$, which is isomorphic 
to $G$. 

The orbit of the point $c(x_0^\pm) \subset M_\pm$ under $G^c$  
is contained in the image of $\xi_\pm: \frakp_\pm \to U_\pm$. 
Thus its preimage $T_\Omega \subset\frakp_\pm$ is well defined and contained in 
$\frakp_\pm$. 

Since $G^c c(x_0^\pm)=  cGc^{-1}c(x_0^\pm) = c G (x_0) \subset M_\pm$ 
is isomorphic to $X= G(x_0)\in M_\pm$, their preimages under $\xi_\pm$ 
are isomorphic, 
hence $T_\Omega = c(\Dd) \cong \Dd$.

\subsection{Other characterizations}
The difference of being of tube type or not is reflected in several geometric properties, for example in the 
structure of the restricted root system.

Recall that $\frakh^-$ is the real span of $ H_j = H_{\xi_j}$. 
If we denote the restriction of a root $\a \in \Psi(\frakg_\CC, \frakh_\CC)$ to $ \frakh^-$ again by $\a$, 
the system $\Xi(\frakg, \frakh^-)$ of $\RR$-roots of $\frakg$ relative to $\frakh^-$ is given by
\begin{enumerate}
\item{tube type:  $\Xi$ of type $C_r=\left\{ \pm\xi_i\pm\xi_j  \,(i\neq j),\,\pm 2\xi_i\right\} $, 
where the roots $\pm\xi_i\pm\xi_j$ have multiplicity $a$, 
and the roots $\pm 2\xi_i$ have multiplicity one.}
\item{non-tube type: $\Xi$ of type $BC_r=\left\{ \pm\xi_i\pm\xi_j \,(i\neq j)\, ,\pm 2\xi_i, \, \pm \xi_i\right\} $, 
where the roots $\pm\xi_i\pm\xi_j$ have multiplicity $a$, 
the roots $\pm 2\xi_i$ have multiplicity one, 
and the roots $\pm \xi_i$ have multiplicity $2b$.}
\end{enumerate}
Note that here $\xi_j$ are the restrictions of the strongly orthogonal roots.

\begin{rem}
From the classification it follows that 
the data $(r,a,b)$ completely determine the associated Hermitian
symmetric space. 
There seems to be no direct proof known.
\end{rem}

We recall some of the equivalent characterizations of tube type symmetric spaces given in 
\cite{Koranyi_Wolf_65_annals} and \cite{Deligne}. 
A new characterization in terms of triples of points on the Shilov boundary 
will be given in Proposition~\ref{connected} in Chapter 4.
\begin{prop}
Let $X$ be an irreducible Hermitian symmetric space. 
The following are equivalent:\\
1) $X$ is of tube type.\\
2) The homogeneous space $\cs$ is a symmetric space.\\
3) $\dim_\RR \cs= \dim_\CC X$.\\
4) $\pi_1(\cs)$ is infinite, in which case $\pi_1(\cs) \cong \ZZ \, 
\text{mod torsion}$.\\
5) The system of $\RR$-roots is of type $C_r$ (see above). \\
6) The generator $Z_\frakg$ of the center of 
$\frakk$ is contained in a three dimensional 
simple Lie subalgebra of $\frakg$.\\
7) The Cayley transform $\Ad(c)$ is of order $4$.\\
8) $\aut(\Dd)$ is disconnected.
\end{prop}

\subsection{Jordan algebras, pairs and triple systems}
Hermitian symmetric domains in their Harish Chandra realization 
or their Cayley transformation 
can be decribed in terms of Jordan pairs and 
Jordan triple systems. 
Detailed descriptions of this approach can be found in 
\cite{Loos_symmetric_spaces, Loos_jordan_pairs, Loos_jordan_triples} and \cite{Bertram}. 
For Hermitian symmetric spaces of tube type, these structures simplify to the structures of 
Euclidean Jordan algebras (see \cite{Faraut_Koranyi} for a treatment of 
Jordan algebras associated to symmetric bounded domains).
We will make use of the structure of a Euclidean Jordan algebra to parametrize 
the Shilov boundary of a Hermitian symmetric space and to define a metric on the 
tangent space of the Shilov boundary.

\section{Examples}
We are going to illustrate the above structures by looking at the examples 
of classical Hermitian symmetric spaces. 
Our guiding example is the space of complex Grassmannians. 
We describe this example in detail. The other examples will be described much shorter.
A calculation of the Bergmann kernels of all classical domains is given in Hua's book \cite{Hua_1}. 
For the concrete realization of the Bergmann kernel of the exceptional domains see \cite{Drucker}.

\begin{notation*}
A $p\times q$ 
matrix consisting of a $p \times p$ matrix and $(q-p)$ zero columns
or of a $q \times q$ matrix and $(p-q)$ zero lines, will be denoted by 
$(A,0)$.
\end{notation*}
\subsection{Complex Grassmannian}
Let $V$ be a $n=p+q$--dimensional $\CC$--vector space, endowed with a Hermitian form $h$ of signature $(p,q)$. 
Assume without loss of generality that $p\leq q$. 
Denote by $(e_1, \dots, e_p, e_{p+1}, \dots, e_n)$ a basis of $V$, then we may suppose that 
$h$ is defined as 
\bqn
h(z,w) = \sum_{i=1}^p z_i \ol{w_i} - \sum_{i=1}^q z_{p+i} \ol{w_{p+i}}.
\eqn
Let $W_+$ be the span of  $(e_1, \dots, e_p)$ and  $W_-$ the span of $( e_{p+1}, \dots, e_n)$ in $V$. 
Then $h|_{W_+}$ is positive definite, and the restriction of $h$ to $W_-= (W_+)^{\perp_h}$ is negative definite. 

We define 
\bqn
X_{p,q}:= \left\{  W^p \subset V \,|\, h|_{W^p}>0\right\} . 
\eqn

We denote by $s_W$ the reflection with respect to $W^p\subset V$ on $V$, 
then $s_W$ induces a map 
\bqn
s_W: X_{p,q} \to X_{p,q}
\eqn
with $s_W \circ s_W = \id$ and $\fix(s_W) = W \in X_{p,q}$, giving 
$X_{p,q}$ the structure of a symmetric space. 

To realize $X_{p,q}$ as homogeneous space $G_{p,q}/K$, 
choose $x_0= W_+\in X_{p,q}$ as base point. 
Let $G_{p,q} = \aut(V,h)$, i.e.
\bqn
G_{p,q} = {\rm SU}(p,q) :=\left\{  g \in {\rm SL}(n,\CC) \,|\, g^* J g = J\right\} , 
\eqn
where $ g^* = \ol{g}^t$ and 
$J = \begin{pmatrix} \id_p & 0 \\ 0 & -\id_q \end{pmatrix}$. 

Since ${\rm SL}(n,\CC)$ acts transitively on the $p$-dimensional subspaces of $V$, 
the group $G_{p,q}$ acts transitively on the  $p$-dimensional subspaces of $V$ on which the 
restriction of $h$ is positive definite. 

Writing an element $g$ as block matrix $g= \begin{pmatrix} A & B \\ C & D \end{pmatrix}$ with 
$A \in \Mat(p,p,\CC)$, $B \in \Mat(p,q,\CC)$, $C \in \Mat(q,p,\CC)$, $D \in \Mat(q,q,\CC)$ we obtain the following conditions: 
\bqn
A^* A-C^* C=\id_p \\
 B^* B- D^*D = -\id_q\\
B^* A-D^* C= 0 \\
\det(g) = 1.
\eqn

The Lie algebra $\frakg_{p,q}$ of $G_{p,q}$ is given by 
\bqn
\frakg_{p,q} = \fraks\fraku(p,q) :=\left\{  l \in \fraks\frakl(n,\CC) \,|\, J l =-l^* J\right\} , 
\eqn
implying the following conditions on $l= \begin{pmatrix} A & B \\ C & D \end{pmatrix}$:
\bqn
A^* =-A, D^* = -D, B^*= C, tr(g) = 1.
\eqn
An element $g= \begin{pmatrix} A & B \\ C & D \end{pmatrix}\in G_{p,q}$ preserves $W^+$ iff $C=0$, 
hence with the above conditions $B=0, A^* A=\id_p , D^*D = \id_q$, i.e. 
\bqn
Stab_{G_{p,q}}(W^+) = {K}(p,q)= {\rm S}({\rm U}(p) \times {\rm U}(q)) 
\eqn
and $X_{p,q}= {\rm SU}(p,q)/ {\rm S}({\rm U}(p) \times {\rm U}(q))$.

The Lie algebra $\frakg_{p,q}$ decomposes as $\frakg= \frakk\oplus \frakp$ 
with 
\bqn
\frakk &=& \left\{  l \in \fraks\fraku(p,q) \,|\,l= \begin{pmatrix} A & 0 \\ 0 & D \end{pmatrix}\right\} ,\\
\frakp &=&\left\{  l \in \fraks\fraku(p,q) \,|\,l= \begin{pmatrix} 0 & B \\ B^* & 0 \end{pmatrix}\right\} .
\eqn 


The subalgebra $\fraka \subset \frakp$, defined by $B=(Diag(\lambda_1, \dots, \lambda_p),0)$ 
is a maximal abelian subalgebra of $\frakg$. 
Thus $F= \exp(\fraka)x_0= A x_0$ is a maximal flat in $X_{p,q}$ and $rk_X= p$.

The flats are the $p$-dimensional subspaces $W\subset V$, which are spanned by a basis 
$e_i+\lambda_i t e_{p+i}$ with $(\lambda_i t)^2 <1$, $i=1,\dots p$.

Furthermore 
\bqn
\frakk_\CC &=& \left\{  l \in \fraks\frakl(n,\CC) \,|\,l=\begin{pmatrix} A & 0 \\ 0 & D \end{pmatrix}\right\} ,\\
\frakp_+ &=& \left\{  l \in \fraks\frakl(n,\CC) \,|\,l= \begin{pmatrix} 0 & B \\ 0 & 0 \end{pmatrix}\right\} ,\\
\frakp_- &=& \left\{  l \in \fraks\frakl(n,\CC) \,|\,l= \begin{pmatrix} 0 & 0 \\ C & 0 \end{pmatrix}\right\} .
\eqn 

The Borel embeddings 
are just the embedding of $X_{p,q}$ into the compact Grassmann manifold of $p$-dimensional respectively $q$-dimensional subspaces  
in $V$ 
\bqn
X_{p,q} &=& \left\{  W^p \subset V \,|\, h|_{W^p}>0\right\}  \subset   \left\{  W^p \subset V \right\} =X^c_{p,q}\\
X_{p,q} &=& \left\{  W^q \subset V \,|\, h|_{W^q}<0\right\}  \subset  \left\{  W^q \subset V \right\} =X^c_{p,q}.
\eqn

We describe the Harish-Chandra embedding 
$\Phi: X_{p,q} \to \frakp_- \cong \Mat(q,p,\CC)$, in geometrical terms. 

The base point $W^+$ induces a direct decomposition of  
$V= W^+\oplus W^-$. 
Given any other point $W \in X_{p,q}$, we may decompose a vector $w\in W$ 
into its parts with respect to the decomposition $V= W^+\oplus W^-$, 
\bqn
w= \pr|_{W^+}(w) + \pr|_{W^-}(w) = v +   \pr|_{W^-}\circ \pr^{-1}|_{W^+}(v),  
\eqn
with $ v=\pr|_{W^+}(w)$.
Thus we can write $W$ as graph of the linear map 
\bqn
T_W= \pr|_{W^-}\circ \pr^{-1}|_{W^+}: W^+ \to W^-. 
\eqn

With respect to our specific choice of basis we identify the space of linear maps 
$L( W^+, W^-)\cong \Mat(q,p,\CC)$.
The condition that $h|_W >0$ translates to $\id_p- T^*_W \circ T_W >0$, thus 
the Harish-Chandra embedding $\Phi: X_{p,q} \to \frakp_- \cong \Mat(q,p,\CC), \, W\mapsto T_W$ 
realizes $X_{p,q}$ as 
\bqn
\Dd_{p,q}= \left\{  Z \in \Mat(q,p,\CC) \,|\, \id_p - Z^*Z>0\right\} .
\eqn  

The action of $G_{p,q}$ on $\Dd_{p,q}$ becomes 
\bqn
g(Z) =  \begin{pmatrix} A & B \\ C & D \end{pmatrix} (Z) = (AZ + B) (CZ+D)^{-1}.
\eqn

The topological compactification of $\Dd_{p,q}$ is given by 
\bqn
\ol{\Dd_{p,q}}= \left\{  Z \in \Mat(q,p,\CC) \,|\, \id_p - Z^*Z\geq 0\right\} 
\eqn  
and the different boundary components are 
given by 
\bqn
G o_i= \left\{  Z \in \Mat(q,p,\CC) \,|\, \id_p - Z^*Z\geq 0, \, \rk(\id_p - Z^*Z)=r-i \right\} .
\eqn 
In particular, the Shilov boundary is 
\bqn
\cs_{p,q}=\left\{  Z \in \Mat(q,p,\CC) \,|\, \id_p = Z^*Z\right\} .
\eqn 
If $p=q$, then $\cs_{p,p}={\rm U}(p)$.

In the special case of $X= \HH^2$, 
\bqn
G_{1,1}=
{\rm SU}(1,1)=\left\{ \begin{pmatrix} a & b \\ \ol{b} & \ol{a} \end{pmatrix}\,|\, a,b\in \CC \text{ and } |a|^2 -|b|^2=1\right\} ,
\eqn
we describe the Harish-Chandra embedding by decomposing 
an element $g\in G$ with respect to the product $P_+\times K_\CC P_-$:
\bqn
g= \begin{pmatrix} a & b \\ \ol{b} & \ol{a} \end{pmatrix}=
\begin{pmatrix} 1 & b \ol{a}^{-1} \\ 0 & 1 \end{pmatrix} 
\begin{pmatrix} \ol{a}^{-1} & 0 \\ 0 & \ol{a} \end{pmatrix}
\begin{pmatrix} 1 & 0 \\ \ol{a}^{-1} \ol{b} & 1 \end{pmatrix},
\eqn
Then  $\Phi: X=G/K \to \frakp_+$ is given by 
$\Phi(g x_0)=b{\overline a} ^{-1} e_+=ze_+,$ with $|z|<1$.
Hence, identifying $\frakp_+$ with $\CC$, 
we get an $K$-equivariant biholomorphic map $\Phi: X\to \DD$. 
The map $\Phi$ sends the base point $x_0$ to $0\in\CC$ and 
the geodesic $\g=\exp(t\lambda x)x_0$ in $X$ 
to a geodesic $\g' = \tanh (t \lambda)e_+$ in $\DD$. 
The corresponding embedding into $\frakp_-$ gives the same embedding as the 
 geometric construction before.

Using the Cartan decomposition of ${\rm SU}(p,q)=G= KA^+ K$ and letting $X_i = \begin{pmatrix} 0 & E_i^* \\ E_i & 0 \end{pmatrix}$ with $E_i \in \Mat(q,p,\CC)$ 
given by $(E_i)_{i,i}=1 $ and all other entries zero be a basis of $\fraka$ corresponding to a set of strongly orthogonal roots. 
We see that a geodesic 
$\g=\exp(t \sum_{i=1}^p \lambda_i X_i)x_0$ in $X_{p,q}$ is mapped under 
the Harish-Chandra embedding to the geodesic $\g' = \sum_{i=1}^p\tanh (t \lambda_i) E_i $ in $\Dd_{p,q}$.
One can see directly that the images of all geodesics with $\lambda_i>0$ for all $i$ 
converge to the point $Z= (\id_p, 0)\in \cs_{p,q}$.

We compute the Bergmann kernel for $\Dd_{1,1}$.
Let $(z,w)\in \Dd_{1,1}\times \Dd_{1,1}$. 
The automorphy kernel $K(z,w)$ is the $K_\CC$-part 
of $\exp(-\ol{w} e_+)\exp (ze_-)$ in the decomposition $P_- K_\CC P_+$. 
Now 
\bqn
\exp(-\ol{w} e_+)\exp (ze_-)= \begin{pmatrix} 1 & -\ol{w} \\ 0 & 1
\end{pmatrix}
\begin{pmatrix} 1 & 0 \\ z & 1\end{pmatrix} = 
\begin{pmatrix} 1-\ol{w}z & -\ol{w} \\ z & 1 \end{pmatrix} = \\
\begin{pmatrix} 1 & 0 \\ \frac{z}{1-\ol{w}z} & 1\end{pmatrix}
\begin{pmatrix} 1-\ol{w}z & 0 \\ 0 & \frac{1}{1-\ol{w}z}\end{pmatrix} 
\begin{pmatrix} 1& \frac{-\ol{w}}{1-\ol{w}z} \\ 0 & 1\end{pmatrix}.
\eqn
 
Hence 
\bqn
K(z,w) = \begin{pmatrix} 1-\ol{w}z & 0 \\ 0 &
  \frac{1}{1-\ol{w}z}\end{pmatrix}.
\eqn
Since 
\bqn
K(z,w)  \begin{pmatrix} 1 & 0 \\ a & 1\end{pmatrix} 
K(z,w)^{-1} = \begin{pmatrix} 1 & 0 \\ \frac{a}{(1-\ol{w}z)^2} &
  1\end{pmatrix},
\eqn 
we get 
\bqn
k(z,w)= \det (\Ad|_{\frakp_-} K(z,w))=(1-\ol{w}z)^{-2}.
\eqn

From this formula for the automorphy kernel (using once again the strong orthogonality of the roots), we obtain 
a formula for 
the automorphy kernel of $\Dd_{p,q}$ for 
points lying in 
\bqn
\fraka_-=\left\{ \sum_{j=1}^{r} \tanh(\lambda_j) E_{j} \,\, \lambda_j\in
\RR\right\} \subset \frakp_-.
\eqn 
Namely if $Z=  \sum_{j=1}^{r} \tanh(\lambda_j) E_{j}$ and $W =\sum_{j=1}^{r} \tanh(\mu_j) E_{j}$, 
then $K(Z,W)=  \id_p - W^* Z$ , furthermore  $\Ad|_{\frakp_-} K(Z,W)=
(\id_p - W^* Z)^{-2\dim \frakp_-}$ 
and the normalized Bergmann kernel satisfies 
\bqn
k(Z,W)= (\det (\id_p - W^* Z))^{-2}= \Pi_{i=1}^{p} (1-\tanh(\mu_i)\tanh(\lambda_i))^{-2}.
\eqn

Since the maximal flats are $K$-conjugate, 
we may after suitable conjugation with some $g\in K$ assume that any two points $Z,W\in \Dd_{p,q}$ 
lie in $\fraka_-$, thus we have determined the Bergmann kernel on $\Dd_{p,q}$.


We finish with some remarks about the upper half space model in the tube type case, i.e. for $\Dd_{p,p}$. 
The Cayley element $c=\exp (\frac{\pi}{4} i \sum_{j=1}^r X_j)$ is given by 
\bqn
c= \frac{1}{\sqrt{2}}  \begin{pmatrix} \id_p & (i(\id_p),0) \\
  (i(\id_p),0) & \id_q\end{pmatrix}.
\eqn
Hence, the image of $Z\in \Dd_{p,p}$ is $\frac{1}{\sqrt{2}} (Z+  i\id) (iZ+\id)^{-1}$. 
In particular, 
\bqn
c: \Dd_{p,p} \to T_\Omega = {\rm Herm}(p, \CC) + i {\rm Herm}^+(p,\CC),
\eqn
where ${\rm Herm}^+(p,\CC)$ denotes the Hermitian $(p\times p)$ matrizes, whose associated quadratic form 
is positive definite.
Furthermore $Z \in \cs_{p,p}$ with $(iZ+\id)$ invertible is mapped to a Hermitian $p\times p$ matrix.

\subsection{Moduli space of complex structures}
(A detailed description can be found in Satake's book \cite[Chapter~II,7]{Satake_book}.)
Let $V$ be a real vector space with a nondegenerate skew-symmetric 
bilinearform $h$. Then $V$ is even dimensional and 
with respect to a {\em symplectic basis} $(e_1, \dots, e_{2n})$, 
the form $h$ is given by
\bqn
h=  \begin{pmatrix} 0 & \id_n\\ -\id_n & 0\end{pmatrix}.
\eqn
The symplectic group ${\rm Sp}(2n, \RR)$ is defined by 
\bqn
{\rm Sp}(2n, \RR)= \aut(V, h) \subset {\rm SL}(2n, \RR).
\eqn 
Denote by $\Jj(V)$ the set of all complex structures on $V$. 
The symmetric space associated to  ${\rm Sp}(2n, \RR)$ can be defined as the 
set 
\bqn
X=\left\{ I\in \Jj(V) \, |\, h(\cdot, I \cdot) \text{ is symmetric and positive definite}\right\} .
\eqn
A basepoint is, with respect to the chosen basis of $V$, 
given by 
\bqn
I_0 = \begin{pmatrix} 0 & -\id_n\\ \id_n & 0\end{pmatrix}.
\eqn
Let $V_\CC$ be the complexification of $V$. Extend $h$ to a 
$\CC$-bilinear form on $V_\CC$. 
For $I\in X$ define $V_\pm (I)$ to be the 
$(\pm i)$-eigenspaces of the action of $I$ on $V_\CC$. 
Then $V_\CC= V_+ \oplus V_-$, and $\ol{V_+}= V_-$. 
Furthermore $V_\pm$ are isotropic subspaces for the extension of $h$. 
The restriction of $i h $ to $V_+\times V_-$ is nondegenerate and 
gives an identification $V_-\cong V_+^*$ of $V_-$ with the dual space of 
$V_+$.
Define a Hermitian form $A_h$ on $V_\CC$ by 
\bqn
A_h(v, w)= i h(\ol{v}, w).
\eqn
The restriction of $A_h$ to $V_+$ is positive definite, 
the restriction to $V_-$ is negative definite and the restriction 
to $V_+\times V_-$ is zero. 
Thus the Hermitian form $A_h$ is a nondegenerate Hermitian form on $V_\CC$ of 
signature $(n,n)$.

Every complex structure $I\in X$ defines a subspace $V_\pm\subset V_\CC$, 
this map identifies $X$ with the subset $X^\pm_{2n}$ of 
the complex Grassmannian of $n$-dimensional subspaces in a complex 
$(2n)$-dimensional vector space: 
\bqn
X^\pm_{2n}:= \left\{  V_\pm \subset V_\CC \,|\, h|_{V_\pm}=0,\,  
{A_h}|_{V_\pm\times V_\mp} >0\right\} . 
\eqn

This also gives an embedding 
\bqn
X_{2n}\subset X_{n,n}, 
\eqn
where $ X_{n,n}$ is the symmetric space associated to 
${\rm SU}(n,n)$ described above. This embedding is holomorphic, 
sending the central element in the Lie algebra of the maximal compact subgroup 
defining the complex structure 
on $X_{2n}$ to the central element defining the complex structure on $X_{n,n}$.

We describe the matrix representation of the Lie algebra 
$\frakg= \fraks\frakp(2n, \RR)$ with respect to this embedding. 
This is a representation with respect to a complex basis of $V_\CC$, 
not with respect to a real basis. 
\bqn
\frakg=\begin{pmatrix} X_1 & X_2 \\ \overline{X_2} & \overline{X_1} \end{pmatrix} =\frakk+\frakp, \\ 
\frakk= \begin{pmatrix} X_1 & 0 \\ 0 & \overline{X_1} \end{pmatrix},\,
\frakp= \begin{pmatrix} 0 & X_2 \\ \overline{X_2} & 0 \end{pmatrix},
\eqn
where $X_i \in \Mat_{n\times n} (\CC)$, $X_1^* =-X_1$, $X_2^T= X_2$.
The element $Z_0$ generating the center of $\frakk$ is given by 
\bqn
Z_0 = \frac{i}{2} \begin{pmatrix} \id_n & 0 \\ 0 & -\id_n \end{pmatrix}.
\eqn
The standard maximal abelian Lie subalgebra $\fraka\subset \frakp$ and the maximal polydisc $\frakr$
are:
\bqn
\fraka&=&\left\{  \begin{pmatrix} 0 & A \\ A & 0 \end{pmatrix}| \, A=\diag(a_1,\dots, a_n), \, a_i\in \RR\right\} , \\
\frakr&=&\left\{ \begin{pmatrix} 0 & Z \\ \ol{Z} & 0 \end{pmatrix}|\, Z=\diag(z_1, \dots, z_n) , \, z_i \in \CC\right\} .
\eqn
Further
\bqn
\frakp_+ =\begin{pmatrix} 0 & Z \\ 0 & 0 \end{pmatrix} \\
\frakp_-=\begin{pmatrix} 0 & 0 \\ Z & 0 \end{pmatrix}, 
\eqn
with $Z \in \sym(n,\CC)$.

Similar to the case of $X_{n,n}$  we obtain a bounded domain model and 
an upper halfspace model:
\bqn
\Dd_{2n}=\left\{  Z \in \sym(n,\CC) \,|\, \id_p - Z^*Z>0\right\}  \subset \Dd_{n,n}\\
T_\Omega= \sym(n,\RR) \oplus i \sym(n, \RR)^+,
\eqn
where $\sym(n, \RR)^+$ are the positive definite symmetric matrices. 
In particular $ T_\Omega$ is the set of all 
complex valued 
symmetric $(n\times n)$-matrices with positive definite imaginary part.

The descriptions of the Bergmann kernel are similar to the case of $X_{n,n}$, 
indeed the structures induced from the embedding $X_{2n} \to X_{n,n}$ 
give the corresponding structures on $X_{2n}$.

\subsection{Quaternion Grassmanians}
Let $V$ be a $n$-dimensional 
(left)-vector space over the quaternions $\HH$, 
endowed with nondegenerate sesquilinear anti-Hermitian form $h$.
Denote by $\sigma$ the conjugation in $\HH$ fixing $\RR$.  
\bqn
h(ax, by)= ah(x,y) \sigma(b)\\
h(x,y)= \sigma(h(y,x)).
\eqn
The group $G={\rm SO}^*(2n)$ is defined as 
$G={\rm SO}^*(2n)= \aut_\HH(V, h)$.

Identifying $\HH= \CC \oplus j \CC$, 
there is a uniquely  determined Hermitian form $A_h$ of signature 
$(n,n)$ on $V^\CC= V$ viewed as complex vector space, 
such that 

\bqn
h(x,y)= i(A_h(x,y) - j A_h(xj, y))\\
A_h(xj,yj)= -\ol{A_h(x,y)}. 
\eqn
In particular, $S(x,y) := -A_h(xj, y)$ is a symmetric bilinear form on $V^\CC$ 
and $G={\rm SO}^*(2n)$ can be defined as 
$G={\rm SO}^*(2n)= \aut_\CC(V^\CC, A_h, S)$. 
The maximal compact subgroup of $G$ is $K = {\rm U}(n)$.

This realization of $G$ induces an embedding
\bqn
X_{*2n}= \left\{  V^n_-\subset V^\CC \, |\, S|_{V_-} =0, \,  
{A_h}|_{V_-}\leq 0 \right\} \subset X_{n,n}.
\eqn
This embedding is holomorphic and induces a corresponding embedding 
\bqn
\Dd_{*2n}=\left\{  Z\in \Mat(n,n,\CC)\, |\, Z^T=-Z, \, \id_n-\ol{Z}Z >0 \right\} 
\subset \Dd_{n,n} 
\eqn
Nevertheless the space $\Dd_{*2n}$ is not always of tube type, it is of 
tube type if and only if $n$ is even. The rank of $X_{*2n}$ is $[n/2]$.
Inducing the structures of $\Dd_{n,n}$ on $\Dd_{*2n}$ is thus a bit more 
complicated. The restriction of the Bergmann kernel on $\Dd_{n,n}$ to 
$\Dd_{*2n}$ gives twice the intrinsically defined Bergmann kernel. 

We describe the Shilov boundary of $\Dd_{*2n}$.
The Shilov boundary $\cs_{*2n}$ is given as the $K$-orbit of the 
matrix 
$Z=\begin{pmatrix} J_{2r_X} & 0 \\ 0 & z \end{pmatrix}$, 
where $|z|^2<1$ and $J_{2r_X}$ is the matrix 
with $r_X$ blocks of the form $j= \begin{pmatrix} 0 & 1 \\ -1 & 0 \end{pmatrix}$ on the diagonal 
and zero entries elsewhere. 
If $n$ is even then $Z= J_{n}$ and the Shilov boundary 
$\cs_{*2n}$ embeds into $\cs_{n,n}$. 
If $n$ is odd, the Shilov boundary  
$\cs_{*2n}$ is not contained in $\cs_{n,n}$

\subsection{The Quadric}
Let $V$ be a $(2+n)$-dimensional vector space endowed with 
a symmetric bilinear form $h$ of signature $(2,n)$.
Let $G={\rm SO}(2,n)= {\rm SO}(V,h)$.
Denote by $V_\CC$ the complexification of $V$ and extend the 
bilinear form $h$ to $V_\CC$. 
The Hermitian form 
$A_h(x,y) := 2 h(\ol{x}, y)$ on $V_\CC$ is of signature $(2,n)$.
Take a decomposition of $V_\CC = V^\CC_+ \oplus V^\CC_-$, 
with $\pm {A_h}|_{V^\CC_\pm} >0$. 
Then $V^\RR_+:= V^\CC_+ \cap V$ is a $2$-dimensional real subspace 
of $V$ admitting two complex structures $\pm I$ leaving 
$h|_{V^\RR_+}$ invariant. 
Put $W_+ :=\left\{  x\in V^\CC_+ \, |\, Ix= ix\right\} $ 
then $V^\CC_+ = W_+ \oplus \ol{W_+}$ and 
$h|_{W_+}=0$, ${A_h}|_{W_+} > 0$.

The space 
\bqn
X_{2+n}:=\left\{  W_+\subset V^\CC \, |\, \dim_\CC(W_+)=1, \,h|_{W_+}=0,\, {A_h}|_{W_+} > 0 \right\} 
\eqn
is a complex manifold with two components. One component 
is the symmetric space associated to $G={\rm SO}(2,n)$.
We describe the Lie algebra: 
\bqn
\frakg=\begin{pmatrix} X_{11} & X_2 \\ X_2^T & X_{12} \end{pmatrix} =\frakk+\frakp, \\
\frakk= \begin{pmatrix} X_{11} & 0 \\ 0 & X_{12} \end{pmatrix}, \,
\frakp= \begin{pmatrix} 0 & X_2 \\ X_2^T & 0 \end{pmatrix}, 
\eqn
where $X_{12} \in \Mat_{n\times n} (\RR)$, $X_{11} \in \Mat_{2\times 2}(\RR)$, 
$X_2 \in \Mat_{2\times n}(\RR)$.
The element $Z_0$ generating the center of $\frakk$ is given by 
$Z_0 = \begin{pmatrix}J & 0 \\ 0 & 0\end{pmatrix},$
where $J=\begin{pmatrix}0 & -1 \\ 1 & 0 \end{pmatrix}$.

The bounded symmetric domain model of $X_{2+n}$ is 
\bqn
D_{2+n}= \left\{  z\in \CC^n\, |\, 1+|z^T z|^2 - 2 z^* z >0 \text{ and } 1-z^*z >0 \right\} ,
\eqn
its Shilov boundary is
\bqn
\cs_{2+n} =\left\{ z\in \CC^n\, |\, z^*z =1, \, e^{i\theta} z \in \RR^n \text{ for some } \theta \in \RR/ 2\pi \ZZ \right\} .
\eqn



\vskip1cm
\chapter{The Shilov boundary}\label{sec:shilov}
A symmetric space of noncompact type can be 
compactified in different ways \cite{Guivarch_Ji_Taylor,Borel_Ji}. 
The different compactifications are in general nonisomorphic.
In some cases the identity map extends 
to a continuous surjective $G$-equivariant map on 
the boundaries (or parts of the boundaries) of different compactifications. 
The Harish-Chandra realization of the Hermitian symmetric space $X$ as a bounded 
symmetric domain $\Dd$ provides a natural compactification 
by taking the topological closure in $\CC^N$. 
The Shilov boundary is not the whole boundary of 
this topological compactification, but only a specific part of it, namely the unique closed $G$-orbit 
in the topological boundary. 
But nevertheless the Shilov boundary captures a lot of the geometry of the 
Hermitian symmetric space and reveals aspects of negative curvature 
in $X$.

\section{A homeomorphism of boundaries}
A Hermitian symmetric space can be viewed as an abstract symmetric space 
$X=G/K$ or via its Harish-Chandra realization as bounded symmetric domain 
$\Dd\subset \CC^N$. 
Since we will work with $X$ as homogeneous space as well as with $X$ realized as bounded symmetric domain $\Dd$, 
an understanding of the relation between the geometric boundary and the topological boundary is important.  

Even though there is no globally defined 
equivariant map $X(\infty)\to \partial \Dd$, 
there is the following relation between 
the geometric and topological boundaries.
It follows from the proof of Lemma \ref{geodesics} in Chapter 3 
that any geodesic with $\g(0)=x_0$ 
converging to a point $\g(\infty)\in X_{\Theta_i}(\infty)$ 
converges to a point 
in $\Ad(K) (o_i)$ and conversely. 
But two geodesics $\g$ and $\eta$ 
which converge to the same point 
$\xi \in X_{\Theta_i}$ do in general not converge 
to the same point in $\partial \Dd$ under the Harish-Chandra map 
if $\g$ and $\eta$ do not meet.
Even though this shows that there is no globally defined 
equivariant map $X(\infty)\to \partial \Dd$, 
there is the following relation between 
the geometric and topological boundaries.
The strata $X_{\Theta_r}(\infty)$, especially $X_{\theta_r}(\infty)$, 
and the Shilov boundary 
are closely related.
\begin{prop}\label{boundaries}
The Harish-Chandra map $\Phi: X\to D$ 
extends to a $G$-equivariant surjective map 
$X_{\Theta_r}(\infty)\to \cs$.
In particular, $\Phi$ extends to a $G$-equivariant homeomorphism
$$
\Phi:X_{\theta_r}(\infty)\to \cs.
$$
\end{prop}
\begin{proof}
We know from Lemma~\ref{geo_shilov} in Chapter 3, that the images under $\Phi$ of two geodesics 
$\g(t), \eta(t)$ representing 
the same point $\xi\in X_{\Theta_r} (\infty)$ 
converge to the same point $\Phi(\xi) \in \cs$. 
Since $\cs$ is the $K$ orbit of a point and the Harish Chandra map is 
$K$-equivariant, this shows that $\Phi$ extends to a $G$-equivariant 
surjective map 
$X_{\Theta_r}(\infty)\to \cs$.
The set $X_{\theta_r}(\infty)$ is the $K$ orbit of a unique point $\xi$ 
which is mapped to $o_r=\sum _{i=1}^r E_i$ under the above map.
The map $\Phi$ thus extends to a $K$- and hence $G$-equivariant
homeomorphism $\Phi:X_{\theta_r}(\infty)\to \cs$.
\end{proof}

\section{The Shilov boundary as homogeneous space}
The embedding $X\to M_\pm$ of $X$ into its compact dual provides another 
compactification of $X$. We want to describe the closed orbit $\cs$ in $\partial \Dd_\pm$ 
in the boundary of $X$ in $M_\pm$ as a 
homogeneous space.

Recall that $X_j$ are real vectors in $\frakp$ associated to the strongly orthogonal roots $\xi_j$.
Define the Cayley element 
\bqn
c=\exp (\frac{\pi}{4} i \sum_{j=1}^r X_j)\in G_\CC.
\eqn

\begin{prop}\label{prop:shilov_M}
The closed boundary orbit $\cs$ can be realized in $M_\pm$ as homogeneous space $G/Q$, 
where $Q=G\cap cK_\CC P_- c^{-1}=G\cap c^{-1}K_\CC P_+ c$ is a real parabolic subgroup of $G$. 
Moreover $G/Q= (G_\CC/\qQ) (\RR)$ with $\qQ=  c^{-1}K_\CC P_+ c\cap
cK_\CC P_- c^{-1}$.
\end{prop}
\begin{proof}
Let $x_0^\pm \in M_\pm$ be the basepoint.
The closed boundary orbit $\cs$ is the $G$-orbit of $\pm i o^\pm_r$. 
Direct computations show that $\xi_\pm(\pm io^\pm_r) = c^{\pm1}(x_0)$.
The stabilizer of $c^{\pm1}(x_0)$ in $G_\CC$ is $\stab_{G_\CC}(c^{\pm1}(x_0))= c^{\mp 1} K_\CC P_\mp c^{\pm 1}$.
Thus the realization of $\cs$ in $M_\pm $  is $G/(c^{\mp 1} K_\CC P_\mp c^{\pm 1}\cap G)$. 
Since the $X_j$ are real vectors, we have that $\sigma(c)= c^{-1}$. Hence 
$\sigma (cK_\CC P_- c^{-1}) = c^{-1} K_\CC P_+ c$, thus $G\cap cK_\CC P_- c^{-1}=G\cap c^{-1}K_\CC P_+ c=Q$. 
Since $Q$ stabilizes a point in $\cs$ it is a parabolic subgroup of type $\theta_r$. 

The group $\qQ=  c^{-1}K_\CC P_+ c\cap  cK_\CC P_- c^{-1}$ is stable under $\sigma$, thus it is defined over $\RR$ 
and $\qQ(\RR) = Q$.
%
\end{proof}

\section{The Bruhat decomposition}
As we have seen in the above section, the Shilov boundary $\cs$ of $\Dd$ 
can be realized as 
\bqn
\cs= G/Q= \gG(\RR)/\qQ(\RR)=(\gG/\qQ)(\RR),
\eqn
where $Q=P_{\theta_r}$ denotes the maximal parabolic subgroup of $G$ 
associated with $\theta = \Delta-\{\g_r\}\subset \Delta$.
By the general theory of cosets of 
$G$ by parabolic subgroups \cite{Borel_lag}, 
the Shilov boundary admits a generalized Bruhat decomposition.

Define $W_\theta$ to be the subgroup 
of the Weyl group $W$ of $\frakg$ relative to $\fraka$, 
which is generated by reflections $s_i$ corresponding to $\g_i\in\theta$.
In each left coset $w W_\theta$ ($w\in W$) 
there exists a unique representative $w^\theta$ of smallest length. 
Let $W^\theta$ be the set of these representatives. 
We define 
\bqn
C(w^\theta)\qQ:= \pP w^\theta\qQ.
\eqn

Then 
\bqn
\gG/\qQ= \bigcup_{w^\theta \in W^\theta} C(w^\theta).
\eqn
If $w_{max}$ denotes the longest element in $W$ 
then 
\bqn
C(w^{\theta}_{max})\qQ= \pP w^{\theta}_{max} \qQ= \qQ w^{\theta}_{max} \qQ
\eqn
 is a Zariski-open $\qQ$-orbit in $\gG/\qQ$.

On the level of $\RR$-points this gives the decomposition
$$G/Q= \bigcup_{w^\theta \in W^\theta} C(w^\theta)(\RR)$$ 
with $C(w^\theta)(\RR)\qQ(\RR)= P w^\theta Q$.

Moreover 
\bqn
C(w^{\theta}_{max}) (\RR)= Pw^{\theta}_{max}= Qw^{\theta}_{max}
\eqn
is open and dense.
All sets $C(w^\theta)(\RR)$ are topological cells.
\begin{lemma}
Define elements $t_k= s_k\cdots s_r \in W$.
The set of representatives $W^\theta$ is given by the elements $(t_1), (t_2 t_1), \dots, (t_r\cdots t_1)\in W$
\end{lemma}
\begin{proof}
The element $t_k$ is reduced and, since $s_r \notin W_\theta$, 
it is a representative of the coset $t_k W_\theta$.
The action of the element $t_k$ on the root system can be calculated 
using the Cartan martix of $C_r$, respectively $BC_r$ 
(for the form of the Cartan matrix see e.g. \cite{Bourbaki}):
$t_k(\xi_i)=\xi_i$ for $i<k$, 
$t_k(\xi_i)=\xi_{i+1}$ for $k\leq i\leq r-1$ and
$t_k(\xi_r)=-\xi_{k}$.
This implies that $W^\theta =\{(t_1),\dots, (t_r\cdots t_1)\}$ is a full set of representatives of smallest length.
\end{proof}

Denote by 
$a$ the multiplicity of roots $\pm \xi_i\pm\xi_j$ $(i\neq j)$ in
$\Xi$, and by $b$ the multiplicity of roots $\pm \xi_i$ in $\Xi$.
\begin{prop}
The cell $C(w^{\theta}_{max})(\RR)= C(t_r \cdots t_1)(\RR)$ is the
unique open and dense orbit of $Q$ in $\cs=G/Q$. 
It is of dimension $(2b+1)r+ {\frac{r(r-1)}{2}} a$. 
The codimension of $C(t_{r-k} \cdots t_1)(\RR)$ in $G/Q$ is $(2b+1)k+{\frac{k(k-1)}{2}} a$.
\end{prop}
\begin{proof}
The first statement follows from the general theory 
since $t_r\cdots t_1$ is a representative of the coset $w_{max} W_\theta$. 
We have to prove the statement about the dimension of the cells $C(t_{r-k} \cdots t_1)(\RR)$. 

Define $\Xi(\theta)$ to be the set of all 
nondivisible positive roots in $\Xi$, 
which do not lie in the span of $\theta$.
To every representative $w^\theta \in W^\theta$ 
we associate 
\bqn
\Xi_{w^\theta}= w^\theta(\Xi(\theta)) \cap \Xi^-.
\eqn
Let $\frakg_{w^\theta}$ be the sum of all root spaces 
$\frakg^\a+\frakg^{2\a}$, 
where $-\a \in \Xi_{w^\theta}$ 
(the root space $\frakg^{2\a}$ is zero, if $2\a$ is not a root), and 
denote by $U_{w^\theta}$ the corresponding unipotent Lie group.
The cell corresponding to $w^\theta$ is then parametrized by $U_{w^\theta}$, 
i.e 
\bqn
C(w^\theta)(\RR)= U_{w^\theta} w^\theta,  \quad
\dim C(w^\theta)(\RR)= \dim U_{w^\theta} w^\theta.
\eqn

The longest element $w^{\theta}_{max}= (t_r\cdots t_1)$ sends 
every positive root into a negative root.
Hence 
\bqn
\Xi_{w^{\theta}_{max}}= \Xi_{(t_r \cdots t_1)}
= \Xi(\theta)^-=\{-(\xi_i+\xi_j)\, (j>i)\,, -\xi_i\}
\eqn 
in the case that $\Xi$ is of type $BC_r$ or 
\bqn 
\Xi_{w^{\theta}_{max}}= \Xi(\theta)^-=\{-(\xi_i+\xi_j)\, (j>i)\,, -2\xi_i\}
\eqn
in the case that $\Xi$ is of type $C_r$, 
with $\dim \frakg_{\xi_i+\xi_j}= a$, 
$\dim \frakg_{2(\xi_i+\xi_j)}= 0$,
$\dim \frakg_{\xi_i}= 2b$ and
$\dim \frakg_{2\xi_i}= 1$.
This gives 
\bqn
\dim C(w_{max}^\theta)(\RR)= (2b+1)r+ {\frac{r(r-1)}{2}} a
\eqn

The above computations of the action of $t_k$ imply
$\Xi_{(t_r \cdots t_1)}-\Xi_{(t_{r-k} \cdots t_1)}= 
\{-\xi_r,\dots -\xi_{r-k+1}, -(\xi_i+\xi_j)\, (j>i\geq r-k+1)\}$. 
As above we have to exchange $-\xi_i$ by $-2\xi_i$ 
in the case that $\Xi$ is of type $C_r$. 
Since we have the same multiplicities as above,
the codimension of $C(t_{r-k} \cdots t_1)(\RR)$ in $G/Q$ is 
\bqn
\dim U_{w^\theta_{max}} - \dim U_{w^{t_{r-k}\cdots t_1}} =
(2b+1)k+{\frac{k(k-1)}{2}} a.
\eqn
\end{proof}

\begin{cor}\label{cor:bruhat}
The Hermitian symmetric space $X$ is of tube type if and only if
the codimension of $C(t_{r-1} \cdots t_1)(\RR)$ in $G/Q$ is $1$.\\
The Hermitian symmetric space $X$ is not of tube type 
if and only if
the codimension of $C(t_{r-1} \cdots t_1)(\RR)$ in $G/Q$ is $\geq3$
\end{cor}
\begin{proof}
By the above proposition $\codim C(t_{r-1} \cdots t_1)(\RR)= 2b+1$, 
where $b\geq 1$ iff $X$ is not of tube type. 
\end{proof}

\begin{cor}\label{cor:codim}
The codimension of proper subcells of 
$\overline{C(t_{r-1} \cdots t_1)}(\RR)$ is at least $2$.
\end{cor}

\begin{defi}
The open dense cell $C(w^\theta_{max})(\RR)$ is called the \emph{Bruhat cell} of $\cs\cong G/Q$ 
associated to $z=eQ$. It is denoted by $B(z)$.
\end{defi}

There is an equivariant version of the Bruhat decomposition describing the $G$-orbits in $G/Q\times G/Q$.
The unique open and dense $G$-orbit in $G/Q\times G/Q$ is called the \emph{equivariant Bruhat cell} and will 
be denoted by $B\subset G/Q\times G/Q$. Fixing a point $z\in G/Q$, the $G$-orbit $B$ is represented 
by $z\times B(z)$. 

The Bruhat decompositions of $G/P_\theta$ are compatible with the natural 
projections $G/P\to G/P_\theta\to G/P_{\theta^\prime}$.

\section{The space of triples of pairwise transverse points}

\begin{defi}\label{def:transversality}
Let $x,y\in G/Q$. We say that $x,y$ are \emph{transverse} if $(x,y)\in B\subset G/Q\times G/Q$.
\end{defi}

We denote the set of transverse pairs in $\cs$ by $\cs^{(2)}$. 
It follows from the definition that transversality is given by a Zariski 
open and dense
condition.

A Hermitian symmetric space $X$ which is not of tube type 
has many maximal subdomains $T$ of tube type of equal rank as $X$.
There is a relation of tranverse pairs in $\cs$ and maximal subdomains of tube type.
\begin{lemma}\label{lem:transverse_tube}
Every two transverse points $(x,y)\in\cs^{(2)}$ determine a unique maximal subdomain 
of tube type $T_{xy}$ containing $x$ and $y$ in its boundary
\end{lemma}
\begin{proof}
If the symmetric space $X$ is of tube type, the claim is obviously true, so we may assume, that $X$ is not of tube type.
Given $(x,y)\in \cs\times \cs = G/Q \times G/Q$ transverse, 
since $G$ acts transitive on $G/Q$, we may assume, that $x = eQ$. 
By transversality $y$ lies in the Bruhat cell $B(x)$ of $G/Q$, 
which is parametrized by the unipotent Lie group 
$U_{w_{max}}$,
with  Lie algebra given by all rootspaces corresponding to the roots  
$\{(\mu_i+\mu_j)\, (j>i),\, \mu_i,\, 2\mu_i\}$. 
The orbit of $y$ under the unipotent Lie group $ U_{T, w_{max}}$, whose Lie algebra is given by all root spaces corresponding to the roots 
$\{(\mu_i+\mu_j)\, (j>i),\, 2\mu_i\}$ is a subspace of $B(x)$. But since the Weyl groups of $BC_r$ and $C_r$ are the same, 
this subspace is the Bruhat cell $B_{T_{xy}}(x)$ of a maximal subdomain of tube type, 
determined by the root system $\{\pm(\mu_i\pm\mu_j)\, (j>i),\, \pm 2\mu_i\}$, hence it determines 
a maximal subdomain of tube type $T_{xy}$ containing $x$ and $y$ in its Shilov boundary.
Assume that there would be another maximal subdomain of tube type $T$ 
containing $x,y$ in its Shilov boundary $\cs_T=G_T/Q_T$.
By the same argument as above, $x=eQ_T$ and $y$ lies in the Bruhat cell $B_T(x)$ of $G_T/Q_T$
parametrized by $ U_{w_{max}}^T$.
This is contained in the Bruhat cell 
of $G/Q$ and therefore coincides with $U_{T, w_{max}}$, which is 
the big cell of $T_{xy}$. Hence $T=T_{xy}$.
\end{proof}

\begin{rem}
In the case of rank one symmetric spaces any two distinct points on the boundary 
are transverse and determine a 
unique geodesic. 
In particular, two transverse points on the boundary 
(which coincides with the Shilov boundary for Hermitian symmetric
spaces of rank one) 
determine a unique complex geodesic. 
In the case of higher rank, there are in general many geodesics and hence 
also many complex geodesics connecting two transverse points in the Shilov boundary.
Lemma \ref{lem:transverse_tube} motivates that maximal subdomains 
of tube type are an appropriate generalization of complex geodesics in rank $1$.
\end{rem}

The set of triples of pairwise transverse points on the Shilov
boundary of a Hermitian symmetric space $X$, 
\bqn
\cs^{(3)}= (G/Q)^{(3)}:=\{(x,y,z)\in (G/Q)^3\,
|\,(x,y)\in B, \, (y,z)\in B, \,(x,z)\in B\}, 
\eqn 
determines whether $x$ is of tube type or not.
From Corollary \ref{cor:codim} we obtain a new characterization of 
non tube type domains.

\begin{thm}\label{connected}
The symmetric space $X$ is not of tube type if and only if $\cs^{(3)}\cong (G/Q)^{(3)}$ is connected.
\end{thm}
\begin{proof}
We prove that if $X$ is not of tube type, then $\cs^{(3)}\cong (G/Q)^{(3)} $ is connected. 
The set $\cs^3\cong (G/Q)^3$ is connected. The codimension of the complement of the equivariant 
Bruhat cell $B$ is greater than $2$ since $X$ is not of tube type. 
Therefore the set $(G/Q)^{(3)}$ is obtained from $(G/Q)^3$ 
by removing finitely many subsets of codimension $\geq 3$, 
hence it stays connected.
It is well known that the set $\cs^{(3)}$ decomposes 
into $r_X +1$ open $G$-orbits if $X$ is of tube type (\cite{Clerc_Orsted_TG}). 
\end{proof}


\begin{lemma}\label{lem:transversality}
Given $x,y\in \cs\cong G/Q$, the following are equivalent:\\
(1) $y \in B(x)\subset G/Q$.\\
(2) $x \in B(y)\subset G/Q$.\\
(3) $(x,y)\in B\subset G/Q\times G/Q$.\\
(4) $h(x,y)\neq 0$, where $h$ is the invariant polynomial related to the Bergmann kernel\\
(5) There exists a geodesic in $\Dd$ joining $x$ and $y$.\\
(6) There exists preimages $x'$ and $y'$ of $x,y$ in the geometric boundary of
$X$ and a geodesic of type $\theta$ with $\theta\subset
\Theta_r$ in $X$ joining $x'$ to $y'$.
\end{lemma}
\begin{proof}
The equivalence of the first three assertions is clear from the above 
definitions. 
The equivalence of (5) and (6) follows from Lemma~\ref{geodesics}.
That (1) implies (6) is  a standard result for 
symmetric spaces of noncompact type (see e.g. \cite{Eberlein}). We may identify 
$G/Q$ with the space of equivalence classes of asymptotic Weyl chambers of type $Q$. Then  
$x\in B(y)$ if and only if there exist two opposite Weyl chambers of type $Q$ joining $x$ to $y$. 
For the 
Shilov boundary $Q$ is a maximal parabolic subgroup, hence Weyl chambers of type $Q$ 
are geodesic rays of type $\theta_r$, which form a geodesic connecting $x$ two $y$ if they are opposite.  
To prove the converse, (6) implies (1), assume that there is a geodesic $\g$ of type $\theta\subset \Theta_r$ in $X$ 
connecting $x$ to $y$. Then $\g$ determines two Weyl chambers of type $\theta$, opposite in $\g(0)$,  
connecting two points $x_\theta, y_\theta \in G/P_\theta$ 
lying in the Bruhat cell of the Bruhat decomposition of $G/P_\theta$, 
which are mapped to $x$ and $y$ under the natural projection $G/P_\theta \to G/P_{\theta_r} = G/Q$. 
Since the natural projections 
respect the Bruhat decompositions (3) follows. 
The implication from (4) to  (5) is shown in \cite{Clerc_Orsted_2}: if $x,y$ are transverse 
we may assume that $x=\sum_{i=1}^r E_i$ and $y=-\sum_{i=1}^r E_i$, hence the 
geodesic $\g(t)= \sum_{i=1}^r \tanh(t) E_i$ connects $x$ to $y$ in $\Dd$.
On the other hand, assume that there is a geodesic joining $x$ to $y$. Since transversality is a $G$-invariant notion we may assume that 
this geodesic passes through $0$ and that $x=\sum_{i=1}^r E_i$, but then (see Lemma \ref{geodesics} in Chapter 3) we have $y=-\sum_{i=1}^r E_i$, 
hence $h(x,y) = 2^r \neq 0$ and (5) implies (4).
\end{proof}

\section{A complex model for the Shilov boundary}
In order to work with the Zariski topology on the Shilov boundary, 
it is convenient to have an explicit description of a complex model
for the Shilov boundary.

Define
\bqn \gG_\CC= G_\CC\times G_\CC
\eqn with real structure  
\bqn
\gG=\{(g,\sigma(g))\,|\, g\in G_\CC\}.
\eqn 
The associated complex conjugation $\Sigma: \gG_\CC\to \gG_\CC$ 
is given by
\bqn
\Sigma(g_1, g_2)= (\sigma(g_2), \sigma(g_1)).
\eqn

Then $\gG_\CC$ is defined over $\RR$. The subgroup $\nN= K_\CC P_-\times K_\CC P_+$ is a complex algebraic subgroup of $\gG_\CC$; it is invariant 
under $\Sigma$ and hence also defined over $\RR$. In particular $\gG_\CC/\nN\cong M_+\times M_-$ is an algebraic variety defined over $\RR$. 
The homomorphism 
\bqn
\Delta: G_\CC\to \gG_\CC, \, \, g\mapsto (g, g)
\eqn
 is defined over $\RR$ and induces an algebraic action of $G_\CC$ on $\gG_\CC/\nN$.
Furthermore $\Delta(G) = Fix_\Sigma(\Delta(G_\CC))$.
Since $\sigma(c) = c^{-1}$, the point $(c,c^{-1}) \nN$ is a real point, $(c,c^{-1}) \nN \in (\gG_\CC/\nN)(\RR)$.

\begin{prop}\label{prop:shilov_complex}
The $G_\CC$-orbit (via $\Delta$) of $(c,c^{-1}) \nN$ is $G_\CC/\qQ$. It is a complex model of the 
Shilov boundary $G_\CC/\qQ(\RR)$ in $\gG_\CC/\nN(\RR)$. 
\end{prop}
\begin{proof}
The stabilizer of  $(c,c^{-1}) N$ in $\gG_\CC$ is
\bqn
\stab_{\gG_\CC}((c,c^{-1}) N) = cK_\CC P_+ c^{-1} \times c^{-1}K_\CC P_-c,
\eqn 
thus the stabilizer in 
$G_\CC$ of $(c,c^{-1}) N$ is $cK_\CC P_+ c^{-1} \cap c^{-1}K_\CC P_-c= \qQ$.
Hence by Proposition~\ref{prop:shilov_M} we get 
\bqn
\cs=G_\CC/\qQ(\RR)\subset \gG_\CC/\nN(\RR),
\eqn
thus
\bqn
\cs = G/(cK_\CC P_-c^{-1}\cap G) \to G/(cK_\CC P_-c^{-1}\cap G)\times
G/(c^{-1}K_\CC P_+c\cap G) \\
 G/(cK_\CC P_-c^{-1}\cap G)\times
G/(c^{-1}K_\CC P_+c\cap G) \subset \gG_\CC/\nN (\RR).
\eqn
\end{proof}

The Bruhat decomposition of $G_\CC/\qQ$ provides a notion of transversality in $G_\CC/\qQ$ which is compatible with the notion 
of transversality defined in Definition~\ref{def:transversality}, because the real points of the complex Bruhat cell are the real Bruhat cell. 
Thus if we define $(G_\CC/\qQ)^{(3)}$ to be the triples of pairwise transverse points in $G_\CC/\qQ$, we have $\cs^{(3)}= (G_\CC/\qQ)^{(3)} (\RR)$.

\section{Parametrizing the Shilov boundary of a tube type domain}
\subsection{The general  case}
A bounded symmetric domain of tube type is, via the Cayley transform, 
$\Tt= c\cdot:\Dd\to T_\Omega$, $c=\exp (\frac{\pi}{4} i \sum_{j=1}^r X_j)\in G_\CC$, biholomorphically equivalent 
to a tube domain
\bqn
T_\Omega:=\{u+iv\,|\, u\in V, v\in \Omega \subset V\},
\eqn
where $\Omega$ is a symmetric convex cone in the real vector space $V$. 
The map $\Tt$ extends to the part of the Shilov boundary consisting of all points $x\in B(o_r)\subset \cs$, 
which are transverse to $o_r \in\cs$. 
The Bruhat cell $B(o_r)$ is mapped to $V$, the complement will be send 
to $\infty$. Hence we may call $(V+i0)\cup \infty$ the Shilov boundary of $T_\Omega$.
The map $\Tt$ respects transversality in the sense that two points $x,y\in \cs$ are transverse 
if and only if its images under $\Tt$ satisfy $\det(\Tt(y)-\Tt(x))\neq 0$.

Given a point $z\in \cs$ we may parametrize a neighbourhood of this point in the following way.
Choose a point $z'\in \cs$ transverse to $z$. 
Then $z$ and a neighbourhood of $z$ in $\cs$ are contained 
in the Bruhat cell $B(z')$. Using the Cayley transformation, we 
identify $B(z')$ with the vector space $V$.
This gives a parametrization of the Bruhat cells. Choosing a cover of $\cs$ by 
$2^{r_X}$ Bruhat cells, we obtain a parametrization of $\cs$ which is
not unique, since it depends on the chosen point $z'$.
 
Once we consider the 
tangent space  of $\cs$ at the point $z$, the parametrization becomes 
independent of the chosen transverse point $z'$. 
We may identify the tangent space 
of $T_z \cs$ at a point $z\in B(z')$ with $V$ by mapping any point 
$w\in B(z') \cong V$ to the vector $w-z$.

The vector space $V$ is furthermore equipped 
with the structure of an Euclidean Jordan algebra 
(see \cite{Faraut_Koranyi} and the Appendix for a treatment of 
Jordan algebras associated to symmetric bounded domains).
Every element $x\in V$ has a spectral decomposition with respect to a so-called Jordan frame, 
that is a set of elements $(c_1, \dots, c_r) \in V$ with $c_i^2=c_i$, $c_i c_j=0$ if $i\neq j$ and 
$e=\sum_{j=1}^{r_X} c_j$ being the identity element of $V$ (see \cite[Theorem~III.1.2]{Faraut_Koranyi}): 
\bqn
x= \sum_{j=1}^{r_X} \lambda_j c_j 
=\sum_{j=1}^{k_+} \lambda_j c_j + \sum_{j=k_+ +1}^{k=k_+ + k_-} \lambda_j c_j + \sum_{j=k+1}^{r_X} \lambda_j c_j,
\eqn 
where $\lambda_j$ are uniquely determined by $x$ 
with $\lambda_j >0$ for $1\leq j\leq k_+$, $\lambda_j <0$ for $k_+ +1\leq j\leq k$ and $\lambda_j = 0$ 
for $j>k$. 
With respect to a spectral decomposition the Jordan algebra determinant of $x$ is $\det(x) = \Pi_{j=1}^r \lambda_j$.
There is a spectral norm on $V$ defined by 
\bqn
s(w) = \tr(w) = \sum_{i=1}^r \lambda_i
\eqn 
for all $w\in V$.
This spectral norm defines a norm on the tangent space of a Bruhat cell 
in $\cs$.

\subsection{The classical cases} 
For the classical domains associated to 
$G= {\rm Sp}(2n,\RR)$, ${\rm SU}(n,n)$, ${\rm SO^*}(4n)$, ${\rm SO}(2,n)$ we give a concrete 
parametrization of the Shilov boundary via quadratic forms and linear maps.
These quadratic forms correspond to the spectral norms in the general case.
They can also be used to define variations of the 
above defined Maslov index (see \cite{Robbin_Salamon, cappell_etal} 
for the usual Maslov index). All this is well known for symplectic case 
and related to Clerc's construction in \cite{Clerc_strict}.

\subsubsection{The Grassmannians}
Assume that $G$ is one of the following 
$G= {\rm Sp}(2n,\RR), {\rm SU}(n,n), {\rm SO^*}(4n)$. 
Let $\KK=\RR,\CC,\HH$ and 
$\iota:\KK \to \KK$ its natural conjugation, let $\sigma: \KK\to \KK$ be the involution 
$\sigma=-\iota$.
Let $V$ be a $2n$-dimensional (right) $\KK$-vector space, endowed with a 
$\sigma$-symmetric nondegenerate bilinear form $h$, i.e. $h(x,y)
=\sigma(h(y,x))$. Then $G={\rm SU}(V, h)$. Denote by $X_G$ the
associated symmetric space.

Denote by $\Lambda(V)$ the space of maximal $h$-isotropic subspaces of $V$, 
i.e. $\cs_{X_G} = \Lambda(V)$. 
Two such subspaces $W_+, W_- \in \L(V)$ are transverse 
if $W_+\cap W_-=\{0\}$.
Any two transverse subspaces  $W_+, W_- \in \L(V)$ span the vector space 
$V= W_+\oplus W_-$. 
Denote by 
\bqn
B(W_-):= \{W\in \L(v) |\, W\cap W_-=\{0\}\}
\eqn
the space of all 
maximal $h$-isotropic subspaces $W\subset V$, which are transverse to $W_-$.
For $W_+ \in B(W_-)$, 
we denote by $\Aa(W_+)$ the space of all $\sigma$-symmetric nondegenerate 
bilinear forms on $W_+$.
\begin{prop}\label{prop:diffeo}
There is a diffeomorphism $a:B(W_-) \to \Aa (W_+)$.
\end{prop}
\begin{proof}
Given $W\in B(W_-)$. In the coordinates given by $V=W_-\oplus W_+$ 
we can write $W = {\rm graph} T_W$, where $T_W\in {\rm Lin}(W_+ , W_-)$ is 
determined by $v+T_Wv\in W$ for all $v\in W_+$.
The map $T$ defines a bilinear form $A_W$ on $W_+$ by 
setting for all $v,w \in W_+$
\bqn
A_W(v,w):= h(v,Tw). 
\eqn

Using that $W_+, W, W_-$ are isotropic, we obtain
\bqn
A_W(v,w)&=& h(v,Tw)=h(v+Tv, w+Tw) -h(Tv, w)= -h(Tv, w)\\
 &=& -\sigma(h(w,Tv))
=-\sigma (A_W(w,v)) = \sigma (A_W(w,v)).
\eqn
The map 
\bqn
a:B(W_-) \to \Aa (W_+), \, W\mapsto A_W
\eqn
 is a diffeomorphism. 
\end{proof}
\begin{rem}
This diffeomorphism corresponds to the extension of the Cayley
transform $\Tt$ from the previous section 
to the Shilov boundary.
\end{rem}
The bilinear form $A_W$ defines a quadratic form $q^+_W: W_+\to \RR$ 
on $W_+$ by setting $q^+_W(v):=  A_W(v,v)$. 
The diffeomorphism $a:B(W_-) \to \Aa (W_+)$ depends on the 
choice of $W_-$. This dependence disappears if we consider the 
tangent space.

\begin{lemma}
Let $W\in \cs$. Then $T_W \cs$ is canonically diffeomorphic to 
$\Aa(W)$. 
\end{lemma}
\begin{proof}
The proof is a direct generalization of the argument in \cite{duistermaat} 
in the symplectic case. For the sake of completeness we repeat it.
Assume that $Z, \tilde{Z} \in B(W)$ and let $W'\in \cs$ be near $W\in \cs$. 
Then  $Z, \tilde{Z} \in B(W')$. We write $\tilde{Z}$ as graph of $F:Z\to W$, 

\bqn
\tilde{Z}=:\{ z+ Fz \, |\, z\in Z\}, 
\eqn
and $W'$ as graph of 
$T:W\to Z$ and $\tilde{T}:W\to \tilde{Z}$ respectively, i.e. 

\bqn
W'=\{ x+ Tx \, |\, x\in W\}
=\{ y+ \tilde{T}y \, |\, y\in W\}.
\eqn

Since $V=W\oplus Z$ we obtain from the equalities $x+Tx=y+\tilde{T}y$ and 
$\tilde{T}y=z+Fz$ by comparison that $z=Tx$ and $y=x-Fz$. 
Hence 
\bqn
x+Tx= (x-FTx)+ \tilde{T}(x-FTx).
\eqn
Taking the product with any
$w\in W$ we get 
\bqn
A_Z(x,w) &=& h(Tx,w)= h(\tilde{T}x,u) -h(\tilde{T}FTx, w)\\
&=&  
A_{\tilde{Z}}(x,w) -h(\tilde{T}FTx, w).
\eqn
But the second term on the right hand side vanishes of second order if 
$W'\to W$, which finishes the proof.
\end{proof}

\subsubsection{The Quadric}
Let $V$ be a $(n+2)$-dimensional vector space with $h$ 
a symmetric bilinear form of signature $(n,2)$ on it.
Let $G={\rm SO}(n,2)= {\rm SO}(V,h)$.
The symmetric space $X_G$ the space of all $2$-dimensional 
subspaces $W\subset V$ such that $h_{|_W}>0$.
The Shilov boundary of the associated 
symmetric space is the set of 
all $h$-isotropic lines in $V$.
Two lines $L, M \in \cs$ are transverse if 
$L\oplus M^\perp= M\oplus L^\perp =V$. 
This implies that the restriction of $h$ to 
$(L\oplus M)$ is nondegenerate 
and of signature $(1,1)$, whereas the 
restriction of $h$ to $(L\oplus M)^\perp$ 
is nondegenerate and of signature $(n-1,1)$.
We denote by $t$ the unique map $t:L\to M$ 
defined by the condition $h(l, tl)=1$ for all 
$l\in L$.

Let $M, L\in \cs$ transverse. 
\begin{lemma}
Suppose that $N\in \cs$ transverse to $L$, 
then 
\bqn
N=\{l+Tl-h(Tl, Tl)tl\, |\, l\in L\},
\eqn 
where $T:L\to (L\oplus M)^\perp$ is a linear map.
\end{lemma} 
\begin{proof}
Decompose an element $n\in N$ into its 
$(L, (L\oplus M)^\perp, M)$-parts, 
$n= l+Tl+atl$. 
Then since $N, L, M$ are isotropic and $Tl\in (L\oplus M)^\perp$, 
we get 
\bqn
0=h(n,n)=h(Tl, Tl) + a h(l,tl), 
\eqn
hence $a= -h(Tl, Tl)$.
\end{proof}

Clearly this parametrization depends on $M$, but 
going to the tangential space we get the following

\begin{lemma}
The tangential space $T_L \cs$ can be identified 
in a canonical way with ${\rm Lin}(\RR, \RR^{n-1,1})$
\end{lemma}

\begin{proof}
Assume $M, \tilde M \in \cs$ are transverse to $L\in \cs$. 
Then 
\bqn
\tilde M=\{m+Tm- h(Tm, Tm) tm\,|\, m\in M\},
\eqn 
where $t:M\to L$, with $h(m, tm)=1$ and $T:M\to (M\oplus L)^\perp$ a linear map.
Let now $L^\prime$ be an isotropic line near $L$ in $\cs$. 
Then $L^\prime$ is still transverse to $M, \tilde M$.
Hence 
\bqn
L^\prime=\{l+Al+h(Al, Al)tl \, |\, l\in L, t:L\to M, A:L\to (L\oplus M)^\perp\}
\eqn
and  
\bqn
L^\prime=\{x+\tilde A x+h(\tilde A x, \tilde Ax)tx \, |\, x\in L, t:L\to \tilde M, \tilde A:L\to (L\oplus \tilde M)^\perp\}.
\eqn 
We want to show that $A$ and $\tilde A$ coincide up to first order.
From 
\bqn
l+Al+h(Al, Al)tl&=&x+\tilde A x+h(\tilde A x, \tilde Ax)tx\\
h(\tilde A x, \tilde Ax)tx&=&m+Tm- h(Tm, Tm) tm
\eqn
we get by comparing to the components in the orthogonal sum decompositions
$V=M\oplus L\oplus (M\oplus L)^\perp$ and $V=\tilde M\oplus L \oplus (\tilde M\oplus L)^\perp$, 
that 
\bqn
x&=& l+ h(Tm, Tm)tl\\
tl&=& -h(Al, Al)tl - (\tilde A x)_M, 
\eqn
where $(\tilde A x)_M$ denotes the $M$-part with respect to the above decomposition.

Hence we obtain
\bqn
\tilde A x = \tilde A(l+ h(Tm, Tm)tl) = \tilde A(l+ h(Tm, Tm)(-h(Al, Al)tl - (\tilde A x)_M)).
\eqn

Taking the product $h(l^\prime,L)$ with $l^\prime \in L^\prime$, 
\bqn
h(Al,L)-h(\tilde A x, L)=\\
h\left( h(Al, Al)tl + \tilde A h(Tm, Tm)(-h(Al, Al)tl - (\tilde A
  x)_M) - h(\tilde A x, \tilde A x)tx, L \right).  
\eqn

But the right side tends to zero in second order as $L^\prime \to L$.
\end{proof}


\vskip1cm
\chapter*{Exkurs: A generalized Maslov index}
\section*{The classical Maslov index}
We want to sketch a generalization of the Maslov index which has been defined on
the space of Lagrangian subspaces  
to the Shilov boundaries of classical 
Hermitian symmetric spaces of tube type. 

Recall the definition of the usual Maslov cycle given by Arnold (\cite{Arnold1}).
Let $V$ be a $2n$ dimensional real vector space. Denote by 
$\Lambda(V)$ the space of Lagrangian subspaces of $V$, $\dim \Lambda(V)=N$. 
The Maslov cycle with respect to a Lagrangian subspace $L\in \Lambda(V)$ 
is defined by 
\bqn
\Sigma(L):=\{ W \in \Lambda(V)\, |\, \dim (W\cap L) >0 \}.
\eqn
The set $\Sigma(L)$ is the space of all Lagrangian subspaces in $  \Lambda(V)$ 
which are not transverse to $L$. It is an algebraic variety 
and a co-oriented codimension one cycle in $\Lambda(V)$, thus it  
determines an element $[\Sigma(L)] \in \h_{N-1}( \Lambda(V), \ZZ)$. 
The class of $\Sigma(L)$ does not depend on $L$.
The Maslov index $\mu \in \h^1( \Lambda(V), \ZZ)$ is defined to be the 
algebraic intersection number of a loop with $\Sigma(L)$.
One problem of this definition of the Maslov index is that given a path, 
not only a loop, in 
$\Lambda(V)$, the 
intersection number of this path with $\Sigma(L)$ is only well defined if the 
endpoints of this path are transverse to $L$.

There are other ways to define a Maslov index being defined for all paths, but
 depending on the choice of a basepoint (see e.g. \cite{Robbin_Salamon}).

\section*{The Maslov cycle in the general case}
We make use of the results about the Bruhat decomposition 
of $\cs$ to construct a generalized Maslov index for Shilov 
boundaries of tube type domains. 

\begin{prop*}
Let $X$ be an irreducible Hermitian symmetric space of tube type with 
Shilov boundary $\cs$. Then there exists a generalized 
Maslov index $\mu \in \h^1(\cs, \ZZ)_{\text{mod torsion}}$.
\end{prop*}

The Maslov index $\mu$ is constructed as in \cite{Arnold1} as intersection 
number of a loop with a codimension $1$ cycle $\Sigma$ in $\cs$. 
The existence of this cycle follows from the Bruhat decomposition of the 
Shilov boundary of $X$. 

Define $\Sigma(z):= \cs \backslash B(z)$. 
\begin{prop*}
The set $\Sigma(z)$ is an algebraic variety and 
a co-oriented codimension one cycle.
\end{prop*}
\begin{proof}
$\Sigma(z)$ is the disjoint union of cells in the Bruhat decomposition of 
$\cs$, giving it the structure of an algebraic variety.
By Corollary~\ref{cor:bruhat} $\Sigma(z)$ has codimension one and is a cycle. 
It inherits it co-orientation from the orientation from $\cs$ and $B(z)$.
\end{proof}

\begin{prop*}
For an element $\a\in \h_1(\cs, \ZZ)_{\text{mod torsion}}$, the intersection 
number 
$\a\cdot \Sigma(z)$ is well defined and independent of $z\in \cs$.
\end{prop*}
\begin{proof}
Since $\cs$ is connected and $\Sigma(z)$ is obtained from $\cs$ by removing a 
contractible set, the cycles $\Sigma(z)$ and $\Sigma(z')$ are homologous for 
$z,z'\in \cs$. Thus the intersection number does not depend on $z\in \cs$. 
\end{proof}

Given a path $\g$ in $\cs$ whose end points lie in $B(z)$, by 
perturbing $\g$ a bit, if necessary, it intersects 
$\Sigma(z)$ transversely in the smooth stratum 
and we may define the Maslov index of $\g$, 
$\mu(\g)$,  to be the geometric oriented intersection number of $\g$ 
with  $\Sigma(z)$.
The equivariant Bruhat cell $B$ gives rise to an equivariant version of 
the Maslov index, which is defined for paths in $\cs\times\cs$.

\subsection*{The classical case}
\begin{prop*}
Let $G= {\rm Sp}(2n,\RR)$.
The generalized Maslov index $\mu \in \h^1(\cs, \ZZ)$ coincides with the 
usual Maslov index $\tilde\mu\in \h^1(\Lambda(V), \ZZ)$. 
\end{prop*}

\begin{proof}
The Shilov boundary is the space of Lagrangian subspaces of $V$. 
The codimension one algebraic varieties defined above are the same.
\end{proof}

For $G= {\rm SU}(n,n), {\rm SO^*}(4n)$ the Shilov boundary 
is also given as the space of maximal isotropic subspaces. 
The Maslov cycle $\Sigma(L)$ consists again of all maximal isotropic subspaces 
having nonzero intersection with $L$. 
\begin{rem}
The parametrization of the Shilov boundary of a tube type domain 
by the vector space with Jordan algebra structure, allows to give a 
generalization of the construction of Robbin and Salamon \cite{Robbin_Salamon} 
of a Maslov index using crossing form. 
Since these things do not lie within the scope of this work, we do not extend 
this idea here. 
\end{rem}
%
%
%
%
%


\vskip1cm
\chapter{Cocycles}\label{sec:cocycle}
\section{The Bergmann and the K\"ahler cocycle}
\subsection{The Bergmann cocycle on an arbitrary bounded domain}
Let $\Dd \subset \CC^n$  be a bounded domain equipped with its 
Bergmann kernel 
\bqn
k: \Dd \times \Dd \to \CC.
\eqn 
Denote by $G=\aut(\Dd)$ the group of holomorphic automorphisms of $\Dd$.
Then the Bergmann kernel has the following invariance property: 
Let $g\in G$, then 
\bqn
k(gx,gy)= j(g,x)^{-1}k(x,y){\overline {j(g,y)}}^{-1},
\eqn
where $j(g,x)$ is the complex Jacobian of $g:\Dd\to\Dd \subset \CC^n$.
Let $\arg k(x,y)$ denote the unique continuous determination 
of the argument of $k$ to $\Dd^2$ such that $\arg k(x,x)=0$ for all
$x\in \Dd$. 
 
We define the \emph{Bergmann cocycle} as 
\bqn
\b(x,y,z):= \arg k(x,y) + \arg k(y,z)+ \arg k(z,x).
\eqn
This is a $G$-invariant continuous alternating cocycle on ${\Dd}^3$

The Bergmann kernel defines, via its K\"ahler form $\omega = \frac{1}{2} d d_\CC \log k(z,z)$, 
another $G$-invariant singular 
$C^1$-cochain given by integration over $C^1$-simplices:
Given such a simplex $\Delta_{x,y,z}$ with vertices $(x,y,z)\in \Dd^{3}$,
we define the \emph{K\"ahler cocycle} as
\bqn
c(\Delta_{x,y,z}) := \int_{\Delta_{x,y,z}} \omega
\eqn 
\begin{rem}
Since the form $\omega$ is exact, the value of $c(\Delta_{x,y,z})$ does not depend on 
the filling of the simplex.
\end{rem}

\begin{lemma}
The cocycles $\b(\Delta_{x,y,z}):= \b(x,y,z)$ and $c(\Delta_{x,y,z})$ define the same 
$G$-invariant singular cohomology class.
\end{lemma}
\begin{proof}
We have to show that $(c-\b)$ is the boundary of a $G$-invariant 
singular cochain. 
Let $\Delta _{x,y,z}$ be a $C^1$-simplex, with sides 
$\g_{xy}, \g_{yz}, \g_{zx}$. Then 
\bqn
\int_\Delta \omega = \frac{1}{2} (\int_{\g_{xy}} d_\CC \log k 
+  \int_{\g_{yz}} d_\CC \log k 
+ \int_{\g_{zx}} d_\CC \log k).
\eqn
Now $ \int_{ g\g_{xy}} d_\CC \log k =  \int_{\g_{xy}} d_\CC g^* \log k$,
but 
\bqn
g^* \log k(x,x) &=& \log k(gx, gx) = \log |j(g,x)|^{-2} k(x,x)\\
&=& \log |j(g,x)|^{-2}+ \log k(x,x).
\eqn
This implies that 
\bqn
\frac{1}{2}(\int_{g\g_{xy}} d_\CC \log k(z,z)& -& \int_{\g_{xy}} d_\CC \log k (z,z))
=\frac{1}{2} \int_{\g_{xy}} d_\CC \log |j(g,z)|^{-2}\\
&=&-\frac{1}{2}\int_{\g_{xy}} 2 d_\CC \Re \log (j(g,z)) \\
&=&-\frac{1}{2}\int_{\g_{xy}} 2 d \Im \log (j(g,z)) 
= (\arg j(g,x)- \arg j(g,y)),
\eqn
where $\arg j(g,x)$ being the continuous determination of the argument. 
On the other hand 
\bqn
\arg k(gx,gy) = -\arg j(g,x) + \arg k(x,y) + \arg j(g,y), 
\eqn
hence the cochain 
$$
\a(\g)=\frac{1}{2}\int_{\g} d_\CC \log k - \arg k(\g(0), \g(1)) 
$$
 is $G$- invariant. 
\end{proof}

\subsection{The Bergmann cocycle on bounded symmetric domains}
If $\Dd$ is a bounded symmetric domain, then the group $G$ acts transitively on $\Dd$ 
and we obtain the equality of the Bergmann cocycle and its associated K\"ahler cocycle 
on the cochain level.
Define $c(x,y,z):=c(\Delta_{x,y,z})$, where $\Delta_{x,y,z}$ is a
geodesic triangle with vertices in $x,y,z$ in $\Dd$.
\begin{lemma}\label{lem:bergmann_kaehler}
If $\Dd\subset \CC^n$ is a bounded symmetric domain then 
$\b(x,y,z) = c (x,y,z).$
\end{lemma}
\begin{proof}
From the above we know that
\bqn
\frac{1}{2}\int_{\g_{xy}} d_\CC \log k -\arg k(x,y)
\eqn
is $G$-invariant. Since $G$ acts transitively we may suppose that $x=0$, 
in which case $ d_\CC \log k =0$ along the geodesic $\g_{0y}$ 
and $\arg k(0,y)=0$. 
\end{proof}
%
\emph{From now on the Bergmann kernel on $\Dd$ will be normalized 
in such a way that the minimal holomorphic curvature of its associated metric is $-1$.}

We would like to extend the continuous cocycle $\b:\Dd^3\to \RR$ 
to a continuous cocycle on $\overline{\Dd}^3$. This is not possible since the 
Bergmann kernel $k$ does not extend continuously to $\overline{\Dd}^2$, but 
only to a part of it. 

Recall that $U' \subset \frakp_+\times \frakp_+$  was defined to be
the subset consisting of $z,w\in \frakp_+\times \frakp_+$ 
such that $\exp(-\sigma(w)) \exp(z)\in P_+ K_\CC P_-$. 
By Lemma~\ref{lem:kernel} in Chapter 3, $k$ extends continuously to $\overline{\Dd}^{(2)}:=\overline{\Dd}^{2}\cap U'$ 
and we may extend the Bergmann cocycle to a 
$G$-invariant continuous alternating cocycle 
\bqn
\beta:\overline{\Dd}^{(3)} \to \RR, 
\eqn
where $\overline{\Dd}^{(3)}:= \{ (x_1,x_2,x_3) \in \overline{\Dd}^3 
\,|\, (x_i,x_j) \in \overline{\Dd}^{(2)}, \forall i\neq j\}$. 
In \cite{Domic_Toledo,Clerc_Orsted_2} it is proven that this cocycle is bounded 
with
\begin{equation}\label{eq:bergmann_bounded}
\max_{(x,y,z)\in {\overline\Dd}^{(3)}} |\beta(x,y,z)|= \pi r_X.
\end{equation}

The properties of triples $(x,y,z) \in \overline{\Dd}^{(3)}$ realizing this maximum are studied in \cite{Clerc_Orsted_2}.
\begin{prop}(\cite{Clerc_Orsted_2})\label{results_c_o}
Let $(x,y,z) \in \overline\Dd^{(3)}$. Assume that $\beta(x,y,z) = \pi r_X$. 
Then
\begin{enumerate}
\item{$x,y,z$ lie on the Shilov boundary of a maximal subdomain of tube type
$T_{xyz}\subset \Dd$.}
\item{The stabilizer $\Stab_G(x,y,z)$ is compact. 
In fact, the set of fixed points of $\Stab_G(x,y,z)$ is a 
holomorphically tightly embedded disc containing $x,y,z$ on its boundary.}
\end{enumerate}
\end{prop}
The Proposition implies in particular that the maximum is obtained on the 
Shilov boundary.

The extension of the Bergmann cocycle is defined on the space of triples of pairwise transverse points in the Shilov boundary 
$\cs^{(3)}$:
\begin{lemma}\label{lem:transversality_shilov}
Let $x,y\in G/Q$ be transverse. Then $(\Phi(x), \Phi(y))\in \overline{\Dd}^{(2)}\cap \cs^2$.
\end{lemma}
\begin{proof}
By definition we have that $(x,y)\in B\subset G/Q\times G/Q$. 
Since $G$ acts transitively on $B$, 
we may assume that $(x,y)= (c(x_0), c^{-1}(x_0))$. 
Thus 
\bqn
\Phi(x) &=& \xi^{-1}\xi(io^+_r)= io^+_r,\\
\Phi(y) &=& \xi^{-1} \xi(-io^+_r)=-io^+_r
\eqn
and 
$(\Phi(x), \sigma(\Phi(y))= (io^+_r, io^-_r)$.

But 
\bqn
\exp(-io^-_r)\exp(io^+_r)= \exp(\frac{1}{2}i o^+_r)k \exp(-\frac{1}{2}io^-_r) \subset P_+ K_\CC P_-,
\eqn 
with $k$ 
being $\exp(4i\sum_{j=1}^{r} H_i)$. Hence $(\Phi(x), \Phi(y))= (io^+_r, -io^+_r)\subset U'\cap \cs^2$.
%
\end{proof}

The restriction of the Bergmann cocycle 
to the Shilov boundary defines a bounded 
$G$-invariant continuous alternating 
cocycle $\b:\cs^{(3)} \to \RR$ on the space of pairwise 
transverse triples in $\cs$. 
From now on \emph{Bergmann cocycle} will refer to this restriction of the Bergmann cocycle to the Shilov boundary.

From Theorem~\ref{connected} in Chapter 4 
we deduce that the Bergmann cocycle behaves differently with 
respect to tube type and non tube type symmetric spaces.

\begin{prop}\label{infinite_values}
If $\Dd$ is not of tube type, then $\b:\cs^{(3)} \to \RR$ takes 
infinitely many values. 
More precisely $\b(\cs^{(3)}) = [-\pi r_X, \pi r_X]$.\\
If $\Dd$ is of tube type, the cocycle $\b$ takes only finitely many values on $\cs^{(3)}$, 
namely $\pi(r_X -2k),$ where $k=0,\dots,r_X$. 
\end{prop}
\begin{proof}
The statement about bounded symmetric domains of tube type is well known \cite{Clerc_Orsted_TG}, 
therefore we give only a proof for the statement about domains, which are not of tube type.
The Bergmann cocycle $\b$ is bounded by $\pm\pi r_X$ and assumes the values $\pm\pi r_X$ for 
three distinct points lying on the boundary of a holomorphic tight
disc in $X$. 
Since $\Dd$ is not of tube type, 
Theorem~\ref{connected} in Chapter 4 implies that $\cs^{(3)}$ is
connected, so its image under the continuous map is a connected set 
containing $\pm\pi r_X$., hence it is $[-\pi r_X, \pi r_X]$.
\end{proof}

\section{The K\"ahler cocycle as bounded cohomology class}
The $G$-invariant normalized K\"ahler form $\omega\in \Omega^2(X)^G$ defines 
a $G$-invariant cocycle $c_G: G^3\to \RR$ by integration over geodesic simplices:
\bqn
c_G (g_0, g_1, g_2) := \int_{\Delta(g_0 x,g_1 x, g_2x)} \omega, 
\eqn
where $x\in X$ is arbitrary and $\Delta(g_0 x,g_1 x, g_2x)$ is the geodesic triangle spanned by 
$(g_0 x,g_1 x, g_2 x)$. 

By Lemma~\ref{lem:bergmann_kaehler} we have $c_G (g_0, g_1, g_2)=\b (g_0 x,g_1 x, g_2x)$. 
In particular 
\bqn
\sup_{g\in G^3} |c_G(g)| = \pi r_X
\eqn 
and thus $c_G\in C_b(G^3)^G$ is a bounded cocycle 
and defines a class $k_X^b = [c_G]_b\in \hcb^2(G)$, called the {\em bounded K\"ahler class}.
with $||k_X^b||=\sup_{g\in G^3} |c_G(g)|$.
Under the isomorphism $\hcb^2(G) \cong \hc^2(G)\cong \Omega^2(X)^G$, 
the bounded K\"ahler class $k^b_X$ corresponds to $\omega$.

\section{The Bergmann cocycle as bounded cohomology class}
The Shilov boundary is a separable compact metrizable space with a continuous $G$-action. 
If the Bergmann cocycle would be defined not only on the space of transverse triples, but on the whole 
$\cs^3$, we could apply Proposition~\ref{prop:borel} of Chapter 1 to conclude that 
$\b$ defines a bounded cohomology class in $\hcb^2(G)$. 
Even though the Bergmann cocycle does not extend to $\cs^3$, 
we can surpass these difficulties using the arguments of 
 \cite[\S~5]{Burger_Iozzi_supq}.

Let 
\bqn
\cs^{(n)}:=\big\{(x_1,x_2,\dots,x_n)\in\cs : x_i, x_j\text{ are 
transverse for all }i\neq j\big\}.
\eqn

Denote by $d_n$ the naturally defined homogeneous coboundary operator. 
Then the complex $(\binfty(\cs^{(n)}),d_n)$ of 
bounded alternating Borel functions on $\cs^{(n)}$, 
endowed with the supremum norm, is a strong $G$-resolution of $\RR$.
Since $G$ acts transitively on $\cs^{(2)}$ there are no $G$-invariant alternating 
Borel functions on $\cs^{(2)}$.  
The Bergman cocycle 
\bqn
\b:\cs^{(3)} \to \RR
\eqn
is an element of 
$\binfty(\cs^{(3)})$.
Applying \cite[Proposition~1.5.2]{Burger_Iozzi_App} there is a 
canonical map 
\bqn
\h^* (\binfty(\cs^{(*+1)})^G) \to \hcb^*(G);
\eqn 
Under this canonical map
the class $[\beta_\Dd]$ corresponds to 
the bounded K\"ahler class $k^b_X\in \hcb^2(G)$. 

%
%
%


\section{The Hermitian triple product}
On a Zariski open subset of the compact dual 
Hermitian symmetric space $M_+$ we define a multiplicative cocycle, 
which we call the \emph{Hermitian triple product}
since it is a generalization of the Hermitian triple product for $G={\rm SU}(p,q)$ defined in \cite{Goldman_book} and \cite{Burger_Iozzi_supq}. 
We relate the restriction of the Hermitian triple product to the
Shilov boundary to the Bergmann cocycle. 
Due to the algebraic nature of the Hermitian triple product we may
later apply arguments in the Zariski topology as in \cite{Burger_Iozzi_supq}. 
We work with the complexification $\gG_\CC$ of $G_\CC$ and the complex model for the Shilov boundary.

Recall that the complexified Bergmann kernel,  $k^\CC:U\to \CC$, was defined on an open subset $U\subset \frakp_+\times \frakp_-$.
We define $U^{(2)}\subset (\frakp_+\times \frakp_-)^2$, 
\bqn
U^{(2)}=\{((z_1, w_1),(z_2, w_2)) \in (\frakp_+\times \frakp_-)^2 \, |\, (z_1, w_2)\in U, \, (z_2, w_1)\in U\}.
\eqn 
Furthermore define 
\bqn
U^{(3)}=\{(z_i, w_i)_{i=1,\dots, 3}\in (\frakp_+\times \frakp_-)^3\,|\, ((z_i, w_i),(z_j, w_j))\in U^{(2)}\, \, \forall i\neq j\}.
\eqn
Note that we also normalize the complexified Bergmann kernel in such a
way that $k(z,w)= k^\CC(z, \sigma(w))$ is the normalized Bergmann kernel. 

Define $A=\CC\times\CC$ with real structure given by $\Delta_A: \CC\to A$, $a\to (a, \overline{a})$.
There is a diagonal action by $\CC$ and a complex conjugation defined on $A$. 
\bqn
\Sigma: A\to A, \, \, (a,b)\mapsto(\overline{b}, \overline{a}).
\eqn

We define a triple product 
\bqn
E:U^{(3)}\to A=\CC\times\CC
\eqn
by
\bqn
E\left( (z_1, w_1), (z_2, w_2), (z_3, w_3)\right) &:=& \\
((k^\CC(z_1, w_2)k^\CC(z_2, w_3)k^\CC(z_3, w_1) &,& k^\CC(z_2, w_1) k^\CC(z_3, w_2)k^\CC(z_1, w_3))
\eqn
\begin{lemma}\label{lem:triple_product}
$E$ satisfies the following relations: 
\begin{enumerate}
\item{Within the domain of definition: \\
$E(\Delta(g)(z_i, w_i)_{i=1,\dots, 3})= \l E((z_i, w_i)_{i=1,\dots,
  3})$ with $g\in G_\CC$, $\l \in \CC^\times$ depending on $g$ and
$z_i, w_i$}
\item{$E\left( (\sigma(z_i), \sigma(w_i))_{i=1,\dots, 3})= \Sigma
      (E((z_i, w_i)_{i=1,\dots, 3})\right) $}
\item{$E\left( (z_i, \sigma(z_i))_{i=1,\dots, 3}\right) = (e^{i\beta(z_1, z_2, z_3)},e^{-i\beta(z_1, z_2, z_3)})$ for all $z_i\in\cs\subset \frakp_+$}. 
\end{enumerate}
\end{lemma}
\begin{proof}
All assertions follow by direct computations using the properties of $k^\CC$ from Lemma~\ref{lem:kernel_prop} in Chapter 3 and the definition of $\b$.
\end{proof}

Define $\uU^{(2)}\subset  (\gG_\CC/\nN)^2$ to be the image of $U^{(2)}$ under the map $\xi\times \xi$, with 
\bqn
\xi=(\xi_+, \xi_-): \frakp_+\times \frakp_-\to \gG_\CC/\nN
\eqn 
and $\uU^{(3)} \subset  (\gG_\CC/\nN)^3$ the image of $U^{(3)}$ under $\xi_3:=\xi\times \xi\times \xi$.

By precomposition with $\xi_3^{-1}: \uU^{(3)} \to U^{(3)}$, the triple product on $U^{(3)}$ defines the Hermitian triple product on 
the subset $\uU^{(3)}\subset (\gG_\CC/\nN)^3$:
\bqn
\HTP_\CC:= E\circ\xi_3^{-1}:\uU^{(3)}\to A^\times/\CC^\times
\eqn
\begin{lemma}
The Hermitian triple product $\HTP_\CC$ is a $\gG_\CC$-invariant multiplicative cocycle. 
\end{lemma}
\begin{proof}
The properties of Lemma~\ref{lem:triple_product} imply the $G_\CC$-invariance of $\HTP_\CC$. 
The cocycle condition follows from the corresponding properties for the 
complexified Bergmann kernel.
\end{proof}

\begin{lemma}
The restriction of the Hermitian triple product 
${\HTP_\CC}|_{\cs^{(3)}}$ is well-defined and takes values in 
$\CC^\times/\RR_+^\times$. 
Furthermore 
\bqn
{\HTP_\CC}|_{\cs^{(3)}}(z_1, z_2, z_3) = e^{i\b(\Phi(z_1), \Phi(z_2), \Phi(z_3))}.
\eqn
%
%
\end{lemma}
\begin{proof}
We realize $\cs$ in $\gG_\CC/\nN$ as in 
Proposition~\ref{prop:shilov_complex} in Chapter 4. Then 
Lemma~\ref{lem:transversality} in Chapter 4 
implies that $\cs^{(3)} \subset \uU^{(3)}$ and 
the Hermitian triple product is well defined on $\cs^{(3)}$. By
property (2) of Lemma~\ref{lem:triple_product} it takes values in the
real points of $A^\times/\CC^\times$, 
which is $\CC^\times/\RR_+^\times$. 
Furthermore with $x_i= \Phi(z_i)$
\bqn
{\HTP_\CC}|_{\cs^{(3)}}(z_1, z_2, z_3) &=&\\
 E((x_1,\sigma(x_1)),( x_2,\sigma(x_2)), (x_3, \sigma(x_3))) &=&\\ 
(e^{i\beta(x_1, x_2, x_3)},e^{-i\beta(x_1, x_2, x_3)})&=& \Delta_A
\left( e^{i\beta(x_1, x_2, x_3)}\right) .
\eqn
\end{proof}

\subsection{Tube type and not tube type}
For $(z_1,z_2)\in G_\CC/\qQ^{(2)}\subset G_\CC/\qQ^2 \cap\uU^{(2)}$, define $\Bb_{12}\subset G_\CC/\qQ $ to be the connected subset
of those elements $z\in G_\CC/\qQ$ such that $(z_1, z_2, z)\in (G_\CC/\qQ)^{(3)}$.
\begin{lemma}
The set $\Bb_{12} \subset  G_\CC/\qQ$ is Zariski open.
The function
\begin{equation*}
P_{12}:\Bb_{12}\to A^\times/\CC^\times
\end{equation*}
given by $P_{12}(z_3):=\HTP_\CC(z_1,z_2,z_3)$ is regular.
\end{lemma}

\begin{lemma}\label{lem:nonconstant}
If $X$ is not of tube type, then for every $l\in\ZZ$, $l\neq0$
$P_{12}^l$ is not constant.
\end{lemma}
\noindent
Here $P_{12}^l(z)=\big(P_{12}(z)\big)^l$, 
where the product is taken in  $A^\times/\CC^\times$.
\begin{proof} 
We have 
\bqn
P_{12}\big(\Bb_{12}(\RR)\big)=
 \operatorname{Image}\big(\HTP_\CC|_{\cs^{(3)}}\big)=
 e^{i\operatorname{Image}(\b|_{\cs^{(3)}})}.
 \eqn
When $X$ is not of tube type the latter is infinite by Proposition~\ref{infinite_values}.
\end{proof}
 
\section{A generalized Maslov cocycle for tube type domains}
We come back to the problem of extending the Bergmann cocycle to a
cocycle defined on all triples, not only on the subset of triples of
pairwise transverse points.

Consider the case where $\Dd=\DD$ and $\cs=S^1$. 
Then $\cs^{(3)}$ is the set of triples of pairwise distinct points  
on $S^1$ and the extension of $\b$ to $(S^1)^{(3)}$ coincides up to a factor $\pi$ 
with the orientation cocycle $c:(S^1)^3\to \RR$, which assigns $\pm 1$ to a positive (resp. negative) oriented 
triple and $0$ to degenerate triples. Hence the cocycle $\b$ can be extended by $\pi c$ 
to a strict cocycle on $(S^1)^3$. Clearly this extension is not
continuous.

For general symmetric domains of tube type there is an extension of the Bergmann cocycle to strict 
cocycle on $\cs^3$  due to Clerc (\cite{Clerc_maslov_tube}). 
Recall that for the symplectic group the Shilov boundary is the space of Lagrangian subspaces 
of a symplectic vector space $V$ with symplectic form $h$. 
The well known (see e.g. \cite{Lion_Vergne}) Maslov cocycle $\tau:\cs^3\to \RR$, 
$\tau(L_1, L_2, L_3)= \sign(B)$, is defined to be the signature of the quadratic 
form $B$ on $L_1\oplus L_2\oplus L_3$ defined by 
$B(l_1, l_2, l_3)= h(l_1, l_2) + h(l_2, l_3)+ h(l_3, l_1)$.

This definition has been generalized to all classical domains 
of tube type \cite{Clerc_strict}.
For a general bounded symmetric domain of tube type Clerc \cite{Clerc_maslov_tube} 
constructed an extension of $\b$ to a strict cocycle on $\cs^3$ using an approach via  
Jordan algebras. He associates to every point $x\in \cs$ a family of curves $\Cc_x$, 
called ``radial'' in $x$, with $g\Cc_x=\Cc_{gx}$ for all $g\in G$. The extension of the argument of the 
Bergmann kernel along curves in $\Cc_x$ 
gives a well-defined function on $\cs^2$. Thus he can extend $\b$ to a $G$-invariant alternating cocycle on $\cs^3$, 
satisfying also $||\b||_\infty= \pi r_x$.
For classical domains the so defined Maslov cocycle $\tau: \cs^3\to\RR$, 
$\tau=\frac{1}{\pi} \b$ coincides with the Maslov cocycle defined via quadratic forms 
and on transverse triples with the Maslov cocycle defined in \cite{Clerc_Orsted_TG}.

We will need a slight variation of the definition of Clerc. 
Let 
\bqn
\cs^{[3]}:= \{(x_1,x_2,x_3)\in \cs^3\, |\, 
x_i \text{ is transverse to } x_j, x_k \text{ for some } i\}
\eqn 
be the space of 
triples, 
where one point is transverse to the other two. We will define a $G$-invariant 
function on $\cs^{[3]}$ whose extension to $\cs^3$, via the cocycle identity, 
coincides with Clerc's Maslov cocycle $\tau$ coming from the extension of $\b$ 
constructed in  \cite{Clerc_maslov_tube}. 
With this variation we prove in particular 
\begin{lemma}\label{strict}
Let $(x_1,x_2,x_3)\in \cs^{[3]}$. If $\beta(x_1,x_2,x_3) = r_X\pi$, then 
$(x_1,x_2,x_3)$ are pairwise transverse.
\end{lemma}
\subsection{The Maslov cocycle for general tube type domains}
Recall that a bounded symmetric domain of tube type is
biholomorphically equivalent, say via $\Tt:\Dd\to T_\Omega$, to a tube domain
$T_\Omega:=\{u+iv\,|\, u\in V, v\in \Omega \subset V\},$ 
where $\Omega$ is a symmetric convex cone in the real vector space $V$. 

We can now define a Maslov cocycle for triples $(x_1,x_2,x_3)\in \cs^{[3]}$. 
Using the transitivity of the $G$-action on $\cs$, 
we may assume that $\Tt(x_3) = \infty$ and 
$y_1=\Tt(x_1), y_2=\Tt(x_2) \in V$, 
then the \emph{Maslov cocycle} is defined to be
\bqn
\tau(x_1, x_2, x_3):= (k_+(y_2-y_1))- k_-(y_2-y_1)),
\eqn
where $k_\pm$ are the numbers of positive respectively 
negative eigenvalues in the spectral decomposition of $(y_2-y_1)$ with respect to a Jordan frame $(c_j)_{j=1,\dots,r}$.

Given any triple $x_1,x_2, x_3\in \cs$, using the operation of $G$ we may assume 
that $x_1,x_2, x_3$ are transverse to $\infty$, hence 
$y_1,y_2,y_3$ in $V$.
The cocycle identity allows us to define $\tau$ on $\cs^3$ as 
\bqn
\tau(x_1,x_2,x_3) &=& \tau (x_1, x_2,\infty) +\tau(x_2,x_3, \infty) -\tau(x_1, x_3, \infty)\\
&=&(k_+(y_2-y_1)- k_-(y_2-y_1))\\
&+&(k_+(y_3-y_2)- k_-(y_3-y_2))\\
&-&(k_+(y_3-y_1)- k_-(y_3-y_1)).
\eqn
This formula for $\tau$ coincides with the formula for the Maslov cocycle  given in 
\cite{Clerc_maslov_tube}. Hence 
we get 

\begin{cor}
Assume that $\Dd$ is of tube type. 
The Bergmann cocycle on  $\cs^{(3)}$ coincides with $\pi \tau$, where $\tau$ 
is the Maslov cocycle defined above. 
\end{cor}
\begin{proof}
In \cite{Clerc_maslov_tube} it is shown that the Maslov cocycle defined there coincides with the Bergmann 
cocycle on triples of pairwise transverse points. 
\end{proof}

We collect some properties of the Maslov cocycle $\tau$.
\begin{lemma}\label{lem:ordering}
If $\tau (x_1, x_2,x_3)= r$, then $k_+(y_2-y_1)=r$.
\end{lemma}
\begin{lemma}\label{lem:max_trans}
Let $(x_1,x_2, x_3)$ be as above with $\tau(x_1,x_2, x_3)= r$, then 
$x_1$ and $x_2$ are transverse.
\end{lemma}
\begin{proof}
From maximality we get $r= k_+(y_2-y_1)-k_-(y_2-y_1)$, where $0\leq k_\pm\leq r$ are the numbers of positive respectively 
negative eigenvalues in the spectral decomposition of $(y_2-y_1)$ with respect to a Jordan frame $(c_j)_{j=1,\dots,r_X}$.
This implies $k_+(y_2-y_1)= r$. Therefore $(y_2-y_1)$ has a spectral decomposition of the form 
$(y_2-y_1)= \sum_{j=1}^r \lambda_j c_j$, where $\lambda_j>0$ for all $j$, hence $\det(y_2-y_1) = \Pi_{j=1}^r \lambda_j>0$, 
in particular $\det(y_2-y_1)\neq 0$ and $x_2$ is transverse to $x_1$.
\end{proof}

From the ordering induced by the Maslov cocycle (see Lemma~\ref{lem:ordering}) it follows
\begin{lemma}\label{lem:converge}
Choose $f_+, f_- \in \cs$ transverse. 
Let $(x_i), (y_i)$ be sequences in $\cs$, 
with $x_i, y_i$ transverse to $f_-$ for all $i$ and 
such that $\tau(f_+, x_i, f_-)$, $\tau(f_+, y_i, f_-)$ and $\tau(y_i, x_i, f_-)$ 
are maximal for all $i$. 
Assume that $x_i \to f_+$ for $i\to \infty$, then 
$y_i \to f_+$ for $i\to \infty$. 
\end{lemma}
\begin{proof}
Indeed, if $x_i \to f_+$, then $\Tt(x_i)\to\Tt(f_+) $. 
By Lemma~\ref{lem:ordering} $k_+(\Tt(x_i) - \Tt(y_i))=r$, 
hence $k_+(\Tt(y_i) - \Tt(f_+)) \leq k_+(\Tt(x_i) - \Tt(f_+)) \to 0$. 
Therefore, since $k_+(\Tt(y_i) - \Tt(f_+))=r$ for all $i$, we get 
$\Tt(y_i)\to \Tt(f_+)$ by the 
uniqueness of the spectral decomposition.
\end{proof}

\begin{rem}
If we assume furthermore that $y_1=0$ we can define $\tau(x_1,x_2,x_3)= k_+(y_2)-k_-(y_2)$. 
\end{rem}
\begin{cor}
The cocycle $\beta$ extends to a strict cocycle on $\cs^{3}$ with 
\begin{enumerate}
\item{$|\beta (x_1,x_2,x_3)|\leq r_X\pi$ for all $x_1,x_2,x_3 \in \cs$}
\item{if $x_3$ is transverse to $x_1$ and $x_2$ and $\beta(x_1,x_2,x_3) = r_X\pi$, then 
$x_1$ and $x_2$ are transverse.}
\end{enumerate}
\end{cor}

\subsection{Maslov cocycle via quadratic forms}\label{sec:maslov_classic}
For the classical domains associated to 
$G= {\rm SU}(V,h)={\rm Sp}(2n,\RR), {\rm SU}(n,n), {\rm SO^*}(4n)$ we give a concrete 
construction of 
the above defined Maslov cocycle via the signature of quadratic forms. 
These quadratic forms constructed 
here are used in section \ref{sec:anosov} to construct metrics on appropriate vector bundles.
Denote by $\Lambda(V)= \cs$ the Shilov boundary, which is isomorphic to the space of maximal $h$-isotropic subspaces of $V$. 
Recall that two such subspaces $W_+, W_- \in \L(V)$ are transverse, 
if $W_+\cap W_-=\{0\}$, in particular $V= W_+\oplus W_-$. 
Denoting by $B(W_-):= \{W\in \L(v) |\, W\cap W_-=\{0\}\}$ the space of all 
maximal $h$ isotropic subspaces $W\subset V$ that are transverse to
$W_-$ 
and by $\Aa(W_+)$, $W_+ \in B(W_-)$, the space of all $\sigma$-symmetric nondegenerate 
bilinear forms on $W_+$. We have constructed a diffeomorphism
$a:B(W_-) \to \Aa (W_+)$ as follows (see Proposition~\ref{prop:diffeo} in Chapter 4): 
Given $W\in B(W_-)$. In the coordinates given by $V=W_-\oplus W_+$ 
we write $W = {\rm graph} T_W$, where $T_W\in {\rm Lin}(W_+ , W_-)$ is 
determined by $v+T_Wv\in W$ for all $v\in W_+$.
The map $T$ defines a bilinear form $A_W$ on $W_+$ for all $v,w \in W_+$
\bqn
A_W(v,w):= h(v,Tw),
\eqn
for all $v,w \in W_+$.
The bilinear form $A_W$ defines a quadratic form $q^+_W: W_+\to \RR$ 
on $W_+$ 
by setting $q^+_W(v):=  A_W(v,v)$. 
\begin{defi}
Given $W_-\in \L(V)$ and $W_+, W\in B( W_-)$, 
define the Maslov cocycle as $\tau (W_+, W,W_-) := \sign (q^+_W)$. 
\end{defi}
\begin{lemma}
If $(W_+, W, W_-) \in \cs^{(3)}$  
then $\tau(W_+, W, W_-) = \frac{1}{\pi} \b (W_+, W, W_-)$.
\end{lemma}
\begin{proof}
Direct computations show that $\tau$ coincides with the Maslov cocycle on 
transverse triples defined in \cite{Clerc_Orsted_TG}, which coincides with 
$\frac{1}{\pi} \beta$.
\end{proof}
\begin{lemma}\label{lem:max_trans_classical}
If $\tau(W_+, W, W_-) =\max_{\cs^{(3)}} |\tau|$, 
then $(W_+, W, W_-)\in \cs^{(3)}$.
\end{lemma}
\begin{proof}
We argue by contradiction. Assume that $W_+, W$ are not transverse, 
then $\dim (W_+\cap W)\geq 1$. Take $0\neq v\in W_+\cap W$. 
Then $q^+_W(v) = A_W(v,v)  = h(v,Tv)= h(v,0)= 0$. Hence $q^+_W$ is degenerate 
and $\sign(q^+_W)< n$.
\end{proof}
\begin{lemma}\label{compare}
Given $W_+, U, W \in B( W_-)$. 
Suppose that $\tau(U,W,L_-)$ is maximal, then $q^+_U\leq q^+_W$.
\end{lemma}
\begin{proof}
By construction  
\bqn
q^+_W (l) - q^+_U(l) 
&=&h( l, T_W(l)) - h(l, T_U(l))
=h(l, T_W(l) - T_U(l))\\
&=&h( l+T_U(l), T_W(l) -T_U(l)) = q(l).
\eqn
But now $l+T_U(l)\in U$ and the sum of the two arguments 
$l+T_U(l) + T_W(l)-T_U(l)= l + T_W(L) \in {\rm graph}(T_W) =W$, 
hence $q(l)= q^U_W(l)$ is the quadratic form on $U$, which is 
associatied to $W$ via the map $U\to W_-$. By assumption this 
form is positive definite, i.e. $q(l) >0$.
\end{proof}
\begin{lemma}\label{lem:converge_classical}
Let $N_i \in \Lambda_0(L_-)$, $M_i \in \Lambda_0(L_-)$ be sequences 
of maximal isotropic 
subspaces 
such that $\tau (L_+, N_i, L_-)$, $\tau (L_+, M_i, L_-)$ and 
$\tau (M_i, N_i, L_-)$ are maximal for all $i$. 
If $N_i \to L_+$ for $i\to \infty$, then 
$M_i \to L_+$ for $i\to \infty$. 
\end{lemma}
\begin{proof}
Indeed, if $N_i \to L_+$, then $q^+_{N_i} \to q_L = 0$. 
By the above Lemma $q^+_{M_i} \leq q^+_{N_i}$, 
hence  $q^+_{M_i} \to q_L = 0$ and therefore $M_i\to L$.
\end{proof}
\begin{rem}
From Lemma~\ref{compare} it 
follows that we could have defined the Maslov cocycle in a different way, 
by fixing $W_-$ and $W_+$ and associating to any $U,V\in B(W_-)$ the quadratic forms $q^+_U, q^+_V$ on $W_+$.
Then $\tau(U,V,W_-)= \sign (q^+_V - q^+_U)$. This coincides with the definition for the general case.
\end{rem}
This description of the cocycle allows us to give a short proof of 
the compactness of the stabilizer of a maximal triple in the 
Shilov boundary.
\begin{cor}
Let $(W_+, W, W_-)\in \cs^{(3)}$ be a maximal triple. 
Then the stabilizer of $(W_+, W, W_-)$ in $G$, 
$\Stab_G(W_+, W, W_-)$, is compact. 
\end{cor}
\begin{proof}
Since $G$ preserves $h$, any element fixing $(W_+, W, W_-)$ 
leaves invariant the definite $\sigma$-symmetric bilinear form $A_W$, and is 
hence contained in the compact subgroup defined by the invariance group of $A_W$ in $G$.
\end{proof}
\begin{rem}
A similar construction can be made for $G={\rm SO}(n,2)$, but the construction of the cocycle is 
more involved.
\end{rem}
%


\vskip1cm
\chapter{Tight homomorphisms and embeddings}
Tight homomorphisms and tight embeddings play an important role in this work. 
They are defined using bounded cohomology in degree two. 
The geometric properties of specific 
(optimal) representatives for these bounded cohomology classes are used to 
determine geometric properties of the associated maps and embeddings. 
The definition can be generalized to every bounded cohomology class, but the interesting 
geometric properties shown in this section depend on the knowledge of 
a nice representative of the bounded cohomology class and on the 
geometric properties of the configuration of points, which maximize 
the norm of this representative. 
Therefore we obtain geometric results only for homomorphisms into 
semisimple Lie groups of Hermitian type. 
Some of the results about totally geodesic tight embeddings suggest that a tight embedding 
is a bounded cohomological avatar of calibrations in deRham cohomology. 
General tight homomorphisms with respect to other classes might also 
inherit nice geometric properties. Unfortunately, for most bounded cohomology classes 
there are neither specific representatives nor any geometric
properties known.

\section{Definition and properties of $\a$-tight homomorphisms}
Let $H, G$ be locally compact 
groups. A continuous homomorphism $\pi: H\to G$ 
induces canonical pull-back maps $\pi^*$ in continuous cohomology and 
$\pi_b^*$ in continuous bounded cohomology, such that the following diagramm commutes:
\bqn
\xymatrix{
\hcb^*(G) \ar[d]^\kappa \ar[r]^{\pi_b^*} & \hcb^*(H) \ar[d]^\kappa\\
\hc^*(G) \ar[r]^{\pi^*} & \hc^*(H)
}
\eqn
where $\kappa$ is the natural comparison map between continuous bounded 
cohomology and continuous cohomology.

The continuous bounded cohomology groups come equipped with a canonical 
seminorm $||\cdot||$, with respect to which $\pi_b^*$ is norm decreasing, 
i.e. 
\bqn
||\pi_b^*(\a)|| \leq ||\a|| \text{ for all } \a\in \hcb^*(G).
\eqn

Requiring that a homomorphism preserves the norm of some classes in 
$\hcb^*(G)$ imposes restrictions on it. 

\begin{defi}\label{def:tight_gen}
Let $H,G$ be 
locally compact groups, $\a \in \hcb^*(G)$. 
A homomorphism $\pi: H\to G$ is called $\a$-tight, if 
$||\pi^*(\a)||=||\a||$.
\end{defi}

\begin{lemma}\label{lem:composition}
Let $\pi:H\to G$, $\psi: G\to K$ continuous homomorphisms of locally compact groups. 
Then $\psi$ is $\a$-tight and $\pi$ is $(\psi^*(\a))$-tight 
if and only if $\phi:= \psi \circ \pi$ is $\a$-tight.
\end{lemma}
Th following properties of tight embeddins are immediate consequences
of results about isometric isomorphisms in bounded 
cohomology from \cite{Burger_Monod_GAFA, Monod_book}.
\begin{lemma}\label{tight_properties} 
Suppose $H,G$ are 
locally compact groups, 
$\a \in \hcb^*(G)$ and $\pi:H\to G$ a $\a$-tight homomorphism.
\begin{enumerate}
\item{Let $L<G$ be a closed subgroup. If $\pi(H)\subset L$ then $\pi$ is $\a_{|_L}$-tight and 
        $||\a_{|_L}||=||\a||$}
\item{Let $H_0<H$ be a closed subgroup. 
        If there exists an invariant mean on $H/H_0$ then $\pi_{|_{H_0}}$ is $\a$-tight 
        and $||\pi_{|_{H_0}}^*\a||= ||\pi^*\a||= ||\a||$.}
\item{Let $R \lhd G$ be a normal closed amenable subgroup  
        and $\overline{\a}$ be the image of $\a$ 
        under the canonical isometric isomorphism $\hcb(G)\cong \hcb(G/R)$ 
        (see Proposition~\ref{prop:subgroups} in Chapter 1). 
        Then $||\overline{\a}|| = ||\a||$ and 
        the induced homomorphism 
        $\overline\pi : H\to G/R$ is $\overline{\a}$-tight.}
\item{Assume that $\a \in \hcb^2(G)$, $G= \Pi_{i=1}^n G_i$ and 
        $\a_i=\a_{|_{G_i}} \in \hcb^2(G_i)$.\\
        Then the induced homomorphisms $\pi_i=\pr_i\circ \pi:H\to G_i$ 
        are $\a_i$-tight.} 
\end{enumerate}
\end{lemma}
\begin{proof}
For $\pi:H\to L<G$ we have 
$||\a_{|_L}||\leq ||\a||$ 
and $\pi^*\a=\pi^*\a_{|_L}$. 
Since $\pi$ is $\a$-tight, we have 
\bqn
||\a_{|_L}||\geq ||\pi^*\a_{|_L}|| 
= ||\pi^*\a||=  ||\a||,
\eqn 
hence 
$||\a_{|_L}||=||\pi^*\a_{|_L}|| $. Moreover 
$||\a_{|_L}||=||\a||$.

The second and the third claim follow immediatly 
from the fact, that the corresponding maps in bounded cohomology 
are norm preserving \cite{Burger_Monod_GAFA}.

Observe that by 
Proposition~\ref{prop:degree_two} in Chapter 1 we have \cite[Corollary 4.4.1.]{Burger_Monod_GAFA} 
$\hcb^2(G)\cong \Pi_{i=1}^{n} \hcb^2(G_i)$ 
with the property that 
$||\a||=\sum_{i=1}^{n} ||\a_i||$, 
hence 
\bqn
\sum ||\pi^* \a_i|| = ||\pi^*\a|| = ||\a|| =\sum ||\a_i||,
\eqn
what together with $||\pi^*\a_i||\leq ||\a_i||$ proves the last claim.
\end{proof}

\begin{rem}
Given a map $f:M\to N$ between two nice enough topological spaces one might say that 
$f$ is $\a$-tight, whenever the induced homomorphism $f:\pi_1 (M) 
\to \pi_1(N)$ is $\a$-tight with 
respect to $\a \in \hb^*(\pi_1(N))$. \\
Tight homomorphisms might also be defined  with respect to bounded
cohomology classes with nontrivial coefficients. 
\end{rem}

\section{K\"ahler tight homomorphisms }
We turn to a specific class of tight homomorphisms which is related to Hermitian symmetric 
spaces.
\begin{defi}
Let $X$ be a Hermitian symmetric space of noncompact type, let $
G$ be the group of isometries of $X$ and $k_X^b\in \hcb^2(G, \RR)$ 
the bounded K\"ahler class.
A continuous homomorphism $\rho:H\to G$  is called \emph{tight}, 
if it is $k_X^b$-tight.
\end{defi}

Let $G=G_1\times \cdots\times G_k$. 
The associated symmetric space $X$ decomposes 
as $X=X_1\times \cdots\times X_k$. 
The choice of a $G$-invariant complex structure on $X$ corresponds to the choice 
of an element $Z_\frakg \in \frakc(\frakk)= \bigoplus_{i=1}^k \frakc(\frakk_i) = \bigoplus_{i=1}^k \RR Z_{\frakg_i}$.
On the simple factors there are only two choices of complex structures corresponding to $\pm Z_{\frakg_i}$.
Hence the choice of a $G$-invariant complex structure corresponds to the choice of an orientation on 
each simple factor.

Fixing a $G$-invariant metric on 
$X$ the choice of a $G$-invariant 
complex structure determines a unique K\"ahler class. 
If we do not fix the $G$-invariant metric, 
the space of K\"ahler classes becomes much bigger, since we have a K\"ahler class 
for each choice of metric and we might scale the metrics on the simple factors independently, getting a $k$-dimensional space 
of invariant metrics and hence of K\"ahler classes.

Fortunately the notion of tightness does not depend on this scaling. 
This is a consequence of the Banach space 
structure of $\hcb^2$ and the following lemma.
\begin{lemma}\label{lem:Banach_norm}
Let $V$ be a Banach space. Let $v_i\in V$, $(i=1,\dots,k)$ be norm one vectors 
and assume that $b=\frac{1}{k} \sum_{i=1}^k v_i$ 
has norm one. Then every convex linear combination 
$v= \sum_{i=1}^k \mu_i v_i$ with $\mu_i >0$ and $\sum_{i=1}^k \mu_i =1$
has norm one.
\end{lemma}
\begin{proof}
The norm of a vector $w\in V$ is  (Hahn-Banach) 
given by $||w||=\sup_{||\lambda||=1} |\lambda(w)|$, where $\lambda:V\to \RR$ is a linear form.
That $b$ is of norm one is equivalent to the fact that for every $\eps>0$ there exists a linear form $\lambda$ such that $\lambda(b)>1-\eps$.
Assume that $\lambda$ is a linear form satisfying 
$\lambda(b)= \lambda(\frac{1}{k} \sum_{i=1}^k v_i) = \frac{1}{k}
\sum_{i=1}^k \lambda(v_i)>1-\eps$. 
Then, since $\lambda(v_i)\leq 1$ for all $i$, 
it follows that $\lambda(v_i)>1-k\eps$ for all $i$. 
Hence $\lambda(v)= \lambda(\sum_{i=1}^k \mu_i v_i) = 
\sum_{i=1}^k \mu_i \lambda(v_i) > (1-k\eps) (\sum_{i=1}^k \mu_i)=(1-k\eps).$
Since $\eps>0$ was arbitrary, this implies that $||v||=1$.
\end{proof}
\begin{lemma}
Fix a $G$-invariant complex structure on the Hermitian symmetric space $X$ associated to $G$.
Then a homomorphism $\rho:H\to G$ is tight if and only if it is tight with respect 
to any K\"ahler class which is compatible with the chosen $G$-invariant complex structure.
\end{lemma} 
\begin{proof}
The second bounded cohomology group of any locally compact group 
$G=G_1\times \cdots\times G_k$, 
$\hcb^2(G,\RR)\cong \bigoplus_{i=1}^k \hcb^2(G_i, \RR)$, is a Banach space, 
on which the canonical norm is an $l_1$-norm. 
For $G$ a connected semisimple group of Hermitian type, 
this $l_1$-norm is defined with respect to the basis given by generators $\a_i\in \hcb^2(G_i, \RR)\cong \RR \a_i$, where $\a_i$ are the  
K\"ahler classes compatible with the induced complex structure on $X_i$, normalized in such a way that $||\a_i||=1$. 
Assume that $\pi:H\to G$ is tight with respect to the K\"ahler class $k_X^b= \frac{1}{k} \sum_{i=1}^k \a_i$. 
Then by Lemma~\ref{tight_properties} $\pi$ is also tight with respect to $\a_i$ for all $i$. 
Hence we have elements $e_i:=\pi^*(\a_i) \in \hcb^2(H)$ of norm $1$ with the property that 
$b=\frac{1}{k} \sum_{i=1}^k e_i$ is also of norm $1$. 
Any other K\"ahler class $\a$ 
compatible with the chosen complex structure is of the form $\a=\sum_{i=1}^k \mu_i \a_i$, 
where $\mu_i>0$. We may assume that $\sum_{i=1}^k \mu_i=1$, hence in particular $||\a||=1$.
The pullback $\pi^*(\a)= \pi^*(\sum_{i=1}^k \mu_i \a_i) = \sum_{i=1}^k \mu_i \pi^*(\a_i)$.
The homomorphism $\pi:H\to G$ is $\a$-tight if and only if the norm of $\pi^*(\a)$ is $1$. 
From Lemma ~\ref{lem:Banach_norm} it follows, that $||\pi^*(\a)||=1$.
\end{proof}

\section{Tight embeddings and their properties}
Let now $Y, X$ be symmetric spaces of noncompact type.
To a totally geodesic embedding $f: Y\to X$ we have associated the corresponding 
homomorphism 
$\pi:H_Y\to G={\rm Is}(X)^\circ$ of some finite extension $H_Y$  
of the connected component of the isometry group ${\rm Is}(Y)^\circ$ of $Y$.
Assume that $X$ is Hermitian.
\begin{defi}
A totally geodesic embedding $f:Y\to X$  
is called \emph{tight} if the corresponding homomorphism $\pi:H_Y\to {\rm Is}(X)^\circ$  (see Section~\ref{sec:totally_geodesic} in Chapter 2) is tight.
\end{defi}
We may give a more geometric definition of tightness of a totally geodesic embedding. 
Let $\omega\in \Omega^2(X)^G$ be a $G$-invariant K\"ahler form on $X$ compatible with the 
$G$-invariant complex structure of $X$. Let $f:Y\to X$ be a totally geodesic embedding, 
then $f^*\omega$ defines an $H_Y$-invariant K\"ahler form on $Y$. In particular if $f^*\omega$ does not vanish, 
$Y$ is Hermitian, but $f^*\omega$ might not 
be compatible with the complex structure on $Y$. 
The embedding $f$ is tight if and only if 
\bqn
\sup_{\Delta\subset Y} \int_\Delta f^*\omega = \sup_{\Delta\subset X} \int_\Delta \omega
\eqn
where the supremum is taken over all geodesic triangles $\Delta$ in
$Y$ or $X$, respectively. 
\begin{rem}
Note that the restriction to Hermitian symmetric spaces of noncompact type 
is natural, 
since the condition of being tight is empty if $\hcb^2(G) = 0$, 
which is the case whenever $G$ is compact or has no factor of 
Hermitian type 
(\cite{Burger_Monod_JEMS}).
\end{rem}

An important tool to study general tight embeddings is the tight holomorphic embedding 
of the disc $\DD$, i.e. the diagonal embedding of $\DD$ into a maximal
polydisc in a Hermitian symmetric space.
\begin{lemma}\label{lem:holo_tight_tight}
The tight holomorphic disc $f(\DD) = C\subset X$ is tightly embedded into $X$.
\end{lemma}
\begin{proof}
Let $\omega$ be the invariant K\"ahler form on $X$. 
Let $\Delta$ be an ideal triangle in $\DD$. Then $\int_\Delta
f^*\omega= \pi r_X$.
\end{proof}

In the case where $H$ is a simple group, $\hcb^2(H, \RR)$ is one dimensional 
and we can relate the pull-back of the K\"ahler class of $G$ to the K\"ahler
class of $H$.

\begin{lemma}\label{tight_max}
Assume that $H$ is simple. Suppose $\pi:H\to G$  is a $k_G$-tight 
homomorphism. Then $\pi^*k_X^b= \frac{r_X}{r_Y} k_Y^b$, where $k_X^b, k_Y^b$ denote 
the bounded K\"ahler classes.
\end{lemma}
\begin{proof}
Since $H$ is simple, the space $\hcb^2(H)$ is one dimensional 
and therefore $\pi^*k_X^b$ is a constant multiple of $k_Y^b$. 
The constant is determined by the ratio $||k_X^b||/||k_Y^b||= \frac{r_X}{r_Y}$.
\end{proof}

We can now use the results of \cite{Clerc_Orsted_2} about the characterization
 of maximal triples (see Proposition~\ref{results_c_o} in Chapter 5)
to exhibit certain properties of tight embeddings:

\begin{prop}\label{tight_embedding}
Let $f:Y\to X$ be a tight embedding then $f$ extends uniquely to a
 continuous equivariant map 
 $\cs_Y\to \cs_X$ of the Shilov boundaries. 
\end{prop}
\begin{proof}
We work with the bounded symmetric domain realizations $\Dd_Y$ and $\Dd_X$ 
of $Y$ and $X$. We may assume that $f(0)=0$.
Assume first that $\Dd_Y=\DD$. 
For any point $x \in \cs_Y=S^1$ we extend the map $f$ 
by mapping $x$ to the endpoint in $\partial\Dd_X$  
of the image of the geodesisc joining $x$ to $0\in \DD$. 
This defines a map 
\bqn
\ol{f} :S^1 \to \partial \Dd_X,
\eqn
which is independent 
of the chosen base point $0\in \DD$. 
Since any distinct points $(x,y,z)\in \cs_1^3$
are joined by geodesics, the corresponding 
image points $( f (x), f (y),  f (z))\in \partial \Dd_X^3$ 
are also joined by geodesics 
and hence (Lemma~\ref{lem:transversality} in Chapter 4) pairwise transverse. 
The assumption that $\pi$ is tight implies that any pairwise distinct points $x,y,z$ 
realize the supremum of $f^*\beta_{\Dd_X}$
\bqn
||f^*\beta_{\Dd_X}(x,y,z)||=r_X||\beta_\DD(x,y,z)||=||\beta_{\Dd_X}||.
\eqn

Hence by Proposition \ref{results_c_o} in Chapter 5 we get that 
$ f (x),  f (y),  f (z)\in \cs$.
This extends $f$ to $ \ol{f}: S^1\to \cs_X$. 
Since  $\ol{f}$ is independent of the base point in $\DD$, it is 
equivariant with respect to $\pi$.

Suppose now that $\Dd_Y$ is arbitrary.
 For any point $x\in \cs_Y$ take the unique tight holomorphic disc $\DD_x$ through $0$ in $\Dd_Y$, containing $x$ in its boundary.
The map $f_x=f_{|_{\DD_x}}$ extends by the previous argument to the
 Shilov boundary. 
Making this construction for all $y\in \cs_Y$ we 
can extend $f$ to the Shilov boundary.
This extension is independent of the chosen point $0\in \Dd_Y$. 
Indeed, let $a\in \Dd_Y$ be another point and $\DD^a_x$ the tight holomorphic disc 
through $a$ containing $x$ in its boundary. Then the two geodesics
 $\g_{0x}$ and $\g_{ax}$ lying in $\DD_x$ and $\DD^a_x$ respectively
 are both maximal singular geodesics which  
lie in bounded distance from each other by
 Proposition~\ref{boundaries} of Chapter 4.
Their images under $f$ remain in bounded distance from each other, 
hence determine the 
same point in the geometric boundary and hence by
 Proposition~\ref{boundaries} of Chapter 4 the same point in the Shilov boundary $\cs_X$. 
This implies that the extension $\ol{f}$ of $f$ is $\pi$-equivariant and thus a 
diffeomorphism onto its image.
\end{proof}

As Corollary we get the analog of Lemma~\ref{tight_max} 
for the Bergmann cocycle
\begin{cor}
Assume that $H$ is simple. Suppose $\pi:H\to G$  is a $k_G$-tight 
homomorphism. Then $\phi^*(\b_G)= \frac{r_X}{r_Y} \b_H$, where $\b_G, \b_H$ denote 
the normalized Bergmann cocycles.
\end{cor}
\begin{proof}
We have the following commuting diagramm 

\bqn
\xymatrix{
\h^* (\binfty(\cs_X^{(*+1)})^G) \ar[d]^{\phi^*} \ar[r] & \hcb^*(G) \ar[d]^{\rho^*}\\
\h^* (\binfty(\cs_Y^{(*+1)})^H) \ar[r] & \hcb^*(H)
}
\eqn
Since $H$ acts transitively on $\cs_Y^{(2)}$, there are no coboundaries in degree two and 
the pull-back $\phi^*(\b_X) \in \binfty(\cs_Y^{(*+1)})^H$ 
of $\b_X$ is a multiple of $\b_Y \in \binfty(\cs_Y^{(*+1)})^H$ since
$H$ is simple.
\end{proof}
Based on the different behaviour of the Bergmann cocycle with respect to (non) tube type 
Hermitian symmetric spaces, we obtain the following obstruction to the
existence of tight totally geodesic embeddings.
\begin{cor}\label{prop:tight_non_tube}
Let $X$ be a Hermitian symmetric space of tube type. 
Suppose that $Y$ is a Hermitian symmetric space, $f:Y\to X$ a tight totally
geodesic embedding. Then $Y$ is of tube type. 
\end{cor}
\begin{proof}
Since $f:Y\to X$ is tight, it extends to $\ol{f}:\cs_Y\to \cs_X$. 
The pull-back $\ol{f}^*(\b_X)$  is a multiple of $\b_Y$. 
It takes only finitely many values, whereas  -- if $Y$ is not of tube
type -- 
$\b_Y$ assumes a continuous interval of values. 
\end{proof}
\begin{rem}
The converse statement of Proposition~\ref{tight_embedding} is not true. 
There are totally geodesic embeddings which extend to a map 
of the Shilov boundaries but which are not tight. 
\end{rem}
\subsection{Tight homomorphisms into Lie groups of Hermitian type}
\begin{prop}\label{centralizer}
Suppose that the inclusion of a semisimple subgroup $H< G$ is a tight homomorphism. Then 
the centralizer $\zent_{G}(H)$ of $H$ in $G$ 
is compact.
\end{prop}
\begin{proof}
Denote by $Y$ and $X$ the Hermitian symmetric spaces associated to $H$ and $G$.
This time we will work first in the homogeneous model of $Y$ and $X$. 
Choose a 
triple of pairwise transverse points $x,y,z\in \cs_Y$ such that 
$f^*\beta(x,y,z)=r_X\pi$. 
Let $\g_{xy}$, $\g_{yz}$, $\g_{zx}$ denote geodesics 
in $f(\Dd_Y)$ connecting $f(x), f(y), f(z)$. 
Since $f(x), f(y), f(z)\in\cs_X$, 
Proposition \ref{boundaries} in Chapter 4 
implies that the geodesics are of some 
(possibly different) type $\theta\in\Theta_r$. 
Each element $g\in \zent_{G}(H)$ 
maps the geodesics to geodesics $g\g_{xy},g\g_{yz},g\g_{zx}$, 
such that $d(\g(t),g\g(t))$ is uniformly bounded. 
Furthermore, since $g\in G$, $g\g$ is of the same type as $\g$. 
By the definition of the geometric boundary 
$g\g(\infty)=\g(\infty)$ 
and hence, 
by Proposition \ref{boundaries} in Chapter 4,
$g\g_{xy},g\g_{yz},g\g_{zx}$ connect 
the same points $x,y,z$ in $\cs_X$.
In particular 
\bqn
 \zent_{G}(H) < \Stab_{G} (f(x),f(y), f(z)),
\eqn 
which is compact by 
Proposition \ref{results_c_o} in Chapter 4.
\end{proof}

\begin{thm}\label{thm:tightness}
Let $\pi:H\to G$ a (continuous) tight homomorphism. 
Then $L=\overline{\pi(H)}^Z$ 
is reductive with compact center. The symmetric space associated to 
the semisimple part of $L$ is Hermitian symmetric.
\end{thm}
\begin{proof}
Let $L:=\overline{\pi(\Gamma)}^{\rm Z}(\RR)$ be the
real points of the Zariski closure of $\pi(H)$.  We may assume 
(Lemma~\ref{tight_properties}) that $L$ is connected and that the 
inclusion $L\to G$ is tight. 
Let $L'=L/R$ be the semisimple part of $L$. The inclusion $L'\to G$ 
is tight by Lemma~\ref{tight_properties}.
Hence by Proposition \ref{centralizer} the centralizer of $L'$ 
in $G$ is compact. 
Since $\zent_G(L) \subset \zent_G(L')$, the centralizer of $L$ in $G$ is
compact. The center of $L$ is contained in $\zent_G(L)$ and hence it is compact.

Assume that $L$ would not be reductive, then $L$ is contained in a 
proper parabolic subgroup of $G$ (\cite{Borel_Tits}). 
But if $L$ is contained in a proper parabolic subgroup of $G$, 
the centralizer of the Levi component of this parabolic group, 
which is noncompact, is contained in the centralizer of $L$, which is
compact. This is a contradiction. 
Therefore 
$L$ is reductive.
Denote by $M$ the semisimple part of $L$ and let 
$M'=M_1\times \dots \times M_\ell$ be the product 
of the simple factors of $M$ such that
$k^b_L|_{M_i}\neq 0$. 
From Lemma \ref{tight_properties} it follows that the induced homomorphism 
of $M' \to G$ is tight. 
This implies that $M=M'\times \Cc$, 
where $\Cc$ are the compact factors of $M$, 
hence the symmetric space 
$Y$ associated to $M$ is Hermitian and tightly embedded into $Y$.
\end{proof}

\subsection{Tight embeddings and tube type}

Tight embeddings extend to equivariant maps of the Shilov boundaries. 
The Shilov boundary reflects the property of being of tube type or not, 
therefore tight embeddings are sensitive to the difference between 
tube type and not tube type. 

\begin{prop}\label{prop:tight_tube}
Assume that $Y$ is of tube type and $f: Y\to X$ is a tight embedding. 
Then $f(Y)$ is contained in a maximal subdomain of tube type $T\subset X$.
\end{prop}

\begin{proof}
Fix two transverse points $x,y \in \cs_Y$. They are mapped under $\ol{f}$ to 
transverse points $\ol{f}(x), \ol{f}(y) \in \cs_X$, hence they determine (Lemma 
\ref{lem:transverse_tube} in Chapter 4) a unique maximal subdomain of tube type
 $T_{\ol{f}(x), \ol{f}(y)} \subset X$.
Consider the set $M_{xy} \subset \cs$ consisting of all points $z\in \cs$ that are 
transverse to $x$ and $y$ and form a maximal triple with $x$ and $y$. Then,
 by Proposition \ref{results_c_o} in Chapter 5,  
$\ol{f}(M_{xy}) \subset \cs( T_{\ol{f}(x), \ol{f}(y)})$. 
The two sets $\ol{f}(\cs_Y)$ and $\cs(T_{\ol{f}(x), \ol{f}(y)})$ are real smooth connected algebraic subvarieties of $\cs$, 
which intersect in $\ol{f}(M_{xy})$, which is an open set in $\ol{f}(\cs_Y)$. 
This implies  $\ol{f}(\cs_Y) \cap \cs(T_{\ol{f}(x), \ol{f}(y)}) = \ol{f}(\cs_Y)$, hence $\ol{f}(\cs_Y) \subset \cs(T_{\ol{f}(x), \ol{f}(y)})$
\end{proof}


\section{The K\"ahler form on Lie algebra level}\label{algebra}
Assume $H, G$ are simple and of 
Hermitian type.
Then $\frakc(\frakk)$ is one dimensional. 
Let $Z_\frakh, Z_\frakg$ be generators inducing the complex
structures on $Y$ and $X$. 
The Lie algebra $\frakk$ decomposes as 
$\frakk= \frakc(\frakk)\oplus [\frakk,\frakk]$ and for any $X\in \frakk$ 
there exists a unique $\lambda (X)  \in\RR$ such that 
\bqn
X= \lambda(X) Z_\frakg +[\frakk, \frakk].
\eqn
Since $[\frakp,\frakp]\subset \frakk$, define for all $V,W\in\frakp$: 
$\Omega(V,W)=\lambda([V,W])$. This map is linear, antisymmetric 
and $K$-invariant. Hence $\Omega \in \Lambda^2(\frakp)^K$ extends under 
the identification of $T_{x_0}X\cong\frakp$ via left translations to a 
$G$-invariant $2$-form on $X$.  This implies that $\Omega$ is a constant 
multiple of the (normalized) K\"ahler form $\omega_X$, where the constant can be explicitly calculated, 
as it is done in \cite{Burger_Iozzi_supq}.

Note that for the holomorphic tight embedding of the disc into a Hermitian symmetric space of tube type, 
the associated homomorphism of the Lie algebras satisfies $\lambda(\pi(\frakh))= 1$, since $\pi(Z_\frakh)=Z_\frakg$. 
For homorphisms between simple groups of tube type the condition $|\lambda(\pi(Z_0))| = 1$  turns out to 
give a criterion for tightness which is easily computable for explicit representations of Lie algebras.
\begin{lemma}\label{tight_tight}
Assume that $H,G$ are simple and of tube type. 
Let $\pi: H\to G$ be a homomorphism and 
$\pi_*:\frakh\to \frakg$ the corresponding homomorphism of the Lie algebras.
Then $\pi$ is tight if and only if $\pi_*$ satisfies $|\lambda(\pi_*(Z_\frakh))|=1$.
\end{lemma}
\begin{proof}
Consider first the case where $\pi_*:\frakh_0 = \fraks\fraku(1,1) \to \frakg$. 
Then the Lie algebra homomorphism  $\pi_0:\fraks\fraku(1,1) \to \frakg$ associated 
to the holomorphic tight embedding of a disc into $X$ satisfies 
$\pi_0(Z_{\frakh_0})= Z_\frakg$ and 
hence $\lambda(\pi_0(Z_{\frakh_0}))=1$. But $|\lambda(\pi_0(Z_{\frakh_0}))|$ determines the ratio of the 
pull back of the K\"ahler form on $\frakg$ and of the K\"ahler form on $\frakh$. Since every other tight 
Lie group homomorphism gives the same ratio (see Lemma~\ref{tight_max}), 
we have that $\pi:H_0\to G$ is tight if and only if 
$\pi_*:\frakh_0\to \frakg$ satisfies $|\lambda(\pi_*(Z_{\frakh_0})|=1$. 
Suppose that $\pi_*:\frakh\to \frakg$ is an injective Lie algebra homomorphism. 
Precomposition of $\pi_*$ with $\pi_0:\frakh_0\to \frakh$, $\pi_0(Z_{\frakh_0})=Z_\frakh$ gives an injective Lie algebra 
homomorphism $\pi_*\circ\pi_0:\frakh_0\to \frakg$. This corresponds to a tight homomorphism of the 
Lie groups ${\rm SU}(1,1) \to G$ 
if and only if $1=|\lambda(\pi_*\circ \pi_0(Z_{\frakh_0}))|= |\lambda(\pi_*(Z_\frakh))|$. 
The Lie group homomorphism $\pi:H\to G$ is by Lemma~\ref{lem:composition} tight if and only if $\pi\circ\pi_0$ is tight, 
hence if and only if  $|\lambda(\pi_*(Z_\frakh))|=1$.
\end{proof}

\begin{rem}
The same result is true, if we do not assume that $H$ and $G$ are
simple, but if we fix a complex structure both on the Hermitian
symmetric space associated to $X$ and to $Y$, thus element $Z_\frakg$
and $Z_\frakh$ in the center of $\frakk_G$ respectively in the center
of $\frakk_H$.
\end{rem}

\begin{rem}
The above criterion applies to tight embeddings of Hermitian symmetric spaces that are not of tube type.
Let $f:Y\to X$ be a tight totally geodesic embedding of Hermitian
symmetric space. 
Then Proposition~\ref{prop:tight_tube} implies that the restriction 
of $f$ to a maximal subdomain of tube type $T_Y$ in $Y$ gives a well-defined 
totally geodesic tight embedding $f_{|_{T_Y}}:T_Y\to T_X$ into some maximal subdomain of tube type $T_X$ in $X$. 
The embedding $f$ is tight if and only if $f_{|_{T_Y}}$ is tight. We
may now apply the 
above criterion to the tight embedding $f_{|_{T_Y}}$.\\
\end{rem}

The tigthness criterion fits nicely together with the different types of holomorphic totally geodesic embeddings.
Recall the following definitions from \cite{Satake_book}:
\begin{defi}
Suppose $\frakh, \frakg$ are simple Lie algebras of Hermitian type 
with $Z_\frakh, Z_\frakg$ the elements of the center of the maximal 
compact Lie subalgebras, which determine the complex structure. 
Then a homomorphism $\pi:\frakh\to \frakg$ is said to be of type 
\begin{enumerate}
\item{$(\h1)$ if $\ad(\pi(Z_\frakh))= \ad(Z_\frakg)$.}
\item{$(\h2)$ if $\pi(Z_\frakh)=Z_\frakg$.}
\item{$({\h2}^\prime)$ if $\pi$ is $(\h1)$ and the induced $\CC$-linear 
map $\Dd_\frakh\to \Dd_\frakg$ maps the Shilov boundary 
into the Shilov boundary.}
\end{enumerate}
\end{defi}

\begin{cor}(of Proposition~\ref{tight_embedding}) 
Assume that $f:Y\to X$ is a holomorphic tight embedding with corresponding 
$(\h1)$ Lie algebra homomorphism $\pi:\frakh\to \frakg$. Then $\pi$ is a $(\h2^\prime)$ homomorphism.
\end{cor}

Lemma~\ref{tight_tight} together with results of Satake \cite[Proposition~10.12]{Satake_book} imply immediatly
\begin{cor}\label{H2_hom}
Suppose $Y, X$ are Hermitian symmetric spaces of tube type and $f:Y\to X$
is a holomorphic embedding.
Then $f$ is tight if and only if the corresponding homomorphism of Lie algebras 
$\pi:\frakh\to \frakg$ is an $(\h2)$-homomorphism, 
i.e. $\pi(Z_\frakh)=Z_\frakg$.
\end{cor}

\subsection{Examples and Counterexamples}
We give several examples of tight embeddings and of embeddings which are not tight in order 
to illustrate this class of totally geodesic embeddings.\\

1) The inclusion $\G<G$ of a lattice is tight. The bounded cohomology 
of $\G$ is given by the same resolution as for $G$, thus the norm of 
the restriction of a cocycle equals the norm of the cocycle. \\

2) By Corollary~\ref{H2_hom} all $(\h2)$-holomorphic embeddings are tight 
embeddings. There exist such embeddings of all classical groups of tube type 
into all other classical groups of tube type provided that the rank conditions 
are satisfied. They are described in \cite{Satake_annals}. 
An example is the embedding corresponding to the Lie algebra homomorphism 
$\fraks\fraku(n,n) \to \fraks\frakp(4n,\RR)$ defined by 
\bqn
g=\begin{pmatrix} X_{11} & X_2 \\ X_2^* & X_{12} \end{pmatrix}
\to \begin{pmatrix} A & B \\ \ol{B} & \ol{A} \end{pmatrix}, 
\eqn 
with 
\bqn
A=\begin{pmatrix} \ol{X_{12}} & 0 \\ 0 & X_{11} \end{pmatrix}\\
B=\begin{pmatrix} 0 & X_2^T \\ \ol{X_2} & 0 \end{pmatrix}.
\eqn\\

3) Other examples of tight embeddings are the natural embeddings of maximal 
subdomains of tube type into the non tube type domain, for example 
$\fraks\fraku(n,n) \to \fraks\fraku(n,m)$ with $m>n$.\\

4) Interesting non-holomorphic tight embeddings of the disc can be obtained 
from the irreducible representations of ${\rm SL}(2,\RR)$ into 
${\rm Sp}(2n,\RR)$. 
\begin{prop}\label{irreducible_rep}
The homomorphism $\pi:\frakg_0=\fraks\frakl(2,\RR)\to \fraks\frakp(2n,\RR)$ which is induced by the 
$2n$-dimensional irreducible representation of $\fraks\frakl(2,\RR)$ is tight.
\end{prop}
\begin{proof}
We only have to determine $\lambda(\pi(Z_{\frakg_0}))$, 
where 
\bqn
\pi(Z_{\frakg_0})= \lambda(\pi(Z_{\frakg_0}))Z_\frakg + [\frakk, \frakk].
\eqn
This can be done by explicit calculations: 
Let $V^{2n}$ the vector space of homogeneous polynomials of degree $n$ in two variables $x,y$. 
A basis is given by $(P_0, \dots P_m)$, $m= 2n-1$, where $P_k(x,y)= x^{m-k}y^k$. 
The irreducible representation of $\fraks\frakl(2,\RR)$ is given by the following action: 
Let $X=\begin{pmatrix} a & b \\ c & -a \end{pmatrix} \in \fraks\frakl(2,\RR)$, then 
\bqn
\pi(X) P_k(x,y)= a(m-2k) P_k + b(m-k) P_{k+1} + c k P_{k-1}.
\eqn
This action preserves the skew symmetric bilinear form $\langle .,.\rangle$ on $V$, 
defined by $\langle P_k, P_{m-k}\rangle = (-1)^k {{m}\choose{k}}^{-1}$. 
Thus this defines the irreducible representation 
\bqn
\pi:\fraks\frakl(2,\RR) \to \fraks\frakp(2n,\RR),
\eqn 
into the Lie algebra of the symplectic group ${\rm Sp}(V,\langle .,.\rangle)$.
The map $J$ defined by $J P_k= (-1)^k P_{m-k}$ gives a complex structure on $V$ and is an element in the center of 
$\frakk\subset \fraks\frakp(2n,\RR)$. The element $Z_\frakg :=\frac{1}{2} J$ induces the complex structure 
on $\frakp\subset \fraks\frakp(2n,\RR)$ via the adjoint action.

The image of the element 
\bqn
Z_0 = \frac{1}{2} \begin{pmatrix} 0 & -1 \\ 1 & 0\end{pmatrix}
\eqn
is given by 
$\pi(Z_{\frakg_0}) P_k = \frac{1}{2}((k-m)P_{k+1} + k P_{k-1})$. 
Decomposing $\pi(Z_{\frakg_0})= \lambda(\pi(Z_{\frakg_0}))Z_\frakg + [\frakk, \frakk]$, 
we have that 
\bqn
Z_\frakg  \pi(Z_{\frakg_0}) = \frac{-\lambda}{4} \id_V + Z_\frakg [\frakk, \frakk].
\eqn
Since $\tr (Z_\frakg [\frakk, \frakk])= 0$, 
we have that $\tr(Z_\frakg  \pi(Z_{\frakg_0})) = \frac{-\lambda}{4} dim(V)$. 
Hence we have to show that 
\bqn
|\tr(Z_\frakg  \pi(Z_{\frakg_0}))| = \frac{1}{4} dim(V)= \frac{n}{2} =|\tr(Z_\frakg  Z_\frakg)|.
\eqn 
Now $Z_\frakg  \pi(Z_{\frakg_0} P_k = \frac{1}{4} (-1)^{k+1} [(m-k) P_{m-k-1}- k P_{m-k+1}]$. 
Thus the diagonal terms are, for $k={n-1}$, $\frac{1}{4} (-1)^{n} (2n-1-n+1)$ 
and, for $k=n$, $\frac{1}{4} (-1)^n n$, 
hence $| \tr(Z_\frakg  \pi(Z_{\frakg_0}))| =  \frac{n}{2}$.
\end{proof}

We give also some counterexamples:
\begin{notation*}
A $p\times q$ 
matrix consisting of a $p \times p$ matrix and $q-p$ zero columns 
of a $q \times q$ matrix and $p-q$ zero lines, will be denoted by 
$(A,0)$.
\end{notation*}
5) Regarding a complex vector space $V_\CC$ of dimension $(1+n)$ with a Hermitian form of signature$(1,n)$ as 
real vector space $V_\RR$ of dimension $(2+n)$ with a quadratic form of signature $(2,n)$ provides a 
natural embedding ${\rm SU}(1,n)\to {\rm SO}(2,2n)$. 
This embedding $\CC\HH^n \to X_{2,2n}$ is holomorphic but it is not tight. 
The image of the tight holomorphic disc 
\bqn
\begin{pmatrix} 0 & (z,0) \\ (z,0) & 0\end{pmatrix}
\eqn  
is 
\bqn
\begin{pmatrix} 0 & (Z,0) \\ (Z^T,0) & 0\end{pmatrix}
\eqn 
with $Z=\begin{pmatrix} x & y\\ -y & x\end{pmatrix}$, where $x+iy=z$ but the tight holomorphic disc in 
${\rm SO}(2,2n)$ is given by 
\bqn
\begin{pmatrix} 0 & (A,0) \\ (A^T,0) & 0\end{pmatrix},
\eqn where $A=Z=\begin{pmatrix} a_1 & b_1\\ 0 & 0\end{pmatrix}$.\\

Since $\CC\HH^n$ is of rank one, the totally geodesic embedding extends 
continuously to a map of topological compactifications. The image 
does not lie in the Shilov boundary, but in 
the $\Ad(K)$-orbit of $o_1$.

6) Tightness really depends 
on the complex structure and cannot be reduced to the property 
that the totally geodesic map extends to the Shilov boundary, which is 
somehow a real object,  
as the following example shows:

In ${\rm SL}(4,\RR)$ 
there are two copies of ${\rm Sp}(4, \RR)$, which are conjugate. 
Namely 
\bqn
{\rm Sp}(4, \RR)_A :=\{ g \in {\rm SL}(4,\RR) \,|\, g^* J g = J\}\\
{\rm Sp}(4, \RR)_B :=\{ g \in {\rm SL}(4,\RR) \,|\, g^* \tilde{J} g = \tilde{J}\}, 
\eqn
where $J= \begin{pmatrix} 0 & \id\\ -\id & 0 \end{pmatrix}$ and 
$\tilde{J}=\begin{pmatrix} 0 & \Lambda \\ -\Lambda & 0\end{pmatrix}$ 
with $\Lambda=\begin{pmatrix} 0 & 1\\ 1 & 0\end{pmatrix}$. 
The element $s= \begin{pmatrix} \id & 0 \\ 0 & \Lambda\end{pmatrix}$ 
satisfies $ \tilde{J}= sJs^{-1}$ 
and thus conjugates $s^{-1} {\rm Sp}(4, \RR)_B s={\rm Sp}(4, \RR)_A$.

There are two embeddings $i_{A,B}: {\rm SL}(2,\RR) \to {\rm SL}(4,\RR)$ 
namely 
\bqn
i_A\left( \begin{pmatrix} a & b\\ c & d\end{pmatrix}\right) = \begin{pmatrix} a\id
  & b\id \\ c\id & d\id \end{pmatrix}
\eqn
 and 
\bqn
i_B\left( \begin{pmatrix} a & b\\ c & d\end{pmatrix}\right) = \begin{pmatrix} a\id
  & b\Lambda \\ c\Lambda & d\id \end{pmatrix}.
\eqn
These two embeddings are conjugate by $s$. 
The image of ${\rm SL}(2,\RR)$ under these two embeddings lie in 
${\rm Sp}(4, \RR)_A \cap {\rm Sp}(4, \RR)_B$. 
The embedding $i_A$ is tight with respect to ${\rm Sp}(4, \RR)_A$ 
and totally real with respect to ${\rm Sp}(4, \RR)_B$, 
whereas the embedding $i_B$ is tight with respect to ${\rm Sp}(4, \RR)_B$ 
and totally real with respect to ${\rm Sp}(4, \RR)_A$. 

The boundary $S^1$ of $\DD$ is mapped under both embeddings 
into the Shilov boundaries 
$\cs_A$ respectively $\cs_B$.

This example can now easily be extended to ${\rm Sp}(2n, \RR)$ 
giving totally geodesic embeddings of the disc $\DD$ which 
which are not tight, 
neither holomorphic nor totally real, 
but extend to the Shilov boundary.

Note that the totally real embeddings 
described above 
extend to embeddings of ${\rm SL}(2,\CC)$ 
whereas the tight embeddings extend to embeddings of ${\rm SO}(2,2)$ into 
${\rm Sp}(4, \RR)$ .\\
%
%
\section{Classification of tight embeddings}
The holomorphic homomorphisms 
$\pi: \fraks\frakl (2,\RR)\to \frakg$ have been classified by 
their multiplicities (see \cite{Satake_book}). 
In general tight homomorphisms $\pi:\fraks\frakl (2,\RR)\to \frakg$
are not holomorphic. We consider the smallest holomorphically embedded totally 
geodesic subspace containing the image of a tight embedding.

\subsection{Hermitian hull}
We introduce the notion of the Hermitian hull for a subset of 
a Hermitian symmetric space.
\begin{defi}
Let $X$ be a Hermitian symmetric space of noncompact type. 
The Hermitian hull of a subset $V\subset X$ is the smallest sub-Hermitian 
symmetric space $Y\subset X$, such that
$V\subset Y$.
\end{defi}

\begin{lemma}
The Hermitian hull $\Hh(V)$ of $V\subset X$ is well-defined.
\end{lemma}
\begin{proof}
Since $V\subset X$ and $X$ is a sub-Hermitian symmetric space of $X$,
it remains to prove 
that there is a unique smallest sub-Hermitian symmetric space 
$Y\subset X$ containing $V$.
 
Assume there were two such subspaces $Y,Y^\prime$ with $V\subset Y$, 
$V\subset Y^\prime$. 
Then $V\subset Y\cap Y^\prime$. 
But since $Y,Y^\prime$ are totally geodesic submanifolds in $X$, 
i.e. every geodesic in $Y$ respectively in $Y^\prime$ is a
 geodesic in $X$, the intersection
$Y\cap Y^\prime$ is a totally geodesic submanifold, hence a 
subsymmetric space. 
Because $Y, Y^\prime$ are holomorphically embedded, 
i.e. for any $v\in T_pY$, $Jv\in T_pY$, where $J$ denotes 
the complex structure of $X$,
the intersection $Y\cap Y^\prime$ is holomorphically embedded, 
and is therefore a sub-Hermitian symmetric space. 
The minimality of $Y$ and $Y^\prime$ forces now $Y=Y\cap Y^\prime= Y^\prime$.
\end{proof}
An immediate consequence of Proposition~\ref{prop:tight_tube} is the following
\begin{cor}
Let $Y$ be a Hermitian symmetric space of tube type. Suppose that $f:Y\to X$ 
is a tight embedding. Then $\Hh(f)$ is of tube type.
\end{cor}
\begin{rem}
For a totally geodesic embedding 
$f: Y\to X$ of Hermitian symmetric spaces 
we call $\Hh(f)=\Hh(f(Y))$ the Hermitian hull of $f$.
The embedding $f$ is holomorphic iff $\Hh(f)= f(Y)$.
\end{rem}
The following Lemma is immediate from the definitions.
\begin{lemma}
Assume that $f: Y\to X$ is a tight totally geodesic embedding. 
Then the map $i:\Hh(f)\to X$ is a tight holomorphic totally geodesic embedding.
\end{lemma}
\begin{lemma}\label{tight_tube}
Suppose 
$X$ is an irreducible Hermitian symmetric space of tube type, 
$H=\Is(Y)^\circ, G=\Is(X)^\circ$ and 
$\pi:H\to G$ is a $k$-tight homomorphism.
 Denote by $f$ the induced totally geodesic embedding $f: Y\to X$, 
by $W=\Hh(f)$ the Hermitian hull of $f$, by 
$L=\Is(W)^\circ$ be the corresponding subgroup of $G$ 
and by $\frakl$ its Lie algebra.

Then
\begin{itemize}
\item{the center of 
        $\frakk_\frakl$ contains an element $Z_\frakl$ defining the
        complex structure on $W$ such that $Z_\frakg =  Z_\frakl.$}
\item{the corresponding central element $Z_G$ of the maximal compact subgroup $K<G$ is contained 
        in the center of a maximal compact subgroup of $L$.}
\end{itemize}
\end{lemma}
\begin{proof}
The embedding $W\to X$ is tight and holomorphic, therefore the 
corresponding Lie algebra homomorphism is $(\h2)$.
Lemma \ref{H2_hom} implies that under the identity homomorphism 
which embeds $\frakl$ as Lie subalgebra into $\frakg$, the element 
$Z_\frakl$ determined by the holomorphic tight embedding of a disc 
is identified with $Z_\frakg$.\\
Taking the exponential map implies the corresponding claim 
on the level of the associated Lie groups.
\end{proof}

\begin{cor}\label{herm_hull}
Suppose $X$ is of tube type,  
$f: Y \to X$ is a tight totally geodesic embedding 
induced by a tight homomorphism $\pi:H\to G$.
Suppose $W$ is a symmetric subspace of $X$ containing $f(Y)$.
Then the following are equivalent:
\begin{itemize}
\item{$W=\Hh(f)$}
\item{$\frakl = \<\pi (\frakh ),Z_\frakg\>$}
\item{$L=\Is(W)^\circ=\pi(H)\cdot Z_G$}
\end{itemize} 
\end{cor}
\begin{proof}
Assume that $W=\Hh(f)$.
The previous Lemma implies that $Z_\frakg \in \frakl$ 
and by minimality $\frakl=\<\pi_\frakh,Z_\frakg\>$. 
Taking the exponential map gives $L=\pi(H)\cdot Z_G$. 
On the other hand the inclusion $\frakl\to \frakg$ is an
$(\h2)$-homomorphism, 
therefore the symmetric subspace $W$ 
associated to $L$ is clearly 
holomorphically embedded into $X$ and 
it contains the image of $f$. 
Since the Lie algebra homomorphism associated to 
any other holomorphically embedded 
symmetric space $W'$ containing the image of $f$ must be a $(\h2)$-homomorphism
(Corollary~\ref{H2_hom}), it follows that 
$W\subset W'$ for all such $W'$, hence $W=\Hh(f)$.
 
\end{proof}

\begin{rem}
These characterizations of the Hermitian hull are not 
true if $X$ is not of tube type. 
Consider for example the canonical embedding of 
$\fraks\fraku (p,p) \to \fraks\fraku(p,q)$. It is holomorphic and tight, 
but $Z_{p,p}\notin \langle Z_{p,q}\rangle$.
\end{rem}

The above Corollary allows us to speak of the Hermitian hull $\Hh(\pi)$ of a 
tight homomorphism $\pi$ of Lie groups or of Lie algebras if the 
target group/algebra 
is associated to a symmetric space of tube type.
\begin{cor}\label{cor:centralizer_hull}
Suppose $G$ is the connected component of the isometry group of a 
Hermitian symmetric space of tube type and $\pi:H\to G$ a tight injective 
homomorphism. Then $\Hh(\pi) < L:= \zent_G(\zent_G(\pi(H))$.\\
Furthermore $\zent_G(\pi(H))= \zent_G(\Hh(\pi))$.
\end{cor}
\begin{proof}
According to Corollary \ref{herm_hull}, $\Hh(\pi)$ is the subgroup of $G$ 
generated by $\pi(H)$ and $Z_G$. 
By definition of $L$ we have $\pi(H)<L$. By Proposition \ref{centralizer} we have $\zent_G(\pi(H))<K_G$, 
therefore $Z_G$, which is in the center of $K_G$, commutes with each element of $\zent_G(\pi(H))$.
This proves the first claim.
 
For the second claim we have to prove that $\zent_G(\pi(H))\subset \zent_G(\Hh(\pi))$. 
The other inclusion $\zent_G(\Hh(\pi))\subset \zent_G(\pi(H))$ is trivial since $\pi(H) \subset \Hh(\pi)$.
But now $\zent_G(\Hh(\pi))= \zent_G(\pi(H)) \cap \zent_G(Z_G)=
\zent_G(\pi(H)) \cap K_G$, 
Since $\zent_G(\pi(H))\subset K_G$, we obtain the desired inclusion.
\end{proof}
\subsection{Irreducibility}
\begin{defi}
Let $\frakg$ be of tube type. 
A tight homomorphism $\pi: \frakh \to \frakg$ is called Hermitian 
irreducible if 
$\frakg = \Hh(\pi(\frakh))$.
%
\end{defi}

\begin{lemma}\label{irreducible}
Let $\frakg$ be a simple Lie algebra of tube type. 
Suppose 
$\pi:\fraks\frakl(2,\RR)\to \frakg$ is a tight Hermitian irreducible homomorphism. 
Then $\frakg=\fraks\frakp(2n,\RR)$ and the representation 
$\pi:\fraks\frakl(2,\RR)\to \frakg$ is irreducible.
\end{lemma}
\noindent
The following proof uses the notion of semiprincipal subalgebras and
of Satake diagrams. The reader who is not familar with these concepts
is advised to read the corresponding section of the 
Appendix before reading the proof.
\begin{proof}
By Corollary~\ref{herm_hull} we have $\frakg=\<\pi(\fraks\frakl(2,\RR)), Z_\frakg\>$.  
The complexification of the representation 
$\pi_\CC: \fraks\frakl(2, \CC) \to \frakg_\CC$ satisfies 
$\frakg=\<\pi_\CC(\fraks\frakl(2, \CC)), Z_\frakg\>=\frakg_\CC$. 
Since $Z_\frakg$ does not lie in the centralizer of 
$\frakh:=\pi_\CC(\fraks\frakl(2, \CC))$, the centralizer 
of $\frakh$ in $\frakg_\CC$ is trivial. 
Thus $\frakh$ is a semiprincipal three dimensional 
simple subalgebra in $\frakg_\CC$. These subalgebras were classified by 
Dynkin (see \cite{Dynkin, Onishchik_Vinberg_III} and the Appendix).\\
We make a case by case consideration:\\
1) When $\frakg= \fraks\frakp(2n,\RR)$, we have 
 $\frakg_\CC= \fraks\frakp(2n,\CC)$ and the semiprincipal subalgebra 
$\frakh$ is given by the irreducible representation of  
$\fraks\frakl(2,\CC) \to \fraks\frakp(2n,\CC)$.\\
2) Assume that $\frakg= \fraks\fraku(n,n)$, then 
$\frakg_\CC= \fraks\frakl(2n, \CC)$. In this case, the semiprincipal 
three dimensional simple subalgebra is principal and given 
by the irreducible representation of $\fraks\frakl(2,\CC)$. Since 
the dimension is $2n$ this irreducible representation is always contained 
in $\fraks\frakp(2n,\CC)$. Thus with   
$Z_{\fraks\fraku(n,n)} = Z_{\fraks\frakp(2n, \RR)}$, 
we have $ \frakg=<\pi(\fraks\frakl(2,\RR)), Z_\frakg>= \fraks\frakp(2n,\RR)$ and the 
representation is irreducible. \\
3) In the case where $\frakg=\fraks\frako(2,2n-1)$, then 
$\frakg_\CC= \fraks\frako(2n+1, \CC)$ with $\frakh$ being a principal 
three dimensional subalgebra defined by the irreducible representation 
of $\fraks\frakl(2,\CC)$. But any real irreducible representation 
of $\fraks\frakl(2,\RR)$ contained in $\fraks\frako(2,2n-1)$ is 
contained either in $\fraks\frako(2,3)\cong \fraks\frakp(4,\RR)$ or 
$\fraks\frako(2,1)\cong \fraks\frakp(2,\RR)$.

The remaining cases are more subtle, since there are semiprincipal 
three dimensional subalgebras which are not principal.\\
4) Assume that $\frakg$ is $\fraks\frako(2,2(n-1))$ or 
$\fraks\frako^*(2n)$. 
Recall how the Satake diagrams of these Lie algebras look like.
The Satake diagram of $\fraks\frako(2,2(n-1))$ is:\\
\smallskip
\bqn
\cnode(-2,0){3pt}{a}
\cnode(-1,0){3pt}{b}
\cnode*(0,0){3pt}{c}
\cnode*(1,0){3pt}{d}
\cnode*(2,0){3pt}{e}
\cnode*(2.7,0.5){3pt}{f}
\cnode*(2.7,-0.5){3pt}{g}
\ncline{-}{a}{b}
\ncline{-}{b}{c}
\ncline{-}{c}{d}
\ncline[linestyle=dotted]{-}{d}{e}
\ncline{-}{e}{f}
\ncline{-}{e}{g}
\eqn
\medskip 
The Satake diagram for 
$\fraks\frako^*(2n)$ is:\\
\smallskip
\bqn
\cnode*(-2,0){3pt}{a}
\cnode(-1,0){3pt}{b}
\cnode*(0,0){3pt}{c}
\cnode(1,0){3pt}{d}
\cnode*(1.7,0.5){3pt}{e}
\cnode(1.7,-0.5){3pt}{f}
\ncline{-}{a}{b}
\ncline{-}{b}{c}
\ncline[linestyle=dotted]{-}{c}{d}
\ncline{-}{d}{e}
\ncline{-}{d}{f}
\eqn 
\medskip

Since the subalgebra $\frakh$ is the image of the complexification of a 
representation of $\fraks\frakl(2,\RR)$ into $\frakg$, the 
image of the characteristic element $h$ is contained in the root spaces 
corresponding to the white circles in the respective 
Satake diagram. In particular  
$\frakh$ cannot be a principal three dimensional subalgebra, since the
characteristic element of the principal three dimensional subalgebra
is given by the Dynkin diagram with all nodes being white. 
On the other hand, the characterstic elements of the 
semiprincipal three dimensional subalgebras of $\fraks\frako(2n,\CC)$ 
are given by the following diagrams with $k$ black circles, 
$k=1, \dots [\frac{n-2}{2}]$:\\
\smallskip 
\bqn
\cnode(-4,0){3pt}{a}
\cnode(-3,0){3pt}{b}
\cnode(-2,0){3pt}{c}
\cnode(-1,0){3pt}{d}
\cnode*(0,0){3pt}{e}
\cnode(1,0){3pt}{f}
\cnode*(2,0){3pt}{g}
\cnode(3,0){3pt}{h}
\cnode*(4,0){3pt}{i}
\cnode(4.7,0.5){3pt}{j}
\cnode(4.7,-0.5){3pt}{k}
\ncline{-}{a}{b}
\ncline[linestyle=dotted]{-}{b}{c}
\ncline{-}{c}{d}
\ncline{-}{d}{e}
\ncline{-}{e}{f}
\ncline{-}{f}{g}
\ncline[linestyle=dotted]{-}{g}{h}
\ncline{-}{h}{i}
\ncline{-}{i}{j}
\ncline{-}{i}{k}
\eqn
\medskip

Comparing this to the Satake diagrams, 
we see that $\frakh$ is not a 
semiprincipal three dimensional subalgebra of $\fraks\frako(2n,\CC)$. 
Therefore $\frakg$ cannot be $\fraks\frako(2,2(n-1))$ or
$\fraks\frako^*(2n)$.

5) The Satake diagram of 
$\frake_{VII}$ looks like:\\
\smallskip
\bqn
\cnode(-3,0.5){3pt}{a}
\cnode(-2,0.5){3pt}{b}
\cnode*(-1,0.5){3pt}{c}
\cnode*(0,0.5){3pt}{d}
\cnode*(1,0.5){3pt}{e}
\cnode(2,0.5){3pt}{f}
\cnode*(0,-0.5){3pt}{g}
\ncline{-}{a}{b}
\ncline{-}{b}{c}
\ncline{-}{c}{d}
\ncline{-}{d}{e}
\ncline{-}{e}{f}
\ncline{-}{d}{g}
\eqn 
\medskip

As before, since 
the subalgebra $\frakh$ is the image of a complexification of a 
representation of $\fraks\frakl(2,\RR)$ into $\frakg$, the 
image of the characteristic element $h$ is contained in the root spaces 
corresponding to the white circles in the Satake diagram. Thus $\frakh$ 
cannot be 
the principal three dimensional subalgebra. 
The characterstic element of the 
semiprincipal three dimensional subalgebras of $\frake_7$ 
are given by the following 
diagrams:

\smallskip
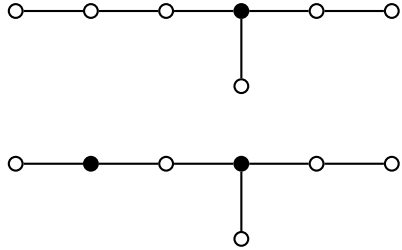
\bqn
\cnode(-3,0.5){3pt}{a}
\cnode(-2,0.5){3pt}{b}
\cnode(-1,0.5){3pt}{c}
\cnode*(0,0.5){3pt}{d}
\cnode(1,0.5){3pt}{e}
\cnode(2,0.5){3pt}{f}
\cnode(0,-0.5){3pt}{g}
\ncline{-}{a}{b}
\ncline{-}{b}{c}
\ncline{-}{c}{d}
\ncline{-}{d}{e}
\ncline{-}{e}{f}
\ncline{-}{d}{g}
\eqn 
\medskip

\smallskip
\bqn
\cnode(-3,0.5){3pt}{a}
\cnode*(-2,0.5){3pt}{b}
\cnode(-1,0.5){3pt}{c}
\cnode*(0,0.5){3pt}{d}
\cnode(1,0.5){3pt}{e}
\cnode(2,0.5){3pt}{f}
\cnode(0,-0.5){3pt}{g}
\ncline{-}{a}{b}
\ncline{-}{b}{c}
\ncline{-}{c}{d}
\ncline{-}{d}{e}
\ncline{-}{e}{f}
\ncline{-}{d}{g}
\eqn 
\medskip
\\
Thus the subalgebra $\frakh$ cannot be a semiprincipal subalgebra of 
$\frake_7$, and hence $\frakg$ cannot be 
$\frake_{VII}$.

Hence we obtain the result: $\frakg=\fraks\frakp(2n,\RR)$ and
$\pi:\fraks\frakl(2,\RR)\to \frakg$ is given by the irreducible
representation of $\fraks\frakl(2,\RR)$.
\end{proof}
\begin{cor}
If $\pi:\fraks\frakl(2,\RR)\to \frakg$ is Hermitian irreducible and $\Hh(\pi)$
simple 
then $\pi$ is unique up to conjugacy.
\end{cor}
\begin{proof}
This follows from Lemma~\ref{irreducible} since there is
 up to equivalence only one irreducible representation 
of $\fraks\frakl(2,\RR)$ in each dimension.
\end{proof}

\begin{prop}\label{hull}
Let $\pi:\fraks\frakl(2,\RR)\to \frakg$ be a tight homomorphism. 
Then $\Hh(\pi)= \oplus_{i=1}^k \fraks\frakp(2n_i,\RR)$, 
with $\sum_{i=1}^k n_i \leq \rk_\RR(\frakg)$.
\end{prop}
\begin{proof}
The Hermitian hull $\Hh(\pi)$ is a semisimple Lie algebra of Hermitian type, 
 $\Hh(\pi)= \oplus_{i=1}^n \frakg_i$, where all $\frakg_i$ are of tube type. 
The representations  $\pi_i:\fraks\frakl(2,\RR)\to \frakg_i$ are again tight homomorphisms 
and $\frakg_i= \Hh(\pi_i)$. Therefore 
by Lemma~\ref{irreducible}, $\frakg_i= \fraks\frakp(2n_i,\RR)$. 
\end{proof}

As Corollary we obtain
\begin{thm}\label{thm:hull}
Suppose that $X$ is a Hermitian symmetric space of noncompact type and
 $f:\DD\to X$ is a tight totally geodesic embedding. 
Then the Hermitian hull $\Hh(f) = \Pi_{i=1}^k Y_i$ 
is a product of Hermitian symmetric subspaces $Y_i$ of $X$, 
where $Y_i$ are Hermitian symmetric spaces 
associated to the symplectic groups ${\rm Sp}(2n_i,\RR)$. 
Moreover, $\sum_{i=1}^k n_i \leq r_X$.
\end{thm}
\begin{proof}
This is the analog of Proposition~\ref{hull} 
for the totally geodesic embeddings associated to the Lie algebra 
homomorphism.
\end{proof}

\subsection{Obstructions to tight embeddings}
Assume that $X,Y$ are of tube type.
\subsubsection{Numerical obstructions}
A tight embedding $f:Y\to X$ satisfies $f^* \b_X = \frac{r_X}{r_Y} \b_Y$.
On the other hand the set of values of $\b_X$ is $r_X, r_X-2, \dots , -r_X$, whereas 
the set of values of $\b_Y$ is $r_Y, r_Y-2, \dots , -r_Y$.
Thus we get that $\frac{r_X}{r_Y} \in \ZZ$ is a necessary condition for the existence 
of a tight embedding $f: Y\to X$.

\subsubsection{Holomorphic tight embeddings}
Since holomorphic tight embeddings correspond to $(\h2)$-homomorphisms, 
they have been completely classified by Satake and Ihara, 
see \cite{Satake_hol_65,Satake_hol_68,Satake_hol_69,Ihara_65,Ihara_66,Ihara_65_suppl}.

This gives also restrictions to the existence of tight embeddings since 
the Hermitian hull of the restriction of a tight embedding of both,
 a tight holomorphic disc  
and an irreducible tight disc in $Y$,
are Hermitian symmetric spaces associated to ${\rm Sp(2n, \RR)}$, for different $n's$, which are 
tight holomorphic and totally geodesic subspaces of $X$.

\section{Open questions}
\subsection{Compact dual}
The criteria for tightness, which we have formulated on the level of
Lie algebra homomorphisms, 
makes also sense for 
Lie algebra homomorphisms of the Lie algebra corresponding to the compact dual Hermitian symmetric space. 
By the duality of compact and noncompact symmetric spaces every tight totally geodesic embedding of a 
Hermitian symmetric space of noncompact type 
into a Hermitian symmetric space of noncompact type corresponds 
to a totally geodesic embedding of the compact duals.
It seems to be of some interest 
to understand the properties of this class of 
totally geodesic embeddings in terms of the 
geometry of compact Hermitian symmetric spaces.   \\

\subsection{Tight embeddings with respect to other classes}
As we already pointed out, the definition of a tight embedding can be
extended to the definition of $\a$-tightness for any bounded cohomology
class $\a$. 
Unfortunately there are not many bounded cohomology classes that are well understood.
To get informations about tight homomorphisms with respect to the bounded K\"ahler class we made not only use of the 
boundedness of this class, but of a geometric characterization of 
the corresponding maximal triples. 
An example of a tight homomorphism is studied by Soma \cite{Soma_third,Soma_tame}. He considers hyperbolic three-manifolds, 
thus representations into ${\rm PSL}(2,\CC)$ and 
proves rigidity results for representations which are tight with respect to the bounded cohomology class 
defined by the volume form on $\HH^3$.

Since the classes in $\hcb^*(G, \RR)$ corresponding to 
characteristic classes of flat $G$-bundles are known to be bounded \cite{Gromov_82,Bucher}, these would 
be the first candidates for which one might study tight embeddings. 
We mention some classes that could be interesting: \\
1) Quaternionic hyperbolic space: The volume form of the 
quaternionic hyperbolic line corresponds to the first Pontrjagin class. 
Calibrations with respect to this class where studied e.g. in \cite{Berger}. 
One might expect that one could obtain rigidity results for respresentations in to the 
isometry group of quaternionic hyperbolic space.  \\
2) Polylogarithms: Goncharov proves a trilogarithm identity 
\cite{Goncharov_config}, which shows that the Borel class in 
$\h_c^5({\rm SL}(3,\CC), \CC)$ is bounded. Gelfand and MacPherson
define \cite{Gelfand_MacPherson} a class of 
generalized simplices in Grassmannian manifolds and relate these simplices to certain polylogarithms and to certain 
characteristic classes. \\
3) Finkelberg gives an explicit bounded representative for a $4$-cocycle on ${\rm O}(n,n)$ in \cite{Finkelberg}.\\
4) The Euler class defined a bounded cohomology class in $\hcb^n({\rm SL}(n,\RR),\RR)$ for $n$ even. A quite explicit respresentative 
is given by Turaev and Ivanov \cite{Turaev_Ivanov}.\\
 
\subsection{Calibrations}
Notice the similarity of the classification 
of holomorphic embeddings of the disc, done by Satake
\cite{Satake_book} and our classification of tight embeddings of the disc. 
This suggests that tight embeddings could be a bounded cohomological 
avatar of calibrations. 
Similar to calibrations, tight embeddings are defined by considering the equality
case of an inequality. And similar to the problem in the study of calibrations, where one has to compute 
the comass, here one encounters the problem of computing 
the sup norm of the bounded cohomology class respectively a cocycle
representing it. 

\subsection{Tight mappings}
Instead of totally geodesic embeddings we might study 
tight embeddings of quotients , i.e. of locally Hermitian symmetric
spaces or even more generally immersion or even more general maps into 
Hermitian symmetric spaces or their quotients. 
This question seems to be related to work of Mok and Eyssidieux and of
Clozel and Ullmo on holomorphic subvarieties of locally Hermitian
symmetric spaces. 


\vskip1cm
\chapter{Group actions on Hermitian symmetric spaces}
The results about Hermitian symmetric spaces and tight embeddings 
may be used to study isometric actions of finitely generated groups $\G$ on 
Hermitian symmetric spaces $X$.
An action of a group $\G$ on $X$ by isometries 
is given by a homomorphism 
\bqn
\rho:\G\to G={\rm Is}(X)^\circ.
\eqn
Let $G$ be a noncompact connected simple Lie group of Hermitian type, 
i.e. the associated symmetric space $X=G/K$ is an irreducible Hermitian 
symmetric space of noncompact type. 
Let $\G$ be a finitely generated group.
We denote by $k_X^b \in \hcb^2(G)$ the (normalized) bounded K\"ahler class. 

Given a representation $\rho: \G\to G$, the canonical pull-back map 
in bounded continuous cohomology 
$\rho^*_b: \hcb^*(G) \to \hb^*(\G)$ 
defines an invariant with values in the second bounded cohomology 
group of $\G$, 
called {\em bounded K\"ahler class of $\rho$}: 
\bqn
\rho^*_b(k_X^b)\in \hb^2(\G).
\eqn

We already used this invariant to define tight embeddings. Now we use it to study representations of 
finitely generated groups $\G$.

\section{Zariski dense representations}\label{sec:zariski_reps}
Let $G$ be a connected simple Lie group with finite center such that 
the associated Hermitian symmetric space $X$ is not of tube type.
Furthermore we let $\G$ denote a finitely generated group. 
An action of $\G$ by isometries on $X$ which does not leave a proper closed 
subspace of $X$ invariant is completely determined by the bounded K\"ahler 
class of its representation. 
This has been proven in \cite{Burger_Iozzi_supq} in the 
case where $G$ is isomorphic to ${\rm SU}(p,q)$, $p\neq q$. 
With the characterization of non-tube type Hermitian symmetric spaces 
(see Chapter~\ref{sec:shilov}) and the 
construction of the Hermitian triple product 
(see Chapter~\ref{sec:cocycle}) we obtain the generalization 
of the results in \cite{Burger_Iozzi_supq} to 
arbitrary simple Lie groups $G$ which are not of tube type. These include 
groups locally isomorphic to ${\rm SU}(p,q)$, $p\neq q$, to ${\rm SO}^*(4k+2)$ and 
$E_{III}$.
\begin{thm}\label{thm:z_dense}
Let $\rho:\G\to G$ be a representation with Zariski dense image. 
Then 
\begin{enumerate}
\item{$\rho^*(k_X^b) \neq 0$.}
\item{$\rho^*(k_X^b)$ determines $\rho$ up to $G$-conjugacy.}
\end{enumerate}
\end{thm}
\begin{cor}
If $\dim \hb^2(\G)< \infty$ then there are up to equivalence only 
finitely many Zariski dense representations 
$\rho:\G \to G$.
In particular, if $\hb^2(\G)=0$, 
there are no Zariski dense representations of $\Gamma$ into $G$.
\end{cor}

\begin{rem}
The results admit a slightly stronger formulation that can be found in \cite{Burger_Iozzi_supq}.
\end{rem}

Most of the arguments of \cite[\S~5 and \S~7]{Burger_Iozzi_supq} for
$G={\rm PSU}(p,q)$ apply to the general case. 
Note that we encounter two difficulties as soon as the rank of $X$ is $r_X\geq 2$. 
First, the Bergmann cocycle is not defined on the doubly ergodic 
Furstenberg boundary $G/P$, but only on $G/Q$, where $Q$ is a maximal parabolic 
subgroup and hence is not amenable. Second, 
the cocycle is not a strict cocycle since it is only defined on 
triples of pairwise transverse points. 
We can surpass these difficulties using the arguments of 
 \cite[\S~5]{Burger_Iozzi_supq}.
For convenience we recall the main steps of the proof. 

\subsection{The boundary map}
An important step is the construction of the boundary map. 
We have already realized the 
bounded K\"ahler class $k^b_X\in \hcb^2(G)$ 
by a class $[\beta_\Dd] \in \h^* (\binfty(\cs^{(*+1)})^G)$ 
defined by the Bergmann cocycle  $\b:\cs^{(3)} \to \RR$.
We would like to realize the bounded K\"ahler invariant as pull-back of 
this Bergmann cocycle via a suitable boundary map.

Let $(S,\nu)$ be a doubly ergodic Poisson boundary for $\G$ 
(for the construction of this boundary for a finitely generated group
$\G$ see for example 
\cite{Burger_Monod_GAFA, Burger_Iozzi_supq}). 
Let $P$ be a minimal parabolic subgroup of $G$. 
\begin{prop}
\cite[Proposition~6.2]{Burger_Iozzi_supq}\label{prop:boundary_map} 
There exists a measurable $\rho$-equivariant 
map $\phi: (S,\nu) \to G/P$. Furthermore, 
for $\nu$-almost every $(x,y)\in S$ the image points are transverse, i.e. 
$(\phi(x), \phi(y))\in B\subset G/P\times G/P$.
\end{prop}
\begin{proof}
Since the actions of $\G$ on $S$ and of $G$ on $G/P$ are amenable, 
there exists a measurable $\rho$-equivariant 
map 
\bqn
\psi:S \to \Mm^1(G/P).
\eqn
Due to the Zariski density of $\rho(\G)$ the map $\phi$ takes values in the space of Dirac 
measures (see \cite[Theorem~6.2]{Burger_Iozzi_supq}), 
which we identify with $G/P \subset \Mm^1(G/P)$. 

It remains to prove the transversality condition.
For a tuple $(z_1, z_2)\in G/P \times G/P$ let 
$\codim(z_1, z_2)$ denote the codimension of the cell of the 
Bruhat decomposition of $G/P \times G/P$ 
containing $(z_1, z_2)$. Then  $\codim(z_1, z_2)= 0$ if and only if $z_1$ and $z_2$ are transverse.
Consider the equivariant bounded measurable map 
\bqn
S\times S\to \NN,\quad (s_1, s_2) \mapsto \codim(\phi(s_1), \phi(s_2)).
\eqn 
Since $\G$ acts ergodically on $S\times S$ this map is essentially constant.
If it is zero, then the map $\phi$ satisfies the transversality condition. 
Assume therefore that it is nonzero. 
Define $\nt(z):=\{ w\in G/P \, |\, w \text{ is not transverse to } z\}$. 
Then $\nt(z)$ is a proper Zariski closed subset of $ G/P$.
For almost all $s_1\in S$ we have that $\phi(s_2)  \in \nt(\phi(s_1))$ 
for almost all $s_2\in S$. Fix one such $s_1$. 
The essential image ${\rm EssIm}(\phi) \subset G/P$ is a closed and $\rho(\G)$-invariant subset of $\cs$.
Assuming that $\codim(\cdot, \cdot)$ is essentially constant and nonzero, 
it follows 
\bqn
{\rm EssIm}(\phi) \subset \bigcap_{\g\in \G} \rho(\g) \nt(\phi(s_1)) =: L\subset G/P.
\eqn
But $L$ is a $\rho(\G)$-invariant proper Zariski closed subset of $G/P$.  
This is a contradiction to the Zariski density of $\rho(\G)$.
Thus  the map $\codim(\cdot, \cdot)$ is essentially the zero 
mapping and 
$\phi:S\to G/P$ satisfies the transversality condition. 
\end{proof}

Since transversality is preserved under the natural projection $G/P \to G/Q$, 
the induced map $\phi: S\to G/Q=\cs$ also satisfies the transversality condition.
This is the map we are going to use to implement the pull-back.
Since $\phi$ maps almost all $x,y,z \in S$ to a triple in $\cs^{(3)}$, 
the pull-back of the Bergmann cocyle $\b_\Dd$ via $\phi$ is well-defined 
and defines a measurable $\G$-invariant alternating essentially bounded 
cocycle 
$\phi^*\b: (S)^3\to \RR$.

Denote by 
$\Zz\la((S)^3)^\Gamma$ the space of measurable $\G$-invariant alternating essentially bounded cocycles on $(S)^3$. 
Then 
there is an isomorphism $\hb^2(\Gamma,\RR) \simeq
\Zz\la((S)^3)^\Gamma$ (see Propsition~\ref{prop:L_alt} in Chapter 1) .
The pull-back $\phi^*\beta\in \Zz\la((S)^3)^\Gamma$ of the Bergmann cocycle via $\phi$ is related to 
the pull-back $\rho^* k_X^b\in\hb^2(\Gamma,\RR)$.
\begin{lemma}\cite[Theorem~8.1.]{Burger_Iozzi_supq}
The cocycle $\phi^* \beta$ corresponds to $\rho^* k_X^b$ under the
canoncical isomorphism $\hb^2(\Gamma,\RR) \simeq \Zz\la((S)^3)^\Gamma$.
\end{lemma}

\subsection{The proof of Theorem~\ref{thm:z_dense}}
\begin{proof}
The main geometric ingredient in the proof is the different behaviour 
of the Bergmann cocycle respectively the Hermitian triple product with respect to tube type and 
non-tube type domains, 
see Proposition~\ref{infinite_values} and Lemma~\ref{lem:nonconstant} in 
Chapter~\ref{sec:cocycle}. 
Since we are going to apply arguments from Zariski topology, the algebraic realization 
of the cocycle given by the Hermitian triple product is crucial.

Let $\rho_i:\G\to G$, $ i=1,\dots, k$ be representations with Zariski dense image. 
Assume that there exists $m_i\in \ZZ\backslash \{0\}$ such that 
$\sum_{i=1}^k m_i \rho_i^*(k_b^X) = 0$. 
Let 
\bqn
\phi_i: S\to \cs = \gG/\qQ(\RR)
\eqn
be the corresponding $\rho_i$-equivariant boundary maps. 
For almost all $x_1,x_2,x_3$ we have that 
\bqn
(\phi_i( x_1),\phi_i(x_2),\phi_i(x_3)) \in \gG/\qQ(\RR)^{(3)}
\eqn
and 
\bqn
\Pi_{i=1}^k \HTP_\CC(\phi_i( x_1),\phi_i(x_2),\phi_i(x_3))^{m_i}= [1] \in \CC^\times \backslash A^\times.
\eqn

Consider the representation 
\bqn
\rho= (\rho_1,\dots, \rho_k) : \G\to \Pi_{i=1}^k \gG_i, \quad
\g\mapsto (\rho_i(\g))_{i=1\dots k}, 
\eqn
where $\gG_i = \gG$ for all $i$ 
and the $\rho$-equivariant measurable map 
\bqn
\phi: S\to  \Pi_{i=1}^k \gG/\qQ(\RR), \quad s\mapsto (\phi_i (s))_{i=1\dots k}.
\eqn

Let $H=\ol{\rho(\G)}^Z$ be the Zariski closure of $\rho(\G)$ in $\Pi_{i=1}^k \gG_i$.
The essential image ${\rm EssIm}(\phi)$ is a closed $\rho(\G)$-invariant subset of $\Pi_{i=1}^k \gG/\qQ(\RR)$.
Fixing $x_1, x_2 \in S$ such that the above equality for the Hermitian triple product holds for almost all $x_3$, 
we define 
$P_i:= P_{\phi_i(x_1), \phi_i(x_2)}$ 
and 
$\Bb_i:= \Bb_{\phi_i(x_1), \phi_i(x_2)}$, 
$\Bb= \Pi_{i=1}^k \Bb_i$. 
 
Then 
\bqn
{\rm EssIm}(\phi)\cap \Bb \subset \{(z_1, \dots, z_k)\in \Bb 
\,|\,\Pi_{i=1}^k P_i^{m_i}  = [1] \} \subset \Pi_{i=1}^k \gG/\qQ(\RR).
\eqn
The set $L:= \{(z_1, \dots, z_k)\in \Bb \,|\,\Pi_{i=1}^k P_i^{m_i}  = [1] \}$ 
is   
a proper Zariski closed subset of $\Bb$ (Lemma~\ref{lem:nonconstant} in Chapter~\ref{sec:cocycle}).
Thus ${\rm EssIm}(\phi)$ is contained in 
a proper Zariski closed subset of $\gG/\qQ(\RR)$ 
and the $\RR$-algebraic subgroup $H$ is 
properly contained in $\Pi_{i=1}^k \gG_i$. 

We can now conclude the statements of the theorem. 
The situation that $\rho^*(k_b^X) = 0$ corresponds to $k=1$. 
But then $H(\RR)=\ol{\rho(\G)}^Z(\RR)$ is  contained in $G$ as a proper Zariski 
closed subgroup. 
This contradicts the Zariski density of $\rho(\G)$, thus 
$\rho^*(k_b^X) \neq 0$

Assuming that there exist two representations with 
$\rho_1^*(k_b^X)=\rho_2^*(k_b^X)$ corresponds to 
the situation that $k=2$, $m_1= -m_2=1$. 
Then  $\pr_i(H(\RR))=G$ , $ i=1,2$, since $\rho_i(\G)$ is Zariski dense in 
$\gG_i$ 
by assumption, 
but  $\rho(\G)$ is not Zariski dense in $\gG$.
Since the groups $\gG_i$  are simple  there exists an isomorphism$\rho_{12}: \gG_1\to \gG_2$ defined over $\RR$ 
such that $\rho_2(\g) = \rho_{12}(\rho_1(\g))$.
This isomorphism induces an isometry 
$T_{12}:X_1\to X_2$ of the symmetric spaces associated to $G_i \cong G$. 
This isometry is holomorphic. 
If not, we would have that $2\rho_1^*(k_b^X)=0$ contradicting the 
first statement of Theorem~\ref{thm:z_dense}. 
\end{proof}

\section{Limit sets in the Shilov boundary}
Let $X$ be a Hermitian symmetric space and $G={\rm Is}^\circ(X)$ the
connected component of its isometry group.
Let $\G<G$ be a (discrete) group. 
We may identify $G$ with the automorphism group $\aut(\Dd)$ of its Harish-Chandra realization $\Dd$.

Define the {\em topological limit set} $L_t(\G, x)\subset \partial\Dd$ in 
the topological boundary of the Harish-Chandra realization $\Dd$ of $X$ by 
taking the topological closure in $\CC^N$ of the $\G$-orbit of a point 
$x\in \Dd$, $\overline{\G x}^t$. The topological limit set 
\bqn
L_t(\G, x)=\overline{\G x}^t\cap\partial \Dd
\eqn
depends on the base point.

\begin{defi}
The \emph{topological Shilov limit set} $\Ll_t(\G)$ is defined to be 
\bqn
\Ll_t(\G)= L_t(\G, x)\cap \cs ,\,x\in \Dd. 
\eqn
\end{defi}

\begin{lemma}
The topological Shilov limit set $\Ll_t(\G)$ is independent of the base 
point $x\in X$.
\end{lemma}
\begin{proof}
Let $(\g_n)_{n\in \NN} \subset \G$ such that $x_n=\g_n x\to \xi \in
\Ll_t(\G, x)$ for $n\to \infty$. 
Let $y\in \Dd$  
and consider $y_n=\g_n y$. Since $\g_n$ acts by isometries 
$\d(x_n, y_n)=\d(x,y)$, which is bounded. Hence by
Lemma~\ref{geodesics} of Chapter 3.
$\xi=\lim_{n\to \infty} x_n$ and $\eta:=\lim_{n\to \infty} y_n$ lie in the same boundary component. 
But since $\xi\in \cs$ its boundary component is one point, hence $\eta=\xi$.
\end{proof}

Recall the definition of the 
{\em Benoist limit set} of $\G$ in $G/P_\theta$ from \cite{Benoist}. 

\begin{defi}
The Benoist limit set $\Lambda_\theta (\G) \subset G/P_\theta$ of $\G<G$ is the set of 
all points $\xi \in G/P_\theta$ for which there exists a sequence $(\g_n) \subset \G$ such that 
$\lim_{n\to \infty} (\g_n)_* \nu_\theta= \delta_\xi$, where $\delta_\xi$ is the Dirac measure of $\xi$ 
and $\nu_\theta$ the unique $K$-invariant probablility 
measure on $G/P_\theta$. 
\end{defi}

\begin{prop}\cite[Lemme~3.6]{Benoist}
Let $\G$ be a Zariski dense subgroup of $G$. 
Then
\begin{enumerate}
\item{$\Lambda_\theta (\G) \subset G/P_\theta$ is Zariski dense.}
\item{Every (nonempty) closed $\G$-invariant subset $F\subset G/P_\theta$ contains 
$\Lambda_\theta (\G)$.  The action of $\G$ on $\Lambda_\theta (\G)$ is minimal, 
$\Lambda_\theta (\G)$ has no isolated points and either $\Lambda_\theta (\G)= G/P_\theta$ or 
$\Lambda_\theta (\G)$ has empty interior.}
\end{enumerate}
\end{prop}

The topological Shilov limit set  is related to the Benoist limit set. 
Recall that there is an isomorphism of the topological
compactification $\ol{\Dd}$ with the 
$\theta_r$-Satake-Furstenberg compactification.
From this we get 

\begin{lemma}\label{lem:limitsets}
The topological Shilov limit set $\Ll_t(\G)$ and the Benoist limit set $\Lambda_{\theta_r}(\G)$ coincide. 
\end{lemma}
\begin{proof}
The isomomorphism of the topological compactification of $X$ by
$\ol{\Dd}$ 
with the $\theta_r$-Furstenberg compactification identifies the 
Shilov boundary $\cs\subset \partial \Dd$ with the space of Dirac
measure on $G/P_{\theta_r}$.
We work with the $\theta_r$-Furstenberg compactification.
Let $\xi\in \Ll_t(\G)$ be a point in the Shilov limit set then, by
definition of the $\theta_r$-Furstenberg compactification, 
there exists a sequence $\g_n \in \G$ such that 
$\g_n 0\to \xi\in \cs$ in $\overline{\Dd}$. Therefore $(\g_n)_* \nu_{\theta_r}\to \delta_\xi \in \Mm(G/P_{\theta_r})$, 
hence $\xi \in \Lambda_{\theta_r}(\G)$.
On the other hand, suppose that $\xi \in \Lambda_{\theta_r}(\G)$, then there exists a sequence $\g_n \in \G$ 
with  $(\g_n)_* \nu_{\theta_r}\to \delta_\xi \in \Mm(G/P_{\theta_r})$, hence $\g_n 0\to \xi\in \cs$.
\end{proof}

%
\begin{lemma}\label{lem:benoist}
Assume that $\rho:\G\to G$ is Zariski dense and discrete.
Then the essential image of the boundary map $\phi: M\to G/P$ coincides with the Benoist limit set 
$\Lambda_{\emptyset}(\G)$.
\end{lemma}
\begin{proof}
The essential image of $\phi$ is a closed $\G$-invariant subset of $G/P$. 
By \cite[Lemme~3.6.]{Benoist} it therefore contains 
$\Lambda_{\emptyset}(\G)$. On the other hand, the essential image is itself minimal, so they coincide. 
\end{proof}

\begin{cor}
If $\rho:\G\to G$ is tight and discrete then the topological Shilov limit set 
$\Ll_t(\rho(\G))$ is a nonempty, closed, connected, $\rho(\G)$-minimal 
subset of $\cs$.
\end{cor}
\begin{proof}
The assertions for Zariski dense representation are proven in
\cite[Lemme~3.6.]{Benoist}. 
The topological Shilov limit set 
$\Ll_t(\rho(\G))$ is the image of the 
topological Shilov limit set of $\rho(\G)$ in the 
Shilov boundary of the 
Hermitian symmetric space associated 
to the connected component of the Zariski closure $L =  \ol{\rho(\G)}^Z (\RR)^0$
under the extension of the 
tight 
totally geodesic embedding associated to $L\to G$.
\end{proof}

\section{Structure of the representation variety}
We recall some structures of the space of homomorphism 
 $\hom(\G, G)$. For proofs see \cite{Goldman_84,Goldman_Millson}.
Let $\G$ be a finitely generated group 
and $G$ the adjoint group of a semisimple real Lie algebra. 
Then, as already mentioned, $G$ 
inherits an algebraic structure. 
We consider the space of all 
homomorphisms $\hom(\G, G)$ equipped with the compact open topology. 
Then a sequence of 
homomorphism $\rho_i \in \hom(\G, G) $ converges to $\rho$ if and only if 
$\rho_i(\g) \to \rho(\g)$ 
for all $\g\in \G$. In particular, convergence can be checked on 
a set of generators of $\G$. 
Choosing a set of generators $\{\g_1, \dots, \g_n\}$ of $\G$, 
the set $\hom(\G, G)$ may be equipped with an algebraic structure, 
which is induced from the algebraic structure from $G$. 
Namely, a homomorphism $\rho$ 
is determined by prescribing 
$n$ elements $g_1,  \dots  , g_n \in G$, 
such that $\rho(\g_i) = g_i$ for all $i= 1, \dots, n$.
This gives an embedding 
 \bqn
\hom(\G, G)\to G^n, 
\eqn
where the latter carries the algebraic structure coming 
from $G$.

The relations $R$ of $\G$ define polynomial equations $\Pp$ in $G^n$ 
And we may identify $\hom(\G, G)$ with $G^n/\Pp$.
The set $\hom(\G, G)$ with the algebraic structure induced from this identification 
is called {\em representation variety} and is denoted by  
\bqn
\Rr(\G, G) :=  G^n/\Pp.
\eqn

The group $G$ acts by conjugation on $\Rr(\G, G)$, 
$ (g\rho)(\g) := g^{-1}\rho(\g) g$. 
The quotient  
\bqn
\Xx(\G, G) := \Rr(\G, G)/G
\eqn
is called {\em character variety}.

This action of $G$ on $\Rr(\G, G)$ is in general not proper 
but we may define the subset of nonparabolic representations 
$\Rr_{np}(\G,G)$ as the set of 
representations $\rho$, such that $\rho(\G)$ 
has no fixed point in the geometric boundary of the symmetric space 
$X$ associated to $G$. This is equivalent to 
requiring that the Zariski closure of the image of $\G$ is a 
reductive subgroup (see \cite{Johnson_Millson}).
The subset $\Rr_{np}(\G,G)\subset \Rr(\G,G)$ is open. The action 
of $G$ on $\Rr_{np}(\G,G)$ is proper and the quotient 
\bqn
\Xx_{np}(\G,G)= \Rr_{np}(\G,G) / G
\eqn
is a separated, locally compact space 
with countable basis. 

Define by $\Rr_{fd}(\G, G)$ the set of all faithful representations with 
discrete image. 
If $\G$ does not contain a normal abelian subgroup of finite index, 
the subset $\Rr_{fd}(\G, G)\subset \Rr(\G, G)$ is closed 
(see \cite{Goldman_Millson}). 
In particular, the subset of the character variety 
\bqn
\Xx_{np, fd} := (\Rr_{np}(\G,G)\cap\Rr_{fd}(\G, G)) /G
\eqn
is a separated, 
locally compact space with countable basis.

The action of $G$ on $\Rr(\G, G)$ is part of an action of 
$\aut(\G)\times G$ on  $\Rr(\G, G)$, defined by
\bqn
((\phi, g)\rho)(\g)= g^{-1}\rho(\phi(\g)) g.
\eqn 
The action of the inner automorphisms ${\rm Inn}(\G)$ 
is absorbed by the action of $G$, 
thus on $\Xx(\G, G)$ the action reduces to the action of the group of outer 
automorphism ${\rm Out}(\G)$. 
This action preserves $\Xx_{np}(\G,G)$ and $\Xx_{np, fd}$.

\subsection{Representations and non tube type}
Let $G$ be a simple Lie group of Hermitian type, which is not of tube type.
Denoting by $\Xx_{Zd}(\G, G)$ 
the equivalence classes of Zariski dense representations in the 
character variety, we obtain from Theorem~\ref{thm:z_dense} 
an injective map 
\bqn
K: \Xx_{Zd}(\G, G) \to \hb^2(\G).
\eqn
In particular, if $\hb^2(\G) =0$, there are no representations of $\G$ 
into $G$ having Zariski dense image. 

\begin{rem}
This is not true for representations into Lie groups $G$ that are of tube type. 
We will later see, that there are for example many different 
Zariski dense representations of surface groups into Lie groups of tube type, 
which have all the same bounded K\"ahler invariant. 
\end{rem}

\subsection{Tight representations}
Let $G$ be a connected semisimple Lie group of Hermitian type. 
Assume 
that the representation $\rho:\G\to G$ is tight, i.e.
\bqn
||\rho_b^*(k_X^b)||=||k_X^b||.
\eqn

Define the subset consisting of tight representations in the 
representation variety 
\bqn
\Rr_t(\G, G) :=\{\rho \in \Rr(\G, G)\,|\, ||\rho_b^*(k_X^b)||=||k_X^b||\}.
\eqn 

As a consequence of Theorem~\ref{thm:tightness} of Chapter 6 
we obtain that the 
set of tight representations conists entirely of 
non parabolic representations.
\begin{cor}
\bqn
\Rr_t(\G, G)\subset \Rr_{np}(\G, G).
\eqn
\end{cor}

On $\Rr_t(\G, G)$ we obtain two different (not necessarily algebraic) 
stratifications.
The first one by considering the subsets $\Rr_t(\G, G_1)$, where $G_1= \zent_G(\zent_G(\rho(\G)))$. 
The other by considering subsets $\Rr_t(\G, H_1)$, where $H_1 = \Hh(\rho)$ is the Hermitian hull of $\rho$.
Since the centralizer in $G$ 
of the image of a tight representation coincides with the centralizer of its Hermitian hull, 
the stratification by the Hermitian hull is 
finer than the stratification of $\Rr(\G, G)$ by $G_1= \zent_G(\zent_G(\rho(\G)))$.


\vskip1cm
\chapter{Representations of surface groups}
We now focus on representations of surface groups. 
Combining the results on Zariski dense and tight representations 
with further techniques from bounded cohomology, we are able to 
obtain geometric information about a certain class of representations of surface groups  
$\rho:\G_g\to G$, $g\geq 2$ into 
semisimple Lie groups of Hermitian type. 

\section{The Toledo invariant}
Let $G$ be a connected semisimple Lie group with finite center 
and associated Hermitian symmetric space $X$. 
Let $\G_g=\pi_1(\Sigma_g)$ be the fundamental group of 
a closed orientable Riemann surface 
and 
\bqn
\rho:\G_g\to G
\eqn
be a representation.
Then $\rho$ induces a canonical pull-back map 
\bqn
\rho_b^*: \hcb^*(G, \RR) \to \hb^*(\G_g, \RR)
\eqn 
in (continuous) bounded cohomology, which fits into the commutative diagram
\bqn
\xymatrix{
\hcb^*(G, \RR) \ar[d]^\kappa \ar[r]^{\rho_b^*} & \hb^*(\G_g, \RR) \ar[d]^\kappa\\
\hc^*(G, \RR) \ar[r]^{\rho^*} & \h^*(\G_g, \RR),
}
\eqn
where $\kappa$ denotes the natural comparison map between bounded and usual 
cohomology. 
Given any class $\a \in \hcb^2(G,\RR)$, we can evaluate the 
pull-back $\rho^*_b(\a) \in \hb^2(\G_g, \RR)$ on the fundamental class 
$[\Sigma_g]\in \h_2^{l_1} (\G_g, \RR)$, given in terms 
of the presentation 
$\G_g=\langle a_i, b_i\, |\, \Pi_{i=1}^g [a_i, b_i] = e\rangle$ 
(see \cite{Goldman_84,Brown}).
This defines the generalized Toledo invariant:
\bqn
{\rm Tol}: \hom(\G_g, G)&\times& \hcb^2(G, \RR) \longrightarrow \RR \\
{\rm Tol}(\rho, \a) &=& \langle \rho^*_b(\a), [\Sigma_g]\rangle.
\eqn

Let $k_X^b\in \hcb^2(G, \RR)$ be the (normalized) bounded K\"ahler class of $G$, 
which is defined by the bounded $G$-invariant $2$-cocycle 
\bqn
c_G(g_1, g_2, g_3)= \int_{\Delta(g_1 x, g_2 x, g_3 x)} \omega_X, \, x \in X
\eqn
 obtained by integrating the K\"ahler form $\omega_X\in \Omega^2(X)^G$ 
of the normalized Bergmann kernel 
over geodesic triangles $\Delta(g_1 x, g_2 x, g_3 x)$.
Then we recover the classical Toledo invariant as ${\rm Tol}(\rho) = {\rm Tol}(\rho, k_X^b)$. 
The Toledo invariant satisfies a Milnor-Wood type inequality

\bqn
|{\rm Tol}(\rho, \a)|=|\langle \rho^*_b(\a), [\Sigma_g]\rangle|\leq ||\rho_b^*(\a)|| |\chi(\Sigma)|\leq ||\a|| |\chi(\Sigma)|, 
\eqn
which gives the inequality (\cite{Domic_Toledo,Clerc_Orsted_2}):
$|{\rm Tol}_\rho|\leq 4\pi(g-1)r_X$ for the classical Toledo invariant.

\subsection{The Toledo invariant as characteristic number}
Suppose that $L\to X$ is a complex Hermitian line bundle on $X$ on which $G$ acts by isomorphisms. 
Let $\Omega$ be the curvature form of $L$. Then $\Omega$ is a $G$-invariant 
$2$-form on X.
In particular $\Omega$ defines a K\"ahler form on $X$ 
and a corresponding bounded K\"ahler class. 
If $X$ is irreducible, $\Omega$ is a multiple of the normalized
K\"ahler form. 
Given a representation $\rho:\G_g\to G$, 
there is an (up to homotopy) 
uniquely defined $\rho$-equivariant map 
$f:\widetilde{\Sigma_g} \to X$, the developing map of $\rho$. 
The pull-back bundle $f^*L\to \widetilde{\Sigma_g}$ 
descends to a Hermitian line bundle 
$L_\rho \to \Sigma$. 
The pull-back  $f^*\Omega$ is the curvature form of $f^*L$, hence 
$\frac{1}{2\pi i} f^*\Omega$ descends to a two form 
$\omega$ on $\Sigma$ representing the first Chern class of $L_\rho$. 
In particular $\int_\Sigma \omega = c_1(L_\rho) \in \ZZ$. 
The Toledo invariant is a constant multiple of $\int_\Sigma \omega$, thus
${\rm Tol}_\rho \in \lambda\ZZ$ for some $\lambda$. 
The explicit constant depends on the choice of the line bundle and
hence on the  group $G$. 
For $G= {\rm PU}(1,n)$ and ${\rm SU}(1,n)$ 
the constants are computed  explicitly in \cite{Goldman_Kapovich_Leeb}.

The Toledo invariant is locally constant and thus invariant under deformations of the representation. 
In particular, considering the Toledo invariant as a function on the representation variety 
${\rm Tol}: \Rr(\G_g, G) \to \RR$, respectively on the character variety ${\rm Tol}: \Xx(\G_g, G) \to \RR$, 
its level sets are unions of connected components.

\subsection{Computing the Toledo invariant}
We follow the computation for the Euler number given by Milnor 
in \cite{Milnor}.
Let 
\bqn
\G_g= \langle a_1, b_1, \dots, a_g, b_g \, |\, \Pi_{i=1}^g [a_i, b_i]=1\rangle
\eqn
be a presentation of the surface group.
Given a representation in terms of the generators 
$\rho(a_i) = A_i\in G$, $\rho(b_i) = B_i \in G$. 

A multiple of the K\"ahler class $k_X\in \hc^2(G)$ corresponds to a central extension 
\bqn
\xymatrix{
1\ar[r] & \pi_1(G) \ar[r]^i  & \tilde{G}\ar[r]^p & G \ar[r]& 1}. 
\eqn
Choose lifts $\tilde{A_i}$ and $\tilde{B_i}$ in $\tilde{G}$ 
with $p(\Pi_{i=1}^g [\tilde{A_i}, \tilde{B_i}])=1$. Since the
extension is central this expression does not depend on the chosen
lifts.  
Then $\Pi_{i=1}^g [\tilde{A_i}, \tilde{B_i}]$ lies in the image of $i$ 
and 
the map $\Sigma_g \mapsto i^{-1}((\Pi_{i=1}^g [\tilde{A_i}, \tilde{B_i}])
\in \pi_1(G)$ corresponds to $\rho^*(k_X) \in \h^2(\Sigma_g, \ZZ)$.
 
To compute the Toledo invariant in terms of $A_i, B_i$, 
we choose a contraction of $G$ to $K$ and a homotopy of $K$ to $S^1$ 
to obtain the diagram 
\bqn
\xymatrix{
\pi_1(G)\ar[d] \ar[r]^i & \tilde{G} \ar[d]^{\tilde{r}} \ar[r]^p & G \ar[d]\\
\pi_1(K) \ar[d] \ar[r] & \tilde{K} \ar[d]^{\tilde{h}} \ar[r] &K \ar[d]\\
\ZZ \ar[r] & \RR \ar[r]^{exp} &S^1
}
\eqn

The Toledo invariant is ${\rm Tol_\rho}$ is up to a 
multiple of $\pi$ equal to 
$- \frac{1}{2\pi} \tilde{h} \circ \tilde{r}(\Pi_{i=1}^g [\tilde{A_i}, \tilde{B_i}]) \in \ZZ \subset \RR$.

\section{Maximal representations}
Since the Toledo invariant is invariant under deformations, 
it can be used to distinguish certain parts of the representation
variety 
which we will characterize in the following.

\begin{defi}
Let $k_X^b \in \hcb^2(G, \RR)$ be the bounded K\"ahler class. 
A representation $\rho:\G_g\to G$ is called \emph{maximal}, if
\bqn 
{\rm Tol}(\rho,k_X^b)=||k_X^b||\, |\chi(\Sigma)|.
\eqn
\noindent
A representation  $\rho:\G_g\to G$ is called \emph{weakly maximal}, if
\bqn
{\rm Tol}(\rho,k_X^b)=||\rho^*(k_X^b)||\, |\chi(\Sigma)|.
\eqn
\end{defi}

\subsection{Examples of maximal representations}

1) Let $G={\rm PSU}(1,1)$. 
A uniformization $\rho_0: \G_g\to G$ is a maximal representation.
This follows from the 
theorem of Gau\ss -Bonnet:
Let $f:\widetilde{\Sigma_g} \to \HH^2$ be the developing map of $\rho$, i.e. 
$f(\Sigma_g)$ has a hyperbolic metric of sectional curvature $K=-1$.
Then 
\bqn
\langle \rho^*(k_X^b), [\Sigma_g]\rangle = \int_{\Sigma_g} f^*(\omega) = \int_{f(\Sigma_g)} \omega= 
- \int_{f(\Sigma_g)} {\rm sec}(M) = 2\pi (2g-2), 
\eqn
where ${\rm sec}(M)$ denotes the sectional curvature of $M$. 
There are several examples of maximal representations into arbitrary $G$ which arise from uniformizations of $\G_g$ in 
${\rm PSU}(1,1)$.\\
2) Let $h_\Delta:{\rm PSU}(1,1)\to G$ be a homomorphism associated with a tight holomorphic embedded disc $C\subset X$. 
Then precomposition of $h_\Delta$ with a uniformization with $\rho_0: \G_g\to {\rm PSU}(1,1)$
gives rise to a maximal representation $\rho= h_\Delta \circ \rho_0: \G_g\to G$. The image of $\G_g$ preserves the tight holomorphic disc 
$C$.\\
3) The precomposition of a homomorphism $h_P:{\rm SU}(1,1)^r\to G$, 
corresponding to a maximal polydisc $P\subset X$, with a 
product $\rho^0=(\rho_1, \dots, \rho_r):\G_g\to {\rm PSU}(1,1)^r$ of $r$ different uniformizations 
of $\G_g$ gives rise to a maximal representation. The image of $\G_g$ does not preserve a tight holomorphic disc in $X$, 
but the maximal polydisc $P\subset X$.\\
4) More generally, given a tight homomorphism $h:H\to G$  with $h^*(k_X^b)= c k_Y^b$ and 
$\rho_0: \G\to H$
a maximal representation, then the 
induced representation $\rho: \G\to G$ is $k$-maximal.\\
5) Not all maximal representations factor through a nontrivial tight homomorphism 
as Theorem~\ref{Zariski_dense} below 
shows. In particular, there exists maximal representations 
into Hermitian symmetric spaces of tube type, which leave no proper 
totally geodesic subspace invariant.

\subsection{Some history}
As explained in the first example, if $G= {\rm PSL}(2, \RR)$, 
all uniformizations 
$\rho:\G_g\to G$ are maximal representations. 
Goldman has shown in \cite{Goldman_thesis} that uniformizations 
can be characterized by the property of being maximal. 
Namely,  given a maximal representation $\rho:\G_g\to G$, 
the action of $\G_g$ via $\rho$ on $\HH^2$ is 
properly discontinuous, cocompact and without fixed points. 
This result gives a computable criterion, namely the Toledo 
invariant, to check whether a representation 
comes from a uniformization or not. 

The next one who studied 
maximal representations was Toledo in \cite{Toledo_89} 
who investigated the case where $G={\rm PSU}(1,n)$. Maximal representations
  $\rho:\G_g\to G$ essentially reduce to uniformizations. More precisely, 
suppose that $\rho:\G_g\to G$ is a maximal representation. Then there exists a 
complex hyperbolic line $\CC\HH^1\cong \DD$ in $\CC\HH^n$, 
which is stabilized by $\rho(\G_g)$. The induced action of $\G_g$ on 
$\CC\HH^1\cong \DD$ is, by Goldman's theorem, properly discontinuous, 
cocompact and without fixed points. In particular, the representation 
$\rho:\G_g\to G$ factors as  
$\rho= (\rho_1, \rho_2):\G_g\to G$, 
where $\rho_1:\G_g\to {\rm PSU}(1,1)$ is a uniformization and 
$\rho_2:\G_g\to {\rm U}(n-1)$ is arbitrary. 
Examples of representations $\rho:\G_g\to {\rm PSU}(1,2)$, whith all possible values of 
the Toledo invariant were constructed in 
\cite{Goldman_Kapovich_Leeb}. Some of these examples are discrete. 

The study of maximal representations into Lie groups of Hermitian type 
that are of rank greater or equal to two turned out to be more complicated.
Partial results were obtained for $G_{p,q}={\rm PSU}(p,q)$, $p\leq q$. 
Hernandez showed in \cite{Hernandez} that any maximal representation 
$\rho:\G_g\to G_{2,q}$ leaves a subspace of tube type associated to 
$ G_{2,2}$ invariant. For the general case $G_{p,q}$ he needed an additional 
assumption to prove the same result. 
Using quite different techniques, namely Higgs-bundles, 
Bradlow, Garc\'{\i}a-Prada and Gothen in \cite{Bradlow_GarciaPrada_Gothen} were  
able to prove the generalization of 
Hernandez' result to arbitrary $p,q$ under the assumption 
that the Zariski closure of the image of the maximal 
representation is reductive. 

As far as we know there were no results known about the 
geometric properties of the actions given by maximal representations 
into higher rank Lie groups before the work of Burger, Iozzi and 
Wienhard in \cite{Burger_Iozzi_Wienhard_ann}. 
Besides the results announced in
\cite{Burger_Iozzi_Wienhard_ann} we present also some further results
obtained since then.

\section{Weakly maximal representations with Zariski dense image}
We prove the following
\begin{prop}\label{prop:main}
Let $G$ be a connected simple Lie group 
with associated Hermitian symmetric space $X$ of rank $r_X\geq 2$.
Suppose  
$\rho : \Gamma_g \rightarrow G$ is a 
weakly 
maximal representation 
with Zariski dense image. 
Then $X$ is a Hermitian symmetric space of tube type on which 
$\G_g$ acts properly discontinuously (via $\rho$). 
If furthermore $\rho$ is maximal, the action of $\G_g$ via $\rho$ has no fixed points.
\end{prop}

\begin{rem}
The (strong) maximality is used to prove the
faithfulness of the representation. 
Probably the injectivity of $\rho$ also holds for weakly maximal
representations. Nevertheless, we do not have a proof for this and 
we are not aware of any nontrivial example of a weakly maximal Zariski dense
representation that is not maximal.
\end{rem}
The proof of the Proposition relies heavily on methods developed in \cite{Burger_Iozzi_App} and arguments 
which were already used in \cite{Iozzi_ern} and \cite{Burger_Iozzi_supq}. 
We make the same integration argument as in \cite{Iozzi_ern} to 
impose the weak maximality of the Toledo invariant on almost all triples.

We assume that $\G_g$ is realized as a 
cocompact lattice in ${\rm PSU}(1,1)$ via some representation $\pi_0:\G_g\to {\rm PSU}(1,1)$. Thus we may take $S^1$ endowed with 
the Lebesgue measure $\Ll$ as a doubly ergodic Poisson boundary 
for $\G_g$ 
(see \cite{Monod_book}).

\subsection{The integral formula}\label{sec:formula}
\begin{prop}\label{prop:formula}
Let $\rho:\G_g\to G$ be a Zariski dense representation and $\phi:S^1\to \cs$ 
the associated equivariant measurable boundary map. Then
\bqn
\int_{\Gamma_g\backslash \mathrm{PSU}(1,1)}
\beta_\Dd(\phi(hx),\phi(hy),\phi(hz))dh=
\frac{{\rm Tol}_\rho}{2\pi|\chi(\Sigma_g)|}\b_\DD(x,y,z)
\eqn
for almost every $x,y,z \in S^1$.
\end{prop}
\begin{proof}
The boundary map $\phi: S^1 \to \cs$ maps almost all $x,y,z \in S^1$
to a triple in $\cs^{(3)}$.  
Therefore the pull-back of the Bergmann cocyle $\b_\Dd$ via $\phi$ is well-defined 
and defines a measurable $\G_g$-invariant alternating essentially bounded cocycle 
$\phi^*\b_\Dd: (S^1)^3\to \RR$ corresponding to 
$\rho^* k_X^b$ under the canonical 
isomorphism $\hb^2(\Gamma_g,\RR) \simeq \Zz\la((S^1)^3)^{\Gamma_g}$.

Since we realized $\G_g$ via $\pi_0$ as a cocompact lattice in ${\rm PSU}(1,1)$, 
we may define the map 
\bqn
t_b: \hb^2(\G_g)\to \hcb^2({\rm PSU}(1,1))
\eqn
by integration with respect to a probability measure $dh$ on $\Gamma_g\backslash \mathrm{PSU}(1,1)$
\bqn
t_b(\alpha)(x,y,z) := \int_{\Gamma_g\backslash \mathrm{PSU}(1,1)} \alpha(hx, hy, hz) dh. 
\eqn
The same map, denoted by $t$, can be defined in ordinary cohomology.

Since $\hcb^2({\rm PSU}(1,1))$ is one dimensional (see
\cite{Guichardet_Wigner}), 
we obtain that $t_b(\phi^*\beta_\Dd)$ is a 
constant multiple of the Bergmann cocycle $\b_\DD$, i.e. 
$t_b(\phi^*\beta_\Dd)=\lambda c$.

We compute $\lambda$ using the following commutative diagram:
\bqn
\xymatrix{k_X^b \in \hcb^2(G) \ar[r]^{\cong} \ar[d]_{\rho^*} & H_c^2(G)\ar[d]^{\rho^*} &\\
\hb^2(\G_g) \ar[r]^{\kappa} \ar[d]_{t_b} & H^2(\G_g) \ar[r]^{\cong} \ar[d]^t & H^2_{dR} (\CC\HH^1/\pi_0(\G_g))\\
\RR \b_\DD=\hcb^2({\rm PSU}(1,1)) \ar[r]^\cong & H_c^2({\rm PSU}(1,1)) \ar[r]^\cong & \Omega^2(X)^{{\rm PSU}(1,1)} \ar[u]_{{\rm res}}
}
\eqn

Under the composition of the isomorphisms 
of the bottom line with the restriction map, 
the class $[\b_\DD]\in \hcb^2({\rm PSU}(1,1))$ is mapped to
$\omega_0={\rm res}(\omega_X) \in H^2_{dR} (\CC\HH^1/\pi_0(\G_g))$
(see \cite{Guichardet_Wigner}). 
Since $\langle \omega_0, [\Sigma_g]\rangle= 2\pi |\chi(\Sigma_g)|$ by Gau\ss-Bonnet,
we get 
\bqn
{\rm Tol}_\rho=\langle\rho^*\kxb, [\Sigma_g]\rangle= 
{\lambda}\langle \omega_0, [\Sigma_g]\rangle
={\lambda} 2\pi |\chi(\Sigma_g)|.
\eqn
Hence $\lambda=\frac{{\rm Tol}_\rho}{2\pi|\chi(\Sigma)|}$ and
\bqn
\int_{\Gamma_g\backslash \mathrm{PSU}(1,1)}
\beta_\Dd(\phi(hx),\phi(hy),\phi(hz))dh=
\frac{{\rm Tol}_\rho}{2\pi|\chi(\Sigma_g)|}\b_\DD(x,y,z).
\eqn
\end{proof}

Proposition~\ref{prop:formula} implies 
\begin{cor}\label{lem:triple_max}
Let $\rho:\G_g\to G$ be 
weakly 
maximal.
For almost all $x,y,z\in S^1$ we get $\phi^*\beta_\Dd(x,y,z) = ||\rho^*(k_x^b)||c(x,y,z) $, 
where $c$ is the orientation cocycle on $S^1$.
In particular, if $\rho$ is maximal then
$\|\rho^*(k_X^b)\|=\|k_X^b\| = \pi r_X$ and for almost every $x,y,z$,
we have 
$\beta_\Dd (\phi(x),\phi(y), \phi(z))=\pi r_X$.
\end{cor}
\begin{proof}
The assumption that
${\rm Tol}_\rho= 2|\chi(\Sigma_g)| ||\rho^*(k_X^b)||$ implies together with 
 $||{\phi^*\beta_\Dd}||=\|\rho^*(k_X^b)\|$ that the above integral equals 
the maximum value of the integrand. 
Thus
for almost every $x,y,z \in S^1$,
$\beta_\Dd(\phi(x),\phi(y),\phi(z))= {\|\rho^*(k_X^b)\|} \frac{1}{\pi}\beta_\DD (x,y,z) ={\|\rho^*(k_X^b)\|} c (x,y,z).$
\end{proof}

\subsubsection{Putting together the arguments}
We need the following dichotomy for Zariski dense subgroups
\begin{lemma}\label{lem:dichotomy}
Let $G$ be a simple Lie group and $\G<G$ a Zariski dense subgroup. 
Then $\G$ is either dense or discrete.
\end{lemma}
\begin{proof}
Let $H=\ol{\G}^\circ<G$ be the connected component of the closure of
$\G$ in the Lie group topology.
By Cartan's theorem, $H$ is a Lie subgroup of $G$. Consider the
adjoint action of $G$ on $\frakg$. Let $\frakh\subset \frakg$ be the
Lie algebra of $H$. Then $\frakh$ is $\Ad(\G)$-invariant, and by the
Zariski density of $\G$ it is $\Ad(G)$-invariant. Hence $\frakh=0$ or 
$\frakh=\frakg$, since $G$ is simple. The first case implies
$H=\{1\}$
and the second case $H=G$.  
\end{proof}

\begin{proof} (of Proposition~\ref{prop:main} for weakly maximal representations)
Assume that $\rho:\G_g\to G$ is a weakly maximal representation with Zariski dense image.
From Corollary~\ref{lem:triple_max} we obtain that 
for almost all $x,y,z \in S^1$ the cocycle $\phi^*\beta(x,y,z)= ||\rho^*(k_X^b)|| c(x,y,z)$ 
is a constant multiple of the 
orientation cocycle $c$ on $S^1$. 
Assume that $X$ is not of tube type. 
Then, fixing $x,y \in S^1$, we have with the notation of 
Lemma~\ref{lem:nonconstant} in Chapter 5 
that 
\bqn
{\rm Ess \,Im}{\phi}\cap \Bb_{\phi(x), \phi(y)} \subset 
Z:=\{ z\in \Bb_{\phi(x), \phi(y)}\,|\, P_{\phi(x), \phi(y)}= const\}.
\eqn
But, by Lemma~\ref{lem:nonconstant} in Chapter 5, 
$Z$ is a proper Zariski closed subset of 
\bqn
\Bb_{\phi(x), \phi(y)}\subset \cs.
\eqn
This contradicts the Zariski density of $\rho(\G_g)$. 
Hence $X$ has to be of tube type.

From now on we may assume that $X$ is of tube type. 
Assume that $r_X\geq 2$, then 
$G$ has at least three open orbits in 
$\cs^{(3)}= \Oo_-\cup \Oo_+ \cup \Oo_l$. 
By Corollary~\ref{lem:triple_max},  ${\rm Ess\,Im}(\phi^3)$ 
is contained in the closure of two open orbits in $\cs^{3}$, namely 
$\Oo_\pm= \{(x,y,z)\in \cs^{(3)}:\beta_\Dd(x,y,z)=\pm ||\rho^*(k_X^b)||\}$. 
Hence ${\rm Ess\,Im}(\phi^3)\neq\cs^{3}$, which implies that 
$\rho(\G_g)$ is not dense, and thus discrete
(Lemma~\ref{lem:dichotomy}). Therefore $\rho(\G_g)$ acts properly
discontinuous on $X$.
\end{proof}

Before we can prove the stronger statement for maximal representations we 
have to establish some additional properties of the boundary map $\phi:S^1\to \cs$.

\subsection{The boundary map}
We would like to have a refinement of the transversality properties of
the $\rho$-equivariant measurable boundary map $\phi: S^1\to \cs$. 
This will be useful to prove the injectivity of maximal representations.
\subsubsection{Transversality property I}
\begin{prop}\label{prop:zero}
Let $\rho:\G_g\to G$ be a representation with Zariski dense image. 
Then there exists a $\rho$-equivariant measurable 
boundary map 
\bqn
\tilde\phi:S^1 \to G/P
\eqn
such that the preimage of any proper $\RR$-algebraic subset $V\subset \GG/\PP$ has measure zero, 
i.e. $\Ll(\phi^{-1} (V))=0$.
\end{prop}

\begin{lemma}\label{lem:full_measure}
 If $A\subset S^1$ is of positive measure, 
then there exists a sequence $(\g_n)_{n\in \NN}$ in $\G_g$ 
such that $\Ll(\g_n A) \to 1$ as $n\to \infty$.
\end{lemma}
\begin{proof}
Note that the round measure on $S^1$ 
is the Patterson-Sullivan measure $\Ll=\mu_0$ with respect to the origin $0$, 
if $S^1$ is realized as the boundary of the Poincar\'e disk.
Let $\xi$ be a density point of $A \subset S^1$. 
Take $(\g_n)_{n\in \NN}$ such that $\g_n^{-1} 0 \to \xi$ converge radially. 
Let $C>0$ be a constant, which will be chosen appropriately later. 
Denote by $b(x,y,C)$ the shadow in $S^1$ of the ball of 
radius $C$ around $y$ seen from $x$, i.e.
\bqn
b(x,y,C):=\{ \xi \in S^1\,|\, \exists \g_{x\xi} \text{ with }
\g_{x\xi}\cap B(y,C) \neq \emptyset\}.
\eqn

For any $\eta \in S^1$ we define the Busemann function as usual by
$$
B_\eta(x,y) = \lim_{z\to \eta} [d(x,z) - d(y,z)]. 
$$
The change of the 
Patterson-Sullivan measures on $S^1$ with respect 
to the base point can be described in terms of 
the Busemann functions: 
\bqn
d\mu_x(\eta)= e^{- B_\eta(x,y)} d\mu_y(\eta).
\eqn
To prove the claim it suffices to show that for all $\e >0$ there exist
constants $C(\e), N(\e,C)$ such that for all $n\geq N$:
\bqn
\frac{\mu_0 (b(\g_n 0, 0, C)\cap \g_n A)}
{\mu_0 (b(\g_n 0, 0, C))}\geq 1-\eps.
\eqn

Modifying the expression, we get 
\bqn
\frac{\mu_0 (b(\g_n 0, 0, C)\cap \g_n A)}{\mu_0 (b(\g_n 0, 0, C))} &=&
\frac{\mu_{\g_n^{-1}0} (b(0,\g_n^{-1}  0, C)\cap A)}{\mu_{\g_n^{-1} 0} (b(0,\g_n^{-1}  0, C))}=\\
1- \frac{\mu_{\g_n^{-1}0} (b(0,\g_n^{-1}  0, C)\cap A^c)}{\mu_{\g_n^{-1} 0} (b(0,\g_n^{-1}  0, C))}
&=& 1- \frac{\int_{b(0,\g_n^{-1} 0, C)\cap A^c}  e^{- B_\xi(\g_n^{-1} 0,0)} d\mu_0(\xi)}
{\int_{b(0,\g_n^{-1} 0, C)}  e^{- B_\xi(\g_n^{-1} 0,0)} d\mu_0(\xi)}
\eqn

Since $\g_n^{-1} 0 \to \xi$ radially, we have that for suitable $N$ and for all 
$n\geq N$, the following holds for $\xi$ belonging to the region of integration:
$$|B_\xi(\g_n^{-1} 0,0) - d(\g_n^{-1} 0,0)| \leq K_1(C,N).$$
Therefore the above expression is greater than or equal to
\bqn
\geq 1-K_2(C,N)\frac{\mu_{0} (b(0,\g_n^{-1}  0, C)\cap A^c)}{\mu_{0} (b(0,\g_n^{-1}  0, C))}.
\eqn
Since $\xi$ was chosen to be a density point of $A$ the second summand 
becomes  arbitrarily small as $n$ tends to infinity.  
Choosing for a given $\eps>0$ the constant $C(\e)$ in such a way that 
$\mu_0 (b(x,0,C)) \geq 1-\eps$ for all $x\notin B(0,C)$ and
 a constant $N(\e,C)$, we obtain that for all $n\geq N$: 
\bqn
\frac{\mu_0 (b(\g_n 0, 0, C)\cap \g_n A)}
{\mu_0 (b(\g_n 0, 0, C))}\geq 1-\eps.
\eqn
\end{proof}

\begin{lemma}\label{lem:proper_sub}
Let $V$ be a proper $\RR$-algebraic subset $V\subset \GG/\PP$.
Let $(\g_n)_{n\in \NN}$ be a sequence in $\G_g$. Then there exists a proper $\RR$-subvariety 
$W\subset \GG/\PP$ 
such that for all $\eps >0$ there exists an 
$N\in \NN$ such that for all $n\geq N$:
$$
\rho (\g_n) V_\RR \subset \Nn_\eps(W_\RR).
$$
\end{lemma}
\begin{proof}
Let $I(V)$ be the ideal in $\RR[\GG/\PP]$ that defines $V$; 
denoting by $I(V)_d= I(V) \cap \RR_d[\GG/\PP]$ the homogeneous 
components, we can write $I(V) = \oplus_{d\geq 1} I_d(V)$. 
Pick now a number $d$ such that $I_d(V) \neq 0$, 
let $l=\dim_\RR I_d(V) \geq 1$ and consider $I_d(V)$ as a point 
in $\Gr_l(\RR_d[\GG/\PP])$ in the Grassmannian of $l$-planes in 
the vector space of homogeneous polynomials of degree $d$ on $\GG/\PP$.
By compactness of $\Gr_l(\RR_d[\GG/\PP])$, after passing to a subsequence, 
$I_d(\rho(\g_n) V) = \rho(\g_n)I_d(V)$ converges to 
$E \in \Gr_l(\RR_d[\GG/\PP])$. 
Given a nonzero element $0\neq P\in E$, there is a sequence 
$P_n\in I_d(\rho(\g_n) V)$ such that $P_n\to P$. 
We define $W$ to be the zero locus of $P$, i.e. 
$W:=\{p\in \GG/\PP \,|\, P(x)=0\}$.

Passing to a subsequence we can assume that $\rho(\g_n) V_\RR$ 
converges in the Gromov-Hausdorff topology 
to a closed set $F\subset \GG/\PP(\RR)$. We want to show that $F \subset W$.
There exists a representation (\cite{Prasad_Raghunathan}) 
$\pi:\GG \to {\rm GL}(U)$, defined over $\RR$, and an element 
$x_o \in \PP(U)(\RR)$, 
such that $\PP= \stab_\GG(x_0)$. 
We identify $\GG/\PP (\RR) = \pi(\GG(\RR))x_o \subset \PP(U)(\RR)$. 
Let $k$ be the dimension of $U$. 
For a nonzero element $[f]\in F$, now represented as $f\in \RR^k$,
there exists a sequence $[v_n]\in \rho(\g_n) V_\RR$ converging to $[f]$, 
such that $v_n \to f$, hence $P_n(v_n) =0$ and $P_n\to P$.
If we realize elements of $\GG/\PP(\RR)$ 
by elements of $\RR^k$, then, in $\RR^k$, $P_n \to P$ 
converge uniformly on compact sets. 
Therefore $P(f)$ is arbitrarily close to $P(v_n)$ for $n$ big enough. 
But $P(v_n)$ is gets arbitrarily close to $P_n(v_n)=0$ for big $n$, i.e.
for all $\eps >0$ there is an $N\in \NN$ such that for all $n\geq N$ we have:
\bqn
|P(f)-P(v_n)|\leq \eps \\
|P(v_n)-P_n(v_n)|\leq \eps.
\eqn
 hence $|P(f)|=|P(f)-P_n(v_n)|\leq 2\eps$ for all $\eps >0$, 
thus $P(f) =0$.
\end{proof}

\begin{proof} (of Proposition~\ref{prop:zero})
The existence of a measurable $\rho$-equivariant boundary map 
$\phi:S^1\to G/P$  
was proven above (see Section~\ref{sec:zariski_reps} in Chapter 7).
Thus it only remains to prove that $\Ll(\phi^{-1} (V))=0$ 
for any proper $\RR$-algebraic subset $V\subset \GG/\PP$.
We argue by contradiction. Let $A= \{ x\in S^1\,|\, \phi (x) \in V\}$ and 
assume that $\mu_0 (A) >0$. Then by Lemma~\ref{lem:full_measure} there is 
a sequence $(\g_n)_{n\in \NN}$ 
such that $\mu_0(\g_n A) \to 1$.
Since $\phi$ is equivariant, $\g_n A=\{ y\in S^1\, |\, \phi (y) \in \rho (\g_n) V \}$,
and the set $E_N = \bigcup_{n=N}^\infty \g_n A$ is of measure $1$ for all $N$. 
Now Lemma~\ref{lem:proper_sub} 
implies that for $\e= \frac{1}{m}$ there exists a $N(m)$ such that 
\bqn
\phi(E_{N(m)}) \subset \Nn_{\frac{1}{m}}(W_\RR).
\eqn
Hence $\bigcap_m E_{N(m)}$ is a set of full measure that is sent into $W_\RR$. 
This implies that ${\rm EssIm} \phi \subset W$ contradicting the Zariski density of $\rho(\G_g)$ 
since $W_\RR$ is a proper $\RR$-subvariety of $\GG/\PP$.
\end{proof}

The map $\phi$ does not take values in the Shilov boundary $\cs\cong G/Q$, 
but in $G/P$, 
where $P$ is a minimal parabolic subgroup of $G$. The boundary map 
$\phi:S^1\to \cs$ is obtained by composing $\phi$ 
with the natural projetion from 
$G/P$ to $G/Q$:

As a Corollary we get:
\begin{cor}\label{cor:zero_measure}
Let $\rho(\G_g)$ be Zariski dense. Then there exists a $\rho$-equivariant measurable 
boundary map 
\bqn
\phi:S^1 \to \cs,
\eqn
such that the preimage of any proper $\RR$-algebraic subset $V\subset G/Q$ has measure zero, 
i.e. $\Ll(\phi^{-1} (V))=0$.
In particular for almost every $x,y \in S^1$, $\phi(x)$ and $\phi(y)$ are transverse.
\end{cor}
\begin{proof}
Transversality of two points on the Shilov boundary is, given in 
terms of the Bruhat decomposition, a Zariski open condition. In particular, 
fixing $x$, the set of the points $y\in\cs$ which are not transverse to $x$
is a proper $\RR$-algebraic subset of $G/Q$.
\end{proof}
%

\subsubsection{Transversality property II}
To prove the injectivity of maximal representations  
we study the essential graph of the boundary map $\phi:S^1\to \cs$. 
The (strong) maximality enters in the proof of Lemma~\ref{lem:graph_transverse}. 
We are assuming that the representation $\rho:\G_g\to G$ is maximal and that $X$ is of tube type. 
Thus we are in the following situation:

The map  $\phi:S^1\to \cs$ is measurable $\rho$-equivariant 
and maps almost every $x,y,z\in S^1$ 
to a triple in $\cs^{(3)}$ 
such that $\beta_\Dd (\phi(x),\phi(y), \phi(z))=\pi r_X$.

Denote by $E={\rm EssIm}\phi$ the essential image of $\phi$.
Then for almost every $x\in S^1$, $\phi(x) \in E$.
We define the essential graph $F$ of $\phi$ to be the 
support of the push-forward of the Lebesgue measure $\Ll$ on $S^1$ 
under the map 
\bqn
\Phi:S^1&\to& S^1\times \cs\\
x &\mapsto& (x,\phi(x)).
\eqn
This is $F:= \supp(\Phi_*(\Ll))\subset S^1\times \cs$.

\begin{lemma}\label{lem:graph_max}
Let $(x_i,f_i)\in F$. Assume that all $x_i$ are pairwise distinct 
and all $f_i$ pairwise transverse, then $\beta(f_1,f_2,f_3) = r_X \pi c(x_1,x_2,x_3)$.
\end{lemma}
\begin{proof}
We may take (Corollary~\ref{cor:zero_measure}) 
small neighbourhoods $U_i\subset F$ of 
positive measure around $(x_i,f_i)$ such that 
for $(y_i, g_i) \in U_i$ the transversality is preserved. 
Let $(x_{i_n},f_{i_n})\in U_i$ be sequences converging to $(x_i, f_i)$ such that 
$\beta(f_{1_n},f_{2_n},f_{3_n}) = r_X \pi c(x_{1_n},x_{2_n},x_{3_n})$.  
By continuity of $\beta$ and $c$ on $\cs^{(3)}$ resp. $(S^1)^{(3)}$, the same holds for the limit points.
\end{proof}

Since $X$ is assumed to be of tube type, 
the Bergmann cocycle $\b$ extends to a strict cocycle on $\cs$ 
(see Lemma~\ref{strict} in Chapter 5) 
that satisfies the following properties:
\begin{enumerate}\label{eq:properties}
\item{$|\beta (f_1,f_2,f_3)|\leq r_X\pi$ for all $f_1,f_2,f_3 \in \cs$}
\item{if $\beta(f_1,f_2,f_3) = r_X\pi$ and $f_3$ is transverse to both 
$f_1$ and $f_2$, then $f_1,f_2,f_3$ are pairwise transverse.}
\end{enumerate}

\begin{lemma}\label{lem:graph_transverse}
Let $(x_1,f_1),(x_2,f_2) \in F$ with $x_1\neq x_2$, then $f_1$ is transverse to $f_2$.
\end{lemma}
\begin{proof}
We may pick (Corollary~\ref{cor:zero_measure}) a point $a \in (x_1,x_2)$ such that $\phi(a)\in E$ and $\phi(a)$ 
is transverse to $f_1$ and $f_2$, and a point $b\in (x_2,x_1)$ 
such that $\phi(b) \in E$ and $\phi(b)$ transverse to $\phi(a), f_1, f_2$.
With Lemma~\ref{lem:graph_max} the cocycle identity for $\beta$ reads:
\bqn
0&=& 
\beta(\phi(a),f_2,\phi(b)) - \beta(f_1,f_2,\phi(b)) + \beta(f_1,\phi(a),\phi(b)) 
- \beta (f_1,\phi(a), f_2)\\
&=&
r_X\pi - \beta(f_1,f_2,\phi(b)) + r_X\pi - \beta (f_1,\phi(a), f_2).
\eqn
Hence 
$\beta(f_1,f_2,\phi(b)) +  \beta (f_1,\phi(a), f_2) = 2r_X\pi$, 
and properties $(1)$ and $(2)$ finish the proof.
\end{proof}

\begin{proof} (of Proposition~\ref{prop:main} in the maximal case)
We may assume that $X$ is of tube type. 

Assume that $\rho$ is maximal, but not injective. Then there exists a $\g\neq e\in \G_g$ 
with $\rho(\g) = e\in G$. Thus for any $(x, f) \in F$ we have that $(\g x,\rho(\g)f)= (\g x,f)\in F$, 
contradicting Lemma~\ref{lem:graph_transverse}, 
since $\g x\neq x$.
\end{proof}

\subsubsection{Left/Right continuity}
Let $\phi:S^1\to \cs$ be the equivariant measurable boundary map,
associated to a maximal representation $\rho:\G_g\to G$. Assume 
that the Hermitian symmetric space associated to $G$ is of tube type. 
The properties of the generalized Maslov cocycle we used to prove the 
injectivity of the representation imply also the following statement.

\begin{lemma}\label{lem:map_left_right}
There is a right continuous map $\phi^+$ and a left continuous map $\phi^-$, 
$\phi^\pm: S^1\to \cs$, such that $\phi= \phi^\pm$ almost everywhere.
\end{lemma}

From Lemma~\ref{lem:converge} in Chapter 5 
we obtain the following property of the essential graph $F$ of $\phi$:
\begin{lemma}\label{left_right}
Let $(x,y]\subset S^1$, 
$E_{(x,y]}:=\{f\in \cs\,|\, \exists t \in (x,y], (t,f) \in F\}$ and
$F_x:=\{x\}\times \cs \cap F$. 
Then $\overline{E_{(x,y]}} \cap F_x$ is reduced to a point.
Let $[x,y)\subset S^1$, 
$E_{[x,y)}:=\{f\in \cs\,|\, \exists t \in [x,y), (t,f) \in F\}$ and
$F_y:=\{y\}\times \cs \cap F$. 
Then $\overline{E_{[x,y)}} \cap F_y$ is reduced to a point.
\end{lemma}
\begin{proof}
Arguing by contradiction we prove the 
first claim, the second follows by a similar argument.
Let $f,f^\prime \in \overline{E_{(x,y]}} \cap F_x$, $f_- \in F_x$
transverse to $f$ and $ f'$.
Let $x_i\to x$ be a sequence in $(x,y]$, $f_i \in F_{x_i}$ a sequence of 
points in $\cs$
 with $f_i\to f$. 
Pick another sequence $y_i \to x$ with $y_i \in (x,x_i)$ and 
$f^\prime_i \in F_{y_i}$ with $f^\prime_i \to f^\prime$. 
We may assume that $f_-$ is transversal to all $f_i, f^\prime_i$. 
By Lemma~\ref{lem:graph_max},  $\tau(f, f^\prime_i, f_i)$ is maximal and 
$f_i\to f$. We deduce from Lemma~\ref{lem:converge} in Chapter 5, that $f^\prime_i\to f$ and 
hence $f= f^\prime$. 
\end{proof}

\begin{proof} (of Lemma~\ref{lem:map_left_right})
For all $x\in S^1$ with $\# F_x =1$ we define $\phi^\pm$ to be $\phi(x)\in F_x$. 
Whenever $\# F_x >1$ we define $\phi^+(x)$ to be the right limit point of $\phi$, i.e. 
$\overline{E_{(x,y]}} \cap F_x$ and $\phi^-(x)$ to be the left limit point of $\phi$, 
i.e. $\overline{E_{[y,x)}} \cap F_x$. 
\end{proof}

\section{Characterizing maximal representations}
We refine the definition of maximal representations, taking into account the 
cohomology class $\a\in \hcb^2(G, \RR)$.
\begin{defi}\label{k_max}
A representation $\rho:\G_g\to G$ is called \emph{$\a$-maximal}, if
${\rm Tol}(\rho,\a)=||\a||\, |\chi(\Sigma)|$. 
\end{defi} 

With the theory of tight embeddings at hand, we are able to 
prove: 
\begin{thm}\label{thm:main}
Let $G$ be a connected semisimple real algebraic group with finite center and assume that its associated 
symmetric space is Hermitian. 
Let $\rho :\G_g \rightarrow \rm G$
be a maximal representation. Then:\\
(1) $\rho(\G_g)$ stabilizes a maximal Hermitian symmetric subspace of tube type $T\subset X$.\\ 
(2) The Zariski closure $L$ of $\rho (\G_g)$ is reductive.\\
(3) The symmetric subspace associated to $L$ is a Hermitian symmetric space of tube type and the inclusion $Y\to T$ is tight. \\
(4) The action of $\G_g$ on $Y$ via $\rho$ is properly 
discontinuous without fixed points.
(5) There exists a Hermitian symmetric subspace of tube type $T_Y$
containing $Y$ on which $\rho(\G_g)$ acts properly discountinuously
without fixed points.
\end{thm}

\begin{rem}
A Hermitian symmetric subspace is a totally geodesically and
holomorphically embedded subspace.
\end{rem}
As in Lemma \ref{tight_properties} in Chapter 6 
we can deduce the following properties  
of $\a$-maximal representations from results on bounded cohomology in 
\cite{Monod_book}.
\begin{lemma}\label{max_properties}
Suppose $\rho: \G_g\to G$ is a representation, 
$\a \in \hcb^2(G)$. If $\rho$ is $\a$-maximal, the 
following holds.
\begin{enumerate}
\item{Suppose $L<G$ is a closed subgroup and $\rho(\G_g) \subset L<G$.\\
        Then $\rho$ is $\a_{|_L}$-maximal and 
        $||\a_{|_L}||=||\a||$}
\item{Suppose $\G_0 \lhd\G_g$ is a normal subgroup of finite index.\\
        Then $\rho_0:=\rho_{|_{\G_0}}$ is $\a$-maximal 
        and $||\rho{|_{\G_0}}^*\a||= ||\rho^*\a||= ||\a||$.}
\item{Suppose $R \lhd G$ is a normal closed amenable subgroup 
        and let $\overline\a$ be the image of $\a$ 
        under the canonical isometric isomorphism $\hcb(G)\cong \hcb(G/R)$ 
        (see \cite[Theorem 4.1.1.]{Burger_Monod_GAFA}). \\
        Then $||\overline\a|| = ||\a||$ and 
        the induced homomorphism 
        $\overline\rho : \G_g\to G/R$ is $\overline\a$-maximal.}
\item{Assume $G= \Pi_{i=1}^n G_i$ and 
        $\a_i=\a_{|_{G_i}} \in \hcb^2(G_i)$.\\
        Then the induced representations $\rho_i=\pr_i\circ \rho:\G_g\to G_i$ 
        are $\a_i$-maximal.} 
\end{enumerate}
\end{lemma}
\begin{proof}
The first, the third and the fourth claim follow immediately from the corresponding properties of tight
homomorphisms.
For the second claim let $\Sigma'$ denote the covering of $\Sigma_g$ associated to the
subgroup $\G_0$. Let $d$ be the degree of the covering. We obtain the
following 
inequality: 
\bqn
{\rm Tol}(\rho_0) \leq 2|\chi(\Sigma')|\, ||\a||= d 2
|\chi(\Sigma_g)|\, ||\a|| = d {\rm Tol} (\rho) = {\rm Tol}(\rho_0).
\eqn
Hence $\rho_0: \G_0\to G$ is maximal. 
\end{proof}

\begin{proof} (of Theorem~\ref{thm:main})
Let $L=\overline{\rho(\G_g)}^Z (\RR)$ be the Zariski closure 
of the image. If $L$ is not connected 
consider $L' = \overline{\rho(\G_g)}^Z (\RR)^\circ$
its connected component of the identity. Then $\G^\prime_g:=
\rho^{-1}(L'\cap \rho(\G_g)$ is a normal subgroup of finite index in
$\G_g$. It follows from Lemma~\ref{max_properties} that the representation $\rho':
\G^\prime_g \to L'$ is still maximal, thus 
passing to a finite index normal subgroup of $\G_g$
(Lemma~\ref{max_properties}), 
we may assume 
that 
$L=\overline{\rho(\G_g)}^Z (\RR)$, the Zariski closure 
of the image, is connected. 
A maximal representation $\rho:\G_g\to G$ is in particular 
a tight homomorphism. Thus  
Theorem~\ref{thm:tightness} of Chapter 6 implies that $L<G$ is reductive 
and its associated symmetric space is Hermitian. 
Denote by $M$ the semisimple part of $L$ and let 
$M'=M_1\times \dots \times M_\ell$ be the product 
of the noncompact simple factors of $M=M'\times \Cc$, 
with associated Hermitian symmetric space $Y=Y_1\times \cdots\times Y_\ell$.
The corresponding representations $\rho_i:\G_g\to M_i$ are $k_i$-maximal 
(Lemma~\ref{max_properties}), 
where $k_i=\pr_i({k_x^b}_{|_L})$ 
are constant multiples of the bounded K\"ahler forms 
$k_{Y_i}^b\in \hcb^2(M_i, \RR)$. 
Hence the representations $\rho_i$ are maximal and 
by Proposition~\ref{prop:main} the symmetric spaces $Y_i$ are 
Hermitian symmetric spaces of tube type and the action of 
$\G_g$ via $\rho_i$ on $Y_i$ is properly discontinuous without fixed points.
Hence the symmetric space $Y= \Pi_{i=1}^{\ell} Y_i$ is a Hermitian symmetric space of tube type 
and the action of $\G_g$ via $\rho$ is properly discontinuous without fixed points. 
Since the inclusion 
$L\to G$ is tight (Lemma~\ref{tight_properties} in Chapter 6), 
the embedding $Y\to X$ 
is tight. 
It follows from Proposition~\ref{prop:tight_tube} of Chapter 6 
that $Y$ is contained in $T\subset X$, a maximal 
subdomain of tube type in $X$, that is hence stabilized by $\rho(\G_g)$. 
Furthermore the Hermitian hull $T_Y$ of $Y$ is preserved by
$\rho(\G_g)$. Since the centralizer of $L$ in ${\rm Is}(T_Y)$ is
trivial, $\rho(\G_g)$ acts still properly discontinuously and without
fixed points on $T_Y$. 
\end{proof}

\subsection{Construction of maximal representations with Zariski dense image}
We construct maximal representations 
with Zariski dense image in the case where $X$ is of tube type, 
showing that Theorem~\ref{thm:main} is optimal.

\begin{thm}\label{Zariski_dense}
Let $X$ be a Hermitian symmetric space of tube type and $P\subset X$
be a maximal polydisc. 
Let $\rho_0: \G_g\to {\rm Is}(X)^\circ$ 
be the maximal representation obtained from composing a diagonal discrete injective embedding 
$\G_g\to {\rm SU}(1,1)^r$ with $h_P$. 

Then $\rho_0$ admits a continuous deformation 
\bqn
\rho_t:\G_g\to {\rm Is}(X)^\circ\, , \, t\geq 0
\eqn
with $\rho_t$ maximal and $\rho_t(\G_g)$ Zariski dense in 
${\rm Is}(X)^\circ$ for $t>0$.
\end{thm}

The proof uses the following Lemmata.
\begin{lemma}
Suppose $\rho_i:\G_g\to {\rm PSU}(1,1)$, 
$i=1,\dots, r$ are representations with 
Zariski dense image which 
are pairwise not conjugate. 
Then $\rho= (\rho_1,\dots,\rho_r): \G_g\to {\rm PSU}(1,1)^r$ 
has Zariski dense image.
\end{lemma}
\begin{proof}
We argue by contradiction. 
Denote $L:={\rm PSU}(1,1)$. 
Assume that the Zariski closure $H$ of $\rho(\G_g)$ is a proper 
subgroup of $L^r$. Since $\rho_i$ has Zariski dense image in $L$, the 
projections $\pr_i$ on the factors of $L^r$ satisfy $\pr_i(H)=L$ for all $i$.
Denote by $pr^i$ the projection on the first $i$ factors.
Let $1\leq s<r$ be the number such that $pr^s(H)=L^s$ and 
$H^\prime:=pr^{s+1}(H)$ is a proper subgroup of $L^{s+1}$.
Then $H^\prime <L^s\times L$ is algebraic and the projections of $H^\prime$
onto the two factors are surjective.
We represent an element $x\in H^\prime$ as $x=(a,b)\in L^s\times L$. 
Since the projection onto the second factor is surjective, we have that 
$H^\prime \cap (e\times L)$ is normal in $(e\times L)$.
Since $L$ is simple we have to consider the following two cases:\\
1) $H^\prime \cap (e\times L)= (e\times L)$, which contradicts the choice of $s$.\\
2) $H^\prime \cap (e\times L)= (e\times e)$, which implies that $H^\prime$ is 
the image of a surjective homomorphism of $L^s\to L$ of a semisimple
group onto a simple group. 
Such a homomorphism can only be the projection onto a factor composed with an automorphism of $L$. 
Hence $H^\prime< L^{s-1}\times \{(l,\a(l)\,|\, l\in L,\, \a\in \aut(L)\}$. 
After reordering this gives $H<L^{r-2}\times \{(l,\a(l))\,|\, l\in L,\, \a\in \aut(L)\}$, 
hence $\a(\rho_{r-1}(a))=\rho_r(a)$, thus $\rho_{r-1}$ is conjugate to 
$\rho_r$.
This contradicts the assumption and proves the claim.
\end{proof}

\begin{lemma}\label{lem:twists}
Let $H, G$ be simple 
Lie groups of Hermitian type with associated symmetric spaces $Y, X$ of tube type and of the same rank.
Assume that $\pi: H\to G$ is a (virtually) injective tight homomorphism auch that the induced totally geodesic embedding 
$Y\to X$ is holomorphic. Let $C\subset Y\subset X$ be a tight disc. Then there exists $k\in \stab_G(C)$ with $k\notin \stab_H(C)$.
%
\end{lemma}
\begin{proof}
It follows from Corollary~\ref{H2_hom} that Lie algebra homomorphism
corresponding to $\pi$ is an
$(\h2)$-homomorphism. These homomorphism were classified in 
\cite{Satake_hol_65, Ihara_65}. From a case by case checking (see Example at the end of this Section) 
it follows that in all situations $\stab_H(C)\subsetneq \stab_G(C)$:\\
1) $G={\rm Sp}(2n,\RR)$: $\stab_{{\rm Sp}(2n,\RR)}(C)\cong {\rm O}(n)$. There are no $(\h2)$-holomorphic embeddings of 
simple Lie algebras of maximal rank.\\
2) $G= {\rm SU}(n,n)$:  $\stab_{{\rm SU}(n,n)}(C)\cong {\rm U}(n)$. The only possible situation is $H ={\rm Sp}(2n,\RR)$, 
i.e. ${\rm O}(n)\subsetneq {\rm U}(n)$.\\ 
3) $G={\rm SO}^*(4n)$:  $\stab_{ {\rm SO}^*(4n)}(C)\cong {\rm SU}(m,\HH)$
The only possible situations are $H= {\rm SU}(n,n)$ or $H= {\rm Sp}(2n,\RR)$. We have to consider only the case $H= {\rm SU}(n,n)$. 
But ${\rm U}(n) \subsetneq {\rm SU}(m,\HH) $.\\
4) $G={\rm SO}(2,n)$: $\stab_{{\rm SO}(2,n)}(C)\cong {\rm O}(n-1)$. The only possible situations are $H={\rm SO}(2,m)$ with $m<n$. 
But ${\rm O}(m-1)\subsetneq {\rm O}(n-1)$ for $n>m$.\\
5) Exceptional case: $\stab_{E_7}(C)\cong F_4$.
\end{proof}
\begin{proof} (of Theorem~\ref{Zariski_dense})
For the construction of a desired representation,
realize the fundamental group as an amalgamated product over 
a separating geodesic, $\G_g=A*_{\langle \gamma\rangle}B$. 
Choose a hyperbolization of $\G_g$, 
$\pi:\G_g\to {\rm PSU}(1,1)$ 
and use the diagonal embedding 
\bqn
\Delta: {\rm PSU}(1,1)\to {\rm PSU}(1,1)^r
\eqn
to define hyperbolizations 
\bqn
\rho_i:=\pr_i\circ\Delta\circ\pi|_A: A\to {\rm PSU}(1,1)
\eqn
and
\bqn
\omega_i:=\pr_i\circ\Delta\circ\pi|_B:B\to  {\rm PSU}(1,1).
\eqn 
Let $t_A, t_B: {\rm PSU}(1,1)^r \to G$ be two different homomorphisms  
which coincide on $\Delta({\rm PSU}(1,1))$. 
Choose now two one-parameter families of deformations $\rho_i^t, \omega_i^t$, such that the
$\rho_i^t$'s, ${i=1,\dots,r}$, respectively the $\omega_i^t$'s, are pairwise not conjugate for all $t>0$ and 
$\rho_i^t(\gamma)=\rho_i(\gamma)$, respectively $\omega_i^t(\gamma)=\omega_i(\gamma)$, for all $t$.
By the above Lemma, the representations of $A$, 
respectively $B$, given by $\rho^t(a)=t_A(\rho_1^t(a),\dots,\rho_r^t(a))$, 
respectively $\omega^t=t_B(\omega_1^t(a),\dots,\omega_r^t(a))$, have Zariski dense image in $t_A({\rm PSU}(1,1)^r)$, 
respectively $t_B({\rm PSU}(1,1)^r)$ for $t>0$. 
They define a representation
\bqn
\pi^t:\G_g\to G
\eqn
by the universal property of amalgamated products. 
Constructed as a deformation of a maximal representation,  
$\pi^t$ has maximal Toledo invariant.
Hence (Theorem~\ref{thm:main}) the Zariski closure 
of its image is reductive, and of maximal rank, 
since it contains the image of $t_A$. 
The symmetric space $Y_t$, $t>0$ 
corresponding to its semisimple part $H_t$ is of tube type. 
It contains a tight holomorphic disc, which is also a tight holomorphic 
 disc in $X$. 
For the corresponding Lie algebra homomorphisms $\pi_0:\fraks\frakl(2,\RR)\to \frakh$, 
$\psi: \frakh\to \frakg$ and $\pi= \psi \circ \pi_0:\fraks\frakl(2, \RR)\to\frakg$ this means 
that $\pi_0$ and $\pi$ are $(\h2)$, hence $\psi$ is also $(\h2)$ and the associated embedding 
$Y\to X$ is holomorphic. By Lemma~\ref{lem:twists} there are elements 
$k\in \stab_G(t_A(\Delta({\rm PSU}(1,1))))=\Stab_G(t_B(\Delta({\rm PSU}(1,1))))$
such that the subgroup of $G$ generated by the image of $t_A$ and $t_B=kt_A$ 
coincides with $G$. 
\end{proof}

\begin{exo}
As illustration, we describe the standard polydisc and the Lie algebra 
$\frakk_0$ of its stabilizer explicitly for all classical cases.

1) Case $\frakg=\fraks\frakp(2n,\RR)$. 
The Lie algebra is given by
\bqn
\frakg=\begin{pmatrix} X_1 & X_2 \\ \overline{X_2} & \overline{X_1} \end{pmatrix} =\frakk+\frakp, \\ 
\frakk= \begin{pmatrix} X_1 & 0 \\ 0 & \overline{X_1} \end{pmatrix},\,
\frakp= \begin{pmatrix} 0 & X_2 \\ \overline{X_2} & 0 \end{pmatrix},
\eqn
where $X_i \in \Mat_{n\times n} (\CC)$, $X_1^* =-X_1$, $X_2^T= X_2$.

The standard maximal abelian Lie subalgebra $\fraka\subset \frakp$ and the maximal polydisc $\frakr$
are:
\bqn
\fraka&=&\left\{  \begin{pmatrix} 0 & A \\ A & 0 \end{pmatrix}| \, A=\diag(a_1,\dots, a_n), \, a_i\in \RR\right\} , \\
\frakr&=&\left\{ \begin{pmatrix} 0 & Z \\ \ol{Z} & 0 \end{pmatrix}|\, Z=\diag(z_1, \dots, z_n) , \, z_i \in \CC\right\} .
\eqn
The tight disc is given as subset of $\frakr$ where $z_i=z$ for all $i$.
The Lie algebra of the subgroup $K_0<K$ fixing the tight disc 
is 
\bqn
\frakk_0=\left\{ \begin{pmatrix} X & 0 \\ 0 & X \end{pmatrix}|\, X\in \frako(n)\right\} .
\eqn
\\ 
2) Case $\frakg= \fraks\fraku(n,n)$: 
The Lie algebra $\frakg$ is given by: 
\bqn
\frakg=\begin{pmatrix} X_{11} & X_2 \\ X_2^* & X_{12} \end{pmatrix} =\frakk+\frakp, \\ 
\frakk= \begin{pmatrix} X_{11} & 0 \\ 0 & X_{22} \end{pmatrix},\,
\frakp= \begin{pmatrix} 0 & X_2 \\ X_2^* & 0 \end{pmatrix}, 
\eqn
where $X_i \in \Mat_{n\times n} (\CC)$, $X_{1i}^* 0=-X_{1i}$, $\tr
X_{11}+ \tr X_{12}=0$.

The standard maximal abelian Lie subalgebra $\fraka\subset \frakp$ and the maximal polydisc $\frakr$
are:
\bqn
\fraka&=&\left\{  \begin{pmatrix} 0 & A \\ A & 0 \end{pmatrix}| \, A=\diag(a_1,\dots, a_n), \, a_i\in \RR\right\} , \\
\frakr&=&\left\{ \begin{pmatrix} 0 & Z \\ \ol{Z} & 0 \end{pmatrix}|\, Z=\diag(z_1, \dots, z_n) , \, z_i \in \CC\right\} .
\eqn
The tight disc is given as subset of $\frakr$ where $z_i=z$ for all $i$.

The Lie algebra of the subgroup $K_0<K$ fixing the tight disc 
is 
\bqn
\frakk_0=\left\{ \begin{pmatrix} X & 0 \\ 0 & X \end{pmatrix}|\, X\in \fraku(n)\right\} .
\eqn
\\
3) Case $\frakg=\fraks\frako^*(4n)$. 
The Lie algebra can be described as:
\bqn
\frakg=\begin{pmatrix} X_1 & X_2 \\ -\ol{X_2} & \ol{X_1} \end{pmatrix} =\frakk+\frakp, \\
\frakk= \begin{pmatrix} X_1 & 0 \\ 0 & \ol{X_1} \end{pmatrix}, \,
\frakp= \begin{pmatrix} 0 & X_2 \\ -\ol{X_2} & 0 \end{pmatrix}, 
\eqn
where $X_i \in \Mat_{2n\times 2n} (\CC)$, $X_1^* 0=-X_1$, $X_2^T= -X_2$.

The standard maximal abelian Lie subalgebra $\fraka\subset \frakp$ and the maximal polydisc $\frakr$
are:
\bqn
\fraka&=&\left\{  \begin{pmatrix} 0 & A \\ -\ol{A} & 0 \end{pmatrix}| \, A=\begin{pmatrix} 0 & B \\ -B & 0 \end{pmatrix},\, B= \diag(b_1,\dots, b_n), \, b_i\in \RR\right\} ,\\
\frakr&=&\left\{ \begin{pmatrix} 0 & Z \\ -\ol{Z} & 0 \end{pmatrix}|\, Z=\begin{pmatrix} 0 & Y \\ -Y & 0 \end{pmatrix}, \, Y= \diag(y_1, \dots, y_n) , \, y_i \in \CC\right\} .
\eqn
The tight disc is given as subset of $\frakr$ where $y_i=y$ for all $i$.
The Lie algebra of the subgroup $K_0<K$ fixing the tight disc 
is 
\bqn
\frakk_0=\fraks\fraku(2n,\HH)=\left\{ \begin{pmatrix} X & 0 \\ 0 & X
  \end{pmatrix}|\, X=\begin{pmatrix} A & B \\ B & \ol{A}
  \end{pmatrix}\right\} ,
\eqn
where $A$ is arbitrary and $B$ is symmetric and purely imaginary.

4) Case $\frakg=\fraks\frako(n,2)$. 
This case is a bit more subtle than the previous cases.
The Lie algebra is: 
\bqn
\frakg=\begin{pmatrix} X_{11} & X_2 \\ X_2^T & X_{12} \end{pmatrix} =\frakk+\frakp, \\
\frakk= \begin{pmatrix} X_{11} & 0 \\ 0 & X_{12} \end{pmatrix}, \,
\frakp= \begin{pmatrix} 0 & X_2 \\ X_2^T & 0 \end{pmatrix}, 
\eqn
where $X_{11} \in \Mat_{n\times n} (\RR)$, $X_{12} \in \Mat_{2\times 2}(\RR)$, 
$X_2 \in \Mat_{n\times 2}(\RR)$.
The element $Z_0$ generating the center of $\frakk$ is given by 
\bqn
Z_0 = \begin{pmatrix}0 & 0 \\ 0 & J \end{pmatrix},
\eqn
where $J=\begin{pmatrix}0 & -1 \\ 1 & 0 \end{pmatrix}$.
The standard maximal abelian Lie subalgebra $\fraka\subset \frakp$ and the maximal polydisc $\frakr$
are:
\bqn
\fraka&=&\left\{ \begin{pmatrix} 0 & A \\ A^T & 0 \end{pmatrix}| \,
A=\begin{pmatrix} a_1 & 0 \\ 0 & a_2 \\ 0 & 0 \end{pmatrix} \, a_i\in \RR\right\} , \\
\frakr&=&\left\{ \begin{pmatrix} 0 & Z \\ Z^T & 0 \end{pmatrix}|\,
Z=\begin{pmatrix} a_1 & b_1 \\ -b_2 & a_2 \\ 0 & 0  \end{pmatrix},\, a_i, b_i \in \RR\right\} .
\eqn
The tight disc is given as the image of the $(\h2)$-holomorphic image of 
$\fraks\frako(1,2)\to \fraks\frako(n,2)$, hence it is 
the subset of $\frakr$ where $-b_2=a_2=0$. 
(It is not the subset with $b_2=b_1$ and $a_1=a_2$.)

The Lie algebra of the subgroup $K_0<K$ fixing the tight disc 
is 
\bqn
\frakk_0=\left\{ \begin{pmatrix}X & 0 \\ 0 & \id \end{pmatrix}|\, X=
\begin{pmatrix}1 & 0 \\ 0 & X_1\end{pmatrix}, \, X_1 \in \frako(n-1)\right\} .
\eqn
\end{exo}

\section{The maximal part of the representation variety}
\subsection{Structure of the representation variety}
Let $\G_g$ be the fundamental group of a closed Riemann surface, $g\geq 2$. 
and $G$ be the adjoint group of a semisimple real Lie algebra. 
Choose a presentation  
\bqn
\G_g= \langle a_1, b_1, \dots, a_g, b_g \, |\, \Pi_{i=1}^g [a_i, b_i]=1\rangle.
\eqn

Consider the {\em representation variety} $\Rr(\G_g, G) =  G^n/\Pp$, 
where $\Pp$ is generated by the one relation $  \Pi_{i=1}^g [a_i, b_i]=1$.

The  Zariski tangent space to $\Rr(\G_g, G)$ at a point $\rho\in \Rr(\G_g, G)$ 
has dimension $(2g-1) \dim(G) + \dim \zent_G(\rho(\G_g))$, 
where $\zent_G(\rho(\G_g))$ denotes the centralizer in $G$ of the image $\rho(\G_g)$. (\cite{Goldman_84})

It follows that $\rho$ is a simple point if $(\zent_G(\rho(\G_g))$ is finite. 
A singular point of $\Rr(\G_g, G)$ is a simple point of a natural submanifold $\Rr(\G_g, G_1)$
of  $\Rr(\G_g, G)$, where $G_1= \zent_G(\zent_G(\rho(\G_g))$.
Denote the set of simple points in $\Rr(\G_g, G)$ by $\Rr(\G_g, G)^-$. 
The natural (algebraic) stratification of  $\Rr(\G_g, G)$ is given by
all $G$-orbits of submanifolds  $\Rr(\G_g, G_1)^-$. 

The action of $G$ is locally free on $\Rr(\G_g, G)^-$.  The quotient 
$\Rr(\G_g, G)^-/G$ is not necessarily Hausdorff. 
Consider the open subset  $\Rr(\G_g, G)$ of  $\Rr(\G_g, G)^-$, 
consisting of representations $\rho$, such that $\rho(\G_g)$ is not contained 
in a proper parabolic subgroup of $G$. 
Then the group $G$ acts properly on  $\Rr(\G_g, G)^{--}$, 
and the quotient $\Rr(\G_g, G)^{--}/G$ is an analytic manifold 
of dimension $(2g-2) \dim(G)$.

The group of outer 
automorphism ${\rm Out}(\G_g)$ can be identified with the mapping class group of $\Sigma_g$, acting thus in $\Xx(\G_g, G)$.
 

\subsection{Structure of the set of maximal representations}
We define the subset of maximal representations
\bqn
\Rr_{max}:=\left\{  \rho \in\Rr(\G_g, G)\, |\, \rho \text{ is maximal } \right\}  \subset \Rr(\G_g, G).
\eqn 
Since the Toledo invariant is locally constant on $\Rr(\G_g, G)$, the set $\Rr_{max}$ is a union of 
connected components of $\Rr(\G_g, G)$. 

The number of connected components can be in certain cases computed using 
Higgs bundles and related methods. 
Of particular interest is the result of Gothen 
(\cite{Gothen}) showing that for $G= {\rm SU}(2,2)$ the subset of maximal 
representations is connected, whereas for $G={\rm Sp}(4,\RR)$ the number of 
connected components grows exponentially in the genus $g$ of the surface 
group. Namely there are $3(2^{2^g}) +2g-4$ connected components in $\Rr_{max}$.
This implies that there are maximal representations which cannot 
be deformed to a representation that extends to a homomorphism of 
${\rm SL}(2,\RR)$.

A maximal representation is in particular a tight homomorphism, hence  
\bqn
\Rr_{max}\subset \Rr_t \subset \Rr_{np}.
\eqn
Denote the simple points of $\Rr_{max}$, i.e. the set of maximal representations 
with $\zent_G(\rho(\G_g) = \zent_G(\Hh(\rho))$ being finite,  by 
$\Rr^-_{max}= \Rr_{max}\cap \Rr^-(\G_g, G)$. 
Then   $\Rr^-_{max} \subset \Rr_{np,fd}$.
At the non simple points we have locally a product structure 
\bqn
\Rr_{max}(\G_g, G) = \Rr_{max, fd}(\G_g, G_1) \times \Rr(\G_g, C),
\eqn 
where $G_1$ is the semisimple noncompact part of the 
Zariski closure of $\rho(\G_g)$ and $C$ is the compact part of the 
Zariski closure of $\rho(\G_g)$. The set $\Rr(\G_g, C)$ is compact.
In particular, the results of Parreau (\cite{Parreau}) on compactifications of 
representation varieties apply to the subset of maximal representations, providing a compactification of 
$\Xx_{max}(\G_g, G):= \Rr_{max}(\G_g, G)/G$.

\subsection{Relation to the Hitchin component}
Starting from the theorem of Goldman (\cite{Goldman_thesis}), that the 
maximal part of the character variety 
\bqn
\Xx_{max}(\G_g, {\rm PSL}(2,\RR)) 
\subset
\Xx(\G_g, {\rm PSL}(2,\RR))
\eqn
are the Teichm\"uller spaces  of $\Sigma_g$, one may view our results 
on maximal rerpresentations as a generalization of Teichm\"uller space in the 
context of Hermitian symmetric spaces. 

In the case when  $G$ is a split real simple Lie group Hitchin \cite{Hitchin}, in a very different way using Higgs bundles, 
singled out one connected component of the 
character variety $\Xx(\G_g, G)$, called {\em Hitchin component}, 
containing a copy of Teichm\"uller space. 
He proves that these components are smooth and homeomorphic to a ball of 
dimension $|\chi(\Sigma_g)| \dim(G)$. Unfortunately he obtains no information 
about geometric properties of the representation contained in this component.
The only split real simple group which is of Hermitian type is the group 
${\rm Sp}(2n,\RR)$. The Hitchin component in this case is the deformation 
space of representation $\rho_0: \G_g \to{\rm PSp}(2n,\RR) $ obtained as 
composition of a Fuchsian representation $\G_g\to {\rm PSL}(2, \RR)$ with 
the irreducible representation ${\rm PSL}(2,\RR) \to  {\rm PSp}(2n,\RR)$. 

From this we obtain the following 
\begin{cor}
The Hitchin component $\Rr_{Hitchin}(\G_g, {\rm PSp}(2n,\RR))$ 
is contained as a (proper) subset in $\Rr_{max}(\G_g, {\rm PSp}(2n,\RR))$. 
In particular all representations in   
$\Rr_{Hitchin}(\G_g, {\rm PSp}(2n,\RR))$ are discrete and faithful.
\end{cor}

\begin{rem}
The number of connected components of 
$\Rr_{max}(\G_g, {\rm Sp}(4,\RR))$ grows exponentially in the genus of 
the surface \cite{Gothen}. Furthermore, also the number of connected
components of 
$\Rr_{max}(\G_g, {\rm Sp}(4,\RR))$ which are 
smooth and homeomorphic to Hitchin's component 
grows exponentially in the genus of the surface. 
There are $2^{2^g}$ such components. 
One might suspect that the action of 
$Out(\G_g)$ permutes different components.
\end{rem}

Recently the Hitchin component for $G= {\rm SL}(n, \RR)$ has been studied by Labourie \cite{Labourie_anosov}. In his very nice work he proves 
that all representations in the Hitchin component are discrete, faithful 
and loxodromic by relating them to a dynamical structure, 
which he calls Anosov 
representation. 
Following his suggestion we obtained a similar structure for all 
maximal representations.

\begin{rem}
Other generalizations of Teichm\"uller space for $G$ being a split
real semisimple Lie group have been recently defined 
by Fock and Goncharov (\cite{Fock_Goncharov}) for punctured surfaces. 
It seems that their positivity condition is similar to 
the condition that positive oriented triples on $S^1$ are mapped to 
maximal triples on the Shilov boundary, which 
we obtained after applying the methods from bounded cohomology. 
\end{rem}

\section{Anosov flows}\label{sec:anosov}
We are going to prove that maximal 
representations $\rho: \G_g\to G$ are holonomy representations of
Anosov structures on the Riemann surface $\Sigma_g$.

These structures are very similar to the structures Labourie obtained for the representations 
in the Hitchin component for $G={\rm SL}(n,\RR)$. 
Indeed, the strategy of our proof was suggested by him.

\subsection{A model for the geodesic flow}
Assume that $\Sigma_g$ is equipped with some hyperbolic metric such
that ${\widetilde\Sigma_g} =\DD$.
The lift of the unit tangent bundle of $\Sigma$ with respect to $\G_g$ is $T^1\DD$.
We denote by $g_t:T^1\Sigma_g\to T^1\Sigma_g$ the geodesic flow, and by $\tilde{g_t}$ 
its lift to $T^1\DD$. 
An element $u=(x,v)\in T^1\DD$ determines uniquely a positively oriented triple
$(u_+,u_0, u_-)$ of
(distinct) points on $S^1$, 
where $u_\pm$ are the endpoints for $t\to \pm \infty$ of the geodesic in $\DD$ that is
determined by $u$ and $u_0$ is the endpoint of
the geodesic determined by the vector orthogonal to the vector $v$ in
$T_x\DD$ such that $(u_+,u_0, u_-)$ is a positively oriented triple in $S^1$. 
Using this, we identify the unit tangent bundle with $(S^1)^{(3)}_+$, the space of distinct positively 
oriented triples in $S^1$, 
\bqn
T^1\DD\cong (S^1)^{(3)}_+.
\eqn 
In this model the lift of the geodesic flow acts as $\tilde{g_t} (u_+,
u_0, u_-) =(u_+, u_t, u_-)$ with $u_t \to u_\pm$ for $t \to \pm \infty$.

\subsection{Anosov representation}
Let $G$ be a connected semisimple Lie group and $B\subset G/Q\times
G/Q$ an open subset of a generalized flag variety defined by some
parabolic subgroup $Q<G$. 

The distributions associated to the product foliation of $B$ give rise
to a decomposition of the tangent bundle $TB=TB^+\oplus TB^-$, 
\begin{defi} (see \cite{Labourie_anosov})
A representation $\rho:\G_g\to G$ is said to be an {\em Anosov
  representation modelled on $(G,B)$} if there exists a continuous $\rho$- and $\tilde{g_t}$
equivariant 
map
\bqn
F:T^1\DD\to B
\eqn
such that the vector bundle $E=F^*(TB)$, defined as the pull-back of
the tangent bundle of $B$ decomposes under the lift of the geodesic
flow $\psi_t$, into two subbundles $E^\pm= F^*(TB^\pm)$ such that
$(\psi_t)_{t\geq 0}$ acts contracting on $E^+$ and 
$(\psi_t)_{t\leq 0}$ 
acts contracting on $E^-$. 
\end{defi}

\subsection{Maximal representations are Anosov}
Let now $\rho: \G_g\to G$ be a maximal representation into the isometry group 
of a Hermitian symmetric space $X$ of tube type. 
Let $B\subset \cs\times\cs$ be the equivariant Bruhat cell in the 
Shilov boundary of $X$.

Denote by $\phi^\pm$ the right- respectively left-continuous $\rho$-equivariant 
boundary maps $\phi^\pm: S^1 \to \cs$.
Define a map 
\bqn
F:T^1\DD\to B, \,\, (u_+, u_0, u_-)\mapsto (\phi^+(u_+), \phi^-(u_-)).
\eqn
This map is well-defined (Lemma~\ref{lem:graph_transverse}), 
$\rho$-equivariant, and clearly 
invariant under the geodesic flow, i.e. $F\circ \tilde{g_t} (u) = F(u)$. 

The equivariant Bruhat cell inherits 
two foliations $\Ee^\pm$ from the product structure 
of $B\subset\cs\times\cs$: 
\bqn
\Ee^+_{(z_+,z_-)} &=& B(z_-)\times z_-,\\
\Ee^-_{(z_+,z_-)} &=& z_+\times B(z_+).
\eqn
The corresponding distributions give a decomposition of the 
tangent bundle of $B$, 
\bqn
TB= E^+\oplus E^-
\eqn
 with 
\bqn
E^+_{(z_+,z_-)} &=& T_{z_+}B(z_-) = T_{z_+} \cs,\\
E^-_{(z_+,z_-)} &=& T_{z_-}B(z_+) = T_{z_-} \cs.
\eqn

We pull back the tangent bundle of $B$ via the map $F$ 
to obtain a (measurable) vector bundle 
\bqn
E= F^* TB \to T^1\DD.
\eqn
Since $F$ is $\rho$-equivariant this vector 
bundle descends to a 
vector bundle over the unit tangent bundle of the surface, 
\bqn
E_\rho\to T^1\Sigma_g.
\eqn
Because the map $F$ is invariant under the geodesic flow, we may 
lift the geodesic flow 
$g_t$ to a flow 
\bqn
\tilde\psi_t:E\to E,
\eqn
defined by $\tilde\psi_t(u,v) = (\tilde g_t u , v)$. 
This descends to a flow $\psi_t: E_\rho\to E_\rho$.

Since $T^1\Sigma_g$ is compact, all continuous scalar products on
$E_\rho$ are equivalent and 
we may fix one. 

\begin{thm}\label{thm:anosov_splitting}
Let $G$ be the isometry group of a Hermitian symmetric space of tube type.  
Assume that $\rho:\G_g\to G$ is a maximal representation. 
Then there is a continuous 
$\tilde{\psi}_t$-invariant splitting 
\bqn
E_\rho= E^+_\rho \oplus E^-_\rho
\eqn
such that \\
(1) The splitting $ E_{\rho, u}= E^+_{\rho , u} \oplus E^-_{\rho, u}$ is compatible with the splitting of $TB$, 
i.e. $ E^\pm_{\rho , u}= F^*(E^\pm)_u$. \\
(2) The flow $(\psi_t)_{ t\geq 0}$ acts contracting on $E^-_\rho$ and 
$(\psi_t)_{t\leq 0}$ acts contracting on $E^+_\rho$.
\end{thm}

\begin{cor}
Maximal representations are Anosov representations modelled on $(G,
B)$, where $B\subset \cs\times\cs$ is the equivariant Bruhat cell.
\end{cor}

We sketch the structure of the proof of
Theorem~\ref{thm:anosov_splitting}:
Given a maximal representation $\rho:\G_g\to G$, we use the
$\rho$-equivariant right-/left-continuous boundary maps $\phi^\pm:S^1\to \cs$ to
pull-back the tangent bundle of $B\subset \cs\times \cs$ to a bundle $E$
 over the unit tangent bundle of the surface $\Sigma_g$. 
Using 
the parametrization of the tangent space $T_x \cs$ of the Shilov
boundary by the Euclidean Jordan algebra $V$, we define fibrewise a norm
on $E$. This norm is a priori not continuous on $E$. But, with respect
to this norm we obtain a measurable splitting of $E=E^+ \oplus E^-$. 
With this splitting and a comparison of the (measurable) norm to a
continuous norm on $E$ on a bounded subset, we can show that the
splitting 
$E=E^+ \oplus E^-$ is infact continuous. 

\subsection{A norm on $E_\rho$}
Recall that we parametrized the Bruhat cell $B(z)\subset \cs$ using the
Euclidean Jordan algebra $V$. 
In particular we may identify the tangent space 
of $B(z)$ at a point $z'\in B(z)$ canonically with the vector space $V$.
Any element in $w\in V$ admits a unique spectral decomposition and hence 
a spectral norm defined by 
$s(w) = \sum_{i=1}^r \lambda_i$.

For the rest of the chapter we work with this parametrization of the
tangent space of the Bruhat cell. Thus for the tangent space of the
equivariant Bruhat cell $B\subset \cs\times \cs$ we obtain:
\bqn
T_{(z_+, z_-)}B= E^+_{(z_+,z_-)} \oplus E^-_{(z_+,z_-)} = 
T_{z_+}B(z_-) \oplus T_{z_-}B(z_+)= V_{z_+}\oplus V_{z_-}
\eqn
We are going to define a fiberwise norm on 
$E= F^*TB$, making use of the generalized Maslov 
cocycle on $\cs$. 

Let $v\in E_u$, i.e. $v\in T_{(z_+, z_-)}B$ where $z_\pm=\phi^\pm(u_\pm)$. 
Let $v= (v^+, v^-)$ be the decomposition of $v$ with respect to 
$E^+\oplus E^- =TB$.
Let $l(v^\pm)$ be the norm of $v^\pm$ with respect to the inner
product in $T_{z_+}\cs$ given by the invariant Riemannian metric on $\cs$.
We define a norm on $(F^*E^\pm)$ by scaling this norm with a factor
depending on the base point $u=(u_+, u_0, u_-)$.
The images of $u_0$ under $\phi^\pm$ are transverse to the images of $u_\pm$ under the maps 
$\phi^\pm$ (Lemma~\ref{lem:graph_transverse}). 
Thus $\phi^\pm(u_0)$ determine 
points in $B(z_\pm)$, and 
hence a vector 
\bqn
(v_0^+, v_0^-)=(\phi^+(u_0)-\phi^+(u_+), \phi^-(u_0)-\phi^-(u_-))
\in V_{\phi^+(u_+)}\oplus V_{\phi^-(u_-)},\\
V_{\phi^+(u_+)}\oplus V_{\phi^-(u_-)}= E_{F(u)}^+\oplus
E_{F(u)}^-=T_{F(u)} B.
\eqn
Since $V^\pm= V_{\phi^\pm(u_\pm)} $ inherit the structures of Euclidean
Jordan algebras, we may associate to 
the elements $v_0^\pm$ their spectral norms  $s_\pm(v_0^\pm)$ 
(see Chapter~\ref{sec:shilov}). 
We define the norm on $(F^*E^\pm)_u$ by 
\bqn
q^\pm_{u}(v^\pm):=  \pm s_\pm(v_0^\pm) l(v^\pm).
\eqn

This induces a norm $q_{u}(v) = q^+_{u}(v^+)+q^-_{u}(v^-)$ on the fibers of 
$F^*TB$.
\begin{lemma}
The map $q_u: (F^*TB)_u \to \RR$ is a norm
\end{lemma}
\begin{proof}
We prove that $q_u^+: (F^*E^+)_u \to \RR$ is a norm. A similar
argument shows that $q_u^-$ is a norm on $(F^*E^-)_u$. Since
$E^+\oplus E^-=TB$ this implies that $q_u$ is a norm on $(F^*TB)_u $.

The Maslov cocycle evaluated on $(\phi^+(u_+), \phi^+(u_0), \phi^-(u_-))$ is 
the signature of $v_0^+=(\phi^+(u_0)-\phi^+(u_+))$. 
The triple $(u_+, u_0, u_-)$ is positively oriented, hence
$(\phi^+(u_+), \phi^+(u_0), \phi^-(u_-))$ is a maximal triple. 
The maximality means that $\sign(v_0^+)= r$. 
Therefore, all eigenvalues $\lambda_i$ of the spectral decomposition are 
positive, in particular $s(v_0^+)> 0$. So we scaled the norm
$l(v^+)$ in the fiber  $(F^*E^+)_u$
with the positive factor $s(v_0^+)$ depending on $u$. 
\end{proof}
\begin{rem}
From Lemma~\ref{lem:converge} it follows that $s_\pm(v_0^\pm) \to 0$
if $u_0\to u_\pm$. 
\end{rem}

We want to determine the ``expanding'' and ``contracting'' 
subbundles of $E_\rho$ for the flow $\psi_t$.
For $e\in E_u$ let $||e||_u=q_u(e)$ be the norm defined by the above 
construction. 

Define 
\bqn
E^+_u &:=& \left\{  e \in E_u \,|\, 
\lim_{t\to \infty} || \psi_t(e)||_{g_t(u)} = 0\right\} \\
E^-_u &:=& \left\{  e \in E_u \,|\, 
\lim_{t\to -\infty} || \psi_t(e)||_{g_t(u)} = 0\right\} .
\eqn

\begin{lemma}\label{lem:subbundle} 
We have that $E^+_u= (F^*E^+)_u$ and $E^-_u = (F^*E^-)_u$.
\end{lemma}
\begin{proof}
We prove only the first assertion since the second follows by 
similar arguments.
Let $e\in (F^*E^+)_u$ then $|| \tilde\psi_t(e)||_{g_t u}= || e||_{g_t u}= 
q^+_{g_t u}(e)$. But, since $g_t u= (u_+, u_t, u_-)$ with $u_t \to u_+$
as $t\to \infty$,  
Lemma~\ref{lem:converge} in Chapter 5 implies that 
$q^+_{g_t u}(e) \to q^+_{u^+}(e)=0$ for $t\to \infty$, where $u^+=(u_+,u_+, u_-)$. Hence 
\bqn
(F^*E^+)_u \subset E^+_u.
\eqn 
To prove the opposite inclusion we may work on $E$ with 
respect to the flow $\tilde\psi_t$.
Let $(u,e) \in E=F^*TB$ and assume that 
$\lim_{t\to \infty} ||e||_{g_t u} = 0$. 
Observe that $q_{g_t u}(e) = q^+_{g_t u} (e_+) + q^-_{g_t u}(e_-)$, 
where $e=e_+ + e_-$ is the decomposition with respect to 
$(F^*TB)_u =(F^*E^+)_u \oplus(F^*E^-)_u $. 
Then we have $q^-_{g_t u}(e_-)\geq q^-_u(e_-)$ by 
Lemma~\ref{lem:ordering} in Chapter 5. 
But since $\lim_{t\to \infty} ||e||_{q_{g_t u}} = 0$, this gives that
$e_-=0$, hence 
$e=e_+ \in (F^*E^-)_u$.
\end{proof}

This Lemma gives a weak splitting $E=E^+\oplus E^-$ that is a priori
not continuous but at least measurable. We want to show, that this splitting is indeed continuous.

Denote by $\Gr_n(TB)$ the Grassmannian bundle, 
where the fibres are Grassmannians of 
$n$-planes in $T_x B$. The splitting 
$E=E^+\oplus E^-$ defines two maps 
\bqn
\Phi^\pm: T^1\DD \to \Gr_n(TB), \quad u\mapsto  E^\pm_{u}.
\eqn

\begin{prop}\label{prop:continuous}
The maps $\Phi^\pm: T^1\DD \to \Gr_n(TB)$, defined by 
$u \mapsto E^\pm_{u}$, are continuous.
\end{prop}
We need the following Lemma:
\begin{lemma}\label{lem:bounded}
Let $K \subset T^1\DD$ be compact. 
Then the sets
\bqn
J_\pm:=\left\{  v_0^\pm=\phi^\pm (u_0) -\phi^\pm (u_\pm) \,|\, u \in K\right\} 
\eqn
are bounded in $V$.
\end{lemma}
\begin{proof}
 If $K\subset (S^1)^{(3)}$ is a compact set, then the subset 
\bqn
P:=\bigcup_{(x_1,x_2,x_3)\in K} (F_{x_1}\times F_{x_2}\times F_{x_3}) \subset 
\cs^{(3)}
\eqn
is compact. 
Indeed,  let $F^3:=\left\{  ((x_1,x_2,x_3),(f_1,f_2,f_3)) \,|\, (x_i,f_i) \in F\right\} $, 
then $F^3 \subset (S^1)^3 \times \cs^3$ is closed, hence 
$F^3\cap (K\times \cs^3)$ is compact and hence $P$, its projection to $\cs^3$ 
is compact. 

Denote by $\Tt:\cs \to
V$ the Cayley transform.\\
Claim: If $K_1\subset \cs^{(3)}$ is a compact subset
then 
\bqn
\left\{ \Tt(f^0)-\Tt(f^\pm)\,|\, f=(f^+, f^0, f^-)\in K_1\right\} \subset V
\eqn 
is bounded. 
The Cayley transform defines a map $\Tt_0:\cs^{(3)} \to V$ by 
determined by mapping a triple $(f^+, f^0, f^-)$  with $\Tt(f^+) =0,
\Tt(f^-)=\infty$ to the point $\Tt(f^0) \in V$.
This map is continuous.
Assume that $(f^+_n, f^0_n, f^-_n)$ is a sequence of 
maximal triples of pairwise transverse points in $\cs$ converging to a
maximal triple of pairwise transverse points
$(f^+,f^0,f^-)$. Using the transitivity of the action of 
$G$ on transverse tuples, we may assume that $f^+_n=f^+$ and 
$f^-_n=f^-$ are fixed  with $\Tt(f^-)=\infty$, $\Tt(f^+) =0$, and 
$f^0_n\to f^0$. But then the images under $\Tt$ converge
$\Tt(f^0_n)\to \Tt(f^0)$ since the Cayley transform is continuos on
$B(f^-)$.
Thus the sets $J_\pm$ are bounded.
\end{proof}
\begin{proof} (of Proposition~\ref{prop:continuous})
We prove that $\Phi^-$ is continuous. The result for $\Phi^+$ follows
by a similar argument.\\
Let $(u_k)\to u$ in $T^1\DD$. Consider the sequence of the images 
\bqn
L^-_k := E^-_{u_k}=\Phi(u_k)
\eqn
 in $\Gr_n(TB)$.
Let $L$ be any accumulation point of $(L^-_k)$. 
Passing to a subsequence 
$(u_l)$, we have $u_l\to u$, $L^-_l\to L$.
We have to show that $L=\Phi(u) = E^-_{u}$.
Now $\left\{ u_l ,u\right\} $ is compact.

It follows form Lemma~\ref{lem:bounded} that the set of norms $\left\{  q_{(u_k)}, q_u\right\} $ is bounded. 
Hence there exists a $\G_g$-invariant continuous scalar product on $E$, with 
induced norms $||\cdot||_u$ on $V$ and  
a constant $\lambda>0$ such that 
\bqn
\lambda^{-1} ||\cdot||_{u_k}\leq q_{u_k}(\cdot)\leq \lambda ||\cdot||_{u_k}\\
\lambda^{-1} ||\cdot||_u\leq q_{u}(\cdot)\leq \lambda ||\cdot||_u
\eqn
We claim that $L \subset L^-_u$. Then for dimension reasons $L=L^-_u$.\\
Given an element $e\in L$ there exists a sequence $e_l\in L^-_l$ such that 
$e_l \to e$ in $TB$.
We have to show that $\lim_{t\to \infty} q_{g_t u}(e)=0$.
We know that $\lim_{t\to \infty} q_{g_t u_l}(e_l) =0$ and by the above
comparison also that $\lim_{t\to \infty} ||\psi_t(u_l,e_l)||_{g_t u_l} =0$.
We want to show that 
$\lim_{t\to \infty} ||\tilde\psi_t(u,e)||_{g_t u} =
0$, since then also $\lim_{t\to \infty} q_{g_t u}(e) =0$ holds.
So far all our considerations where about the bundle $E$, now we are 
using the fact, that $\G_g$ acts cocompactly on $T^1\DD$. \\
Since the flow $\tilde\psi_t$ is $\G_g$-invariant, we have that 
\bqn
||\tilde\psi_t(u_l,e_l)||_{g_t u_l}= ||\g_*
\tilde\psi_t(u_l,e_l)||_{\g g_t u_l}.
\eqn
Furthermore, since $\Gamma_g$ acts cocompactly, for any $\eps >0$
there exist $N\in \NN$ and $\eta\in \G_g$  such for all $l\geq N$ we have 
\bqn
||\g_* \tilde\psi_t(u_l,e_l) - {\eta}_* \tilde\psi_t(u,e)||\leq \eps.
\eqn
Hence 
\bqn
||\g_* \tilde\psi_t(u_l,e_l)||_{\g g_t u_l} -
||{\eta}_* \tilde\psi_t(u,e)||_{\eta g_t u} \in [-\eps, \eps].
\eqn
Since this holds for every $\eps >0$, 
we have that  $\lim_{t\to \infty} ||\tilde\psi_t(u,e)||_{g_t u} =
0$ and the above comparison of $||\cdot||$ and the norm 
$q_u$ implies the claim.
\end{proof}
\begin{proof}(Proof of Theorem~\ref{thm:anosov_splitting})
It is a direct consequence from Proposition~\ref{prop:continuous}, that the 
splitting $E_\rho= E^+_\rho \oplus E^-_\rho$ with 
$E^\pm_\rho:=\left\{ (u,e)\in E\,|\,v\in  E^\pm_{\rho , u}\right\} $ is continuous.
The statements (1) and (2) follow from Lemma~\ref{lem:subbundle}.
\end{proof}
\subsection{The boundary map} 
The continuity of the maps $\Phi^\pm$ imply the continuity of the
boundary maps $\phi^\pm: S^1 \to \cs$.

\begin{cor}\label{cor:map_continuous}
The boundary map $\phi: S^1\to \cs$ is continuous. 
The essential graph $F\subset S^1\times \cs$ is the graph of $\phi$.
\end{cor}
\begin{proof}
The boundary maps $\phi^\pm$ are the compositions of the maps $\Phi^\pm$ from Lemma~\ref{prop:continuous} with the 
projection onto the base point. Since $\Phi^\pm$ are continuous, the maps $\phi^\pm$ are continuous.
In particular $\phi=\phi^-=\phi^+$ and $F= graph (\phi)$.
\end{proof}

\begin{cor}\label{cor:map_rectifiable}
$F\subset S^1\times \cs$ is rectifiable of uniformly bounded length.
\end{cor}
\begin{proof}
We cover $\cs$ by $2^{r_X}$ Bruhat cells $B(z)\subset \cs$ with
respect to well-chosen different base points $z\in \cs$. 
The parametrization of the Bruhat cell by the Jordan algebra $V$ gives local 
coordinates on $\cs$ and we can parametrize $F= {\rm graph}(\phi)$ locally by 
maps $I \to V\cong B(z), \, t\to z_t\in V\cong \RR^N$. 
Since for every triple $t_1<t_2<t_3$ in $I$ the generalized Maslov
cocycle $\tau(z_{t_1}, z_{t_2}, z_{t_3})$ is maximal, 
these maps are monotone increasing. Hence $F$ is rectifiable and the
length is bounded by constants depending only on the dimension of $\cs$.
\end{proof}
\subsection{Limit sets of maximal representations} 
\begin{thm}\label{thm:limit_set}
Let $\Dd$ be a bounded symmetric domain and $\rho: \G_g\to \aut(\Dd)$ 
a maximal representation. Then, the Shilov limit set 
\bqn
\Ll_\rho:=\Ll_t(\rho(\G_g)) \subset \cs
\eqn 
of $\rho(\G_g)$ is a rectifiable circle, and there is a $\G_g$-equivariant homeomorphism 
$\phi:S^1 \to \Ll_\rho$.
\end{thm}
\begin{proof}
The image of $\phi:\partial\G_g\cong S^1\to \cs$ is the projection of the image of $\phi:S^1\to G/P$ to $G/P_{\theta_r}$. 
The Benoist limit set $\Lambda_{\theta_r}(\G_g)$ is the projection of the Benoist limit set $\Lambda_{\emptyset}(\G_g)\subset G/P$
to  $G/P_{\theta_r}$ By Lemma~\ref{lem:benoist} in Chapter 7, they coincide. But the Benoist limit set $\Lambda_{\theta_r}(\G_g)$ coincides with the 
Shilov limit set  $\Ll_\rho:=\Ll_t(\rho(\G_g))$. The map $\phi$ 
is an homeomorphism onto its image by
Corollary~\ref{cor:map_continuous} and by
Corollary~\ref{cor:map_rectifiable} the image of $\phi$ in $\cs$ is
rectifiable.
\end{proof}

\begin{rem}
The rectifiability of the limit set for any maximal representation 
contrasts the behavious of limit sets of 
quasifuchsian representation, i.e. deformations of 
Fuchsian representations in ${\rm PSL}(2,\CC)$, where the limit set is 
only rectifiable if the representation is Fuchsian \cite{Bowen}. 
This might be considered as a higher rank phenomenon 
since the limit curves obtained by 
Labourie for representations in the Hitchin component of ${\rm PSL}(n,\RR)$ 
\cite{Labourie_anosov} are also rectifiable.
\end{rem}

\subsection{Tight homomorphisms and Anosov flows}
Let $H={\rm Is}(Y)^\circ$ be the connected component of the 
isometry group of a Hermitian symmetric space of tube type admitting a tight 
embedding $f:Y\to X$ into another Hermitian symmetric space of tube type $X$ with 
connected component of the isometry group being $G$. 

Associated to the tight embedding $f: Y\to X$ we consider the 
tight homomorphism 
$\pi:H_Y\to G$, with $H_Y$ being a suitable finite extension 
of $H$ by $\Lambda$, and the 
$\pi$-invariant extension (Proposition~\ref{tight_embedding} in Chapter 6) 
$\ol{f}:\cs_Y\to \cs_{X}$.

Let $\rho: \G_g\to H$ be a maximal representation. Via $\rho$ we construct a finite 
extension $\tilde{\G_g}$ of $\G_g$ corresponding to $H_Y$, such that we 
obtain a representation $\tilde\rho:\tilde\G_g\to H_Y$, which fits in 
the following commutative diagram:
\bqn
\xymatrix{
&\Lambda\ar[d] &\\
\tilde{\G_g} \ar[d] \ar[r]^{\tilde\rho} & H_Y \ar[d] \ar[r]^\pi 
& G\\
\G_g \ar[r]^\rho &  H &
}
\eqn
Since the respective actions of $\tilde\G_g$ and of $H_Y$ on $S^1$
 and on $\cs_Y$ coincide with the corresponding actions of $\G_g$ respectively of $H$, the 
$\rho$-equivariant boundary map $\phi:S^1 \to \cs_Y$ lifts to a 
$\tilde\rho$-equivariant map $\tilde\phi:S^1\to \cs_Y$.

The tight embedding $\pi$ induces (Proposition~\ref{tight_embedding} 
in Chapter 6) a 
$\pi$-equivariant map 
\bqn
\Phi: \cs_Y \to \cs_X.
\eqn
The composition of $\tilde\phi$ 
with $\Phi$ is a $\pi\circ\tilde\rho$-equivariant boundary map 
$\phi_H:S^1\to \cs_X$, which maps almost all triples 
of distinct point in $S^1$ to maximal triples in $\cs_X$. 
We thus may pull back the Anosov-structure modelled on $(G, B_X)$ with $B_X\subset \cs_X\times\cs_X$ to 
an Anosov-structure modelled on $(H_Y, B_Y)$ with $B_Y\subset \cs_Y\times\cs_Y$.
\begin{rem}
All classical domains admit tight embeddings into the Siegel upper half space.
It is not known to the author whether 
there exists a tight embedding of the exceptional domain of tube type into a classical domain of tube type. 
It is well known that there do not exist holomorphic tight
embeddings of this kind (\cite{Satake_hol_65}).
\end{rem}

\subsection{The Anosov flow in the classical case}\label{sec:anosov_class}
Assume that $G$ is one of the classical groups $ {\rm Sp}(2n,\RR)$, ${\rm SU}(n,n)$ or ${\rm SO}^*(4n)$. 
In these cases, we can construct an Anosov structure on the 
associated vector bundle.

Recall from Section~\ref{sec:maslov_classic} 
that $G$ is the invariance group  of an $2n$ dimensional 
$\KK$-vectorspace $V$, $\KK=\RR,\CC,\HH$, equipped with a $\sigma$-symmetric nondegenerate linear 
form $h$. The Shilov boundary $\cs=G/Q$ is the space of 
maximal $h$-isotropic 
subspaces, where two of them, $L_1, L_2$, are transverse if and only if 
$V = L_1\oplus L_2$. 
Recall from Section~\ref{sec:maslov_classic} in Chapter 5 that given $L_-$ transverse to $L_+,L_1, L_2$, we can 
associate to $L_i$ a quadratic form $q_i^+$ on $L_+$ that is positive definite if and only if 
$L_+, L_i, L_-$ is a maximal triple. If the triple $L_1, L_2, L_-$ is maximal, 
then $q^+_1<q^+_2$.

\begin{rem}
Changing the roles of $L_+$ and $L_-$ we can define a quadratic form 
$q^-_i$ on $L_-$, satisfying up to a sign change the same properties as $q^+_i$.
\end{rem}

Given a representation $\rho:\G_g\to G$, we fix a hyperbolic metric on 
$\Sigma_g$
and consider the associated flat bundle over the unit tangent 
bundle of $\Sigma$:
\bqn
E_\rho:= (T^1\DD \times V)/\G_g \to T^1\Sigma,
\eqn
where $\G_g$ acts diagonally on $T^1\DD \times V$, 
$\g(u,e)=(\g u, \rho(\g) e)$.
We can lift the geodesic flow $g_t$ on $T^1\Sigma$ to a flow 
$\tilde\psi_t$ on $\tilde E$ by 
\bqn
\tilde\psi_t(u,v) = (\tilde g_t (u),v),
\eqn 
where $\tilde g_t$ is the lift of the geodesic flow to $T^1\DD$. 
The flow $\tilde\psi_t$ is $\G_g$-invariant and descends to a flow 
$\psi_t$ on $E_\rho$.

All continuous scalar products on $E_\rho$ being equivalent, 
we fix one. 

\begin{thm}\label{thm:anosov_splitting_class}
Let $G={\rm SU}(V, h)$. 
Assume that $\rho:\G_g\to G$ is a maximal representation. 
Then there is a continuous 
$\tilde{\psi}_t$-invariant splitting 
\bqn
E_\rho= E^+_\rho \oplus E^-_\rho
\eqn
such that \\
(1) The splitting $ E_{\rho, u}= E^+_{\rho , u} \oplus E^-_{\rho, u}$ is a decomposition into complementary maximal $h$-isotropic subspaces.\\
(2) The flow $(\psi_t)_{ t\geq 0}$ acts contracting on $E^-_\rho$ and 
$(\psi_t)_{t\leq 0}$ acts contracting on $E^+_\rho$.
\end{thm}
\begin{proof}
We define a metric on the vector bundle $E$ similarly to the
 definition of the metric on the tangent bundle of $B$ above
 using this time the quadratic forms $q_0^\pm$ instead of the spectral norms.
The arguments to prove the Theorem are then analogous to 
the ones used to prove Theorem~\ref{thm:anosov_splitting}.
\end{proof}


\vskip1cm
\chapter*{Appendix}

\section*{Some facts about Lie groups and Lie algebras}
\subsection*{Satake diagrams}
Complex simple Lie algebras are uniquely determined by their 
root systems, which is described by its Dynkin 
diagram. 
Since a complex simple Lie algebra has different noncompact real
forms,   
the noncompact real forms are not determined by the Dynkin diagram of their 
complexification alone. They are determined by their {\em Satake diagram}, which is the 
Dynkin  diagram of their complexification endowed with some additional 
information.
We give a short description how the Satake diagram is defined. (A  
more detailed description can be found in \cite{Onishchik_Vinberg_III}).
Consider a real simple Lie algebra $\frakg=\frakk\oplus \frakp$ and a maximal 
abelian subalgebra $\fraka\subset \frakp$. Denote by $\frakm$ the Lie algebra of 
the centralizer of $\fraka$ 
in $K$. Let $\Xi=\Xi(\frakg, \fraka)$ be the real root system of 
$\frakg$ with respect to $\fraka$. 
The Cartan subalgebras of $\frakg$ containing $\fraka$ are of the form 
\bqn
\frakh=\frakh^+ \oplus \fraka, 
\eqn
where $\frakh^+$ is a Cartan subalgebra of $\frakm$.
The complexification $\frakh_\CC$ of $\frakh$ is a Cartan subalgebra 
of the complexification $\frakg_\CC$ of $\frakg$. 
We consider the connection between the complex root system 
$\Psi=\Psi(\frakg_\CC, \frakh_\CC)$ of  $\frakg_\CC$ with respect to $\frakh_\CC$
and the real root system $\Xi(\frakg, \fraka)$.

The restriction map $r: \frakh_\CC(\RR)^* \to \fraka^*$ is surjective, 
but it is not injective. Denote by $\Psi_0$ the set of 
roots $\a\in \Psi$ such that $r(\a) =0$. 
In $\Psi_1=\Psi\backslash\Psi_0$ there are roots $\a_1,\a_2$ such 
that $r(\a_2)=r(\a_1)$. 

Choose a basis of $\frakh_\CC(\RR)$, such that the first basis vectors 
form a basis of $\fraka$, and consider the corresponding lexicographical orders on 
$\frakh_\CC(\RR)^*$, respectively on $\fraka^*$. 
The restriction map sends positive to positive and negative to 
negative roots. 
The fundamental set of simple positive roots in $\Xi$ is obtained from 
the set of those simple positive roots of $\Psi$, which are contained in 
$\Psi_1$.

The Satake diagram of $\frakg$ 
is now obtained from the Dynkin diagram of the complexification $\frakg_\CC$ 
by blackening all vertices corresponding to simple positive roots of $\Psi$, 
which are contained in $\Psi_0$, and by connecting all 
simple positive roots contained in $\Psi_1$ by arrows  
 whenever they are mapped to the same simple positive root in $\Xi$.
Denote by $b_0$ the number of black vertices and by $a$ the number of arrows. 
Then 
\bqn
\rk(\frakg_\CC) = \rk_\RR(\frakg) + b_0+a.
\eqn
In particular a real Lie algebra is split if and only if the 
Satake diagram coincides with the Dynkin diagram of its complexification.

\subsection*{Principal and semiprincipal subalgebras}
For the proof of Lemma~\ref{irreducible} in Chapter 6 we use the classification 
of special three dimesional simple subalgebras of complex semisimple
Lie algebras which are called semiprincipal and principal. This classification is due to 
Dynkin \cite{Dynkin}.
We just recall the basic facts, to help the reader to understand the 
arguments in Chapter 6. For a more detailed description we refer to the book 
of Onishchik and Vinberg \cite{Onishchik_Vinberg_III} and the translation 
of Dynkin's original article \cite{Dynkin}.

Assume that $\frakg$ is a simple complex Lie algebra.
A three dimensional simple subalgebra of $\frakg$ is given by a 
$\fraks\frakl_2$-triple in $\frakg$, i.e. a triple of 
elements $\{e,h,f\}$ such that the relations 
$[h,e]=2e$, $[h,f]=-2f$, $[e,f]=h$ hold. 
Every $\fraks\frakl_2$-triple is determined by the elements $e$ or $h$. 
The element $h$ is called {\em characteristic} of the $\fraks\frakl_2$-triple.

\begin{prop*}\cite{Dynkin} 
Fix a Cartan subalgebra $\frakh\subset\frakg$. Let 
$\Delta= \{ \a_1, \dots \a_k\}$ be a fundamental system of simple positive 
roots of $\Psi(\frakg,\frakh)$, and let $C$ be the corresponding Weyl chamber.
Then any characteristic of a $\fraks\frakl_2$-triple in $\frakg$ is 
conjugate to a unique element $h\in C$ with ``labels'' 
$\a_i(h)$ being $0$, $1$ or $2$.
\end{prop*}
The $\fraks\frakl_2$-triple can be described up to conjugacy 
by the charactersitc element $h$. Thus it is determined by the Dynkin diagram of $\frakg$ with 
labels $0$, $1$ or $2$ corresponding to the values of $\a_i(h)$ at the nodes.
The element $e$ of an $\fraks\frakl_2$-triple is always a nilpotent element. 
This element is called regular if $e$ is the sum of some  
roots vectors corresponding to the simple positive roots $\a_1, \dots a_k$, 
i.e. $e= e_{\a_1} + \dots + e_{\a_k}$.
The charactersitic of a regular element is conjugate to the element $h\in C$, 
where all labels $\a_i(h)$ are equal to $2$.

A three-dimensional simple subalgebra of $\frakg$ given by a  
$\fraks\frakl_2$-triple $\{e,h,f\}$ is called {\em principal} if 
the nilpotent element $e$ is regular. 
A three-dimensional simple subalgebra of $\frakg$ given by a  
$\fraks\frakl_2$-triple $\{e,h,f\}$ is called {\em semiprincipal} if 
the centralizer of this subalgebra in $G$ is trivial.
A principal three-dimensional simple subalgebra is semiprincipal.
Dynkin showed that the labels of the characteristic elements of a 
semiprincipal  three-dimensional simple subalgebra are $0$ or $2$. 
Therefore, the characteristics of semi-principal three-dimensional
simple subalgebras may be represented by the 
Dynkin diagram, where a vertex is blackened if the corresponding label
$\a_i(h)$ is $0$ 
and is white if the corresponding label $\a_i(h)$ is $2$. 
The characteristic of the principal three-dimensional simple subalgebra 
is thus just represented by the Dynkin diagram with all nodes being white.
For the classification of the semiprincipal three-dimensional simple
subalgebras see \cite{Dynkin} or \cite{Onishchik_Vinberg_III}.

\subsection*{A simple fact}
\begin{prop*}
The Lie algebra $\frakn$ is nilpotent if and only if $N = \exp(\frakn)$ is unipotent. 
Then $\exp: \frakn \to N$ is a birational map, 
the inverse is given by $\log:N\to \frakn$.
\end{prop*}
\begin{proof}
We realize $\frakn$ and $N$ by matrices via its adjoint representations on $\frakn$. 
Then $\exp(\ad(X)) = \sum_{k=0}^\infty \frac{1}{k!} \ad(X)^k$. The Lie algebra $\frakn$ is 
nilpotent, thus for all $X\in \frakn$ there exists a $k\in \NN$ such that $\ad(X)^k = 0$. 
Thus the above map is a finite sum, hence a rational map. 
The logarithm has the following expansion 
$\log (\id +\Ad(Y)) = \sum_{k=0}^\infty \frac{(-1)^n}{k+1} \Ad(Y)^{k+1}$.
Thus if $Y= \exp(X) \in N$, we have 
$\ad(X)= \log(Y)= \log (\exp(\ad(X))= \log(\id + \Ad(\exp(X)) -\id) = \sum_{k=0}^\infty \frac{(-1)^n}{k+1} (\Ad(Y)-\id)^{k+1}$. 
But if $N$ is unipotent, then for any $Y\in N$ there is a $k\in\NN$, such that $(\Ad(Y)-\id)^{k+1} = 0$. 
Thus the sum is finite and hence $\log$ is a rational map.
\end{proof}

\section*{Some facts about Jordan algebras}
For an introduction into Jordan algebras we refer to the book of Faraut and 
Kor\'anyi \cite{Faraut_Koranyi}. We recall the basic
definitions and 
facts from there which are used in this work. 

A vector space $V$ over $\RR$ is said to be a {\em Jordan algebra} 
if there exists a bilinear mapping $V\times V \to V$, $(x,y)\mapsto xy$, the 
multiplication, which is commutative, i.e. $xy=yx$ for all $x,y\in V$ and 
satisfies for all $x,y \in V$:
\bqn
x(x^2y)= x^2(xy).
\eqn
A Jordan algebra is in general not associative. We always assume that $V$ has 
an identity element.

\subsection*{The minimal polynomial}
Denote by $\RR[X]$ the algebra of polynomials in one variable with 
coefficients in $\RR$.
For any $x\in V$ the ideal $\Jj(x):= \{p\in \RR[X]\,|\,  p(x)=0\}$ 
is generated by a monic polynomial $f$, called the minimal polynomial of $x$.
Its degree is called the rank of $x$. An element $x\in V$ is called regular 
it its rank is maximal, i.e. 
$\rk(x) = \max_{y\in V} \rk(y) = \rk(V)=r$.

\begin{prop*}\cite[Proposition~II.2.1]{Faraut_Koranyi}
The set of regular elements is open and dense in $V$. 
There exist polynomials $a_1, \dots, a_r$ on $V$ such that the minimal 
polynomial of any regular element $x\in V $ is given by 
\bqn
f(\lambda, x)= \lambda^r -a_1(x)\lambda^{r-1} + \cdots +(-1)^r a_r(x).
\eqn
The polynomials $a_i$ are unique and homogeneous of degree $i$. 
The coefficient $a_r(x)=\det(x)$ is called the {\em determinant of $x$}, 
the coefficient $a_1(x)=\tr(x) $ is called the {\em trace of $x$}.
\end{prop*}

\subsection*{Spectral decomposition}
A Jordan algebra is said to be Euclidean if there exists a positive definite 
symmetric bilinear form $\langle\cdot,\cdot\rangle:V\times V \to \RR$ such that 
$\langle xu,y\rangle= \langle u,xy\rangle$ for all $x,y,u\in V$.

An idempotent in $V$ is an element $c\in V$ satisfying $c^2=c$.
An idempotent is called primitive if it is non-zero and cannot be written 
as sum of two non-zero idempotents. 
A Jordan frame is a system of primitive idempotents $c_1, \dots, c_m$ 
satisfying 
\bqn
c_ic_k =0 \quad i\neq k\\
\sum_{i=1}^m c_i=e.
\eqn

\begin{prop*}\cite[Theorem~III.1.2.]{Faraut_Koranyi}
Suppose $V$ has rank r. Then for every $x \in V$ there exist a Jordan frame 
$c_1, \dots, c_r$ and real numbers $\lambda_1, \dots, \lambda_r$ such that 
\bqn
x=\sum_{i=1}^r \lambda_i c_i.
\eqn
The numbers $\lambda_i$ (with their multiplicities) are uniquely determined by 
$x$. 
Furthermore 
\bqn
\det(x) = \Pi_{i=1}^r \lambda_i, \quad \tr(x)= \sum_{i=1}^r \lambda_i. 
\eqn
Moreover, $a_i(x)$ is the $
i$-th elementary symmetric function in the $\lambda_i$'s.
\end{prop*}


\bibliographystyle{amsplain}
\bibliography{refs}
\providecommand{\bysame}{\leavevmode\hbox to3em{\hrulefill}\thinspace}
\providecommand{\MR}{\relax\ifhmode\unskip\space\fi MR }
\providecommand{\MRhref}[2]{%
  \href{http://www.ams.org/mathscinet-getitem?mr=#1}{#2}
}
\providecommand{\href}[2]{#2}

\vskip1cm
\end{document}